\documentclass[12pt]{article}

\usepackage{amsmath}
\usepackage{amssymb}
\usepackage{amstext}
\usepackage{graphicx}
\usepackage[all]{xy}

\begin{document}

\begin{titlepage}
\begin{center}
\textsc{\Large An Introduction to Exotic 4-manifolds}\\[0.5cm]
\textsc{\large Dean Bodenham} \\[4.5cm]
\begin{minipage}{10cm}
\begin{center} Abstract \end{center}
This article intends to provide an introduction to the construction of small exotic 4-manifolds. Some of the necessary background is covered. An exposition is given of J. Park's construction in \cite{P1} of an exotic $\mathbb{CP}^{2} \# 7 \overline{\mathbb{CP}^{2}}$. This article does not intend to present any new results. It was originally a Master's thesis, and its aim is merely to provide a leisurely introduction to exotic 4-manifolds that might be of use to interested graduate students.
\end{minipage}


\end{center}
\end{titlepage}

\pagebreak

\begin{flushleft}
\pagenumbering{roman}
\tableofcontents

\pagebreak

\pagenumbering{arabic}
\pagebreak

\newcommand{\mz}{\mathbb{Z}}
\newcommand{\cpt}{\mathbb{CP}^{2}}
\newcommand{\cpo}{\mathbb{CP}^{1}}
\newcommand{\cptbar}{\overline{\mathbb{CP}^{2}}}
\newcommand{\cpobar}{\overline{\mathbb{CP}^{1}}}
\newcommand{\st}{| \;}
\section{Introduction}

A manifold $X$ is \emph{exotic} if it is homeomorphic to another manifold $Y$, but is not diffeomorphic to it. Usually, $Y$ is a well-known manifold and then we say ``$X$ is an exotic $Y$''. \newline

A natural question to ask is ``How easy it is to find an exotic manifold?''. Or, more specifically, ``For a fixed manifold $X$, how many exotic $X$'s are there?'' To say there is an exotic $X$ is the same as saying that $X$ has more than one smooth structure.\newline 

Suppose $X$ is a closed, topological $n$-manifold. If $n \leq 3$, then $X$ has a unique smooth structure. If $n \geq 5$, then $X$ has at most finitely many smooth structures. This much is known. However, if $n=4$, then (as far as we know) $X$ must have infinitely many smooth structures or none at all. \newline

In other words, currently there are no closed smooth 4-manifolds known to have only finitely-many smooth structures. \newline

Perhaps the ``simplest'' smooth 4-manifolds, $S^4$, $S^2 \times S^2$ and $\cpt$, will provide us with an example of a 4-manifold with a unique smooth structure, but proving or disproving that this is the case has turned out to be quite difficult. \newline

One way mathematicians have tackled this problem in recent years is to try and find exotic $\cpt \# n \cptbar$s for as small an $n \in \mathbb{N}$ as possible. Exotic manifolds for the cases $n=9$ and $n=8$ were discovered by S. Donaldson in \cite{Do} and D. Kotschick in \cite{Ko3}, respectively, in the late 1980s. Then, for over 15 years, the next case of $n=7$ lay unsolved. \newline

In 2004, J. Park in \cite{P1} found a symplectic manifold that was homeomorphic, but not diffeomorphic, to $\cpt \# 7 \cptbar$, and since then infinite families of exotic $\cpt \# n \cptbar$s have been found for as low as $n=3$. \newline

\pagebreak

It is the purpose of this article to give an exposition of J. Park's construction. The principal sources are J. Park's paper \cite{P1}, and the papers of R. Fintushel and R. Stern \cite{FS1}, \cite{FS2} and \cite{FS3}. In order for the construction to be followed, some of the necessary background material is covered, most of which can be found in \cite{GS}. It is assumed the reader is familar with algebraic topology (see \cite{H1}) and knot theory (see \cite{Ro}). Kirby Calculus (see \cite{GS}) is used in section \ref{rbdsec}. \newline 

I would like to thank my supervisor, David Gay, for his continued patience and encouragement over the last two years. \newline

I would also like to thank Andr\'as Stipsicz for his excellent lectures and his patience in answering my many questions. \newline

\subsection*{A Note About This Article}
This article does not intend to present any new results, but merely provide some background and give an introduction to exotic 4-manifolds. It is an updated version of a Master's thesis that was submitted to the University of Cape Town in August 2008. It was suggested that it might be of use to other graduate students starting out in the field and so, after some encouragement, it was posted on the arXiv. \newline

Comments\footnote{I am not sure why, but downloading the postscript file and then turning it into a pdf gives diagrams of a slightly better quality than in the downloadable pdf.}, suggestions and corrections are welcomed and can be emailed to deanab17@yahoo.com. Any errors $-$ typographical, grammatical or mathematical $-$ are entirely my own. \newline 

Dean Bodenham \linebreak
December 2008\newline 

\pagebreak


\section{A Quick Review of Manifolds}

The definitions of a manifold, a differentiable structure, a homeomorphism and a diffeomorphism can be found in \cite{GP}, \cite{GS} and other books, but we review them here. For other definitions such as embeddings, isotopies, orientations, tangent spaces, the reader is referred to \cite{GP}. A useful introduction to differential geometry is \cite{I1}. For knot theory, see \cite{Ro}. For algebraic topology, see \cite{H1}. For fibre bundles, see \cite{St}. For characterstic classes, see \cite{MiSt} or \cite{GS}. We mainly follow the definitions presented in \cite{GS}. 

\newtheorem{defh}{Definition}[section]
\begin{defh}
\upshape
A \emph{homeomorphism} is a bijective map $\phi:X \longrightarrow Y$ between two topological spaces $X$ and $Y$ such that both $\phi$ and $\phi^{-1}$ are continuous.
\end{defh}

\newtheorem{defrp}[defh]{Definition}
\begin{defrp}
\upshape
We define $\mathbb{R}^{n}_{+}$ to be the upper half space of $\mathbb{R}^{n}$, i.e. $\mathbb{R}^{n}_{+} = \{(x_1,x_2,\dots, x_{n}) \in \mathbb{R}^{n}| \; x_{n} \geq 0 \}$.
\end{defrp}

\newtheorem{defm}[defh]{Definition}
\begin{defm}
\upshape
An \emph{$n$-dimensional topological manifold} is a separable Hausdorff topological space $X$, such that for every point $p \in X$ there is an open neighbourhood $U$ of $p$ that is homeomorphic to an open subet of $\mathbb{R}^{n}_{+}$.
\end{defm}

\newtheorem{defmr1}[defh]{Remark}
\begin{defmr1}
\upshape
We usually abbreviate ``$n$-dimensional manifold'' simply to ``$n$-manifold''.
\end{defmr1}

\newtheorem{egm0}[defh]{Example}
\begin{egm0}
\upshape
$\mathbb{R}^{n}$ and $\mathbb{R}^{n}_{+}$ are trivial examples of $n$-manifolds. Note that they are not compact.
\end{egm0}

\newtheorem{defm2}[defh]{Definition}
\begin{defm2}
\upshape
Let $X$ be a topological $n$-manifold. \newline
A pair $(U, \phi)$, where $U$ is an open subset of $X$ and $\phi: U \longrightarrow \mathbb{R}^{n}_{+}$ is a homeomorphism of $U$ onto an open subset of $\mathbb{R}^{n}_{+}$, is called a \emph{chart}. \newline

A collection of charts $\{ (U_{\alpha}, \phi_{\alpha}) | \; \alpha \in A \}$ is called an \emph{atlas} if it is a cover of $X$, i.e. $\cup_{\alpha \in A} U_{\alpha} = X$. \newline

The map $\phi_{\beta} \circ \phi_{\alpha}^{-1}$ from the open subset $\phi_{\alpha}(U_{\alpha} \cap U_{\beta}) \subset \mathbb{R}^{n}_{+}$ to the open subset $\phi_{\beta}(U_{\alpha} \cap U_{\beta}) \subset \mathbb{R}^{n}_{+}$ is called the \emph{transition function} between the charts $(U_{\alpha}, \phi_{\alpha})$ and $(U_{\beta}, \phi_{\beta})$. \newline

A topological manifold $X$ with an atlas $\{(U_{\alpha}, \phi_{\alpha})| \; \alpha \in A \}$ is called a \emph{$C^{r}$-manifold} $(r = 1,2, \dots \infty)$ if the transition functions are $C^{r}$-maps. In the case $r=\infty$, $X$ is called a \emph{smooth manifold}.
\end{defm2}

\newtheorem{defm2r}[defh]{Remarks}
\begin{defm2r}
\upshape
\begin{itemize}
\item[(i)] Often a chart is called a \emph{coordinate chart}.
\item[(ii)] If $X$ is a $C^{r}$-manifold ($r > 0$), $X$ is often simply called a \emph{differentiable manifold}. 
\item[(iii)] An atlas on a manifold $X$ with $C^{r}$ transition functions ($r > 0$) is called a \emph{differentiable structure on $X$}.
\end{itemize}
\end{defm2r}

\newtheorem{defm3}[defh]{Definition}
\begin{defm3}
\upshape
Let $X$ be a topological $n$-manifold. The points of $X$ corresponding to the points in $\{(x_1, x_2, \dots, x_{n}) \in \mathbb{R}^{n}_{+} | \; x_{n}=0 \} \cong \mathbb{R}^{n-1}$ form an $(n-1)$-dimensional submanifold of $X$, denoted $\partial X$ and called the \emph{boundary of $X$}.
\end{defm3}

\newtheorem{defm3r}[defh]{Remark}
\begin{defm3r}
\upshape
If we required the homeomorphisms $\phi_{\alpha}$ of the charts $(U_{\alpha}, \phi_{\alpha})$ to be maps into $\mathbb{R}^{n}$ as opposed to $\mathbb{R}^{n}_{+}$, $X$ would have empty boundary, i.e. $\partial X = \emptyset$.
\end{defm3r}

\newtheorem{defm33}[defh]{Definition}
\begin{defm33}
\upshape
We say that a manifold $X$ is \emph{closed} if it is compact and $\partial X = \emptyset$.
\end{defm33}

\newtheorem{egm6}[defh]{Example}
\begin{egm6}
\upshape
The $n$-dimensional sphere $S^n = \{ \mathbf{x} \in \mathbb{R}^{n+1} | \parallel \mathbf{x} \parallel = 1 \}$ is a closed $n$-manifold.
\end{egm6}

\newtheorem{egm7}[defh]{Example}
\begin{egm7}
\upshape
The $n$-dimensional disk $D^n = \{ \mathbf{x} \in \mathbb{R}^{n} | \parallel \mathbf{x} \parallel \leq 1 \}$ is a compact $n$-manifold with boundary $\partial D^{n} = S^{n-1}$.
\end{egm7}

\newtheorem{defm4}[defh]{Definition}
\begin{defm4}
\upshape
(from \cite{I1}) Let $X$ and $Y$ be two $C^{r}$-manifolds, and suppose $X$ is $n$-dimensional and $Y$ is $m$-dimensional. The \emph{local representative} of a map $f: X \longrightarrow Y$ with respect to the charts $(U, \phi)$ and $(V, \psi)$ on $X$ and $Y$, respectively, is the map 
\begin{equation*}
\psi \circ f \circ \phi^{-1}: \phi(U) \subset \mathbb{R}^{n}_{+} \longrightarrow \mathbb{R}^{m}_{+}
\end{equation*}

A map $f: X \longrightarrow Y$ is a \emph{$C^{r}$-map} between two $C^{r}$-manifolds $X$ and $Y$ if the local representatives of $f$ are $C^{r}$ with respect to every chart of the atlases of $X$ and $Y$.
\end{defm4}

\newtheorem{defm5}[defh]{Definition}
\begin{defm5}
\upshape
Let $X$ and $Y$ be two $C^{r}$-manifolds. A homeomorphism $f: X \longrightarrow Y$ is called a \emph{$C^{r}$-diffeomorphism} if both $f$ and $f^{-1}$ are $C^{r}$-maps.
\end{defm5}

\newtheorem{defm5r}[defh]{Remark}
\begin{defm5r}
\upshape
In the case $r=\infty$ we usually call such a map a \emph{diffeomorphism}.
\end{defm5r}


\newtheorem{defm7}[defh]{Definition}
\begin{defm7}
\upshape
Let $W$ be an open neighbourhood of $\mathbb{C}$. A function $f: W \longrightarrow \mathbb{C}$ is called \emph{holomorphic} if it is complex-differentiable at every point in $W$.
\end{defm7}

\newtheorem{defm7r}[defh]{Remark}
\begin{defm7r}
\upshape
Recall from complex analysis that a holomorphic function is also \emph{analytic}, i.e. equals its Taylor series in a neighbourhood of each point of its domain (and is therefore also a smooth function). 
\end{defm7r}

\newtheorem{defm8}[defh]{Definition}
\begin{defm8}
\upshape
An atlas $\{(U_{\alpha}, \phi_{\alpha}) | \; \alpha \in A \}$ on a (real) $2n$-dimensional manifold $X$ is called a \emph{complex structure} if each $\phi_{\alpha}$ is a homeomorphism between $U_{\alpha}$ and an open subset of $\mathbb{C}^{n}$ (identified with $\mathbb{R}^{2n}$), and the transition functions $\phi_{\beta} \circ \phi_{\alpha}^{-1}$ are holomorphic.
\end{defm8}

\newtheorem{defm8r}[defh]{Remark}
\begin{defm8r} \label{compcanor}
\upshape
Complex manifolds are canonically oriented. The following argument comes from \cite{GS}: Firstly, $\mathbb{C}$ is oriented as a real vector space by the ordered basis $(1,i)$. Secondly, the connected group $GL(n; \mathbb{C})$ lies in $GL^{+}(2n; \mathbb{R})$, and so by choosing a complex isomorphism with $\mathbb{C}^{n}$, any $n$-dimensional complex vector space is canonically oriented.
\end{defm8r}

The next theorem from \cite{Mu} (quoted from \cite{GS}) shows that every $C^{r}$-manifold essentially has a smooth structure (for $r > 0$). 

\newtheorem{muthm}[defh]{Theorem}
\begin{muthm}
Suppose that $X$ is a $C^{r}$-manifold and $1 \leq r \leq s$ (including $s=\infty$). Then there is a $C^{s}$-atlas of $X$ for which the induced $C^{r}$-structure is isotopic to the original $C^{r}$-structure on $X$. Moreover, this $C^{s}$-structure is unique up to isotopy (through $C^{r}$-diffeomorphisms); consequently the $C^{r}$-manifold $X$ admits a unique induced $C^{s}$-structure for every $s \geq r$.
\end{muthm}

\newtheorem{muthmr}[defh]{Remark}
\begin{muthmr}
\upshape
Therefore, we only need to focus on classifying the topological manifolds and the smooth manifolds. While every smooth manifold is a topological manifold, we shall see later that there are some topological manifolds that do not admit any smooth structures.
\end{muthmr}

\pagebreak


\section{A Brief Account of Fibre Bundles}

What follows below is the ``provisional definition'' given in \cite{St}, which will be sufficient for our purposes.
 
\newtheorem{fbd}{Definition}[section]
\begin{fbd}
\upshape
A fibre bundle $\mathfrak{B} = (E,B,p,F)$ is a quadruple consisting of
\begin{itemize}
\item[(1)] a topological space $E$ called the \emph{bundle space}
\item[(2)] a topological space $B$ called the \emph{base space}
\item[(3)] a continuous map $p: E \longrightarrow B$ of $E$ onto $B$ called the \emph{projection map}
\item[(4)] a space $F$ called the \emph{fibre}
\end{itemize}

This quadruple must satisfy two conditions. For each $x \in B$,
\begin{itemize}
\item[(i)] the set $F_{x}$ defined by $F_{x} = p^{-1}(x)$, called the \emph{fibre over the point $x$}, must be homeomorphic to the fibre $F$
\item[(ii)] there is a neighbourhood $U$ of $x$ and a homeomorphism $\phi : U \times F \longrightarrow p^{-1}(U)$ such that for all $u \in U$ and all $f \in F$
\begin{equation*}
p \circ \phi (u,f) = u
\end{equation*}
\end{itemize}
\end{fbd}

\newtheorem{fbdr}[fbd]{Remarks}
\begin{fbdr}
\upshape
\begin{itemize}
\item[(i)] Often, we just call the bundle space the \emph{bundle}, the base space the \emph{base} and the projection map the \emph{projection}. 
\item[(ii)] Condition (ii) above is called \emph{local triviality}, and is equivalent to saying that the following diagram commutes
\begin{displaymath}
\xymatrix{ 
U \times F \ar[dr]_{\pi_{1}} \ar[rr]^{\phi} & & p^{-1}(U) \ar[dl]^{p} \\
 & U & }
\end{displaymath}
where $\pi_1$ is the projection of $U \times F$ onto the first factor $U$.
\item[(iii)] It should be noted that a fibre bundle is actually a quintuple $\mathfrak{B} = (E,B,p,F,G)$, where $G$ is a topological tranformation group, called the \emph{structure group}, satisfying certain conditions. However, we shall not really need this additional structure, and shall just consider a fibre bundle to be as described in the definition above.
\end{itemize}
\end{fbdr}

\newtheorem{fbdseq}[fbd]{Remark}
\begin{fbdseq} \label{fbrseqrem}
\upshape
In \cite{St} and \cite{H1} it is proved that a fibre bundle $\mathfrak{B} = (E,B,p,F)$ gives rise to the following long exact sequence of homotopy groups:
\begin{equation*}
\dots \longrightarrow \pi_{n}(F) \longrightarrow \pi_{n}(E) \longrightarrow \pi_{n}(B) \longrightarrow \pi_{n-1}(F) \longrightarrow \dots \longrightarrow \pi_{0}(E) \longrightarrow 0
\end{equation*}
\end{fbdseq}

\newtheorem{fbtr}[fbd]{Definiton}
\begin{fbtr}
\upshape
We define a fibre bundle $\mathfrak{B} = (E,B,p,F)$ to be a trivial bundle if there is a homeomorphism $h: E \longrightarrow B \times F$ which commutes with the projection maps, as above.
\end{fbtr}

With this definition in mind, we have the following important theorem:

\newtheorem{fbtrthm}[fbd]{Theorem}
\begin{fbtrthm}
If the base space $B$ of a fibre bundle $\mathfrak{B} = (E,B,p,F)$ is contractible, then $\mathfrak{B}$ is a trivial bundle.
\end{fbtrthm}

There is a proof of this theorem in \cite{St}, but I prefer the proof in \cite{O1}.

\newtheorem{fbsecdef}[fbd]{Definition}
\begin{fbsecdef}
\upshape
Let $\mathfrak{B} = (E,B,p,F)$ be a fibre bundle. We call the (continuous) map $s:B \longrightarrow E$ a \emph{section} of the fibre bundle if for all $b \in B$ we have $p \circ s (b) = b$.
\end{fbsecdef}


\pagebreak


\section{The Intersection Form}

In this section, let $X$ be a compact, oriented, topological 4-manifold. We follow the notations and definitions presented in \cite{GS}. \newline

Since $X$ is oriented, it admits a fundamental class $[X] \in H_{4}(X, \partial X, \mathbb{Z})$. See \cite{H1} for details.
We now define the intersection form of a 4-manifold $X$.

\newtheorem{dif}{Definition}[section]
\begin{dif}
\upshape
The symmetric bilinear form 
\begin{equation*}
Q_{X}: H^{2}(X, \partial X; \mz) \times H^{2}(X, \partial X; \mz) \longrightarrow \mz
\end{equation*}
is defined by $Q_{X}(a,b) = <a \cup b, [X]>$.
\end{dif}

\newtheorem{difr1}[dif]{Remark}
\begin{difr1}
\upshape
In the definition above $a,b \in H^{2}(X, \partial X; \mz)$ and `$a \cup b$' denotes the cup product between $a$ and $b$, $[X] \in H_{4}(X, \partial X, \mz)$ is the fundamental class of the manifold $X$, and $< . , . >: H^{4}(X, \partial X, \mz) \times H_{4}(X, \partial X, \mz) \longrightarrow \mz $ is the bilinear form where the cohomology class is evaluated on the homology class. 
\end{difr1}

\newtheorem{difr2}[dif]{Remark}
\begin{difr2}
\upshape
\begin{itemize}
\item[(i)] We often denote $Q_{X}(a, b)$ by $a \cdot b$.
\item[(ii)] Since by Poincar\'e duality $H_2(X; \mz) \cong H^2(X, \partial X; \mz)$, $Q_{X}$ is also defined on $H_2(X; \mz) \times H_2(X; \mz)$.
\item[(iii)] By the definition of the intersection form, we have $Q_{\overline{X}} = - Q_{X}$, where $\overline{X}$ is the manifold $X$ with the opposite orientation. \label{difr25}
\end{itemize}
\end{difr2}

If $a$ or $b$ is a torsion element of $H_2(X; \mz)$ then $Q_{X}(a,b)=0$. Therefore, we could consider intersection forms as just being defined on $H_2(X; \mz)$/torsion, which is a finitely-generated free abelian group. \newline

 By choosing a basis $\{b_1,b_2, \dots ,b_{n} \}$ of $H_2(X; \mz)$/torsion, we can represent $Q_{X}$ by a matrix $M$. Since $M$ depends on our choice of basis, using another basis $\{\tilde{b_1},\tilde{b_2}, \dots , \tilde{b_{n}} \}$  could result in $Q_{X}$ being represented by another matrix $\tilde{M}$. However, if $B$ is the basis tranformation between the bases $\{b_{i} \}$ and $\{ \tilde{b_{i}} \}$, and $B^{t}$ is the transpose matrix of $B$, then $\tilde{M} = B M B^{t}$, and $B$ is such that $\mathrm{det}(B) = \pm 1$. Then
\begin{align*}
\mathrm{det} (\tilde{M}) &= \mathrm{det} (B) \mathrm{det} (M) \mathrm{det}(B^{t}) = \mathrm{det} (B)  \mathrm{det}(B^{t}) \mathrm{det} (M)\\ 
&= \mathrm{det} (B)  \mathrm{det}(B) \mathrm{det}(M) = (\mathrm{det}(B))^2 \mathrm{det} (M) = 1 \cdot \mathrm{det} (M)\\ 
&= \mathrm{det} (M)
\end{align*}

This shows that $\mathrm{det}(M)$ is independent of the basis we choose, and sometimes we denote this by $\mathrm{det}(Q_{X})$. We shall also usually identify the intersection form $Q_{X}$ with the matrix representing it.

\newtheorem{def2}[dif]{Definition}
\begin{def2}
\upshape
If $M$ and $\tilde{M}$ are two $n \times n$ matrices over $\mz$ and there is a matrix $B$ (also over $\mz$) such that 
\begin{equation*}
\tilde{M} = B M B^{t}
\end{equation*}
then we say that $M$ and $\tilde{M}$ are equivalent.
\end{def2}

It is easy to check that this definition of \emph{equivalent} is an equivalence relation. \newline

It is natural to ask where the name \emph{intersection form} originates. If $X$ is a smooth manifold, then $Q_{X}(a,b)$ can be interpreted as the (signed) number of intersections of two submanifolds of $X$. In order to make this statment clear, we shall need a little background (which we take from \cite{GS}). \newline

In what follows, $X$ is a closed, oriented, smooth 4-manifold. Similar results hold for cases when $X$ has boundary, is non-compact or is non-orientable, but we shall not consider these cases here.

\newtheorem{def3}[dif]{Definition}
\begin{def3}
\upshape
Let $X^{n}$ be a smooth $n$-dimensional manifold. We say that a class $\alpha \in H_2(X^{n}; \mz)$ is represented by a closed, oriented surface $\Sigma_{\alpha}$ if there is an embedding $i: \Sigma_{\alpha} \hookrightarrow X^{n}$ such that $i_{*}([\Sigma_{\alpha}]) = \alpha$, where $[\Sigma_{\alpha}] \in H_2(\Sigma_{\alpha}; \mz)$ is the fundamental class of $\Sigma_{\alpha}$.
\end{def3}

With this definition in mind, we have the following proposition:

\newtheorem{pro1}[dif]{Proposition}
\begin{pro1}
Let $X$ be a closed, oriented, smooth 4-manifold. Then every element of $H_{2}(X; \mz)$ can be represented by an embedded surface.
\end{pro1}

The proof in \cite{GS} uses results that are beyond the scope of this article.

\pagebreak
\newtheorem{pro1r1}[dif]{Remark}
\begin{pro1r1}
\upshape
If $X$ is simply connected, by the Hurewicz Theorem (see \cite{H1}) $\pi_2(X) \cong H_2(X; \mz)$, which implies that every $\alpha \in H_2(X; \mz)$ can be represented by an immersed sphere. Note the difference: the embedded \emph{surface} above need not have been a \emph{sphere}; for example, it could have been a torus. \newline

Further note that although this immersion is not an embedding in general, one can assume that an immersion $S^{2} \longrightarrow X^4$ intersects itself only in transverse double points (see \cite{GP}).
\end{pro1r1}

Again, suppose that $X$ is a closed, oriented, smooth 4-manifold. Let $a, b \in H^2(X; \mz)$ and let their Poincar\'e duals be $\alpha = PD(a), \beta = PD(b)$, respectively. For the following, see \cite{GP}. \newline

Let $\Sigma_{\alpha}$ and $\Sigma_{\beta}$ be the surface representatives of $\alpha$ and $\beta$ (and therefore of $a$ and $b$), respectively, and suppose that the surfaces $\Sigma_{\alpha}$ and $\Sigma_{\beta}$ have been chosen generically, so that all their intersections are transverse. \newline

If $p \in \Sigma_{\alpha} \cap \Sigma_{\beta}$, the tangent spaces at the point $p$, denoted $T_{p}\Sigma_{\alpha}$ and $T_{p}\Sigma_{\beta}$ are orthogonal to each other (since the surfaces intersect transversely at $p$). \newline

If we concatenate the basis $\{x_1, x_2 \}$ of $T_{p}\Sigma_{\alpha}$ and the basis $\{y_1, y_2 \}$ of $T_{p}\Sigma_{\beta}$, we get a basis $\{x_1, x_2, y_1, y_2 \}$ for $T_{p}X$. \newline

If this basis $\{x_1, x_2, y_1, y_2 \}$ is positive (i.e. defines a positive orientation on $X$), we define the sign of the intersection at $p$ to be positive, otherwise it is a negative intersection. \newline

Note that the sign will not depend on the order of $\{\alpha, \beta \}$, but will depend on the orientations of the embedded surfaces $\Sigma_{\alpha}$ and $\Sigma_{\beta}$. \newline

This leads to the geometric interpretation of $Q_{X}$ (for the proof, see \cite{GS}):

\newtheorem{pro2}[dif]{Proposition}
\begin{pro2}
Let $X$ be a closed, oriented, smooth 4-manifold. Let $a, b \in H^{2}(X; \mz)$. Let $\alpha, \beta \in H_2(X;\mz)$ be their Poincar\'e duals, respectively, and let $\Sigma_{\alpha}$ and $\Sigma_{\beta}$ be their surface representatives, respectively. Let $Q_{X}$ be the intersection form of $X$. \newline

Then $Q_{X}(a,b)$ is the number of points in $\Sigma_{\alpha} \cap \Sigma_{\beta}$, counted with sign.
\end{pro2}

\newtheorem{pro2r1}[dif]{Remark}
\begin{pro2r1}
\upshape
For example, if there are three intersection points in $\Sigma_{\alpha} \cap \Sigma_{\beta}$, 2 negative and 1 positive, then $Q_{X}(a,b) = -1$.
\end{pro2r1}

\newtheorem{pro2r2}[dif]{Remarks}
\begin{pro2r2} \label{comrem1}
\upshape
A complex structure on a manifold $X$ defines an orientation on $X$, and so any complex submanifolds of $X$ are canonically oriented (see Remark \ref{compcanor}). \newline
In particular, if $S$ is a complex surface, and $C_1$ and $C_2$ are complex curves in $S$ that intersect each other transversely, then $Q_{S}(C_1, C_2) \geq 0$. In other words, the transverse intersection of complex submanifolds is always positive. This important fact will be used later. \newline
However, it is worth noting that the \emph{self-intersection} $Q_X(C,C)$ of a complex curve $C$ can be negative, as we shall see later.
\end{pro2r2}

The intersection form $Q_{X}$ of a closed, oriented, topological 4-manifold $X$ is a symmetric, bilinear form on a finitely-generated free abelian group. However, there is another property that the matrices representing $Q_{X}$ have, one that is not immediately apparent.

\newtheorem{def4}[dif]{Definition}
\begin{def4}
A matrix $Q$ is called unimodular if $\mathrm{det}(Q) = \pm1$.
\end{def4}

This is equivalent to saying that the matrix $Q$ is invertible over $\mz$. We can now state the following proposition:

\newtheorem{pro3}[dif]{Proposition}
\begin{pro3}
The intersection form $Q_{X}$ of a closed, oriented, topological 4-manifold $X$ is unimodular.
\end{pro3}

The proof of this proposition can be found in both \cite{GS} and \cite{Sc}. \newline

Now, let us forget about the 4-manifolds for a while, and just look at properties of symmetric, bilinear, unimodular forms defined on a finitely-generated free abelian group.

\pagebreak

\section{Classification of Integral Forms}

Actually, the title of this section is a bit of an abbreviation. We shall just be considering symmetric, bilinear, unimodular forms. Again, we follow \cite{GS}. Another good reference is \cite{MH}. \newline

Let $Q: A \times A \longrightarrow \mz$ be a symmetric, bilinear, unimodular form defined on a finitely-generated free abelian group $A$. We define the following invariants of $Q$: \emph{rank}, \emph{signature} and \emph{parity}. \newline

The \emph{rank} of $Q$ is the dimension of $A$ and is denoted $rk(Q)$. \newline

The \emph{signature} of $Q$ is defined as follows: consider $Q$ as an $n \times n$ matrix with entries in $\mz$, and diagonalize it over $\mathbb{R}$. Denote the number of positive eigenvalues on the diagonal by $b_{2}^{+}(Q)$ and the number of negative eigenvalues on the diagonal by $b_{2}^{-}(Q)$. We finally define the signature of $Q$ as $\sigma(Q) = b_{2}^{+}(Q) - b_{2}^{-}(Q)$.

\newtheorem{r21}{Remark}[section]
\begin{r21}
\upshape
We define $b_{2}(Q) = b_{2}^{+}(Q) + b_{2}^{-}(Q)$, and it is called the \emph{second Betti number}. Clearly, $b_{2}(Q) = rk(Q)$. Sometimes we shall write $b_{2}(Q)$ simply as $b_{2}$ when there is no chance of ambiguity.
\end{r21}

We define the \emph{parity} of $Q$ to be either \emph{odd} or \emph{even}. If for all $a \in A$ $Q(a,a) \equiv 0$ (mod 2), we say that $Q$ is even. Otherwise, we say that $Q$ is odd.

\newtheorem{r22}[r21]{Remark}
\begin{r22}
\upshape
Note that if there is just one element $a' \in A$ such that $Q(a', a') \equiv 1$ (mod 2), it is enough to make $Q$ odd.
\end{r22}

It is worthwhile checking that rank, signature and parity are indeed invariants of a symmetric, bilinear, unimodular form. When we say they are invariants, we mean that two matrices representing the same form should have the same rank, signature and parity. This is thre same as saying that if $X$ and $Y$ are two equivalent matrices with entries in $\mz$, i.e. there is a basis tranformation tranformation matrix $B$ such that $X = B Y B^{t}$, then $rk(X) = rk(Y)$, $\sigma(X) = \sigma(Y)$ and the parity of $X$ is the same as the parity of $Y$. \newline

Suppose $X$ and $Y$ are equivalent matrices representing the same symmetric, bilinear, unimodular form $Q$. The fact that rank and signature are invariants is due to the following theorem from linear algebra (Theorem $6.z_{3}$ in \cite{He}):

\newtheorem{hers1}[r21]{Theorem}
\begin{hers1}
Given the real symmetric matrix $A$ there is an invertible matrix $T$ such that 

\begin{displaymath}
TAT^{t} = 
\left( \begin{array} {ccc}
I_{r} &  & \\
& -I_{s} & \\
&  &  0_{t} 
\end{array} \right)
\end{displaymath}
where $I_{r}$ and $I_{s}$ are respectively the $r \times r$ and $s \times s$ unit matrices and where $0_{t}$ is the $t \times t$ $0$-matrix. The integers $r+s$, which is the rank of $A$, and $r-s$, which is the signature of $A$, characterize the congruence class of $A$. That is, two real symmetric matrices are congruent if and only if they have the same rank and signature.
\end{hers1}

\newtheorem{hersr1}[r21]{Remark}
\begin{hersr1}
\upshape
Herstein's notion of two matrices $A$ and $B$ being \emph{congruent} over $\mathbb{R}$ means there is a non-singular real matrix $T$ such that $B = TAT^{t}$, and the theorem above then proves that if $A$ and $B$ are congruent, they have the same rank and signature. If two matrices are \emph{equivalent} over $\mz$ (as defined in the previous section), then they are clearly \emph{congruent} over $\mathbb{R}$ (as Herstein defines it), and so two equivalent matrices have the same rank and signature. Note that since we also consider our matrices in this section to be unimodular, and therefore invertible, so $t = \mathrm{dim}(0_{t}) = 0$ above.
\end{hersr1}

Finally, a short lemma below proves that parity is also an invariant. We shall use the shorthand `$X \equiv_{2} Y$' to denote `$X \equiv Y$ (mod 2)'.

\newtheorem{parlem1}[r21]{Lemma}
\begin{parlem1} \label{L25}
Let $Q$ be a symmetric, bilinear, unimodular form over $\mz$. Then $Q$ is even if and only if $Q(\alpha_{i}, \alpha_{i}) \equiv_{2} 0$ for each $i = 1, 2, \dots, n$, where $\{ \alpha_{1}, \alpha_{2}, \dots, \alpha_{n} \}$ is a basis for $A$.
\end{parlem1}

Proof: \newline
($\Rightarrow$): By the definition of $Q$ being even. \newline
($\Leftarrow$): Let $A_{0} \subset A$ be the subset of $A$ such that for $a \in A_{0}$, $Q(a,a) \equiv_{2} 0$. If $a, b \in A_{0}$, then 

\begin{align*}
Q(a+b, a+b) &= Q(a, a) + Q(a, b) + Q(b, a) + Q(b,b) \\
&= Q(a,a) + 2Q(a,b) + Q(b,b) \\
\Rightarrow Q(a+b, a+b) &\equiv_{2} 0
\end{align*}

since $Q(a,a) \equiv_2 0$, $Q(b,b) \equiv_2 0$ and clearly $2Q(a,b) \equiv_2 0$. A similar argument proves that $Q(a-b, a-b) \equiv_{2} 0$. Therefore, $a, b \in A_{0}$ implies that $a + b \in A_{0}$ and $a - b \in A_{0}$. \newline

Since the basis $\{ \alpha_{1}, \alpha_{2}, \dots, \alpha_{n} \}$ is contained in $A_{0}$, by applying this argument repeatedly we have for any $\lambda_{1}, \lambda_{2}, \dots, \lambda_{n} \in \mz$ that 

\begin{equation*}
Q(\lambda_{1}\alpha_{1} + \lambda_{2}\alpha_{2} + \dots + \lambda_{n}\alpha_{n}, \lambda_{1}\alpha_{1} + \lambda_{2}\alpha_{2} + \dots + \lambda_{n}\alpha_{n}) \equiv_2 0
\end{equation*}

This implies that $A \subset A_{0}$, which implies $A = A_{0}$, which proves that $Q$ is even. $\Box$

%
%
%
%
%
%
\newtheorem{parlem2r1}[r21]{Remark}
\begin{parlem2r1}
\upshape
The lemma above proves that if $Q$ is even in one basis of $A$, then it is even in any other basis of $A$. Since two equivalent matrices are just the same symmetric, bilinear, unimodular form represented in two different bases, if one matrix is even then so is the other.
\end{parlem2r1}

We define a further notion known as the \emph{definiteness} of the intersection form as follows:
\begin{itemize}
\item[(i)] If $rk(Q) = \sigma(Q)$, $Q$ is called \emph{positive-definite}. 
\item[(ii)] If $rk(Q) = -\sigma(Q)$, $Q$ is called \emph{negative-definite}. 
\item[(iii)] Otherwise, $Q$ is called \emph{indefinite}.
\end{itemize}

If $Q$ is not indefinite, it could simply be called \emph{definite}. \newline


Due to a theorem of Serre (quoted from \cite{Sc}, but the original source is \cite{Se}), these three invariants rank, signature and parity are enough to classify all indefinite intersection forms.

\newtheorem{serre}[r21]{Theorem}
\begin{serre}
Let $Q_1$ and $Q_2$ be two indefinite, symmetric, bilinear, unimodular forms. If $Q_1$ and $Q_2$ have the same rank, signature and parity, then they are equivalent.
\end{serre}

This classification of indefinite forms will be very useful later. However, there is no `nice' classification of definite forms; in fact there are positive-definite forms that are not equivalent, even though they have the same rank, signature and parity (refer to page 14, \cite{GS}). \newline

Before we return to the world of 4-manifolds, there is one more definition.

\newtheorem{cd}[r21]{Definition}
\begin{cd}
An element $x \in A$ is called a characteristic element if $Q(x,\alpha) \equiv Q(\alpha,\alpha)$ (mod 2) for all $\alpha \in A$.
\end{cd}

This leads to an interesting result (for the proof, see \cite{GS}):

\newtheorem{l1220}[r21]{Lemma}
\begin{l1220} \label{L210}
If $x \in A$ is characteristic, then $Q(x,x) \equiv \sigma(Q)$ (mod 8). In particular, if $Q$ is even, then the signature $\sigma(Q)$ is divisible by 8.
\end{l1220}

\pagebreak

\section{The Classification of Topological 4-manifolds} 
\label{class4msec}
Let $X$ be a simply-connected, closed, oriented, 4-manifold.\newline

Recall that $\pi_1(X)=0$ if $X$ is simply-connected. Then since $H_1(X)$ is just the abelianization of $\pi_1(X)$, we also have $H_1(X; \mz) = 0$. We then have $H^1(X; \mz) \cong \mathrm{Hom}(H_1(X;\mz); \mz) = 0$ and by Poincar\'e duality, we then have $H_{3}(X;\mz) = 0$ and $H^3(X;\mz) =0$. This also implies that $H_2(X;\mz) \cong H^{2}(X; \mz)$ has no torsion. Therefore, the intersection form $Q_{X}$ contains all the homological information of $X$. Whitehead first showed that $Q_{X}$ classifies topological 4-manifolds up to homotopy:

\newtheorem{w1}{Theorem}[section]
\begin{w1}
The simply-connected, closed, topological 4-manifolds $X_1$ and $X_2$ are homotopy equivalent if and only if $Q_{X_1} \cong Q_{X_2}$.
\end{w1}

Then, in 1982 M. Freedman proved the following theorem in \cite{F} which shows that $Q_{X}$ actually classifies $X$ up to homeomorphism:

\newtheorem{f1}[w1]{Theorem}
\begin{f1}
For every symmetric, bilinear, unimodular form $Q$ there exists a simply-connected, closed, topological 4-manifold $X$ such that $Q_{X} \cong Q$. \newline

Furthermore, if $Q$ is even, then this manifold is unique up to homeomorphsim. If $Q$ is odd, there are exactly two different homeomorphism types of manifolds with intersection form $Q$, and at most one of these homeomorphism types carries a smooth structure.
\end{f1}

If we restrict our attention to smooth manifolds, this leads to an important corollary:

\newtheorem{f2}[w1]{Corollary}
\begin{f2} \label{frcor1}
If $X_1$ and $X_2$ are smooth, simply-connected 4-manifolds with equivalent intersection forms, then $X_1$ and $X_2$ must be homeomorphic.
\end{f2}

A special case of Freedman's theorem is the topological 4-dimensional Poincar\'e Conjecture:

\newtheorem{f3}[w1]{Corollary}
\begin{f3}
If $X$ is a topological 4-manifold homotopy equivalent to $S^4$, then $X$ is homeomorphic to $S^4$.
\end{f3}

From now on, we shall usually write the invariants rank and signature of the intersection form $Q_{X}$ as $b_2(X) = b_2^{+}(X) + b_2^{-}(X)$ and $\sigma(X) = b_2^{+}(X) - b_2^{-}(X)$, respectively. \newline

Below is an interesting result that Rohlin proved in \cite{R2}:

\newtheorem{r2t1}[w1]{Theorem}
\begin{r2t1} \label{T35}
Let $X$ be a simply-connected, closed, oriented, smooth 4-manifold. If $Q_X$ is even, then the signature $\sigma(X)$ is divisible by 16.
\end{r2t1}

We shall now look at a few examples of 4-manifolds and their intersection forms.

\pagebreak


\section{Examples}

\newtheorem{eg1}{Example}[section]
\begin{eg1}
\upshape
The simplest example is $S^4 = \{\mathbf{x} \in \mathbb{R}^4 | \; \parallel \mathbf{x} \parallel = 1  \}$. Since $H_2(S^4; \mz) = 0$, the intersection form $Q_{S^{4}} = <.> $, where $<n>$ denotes the $1 \times 1$ matrix with the single entry $n \in \mz$, and $<.>$ denotes the ``empty'' intersection form (there are no homology classes to ``intersect'' each other in the case of $S^4$; note this is not standard notation). 
\end{eg1}

\newtheorem{eg2}[eg1]{Example}
\begin{eg2}
\upshape
The next examples are the complex projective spaces. We define $\mathbb{CP}^n = \{\mathbf{z} \in \mathbb{C}^{n+1} | \; \mathbf{z} \neq \mathbf{0} \} / \sim$, where $\mathbf{0} = (0, 0, \dots, 0)$ and the relation $\sim$ is defined as: \newline 
for all $\lambda \in \mathbb{C} \setminus \{0 \}$,  $(\lambda z_0, \lambda z_1, \dots, \lambda z_{n}) \sim (z_0, z_1, \dots, z_{n})$. More compactly, denoting $\mathbb{C} \setminus \{0 \}$ by $\mathbb{C}^{*}$, for all $\lambda \in \mathbb{C}^{*}$, $\lambda \mathbf{z} \sim \mathbf{z}$. 
\end{eg2}

\newtheorem{eg2r}[eg1]{Remark}
\begin{eg2r}
\upshape
Note that a point $P \in \mathbb{CP}^{n}$ is an equivalence class of points, so if $(z_0, z_1, \dots, z_{n}) \in P$, we usually denote $P$ by its \emph{homogeneous coordinates} $[z_0:z_1: \dots : z_{n}]$. For example, all the points $(0,\dots,0,1)$, $(0,\dots,0,2)$, etc. are in the equivalence class $[0:\dots:0:1]$. We call $\cpo$ the \emph{complex projective line} and $\cpt$ the \emph{complex projective plane}.
\end{eg2r}

One can similarly define the real projective spaces, and it is worthwhile looking at $\mathbb{RP}^1$ and $\mathbb{RP}^2$ in order to get a better idea of what these complex projective spaces actually are. \newline

We define $\mathbb{RP}^1 = \{(x,y) \in \mathbb{R}^2 | \; (x,y) \neq (0,0) \} / \sim$, where now $\sim$ is the relation that for all $\lambda \in \mathbb{R} \setminus \{0 \}$, $(\lambda x, \lambda y) \sim (x,y)$. Let us pick an element $(x,y) \in \mathbb{R}^2$ that lies on the unit circle. This element also defines a line through the origin in the direction $(x,y)$, and then $(\lambda x, \lambda y)$ for $\lambda \neq 0$ is just any other element on this line, except the origin, and we identify $(\lambda x, \lambda y)$ with $(x,y)$. \newline

So, we could picture $\mathbb{RP}^{1}$ as follows: we start with $\mathbb{R}^{2} \setminus \{(0,0)\}$, then quotienting out by the relation $\sim$ retracts $\mathbb{R}^{2} \setminus \{(0,0)\}$ onto the unit circle $S^1$, and then identifies antipodal points of the unit circle (since $(x,y) \sim (-x, -y)$). Using this approach for higher dimensions, we can consider $\mathbb{RP}^{n}$ as the unit sphere $S^{n}$ with antipodal points identified. \newline

\pagebreak
Although this `picture' of the real projective spaces doesn't quite extend to the complex case, since the scalars $\lambda$ are then complex, it can be used to show that $\mathbb{CP}^{n}$ is both compact and simply-connected, as the following lemma shows:

\newtheorem{cplem}[eg1]{Lemma}
\begin{cplem}
The complex projective spaces $\mathbb{CP}^{n}$ are compact and simply connected.
\end{cplem}

Proof: (from an exercise in \cite{GS})  \newline
Using the idea of the equivalent formulation of real projective spaces, we have the following equivalent definition of $\mathbb{CP}^{n}$:

\begin{equation*}
\mathbb{CP}^{n} = \{ \mathbf{x} \in S^{2n+1} \subset \mathbb{R}^{2n+2} \cong \mathbb{C}^{n+1} | \; \mathbf{x} \neq \mathbf{0} \} / \sim
\end{equation*}

where $\sim$ is defined as: for all $\lambda \in S^{1}$, $\lambda \mathbf{x} \sim \mathbf{x}$ (recall that the set $\{z \in \mathbb{C} | \; |z| = 1 \} \cong S^1$). \newline

This definition makes $S^{2n+1}$ into an $S^1$-fibration over $\mathbb{CP}^{n}$, and from Remark \ref{fbrseqrem} above we have the following long exact sequence of homotopy groups:
\begin{equation*}
\dots \rightarrow \pi_1(S^1) \rightarrow \pi_1(S^{2n+1}) \rightarrow \pi_1(\mathbb{CP}^{n}) \rightarrow \pi_0(S^1) \rightarrow \pi_0(S^{2n+1}) \rightarrow \dots
\end{equation*}

Since $S^{2n+1}$ is simply connected, $\pi_1(S^{2n+1})= 0$. Recall from \cite{H1} (page 346) that a path-connected space $X$ has $\pi_0(X) = 0$. Therefore, since $S^1$ and $S^{2n+1}$ are path-connected, we have $\pi_0(S^1) = 0$ and $\pi_0(S^{2n+1} = 0)$. So, a portion of our exact sequence becomes
\begin{equation*}
\dots \rightarrow \pi_1(S^1) \rightarrow 0 \rightarrow \pi_1(\mathbb{CP}^{n}) \rightarrow 0 \rightarrow 0 \rightarrow \dots
\end{equation*}

and so $\pi_1({\mathbb{CP}^{n}}) = 0$, and therefore $\mathbb{CP}^{n}$ is simply-connected. Since $S^{n+1}$ is compact, and the projection map $p: S^{n+1} \longrightarrow \mathbb{CP}^{n}$ is continuous, we have that $\mathbb{CP}^{n}$ is compact.

\newtheorem{lemrem1}[eg1]{Remark}
\begin{lemrem1} \label{cpnhomref}
\upshape
The homology groups of $\mathbb{CP}^{n}$ follow an interesting pattern: $H_{i}(\mathbb{CP}^{n}; \mz) \cong \mz$ if $i=2d$, where $d=0,1,\dots,n$, and $H_{i}(\mathbb{CP}^{n}; \mz) = 0$ otherwise. For a proof, see \cite{H1} or \cite{GS} Example 4.2.4.
\end{lemrem1}

Let us now focus our attention on $\cpt$, which is a 4-manifold that is important later, and let us calculate its intersection form (from an exercise in \cite{GS}). \newline

Let $h \in H_2(\cpt; \mz)$ be the fundamental class of the submanifold $H = \{ [x:y:z] \in \cpt | \; x = 0 \}$, and let $h' \in H_2(\cpt; \mz)$ be the fundamental class of the submanifold $H' = \{ [x:y:z] \in \cpt | \; y = 0 \}$. Clearly $H \cap H' = \{[0:0:1] \}$ is a transverse intersection. Since both submanifolds are complex, the intersection is also positive (see Remarks \ref{comrem1}). Therefore, $Q_{\cpt}(h,h') =1$. \newline

Claim: $h$ cannot be the multiple of any other class in $H_2(\cpt; \mz)$, and so it generates $H_2(\cpt; \mz) \cong \mz$.\newline

Proof of claim: Suppose $h$ is a multiple of a class $g$, i.e. $h = mg$ for $m \in \mz, |m| > 1$ (otherwise $g = \pm h$). Then 
\begin{align*}
Q_{\cpt}(h,h') &= Q_{\cpt}(mg,h') \\
&= m Q_{\cpt}(g,h') \\
&= mk
\end{align*}

where $k \in \mz$. Therefore $mk = 1$, which is clearly a contradiction (since $|m| > 1$). So $h$ is not the multiple of any other class in $H_2(\cpt; \mz)$. \newline 

Then, since $H_2(\cpt; \mz) \cong \mz$, $h$ must be a generator of $H_2(\cpt; \mz) \cong \mz$. Furthermore, $H_2(\cpt; \mz) \cong \mz$ implies $rk(Q_{\cpt}) = 1$, and since $Q_{\cpt}(h,h')=1$, we must have $Q_{\cpt} = <1>$. 


\newtheorem{eg22}[eg1]{Example}
\begin{eg22}
\upshape
We define $\cptbar$ to be the manifold $\cpt$ with the opposite orientation. Therefore, by Remark \ref{difr25} (iv), we have $Q_{\cptbar} = -Q_{\cpt} = <-1>$.
\end{eg22}

\newtheorem{eg3}[eg1]{Example}
\begin{eg3}
\upshape
Consider the manifold $\mathbb{CP}^1 \times \mathbb{CP}^1$. Looking at $\mathbb{CP}^1$ more closely, we see that it is actually a real 2-manifold. Moreover, by the lemma above, we know that $\mathbb{CP}^1$ is closed and simply-connected, so by the classification of compact 2-manifolds (see \cite{GP}), it must be homeomorphic to $S^2$. Therefore, $\mathbb{CP}^1 \times \mathbb{CP}^1$ is homeomorphic to $S^2 \times S^2$. \newline

Since $\pi_1(S^2 \times S^2) \cong \pi_1(S^2) \times \pi_1(S^2) \cong 0$, $S^2 \times S^2$ is a simply-connected, closed 4-manifold. Since $H_2(S^2 \times S^2; \mz) \cong \mz \oplus \mz$ (see \cite{H1}), we know the intersection form $Q_{S^2 \times S^2}$ has rank 2. If we choose as a basis for $H_2(S^2 \times S^2; \mz)$ the homology elements $\alpha_1 = [S^2 \times \mathrm{pt}]$ and $\alpha_2 = [\mathrm{pt} \times S^2]$, we can see that $\alpha_1 \cdot \alpha_1 = 0$ and $\alpha_2 \cdot \alpha_2 = 0$, and $\alpha_1 \cdot \alpha_2 = \alpha_2 \cdot \alpha_1 = 1$. Therefore, the intersection form is

\begin{displaymath}
Q_{S^2 \times S^2} = 
\left( \begin{array} {cc}
0 & 1 \\
1 & 0 
\end{array} \right)
\end{displaymath}

We usually denote this particular matrix by $H$.
\end{eg3}

We recall the definition of the \emph{connected sum} of two $n$-manifolds $X_1$ and $X_2$ from \cite{GS}:

\newtheorem{consumdef}[eg1]{Definition}
\begin{consumdef}
\upshape
Let $X_1$ and $X_2$ be two smooth $n$-dimensional manifolds. Let $D_1 \subset X_1$ and $D_2 \subset X_2$ be two embedded $n$-disks, and let $\phi: D_1 \longrightarrow D_2$ be an orientation-reversing diffeomorphism. The \emph{connected sum} $X_1 \# X_2$ of $X_1$ and $X_2$ is defined to be the smooth manifold $(X_1  \setminus D_1) \cup_{\phi|_{\partial D_1}} (X_2 \setminus D_2)$.
\end{consumdef}

\newtheorem{consumdefr}[eg1]{Remarks}
\begin{consumdefr}
\upshape
Note that the connected sum operation is well-defined, in the sense that it does not depend on our choice of disks $D_1, D_2$ or on our choice of homeomorphism $\phi$. We sometimes denote by $\#m X$ the connect sum of $m$ copies of the manifold $X$ (where $m \geq 0$, and if $m=0$ then $\#m X = S^{n}$). Note that $X \# S^n$ is simply $X$.
\end{consumdefr}

\newtheorem{consumdefe1}[eg1]{Remarks}
\begin{consumdefe1}
\upshape
Let $T^2$ denote the familiar genus-1 surface, the torus. Then $T^2 \# T^2$ is the surface of genus 2, and $\#m T^2$ is the surface of genus $m$.
\end{consumdefe1}

Now, there is a simple equation relating the intersection forms of $X_1$ and $X_2$ and their connect sum $X_1 \# X_2$:

\newtheorem{eg4}[eg1]{Lemma}
\begin{eg4}
\upshape
Let $X_1$ and $X_2$ be 4-manifolds with intersections forms $Q_{X_1}$ and $Q_{X_2}$, respectively. Then the connect sum $X_1 \# X_2$ has intersection form
\begin{equation}
Q_{X_1 \# X_2} = Q_{X_1} \oplus Q_{X_2} 
\label{eg4eq} 
\end{equation}
\end{eg4}

\newtheorem{eg4r}[eg1]{Remark}
\begin{eg4r}
\upshape
This is an important lemma that we shall use many times in later sections. Equation $\eqref{eg4eq}$ is proved using a Mayer-Vietoris sequence. See \cite{GS} for the proof.
\end{eg4r}

\newtheorem{eg5}[eg1]{Example}
\begin{eg5}
\upshape
Consider the manifold with intersection form given by the matrix:

\begin{displaymath}
-E_8 = 
\left( \begin{array} {cccccccc}
-2 & 1 & 0 & 0 & 0 & 0 & 0 & 0 \\
1 & -2 & 1 & 0 & 0 & 0 & 0 & 0 \\
0 & 1 & -2 & 1 & 0 & 0 & 0 & 0 \\
0 & 0 & 1 & -2 & 1 & 0 & 0 & 0 \\
0 & 0 & 0 & 1 & -2 & 1 & 0 & 1 \\
0 & 0 & 0 & 0 & 1 & -2 & 1 & 0 \\
0 & 0 & 0 & 0 & 0 & 1 & -2 & 0 \\
0 & 0 & 0 & 0 & 1 & 0 & 0 & -2
\end{array} \right)
\end{displaymath}

If one were to denote the basis used to represent $-E_8$ as it appears above by $\{\alpha_1, \alpha_2, \dots, \alpha_8 \}$, the matrix shows us that $Q_{-E_8}(\alpha_{i}, \alpha_{i}) \equiv 0$ (mod 2) for $i = 1,2, \dots, 8$. By Lemma \ref{L25}, this shows us that $-E_8$ is even. By diagonalizing $-E_8$ over $\mathbb{R}$, one finds that $-E_8$ is negative-definite, since $\sigma(-E_8) = -8 = -rk(-E_8)$. Note that we should have expected that $\sigma(-E_8) \equiv 0$ (mod 8) by Lemma \ref{L210}. Finally, since $\sigma(-E_8)$ is not divisible by 16, by (Rohlin's) Theorem \ref{T35}, the 4-manifold with intersection form $-E_8$ cannot admit any smooth structures. \newline

So, we have our first concrete example of a topological 4-manifold that does not admit a smooth structure.

\end{eg5}

\pagebreak



\section{Symplectic Structures, Almost Complex Structures, Complex Structures, and Characteristic Classes} \label{canclasssection}

In this section we recall the definitions of \emph{symplectic structures} and \emph{almost complex structures}. Almost everything presented in this section is from \cite{MS2}, which goes into far more detail. Another good reference, from where we borrow a few definitions, is \cite{LM}. The reason for this discussion is that we shall need the fact that every symplectic manifold has a \emph{canonical class} associated to it that is \emph{compatible with the symplectic structure}. Theorems presented later on in this article use the canonical class in order to distinguish between smooth structures on a manifold. \newline

The topic of characteristic classes is beyond the scope of this article, yet the canonical class is defined to be a ``certain'' Chern class. The reader can just take this definition at face value, and refer to \cite{MS2}, \cite{GS} or \cite{MiSt} for more on characteristic classes.


\newtheorem{d31}{Definition}[section]
\begin{d31}
\upshape
A \emph{symplectic vector space} is a pair $(V, \omega)$ consisting of a finite-dimensional real vector space $V$ and a bilinear form $\omega: V \times V \longrightarrow \mathbb{R}$ satisfying the following two conditions:
\begin{itemize}
\item[(1)] For all $v,w \in V$, $\omega(v,w) = -\omega(w,v)$.
\item[(2)] For every $v \in V$, if $\omega(v,w) = 0$ for all $w \in V$, then $v = 0$.
\end{itemize}
\end{d31}

\newtheorem{d31rr}[d31]{Remark}
\begin{d31rr}
\upshape
Condition $(1)$ is called \emph{skew-symmetry} and condition $(2)$ is called \emph{non-degeneracy}. Therefore, a symplectic form is a skew-symmetric, non-degenerate form.
\end{d31rr}

\newtheorem{d31rr2}[d31]{Remark}
\begin{d31rr2}
\upshape
A symplectic vector space must be even-dimensional, otherwise condition (2) will not be satisfied.
\end{d31rr2}

\newtheorem{d311}[d31]{Definition}
\begin{d311}
\upshape
Let $X$ be a manifold. We say $\omega$ is a \emph{2-form on $X$}, if for each $p \in X$, $\omega_{p}$ is a skew-symmetric bilinear map on the tangent space of $X$ at $p$, i.e. $\omega_{p}: T_{p}X \times T_{p}X \longrightarrow \mathbb{R}$. Furthermore, $\omega_{p}$ varies smoothly in $p$.
\end{d311}

\newtheorem{d312}[d31]{Definition}
\begin{d312}
\upshape
Let $\omega$ be a 2-form on a manifold $X$. We say that $\omega$ is a \emph{symplectic form} if $\omega$ is closed and $\omega_{p}$ is symplectic on $T_{p}X$ for all $p \in X$.
\end{d312}

\newtheorem{d313}[d31]{Definition}
\begin{d313}
\upshape
A \emph{symplectic manifold} is a pair $(X, \omega)$ where $X$ is a manifold and $\omega$ is a symplectic form.
\end{d313}

\newtheorem{d32}[d31]{Remarks}
\begin{d32}
\upshape
See \cite{MS2} for more on the following remarks:
\begin{itemize}
\item[(i)] If $X$ is a symplectic manifold, then $X$ must be even-dimensional.
\item[(ii)] If $(X, \omega)$ is a symplectic manifold of dimension $2n$, then the $n$-fold wedge product $\omega \wedge \dots \wedge \omega$ is never zero. This implies that a symplectic manifold $(X, \omega)$ is orientable.
\end{itemize}
\end{d32}

If we want to define two symplectic manifolds $(X_1, \omega_1)$ and $(X_2, \omega_2)$ to be equivalent, not only do we need the underlying smooth manifolds to be diffeomorphic, but we also need the two symplectic forms to be related in some way.  The phrase we use for such equivalence is ``$X_1$ is \emph{symplectomorphic} to $X_2$'', and the map (the diffeomorphism that ``preserves'' the symplectic structure) is called a \emph{symplectomorphism}.

\newtheorem{dsplmphsm}[d31]{Definition}
\begin{dsplmphsm}
\upshape
Let $(X_1, \omega_1)$ and $(X_2, \omega_2)$ be two symplectic manifolds, both of dimension $2n$, and let $f: X_1 \longrightarrow X_2$ be a diffeomorphism. Then $f$ is a \emph{symplectomorphism} if $f^{*}\omega_2 = \omega_1$.
\end{dsplmphsm}

\newtheorem{dsplmphsmr}[d31]{Remark}
\begin{dsplmphsmr}
\upshape
$f^{*}\omega_2$ is the pullback of $\omega_2$ by $f$. See \cite{I1}, \cite{MS2}, or \cite{LM} for details.
\end{dsplmphsmr}

\newtheorem{dcsvs}[d31]{Definition}
\begin{dcsvs}
\upshape
Let $V$ be a vector space. A \emph{complex structure} on $V$ is an automorphism $J: V \longrightarrow V$ such that $J^{2} = -I$, where $I$ is the identity automorphism on $V$. With such a structure, $V$ becomes a complex vector space with multiplication $i = \sqrt{-1}$ corresponding to $J$, by the map
\begin{equation*}
\mathbb{C} \times V \longrightarrow V: (s + it, v) \longmapsto sv + tJv
\end{equation*}
$V$ must be even-dimensional over $\mathbb{R}$.
\end{dcsvs}

\newtheorem{dcsvscom}[d31]{Definition}
\begin{dcsvscom}
\upshape
Let $(V,\omega)$ be a symplectic vector space. A complex structure $J$ on $V$ is said to be \emph{compatible with $\omega$} if for all $v,w \in V$
\begin{equation*}
\omega(Jv,Jw) = \omega(v,w)
\end{equation*}
and if for all $v,w \in V$, with $v$ nonzero,
\begin{equation*}
\omega(v,Jw) > 0
\end{equation*}
\end{dcsvscom}

\newtheorem{d33}[d31]{Definition}
\begin{d33}
\upshape
Let $X$ be a $2n$-dimensional manifold. An \emph{almost complex structure} on $X$ is a complex structure $J$ on the tangent bundle $TM$. A non-degenerate 2-form $\omega$ on $X$ is called \emph{compatible} with $J$ if the bilinear form $g$ defined by
\begin{equation*}
g(v,w) = \omega(v, Jw)
\end{equation*} 
defines a Riemannian metric on $X$.
\end{d33}

\newtheorem{d34}[d31]{Definition}
\begin{d34}
\upshape
A Riemannian metric $g$ on $M$ is called \emph{compatible with $J$} if for all $v,w \in T_{p}X$,
\begin{equation*}
g(Jv,Jw) = g(v,w)
\end{equation*}
\end{d34}

All of the above was done to make the following proposition intelligible:

\newtheorem{prop3}[d31]{Proposition}
\begin{prop3}
\upshape
\cite{MS2} \itshape Let $X$ be a $2n$-dimensional manifold. Then
\begin{itemize}
\item[(i)] for each non-degenerate 2-form $\omega$ on $X$, there exists an almost-complex structure $J$ which is compatible with $\omega$.
\item[(ii)] for each almost complex structure $J$ on $X$ there exists a non-degenerate 2-form $\omega$ which is compatible with $J$.
\end{itemize}
\end{prop3}

Now for the canonical class: 
\newtheorem{d35}[d31]{Definition}
\begin{d35}
\upshape
If $X$ has an almost complex structure $J$, its tangent bundle $TX$ and its cotangent bundle $T^{*}X$ are complex rank 2 bundles. The \emph{canonical class} $K = K(X)$ is defined to be the first Chern class of the cotangent bundle, i.e.
\begin{equation*}
K = c_{1}(T^{*}X, J) = -c_{1}(TX, J) \in H^{2}(X, \mz)
\end{equation*}
\end{d35}

\newtheorem{d35r}[d31]{Remark}
\begin{d35r}
\upshape
Note that we shall often write $Q_{X}(K,K) = K \cdot K = K^2$.
\end{d35r}


We then have the following result from \cite{Wu}, quoted from \cite{GS}:

\newtheorem{wut1}[d31]{Theorem}
\begin{wut1} \label{wuthm}
For a given 4-manifold $X$ and an almost-complex structure $J$ on $X$, Let $K = c_{1}(T^{*}X, J)$ be the canonical class. Then  $K^2 = 3 \sigma(X) + 2 \chi(X)$.
\end{wut1}

\newtheorem{wut1r}[d31]{Remark}
\begin{wut1r}
\upshape
In the theorem above, $\sigma(X)$ denotes the signature of the intersection form of $X$ and $\chi(X)$ denotes the Euler characteristic of $X$. Note that we are presenting a very ``watered-down'' version of the original result; an almost-complex structure on $X$ provides two further identities (one concerning another characteristic class, the first \emph{Stiefel-Whitney class} of $X$), and there is an appropriate converse. Since we shall not need these extra identities, or the converse, the current version of the theorem is sufficient for our purposes.
\end{wut1r}


\pagebreak
\section{Blowing Up and Blowing Down}

\newtheorem{s41}{Definition}[section]
\begin{s41}
\upshape
Let $X$ be a smooth, oriented manifold. The connected sum $X' = X \# \cptbar$ is called the \emph{blow-up} of $X$ at a point. We obtain a map $\pi: X' \longrightarrow X$ with the following properties: For a point $P \in X$,
\begin{itemize}
\item[(1)] $\pi|_{X' \setminus \cpobar}: X' \setminus \cpobar \longrightarrow X \setminus \{P\}$ is a diffeomorphsim
\item[(2)] $\pi^{-1}(P) = \cpobar$
\end{itemize}
\end{s41}

\newtheorem{bur}[s41]{Remarks}
\begin{bur}
\upshape
The sphere $\cpobar$ in (ii) is called the \emph{exceptional sphere}. Note that 
\begin{itemize}
\item[(i)] it is contained in the $\cptbar$ summand of $X'$. 
\item[(ii)] its homology class $[\cpobar]$ is usually denoted by $e = [\cpobar] \in H_2(X'; \mz) = H_2(X; \mz) \oplus H_2(\cptbar; \mz)$ \linebreak (so, actually $e \in H_{2}(\cptbar; \mz)$).
\item[(iii)] $Q_X(e,e) = -1$. 
\item[(iv)] We usually call the map $\pi: X' \longrightarrow X$ the projection map.
\end{itemize}
\end{bur}

\newtheorem{bur2}[s41]{Remark}
\begin{bur2}
\upshape
Informally, when we blow up a manifold at a point $P$, we replace the point $P$ with the space of all lines going through $P$, which is a copy of $\cpo$. 
\end{bur2}

\newtheorem{bupt}[s41]{Definition}
\begin{bupt}
\upshape
Let $X$ be a smooth 4-manifold and let $\Sigma$ be a smooth surface in $X$. Suppose we blow up $X$ at a point $P \in \Sigma$, and we denote the projection by $\pi: X' \longrightarrow X$. We define
\begin{itemize}
\item[(i)] the \emph{total transform} of $\Sigma$ to be the inverse image $\Sigma' = \pi^{-1}(\Sigma) \subset X'$.
\item[(ii)] the \emph{proper transform} of $\Sigma$ to be the closure $\tilde{\Sigma} = cl(\pi^{-1}(\Sigma \setminus \{P\}))$.
\end{itemize}
\end{bupt}

\newtheorem{bur22}[s41]{Remark}
\begin{bur22}
\upshape
So, suppose $\Sigma_1$ and $\Sigma_2$ are smooth surfaces in $X$ intersecting each other transversally only in the point $P$. Let $X'$ be the blow-up of $X$ at $P$ and let $\pi:X' \longrightarrow X$ be the projection map. Then, the proper transforms $\tilde{\Sigma}_{1}$, $\tilde{\Sigma}_{2}$ will be disjoint in the blow-up $X'$.
\end{bur22}

There is an inverse operation, called a \emph{blow down}, which can be performed under certain conditions.

\newtheorem{bd1}[s41]{Definition}
\begin{bd1}
\upshape
If the 4-manifold $X$ contains a sphere $\Sigma_{-}$ with $[\Sigma_{-}]^2 = -1$, then $X = Y \# \cptbar$ for some manifold $Y$, and $Y$ is called the \emph{blow down} of $X$.
\end{bd1}

\newtheorem{bd1r}[s41]{Remark}
\begin{bd1r}
\upshape
There are corresponding blow-up and blow-down operations using $\cpt$ instead of $\cptbar$, but then the manifolds cease to be complex, as the following theorem (quoted from \cite{GS}), called the Noether formula, shows. We shall always use $\cptbar$, unless we specify otherwise.
\end{bd1r}

\newtheorem{nf1}[s41]{Theorem}
\begin{nf1}
For a complex surface $S$, the integer $c_1^2(S) + c_2(S) = 3(\sigma(S)+ \chi(S))$ is divisible by 12, or equivalently, $1-b_1(S)+b_2^{+}(S)$ is even. In particular, if $S$ is a simply-connected complex surface, then $b_2^{+}(S)$ is odd.
\end{nf1}

\newtheorem{bur23}[s41]{Remark}
\begin{bur23}
\upshape
Note that the blow-up can be defined holomorphically for complex manifolds (see \cite{GS}) and symplectically for symplectic manifolds (see \cite{MS2}).
\end{bur23}

\newtheorem{bur24}[s41]{Remark}
\begin{bur24}
\upshape
If $X' = X \# \cptbar$ we have
\begin{itemize}
\item[(i)] $b_2^{+}(X') = b_2^{+}(X)$
\item[(ii)] $b_2^{-}(X') = b_2^{-}(X) + 1$
\item[(iii)] $b_2(X') = b_2(X) + 1$
\item[(iv)] $\sigma(X') = \sigma(X) -1 $
\item[(v)] $\chi(X') = \chi(X) + 1$
\end{itemize}
\end{bur24}

\pagebreak
\section{The Rational Blowdown Technique}
\label{rbdsec}

The rational blowdown technique was first discovered by R. Fintushel and R. Stern (\cite{FS1}). It can be thought of as a generalization of the usual blowdown process, in that a certain \emph{configuration} of spheres $C_p$ (basically, just a special collection of spheres intersecting each other in a certain way and with certain self-intersection numbers) is removed from a manifold and is replaced by a \emph{rational 4-ball} $B_p$ which has the same boundary, i.e. $\partial C_p = \partial B_p$. The reason the rational blowdown technique is useful, is that it is relatively easy to calculate the Seiberg-Witten invariants of manifolds constructed with this technique. \newline

For a fixed prime $p$ , a rational 4-ball $B_p$ is a 4-manifold that has the same rational homology as a ball, i.e. $H_k(B_p; \mathbb{Q}) \cong 0$ for $k>0$. However, this does not mean that its integral homology groups $H_k(B_p; \mz)$ are also trivial, but that if they are not trivial, then they are just finite (torsion) groups. \newline

For example, we consider the rational 4-ball $B_p$, given by the Figure \ref{rbdsec}.1.




\newtheorem{bp}{Lemma}[section]
\begin{bp}
$B_p$ has trivial rational homology, and so has the same rational homology as $D^4$.
\end{bp}
Proof: \newline
Using techniques explained in \cite{OzSt}, pages 42 and 43 (and discussed in section \ref{x7proofsection}), we can calculate that $H_1(B_p; \mz) \cong \mz_{p}$, and $H_k(B_p; \mz)$ are trivial for $k \geq 2$ (and $H_{0}(B_p; \mz) \cong 0$ since it is path-connected). So, $B_p$ has the same rational homology as $D^4$ (all homology groups are trivial). $\Box$  

\newtheorem{bprem1}[bp]{Remark}
\begin{bprem1}
\upshape
Using these techniques, it can also be shown that $H_1(\partial B_p; \mz) \cong \mz_{p^2}$.
\end{bprem1}

We define $C_p$ to be the 4-manifold that is the plumbing according to the graph Figure \ref{rbdsec}.2, where $p \geq 2$ and the number of $-2$'s is $p-2$. So, its Kirby diagram is given by Figure \ref{rbdsec}.3. We now need a short lemma.

\newtheorem{cpbound1}[bp]{Lemma} 
\begin{cpbound1}\label{cpzp2lemma}
The boundary $\partial C_p$ is the lens space $L(p^2, p-1)$, and so $\pi_1( \partial C_p) \cong \mz_{p^2}$.
\end{cpbound1}
Proof: \newline
First we note that the continued fraction expansion of $p^2/p-1$ is
\begin{equation*}
\frac{p^2}{p-1} = p+2 - \frac{1}{2 - \frac{1}{2 - \dots}}
\end{equation*}

or, $[p+2, 2, \dots, 2]$, where there are $p-2$ 2's. By using the slam-dunk technique in \cite{GS} from right to left on Figure \ref{rbdsec}.3, we get a single unknot with coefficient $p^2/p-1$, which shows that $\partial C_p \cong L(p^2, p-1)$. \newline

It is well-known (see for example \cite{Ro}) that $\pi_1(L(a,b)) \cong \mz_a$ (whatever the value of $b$ is), and so 
\begin{equation*}
\pi_1(\partial C_p) \cong \pi_1(L(p^2,p-1)) \cong \mz_{p^2}.
\end{equation*}
\begin{flushright}
$\Box$
\end{flushright}

We shall now show that $\partial C_p \cong \partial B_p$. We do this indirectly by showing that $B_p \cup_{\partial} \overline{C_p} \cong \# (p-1) \cpt$, which also shows that $C_p$ embeds in $\# (p-1) \cptbar$. First, we need another (equivalent) Kirby diagram for $C_p$.

\newtheorem{cpequiv1}[bp]{Lemma}
\begin{cpequiv1}
Figure \ref{rbdsec}.4 is also a Kirby diagram for $C_p$, given in Figure \ref{rbdsec}.3.
\end{cpequiv1}
Proof: \newline
We start with Figure \ref{rbdsec}.4, and by a sequence of handleslides and a handle-cancellation, the diagram shall become Figure \ref{rbdsec}.3. We start by sliding a $-1$-framed meridian along the dotted circle to the place where the dotted circle and the $0$-framed circle twist, and then slide the $0$-framed circle over the meridian, which ``removes'' a twist. See Figure \ref{rbdsec}.5.\newline

We slide this ``used-up'' meridan along both strands, below a twist, and then use another $-1$-framed meridian to remove another twist. Note that this handleslide decreases the framing of the 0-framed circle by 1 each time. We do this with each of the meridians, and obtain Figure \ref{rbdsec}.6. \newline

We then slide the $-1$-framed meridians over each other, from top to bottom, as in Figure \ref{rbdsec}.7, which gives us Figure \ref{rbdsec}.8. Note that $p-2$ of the meridians become $-2$-framed meridians. \newline

We then slide the $-(p-1)$-framed circle over the only $-1$-framed meridian, using handle-subraction to get \ref{rbdsec}.18. We do this by looking at the part of the diagram indicated by Figures \ref{rbdsec}.9 and \ref{rbdsec}.10. First, we assign orientations to the $-(p-1)$-framed circle and the $-1$-framed meridian as in Figure \ref{rbdsec}.11 (note that we have chosen orientations so that the linking number of the two circles is $1$). Then, we draw a parallel copy of the $-1$-framed meridian, as in Figure \ref{rbdsec}.12, and perform the handle subtraction to get Figure \ref{rbdsec}.13. \newline

The framing of the $-(p-1)$-framed circle becomes
\begin{equation*}
-(p-1) - 1 - 2(1) = -p + 1 -1 -2 = -p -2
\end{equation*}
and note that it is linked with both the $-1$-framed meridian and the ``lowest'' $-2$-framed circle. \newline

We can drop the orientations and perform a few Reidemeister moves, as in Figures \ref{rbdsec}.14 to \ref{rbdsec}.17, to finally get \ref{rbdsec}.18. \newline

Now, since the dotted circle is only linked with the $-1$-framed meridian, we can perform a handle-cancellation, and so we are left with a link diagram as in Figure \ref{rbdsec}.3, as required. $\Box$ \newline

We now prove the main result.

\newtheorem{rbmain}[bp]{Proposition}
\begin{rbmain}
$B_p \cup_{\partial C_p} \overline{C_p}$ is diffeomorphic to $\# (p-1) \cpt$, and so $\partial C_p \cong \partial B_p$.
\end{rbmain}

Proof: \newline
We start with $DC_p$, the double of $C_p$. As we know from the section on Kirby Calculus, $DC_p$ is formed simply by attaching $0$-framed meridians to each link component, and then adding a 3-handle and a 4-handle. So, taking the double of $C_p$ in Figure \ref{rbdsec}.4 we get Figure \ref{rbdsec}.19. \newline

We perform surgery inside $C_p \subset DC_p$ twice, first to change the dotted circle into a 0-framed circle, and then to change the original 0-framed circle into a dotted circle. Although this surgery changes the 4-manifold, it does not change its boundary. So, although the surgery might change something ``inside'' $DC_p$, its boundary $C_p \cup_{\partial C_p} \overline{C_p}$ will be unchanged.

 So, we now have Figure \ref{rbdsec}.20. \newline

We slide each $0$-framed meridians over the $-1$-framed meridian it is linked with, so that the framing of the $0$-framed meridian becomes $+1$ and the meridians become unlinked, as shown in Figure \ref{rbdsec}.21. Blowing down the $-1$-framed circles in $C_p$ gives us Figure \ref{rbdsec}.22. Note the change in the framing of the ``large'' $0$-framed circle. Note that we can see this as a diagram of $B_p \cup_{\partial} \overline{C_p}$, since we did not do anything in the $\overline{C_p}$ half of $DC_p$, which now consists of the 1-framed meridians, the 0-framed meridian, the 3-handle and the 4-handle, while $B_p$ consists of the $(p-1)$-framed circle twisted around the dotted circle, as in Figure \ref{rbdsec}.1. \newline

Sliding the $(p-1)$-framed circle over its 1-framed meridians, as in Figure \ref{rbdsec}.23, we get Figure \ref{rbdsec}.24. Note that each handle-slide in Figure \ref{rbdsec}.23 decreases the framing of the $(p-1)$-framed circle circle by 1. \newline

We now use the 0-framed meridian to unlink the dotted circle and the 0-framed circle (Figure \ref{rbdsec}.25, repeated $p$ times) to get Figure \ref{rbdsec}.26. We now perform two handle-cancellations: first, the dotted circle and its 0-framed meridian cancel as a 1-handle/2-handle cancelling pair, and then the unlinked 0-framed circle cancels with the 3-handle, as a 2-handle/3-handle cancelling pair. This leaves us with Figure \ref{rbdsec}.27, which is simply $\# (p-1) \cpt$, as required. $\Box$.

\newtheorem{cporrem1}[bp]{Remark}
\begin{cporrem1}
\upshape
Sometimes we shall denote the boundary of $C_p$ as $L(p^2,1-p)$, as other authors do, instead of as $L(p^2,p-1)$. This is fine, since $\overline{L(a,b)} = L(a,-b)$, and so we are just considering the opposite orientation.
\end{cporrem1}

\newtheorem{cporrem2}[bp]{Remark}
\begin{cporrem2}
\upshape
It should be noted that we have only shown that there is a diffeomorphism $\phi:\partial B_{p} \longrightarrow \partial C_p$ (i.e. it is a self-diffeomorphism of $\partial B_p = L(p^2, 1-p)$). In order for this operation of rationally blowing down to be well-defined, we need that a self-diffeomorphism of the boundary $\partial B_p$ always extends to a diffeomorphism over the whole rational ball $B_{p}$. Fortunately, the following theorem due to Bonahon in \cite{Bo} shows that there are not too many self-diffeomorphisms of $\partial B_p$ to consider.
\end{cporrem2}

\newtheorem{bon1}[bp]{Theorem}
\begin{bon1}
$\pi_{0}(\mathrm{Diff}(L(p^2, 1-p))) \cong \mathbb{Z}_2.$
\end{bon1}

\newtheorem{cporrem3}[bp]{Remark}
\begin{cporrem3}
\upshape
So, this theorem is saying that, up to homotopy, there are exactly two non-homotopic maps that are self-diffeomorphisms of $\partial B_p$. The identity map is one of these diffeomorphisms, and it clearly extends to a diffeomorphism over $B_p$. As noted in \cite{GS}, if we consider Figure \ref{rbdsec}.1 as in Figure \ref{rbdsec}.28 below, we see that it is symmetric (by a $180^{\circ}$ rotation about the $y$-axis).\newline

Let us denote by $R$ the map that performs this roation. Clearly, $R$ is a self-diffeomorphism of $\partial B_p$, and it is also fairly clear that this self-diffeomorphism extends to $B_p$. If we consider what $R$ does to a meridian $m_1$ of the dotted circle $c_1$, we see that it inverts the meridian (gives it the opposite orientation). This shows that $R$ is not homotopic to the identity map. \newline


We have therefore found two non-homotopic self-diffeomorphisms of $\partial B_p$ which extend to $B_p$, and so we have the following theorem, as in \cite{GS}.
\end{cporrem3}

\newtheorem{bon2}[bp]{Theorem}
\begin{bon2} \label{bon2label}
Any self-diffeomorphism of $\partial B_p$ extends to $B_p$.
\end{bon2}

This finally allows us to give the definition of the rational blowdown of a 4-manifold $X$, as in \cite{GS}, which by Theorem \ref{bon2label} is well-defined up to diffeomorphism for a fixed $X$ and a fixed $C_p$ embedded in $X$.

\newtheorem{rbdef1}[bp]{Definition}
\begin{rbdef1}
Assume that $C_p$ embeds in the 4-manifold $X$, and write $X$ as $X= X_{0} \cup_{L(p^2, p-1)} C_{p}$. The 4-manifold $X_p = X_0 \cup_{L(p^2,p-1)}B_p$ is by definition the rational blowdown of $X$ along the given copy of $C_p$.
\end{rbdef1}

\newtheorem{rbrem1}[bp]{Remark}
\begin{rbrem1}
\upshape
We shall need the observation that if $X$ and $X \setminus C_p$ are simply connected, then $X_p$ is simply connected. This will be shown in section \ref{x7proofsection} for the case $p=7$.
\end{rbrem1}

\pagebreak

\begin{center}
\begin{minipage}{8cm}
\includegraphics[width=8cm]{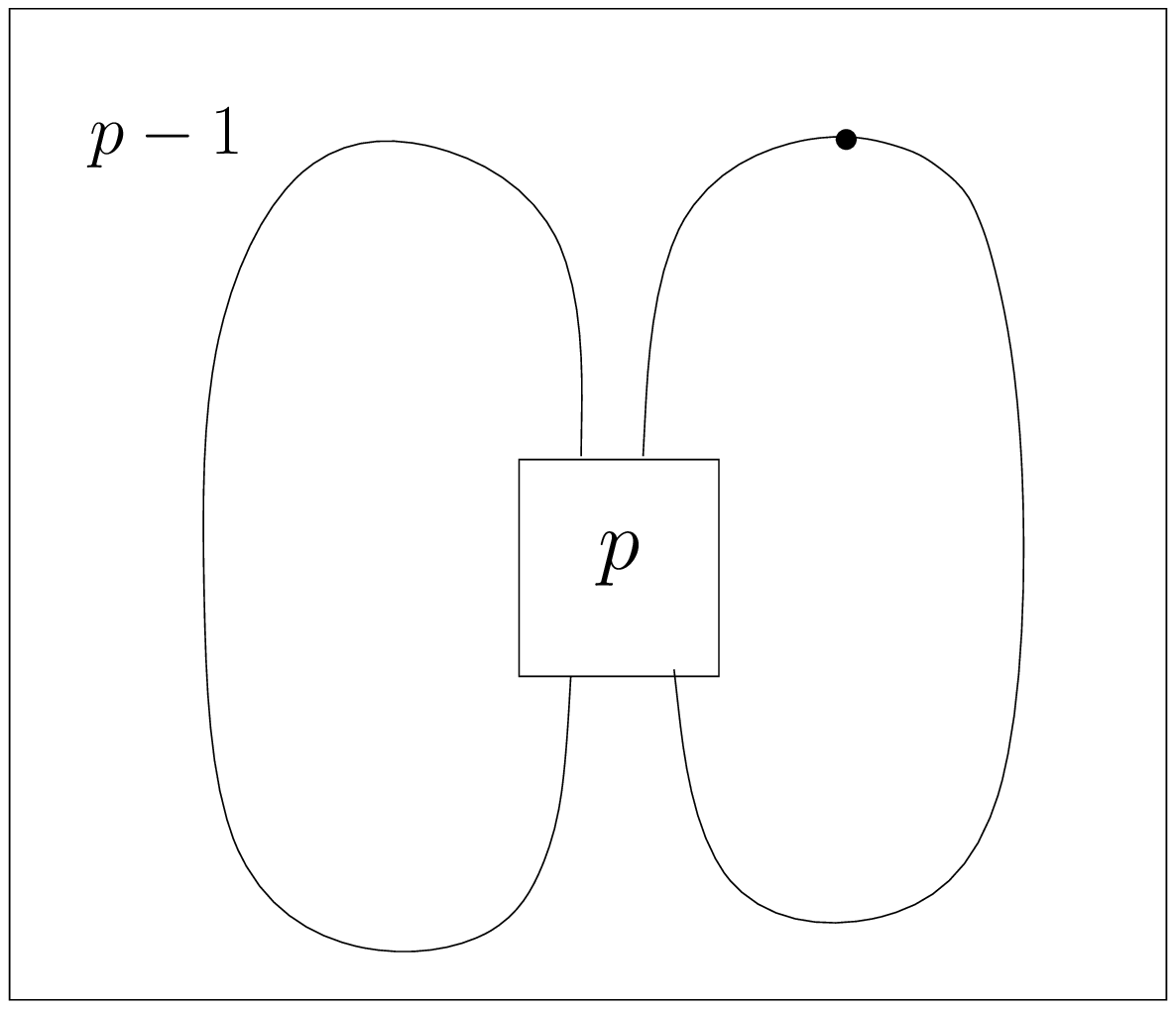}
\end{minipage}
\end{center}
\begin{center}
Figure \ref{rbdsec}.1
\end{center}

\begin{center}
\begin{minipage}{8cm}
\includegraphics[width=8cm]{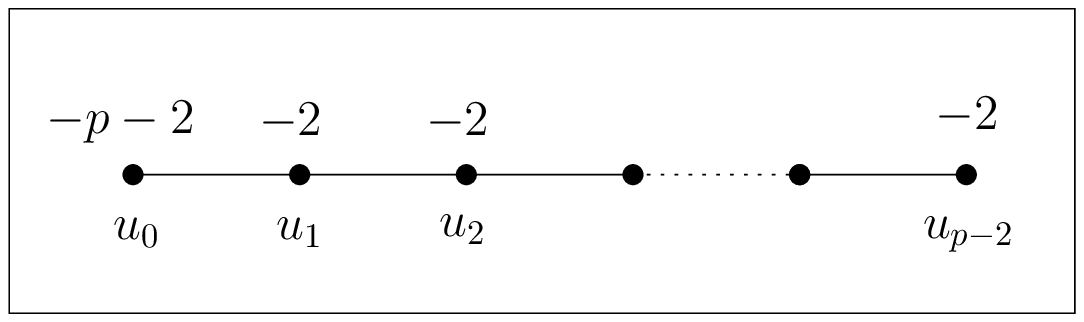}
\end{minipage}
\end{center}
\begin{center}
Figure \ref{rbdsec}.2
\end{center}

\begin{center}
\begin{minipage}{8cm}
\includegraphics[width=8cm]{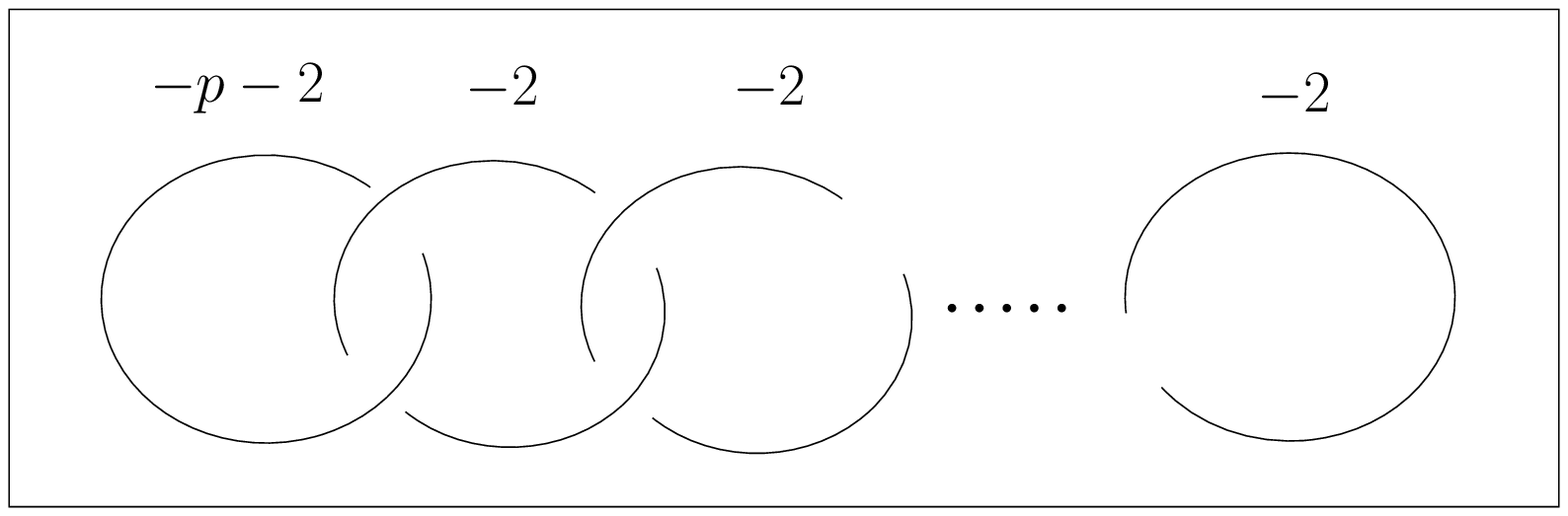}
\end{minipage}
\end{center}
\begin{center}
Figure \ref{rbdsec}.3
\end{center}

\begin{center}
\begin{minipage}{8cm}
\includegraphics[width=8cm]{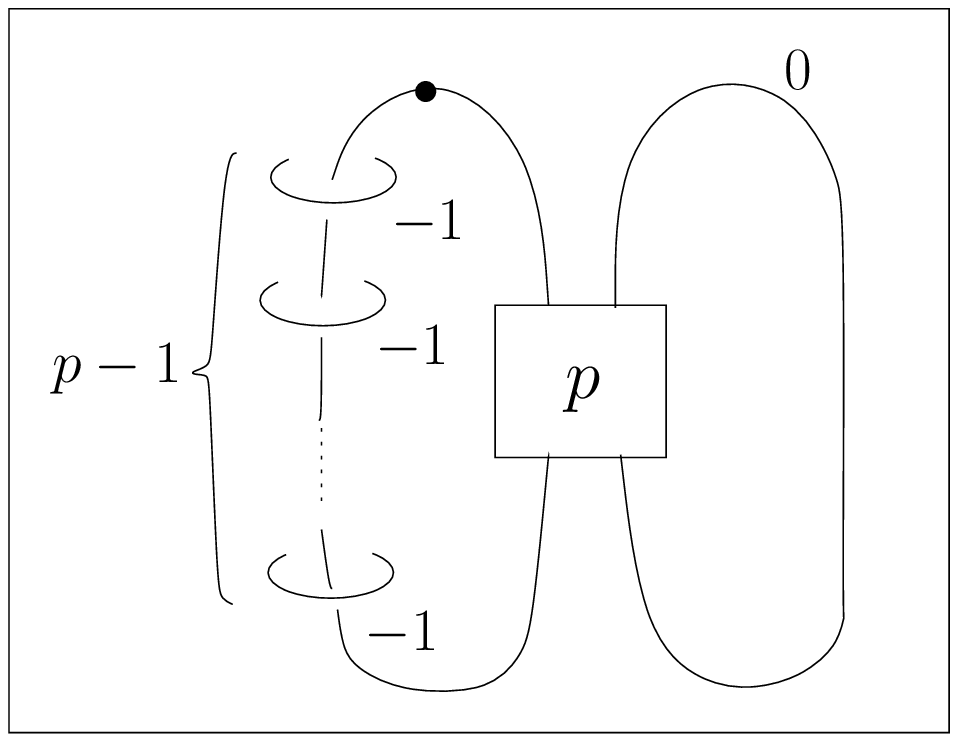}
\end{minipage}
\end{center}
\begin{center}
Figure \ref{rbdsec}.4
\end{center}

\begin{center}
\begin{minipage}{13cm}
\includegraphics[height=4cm]{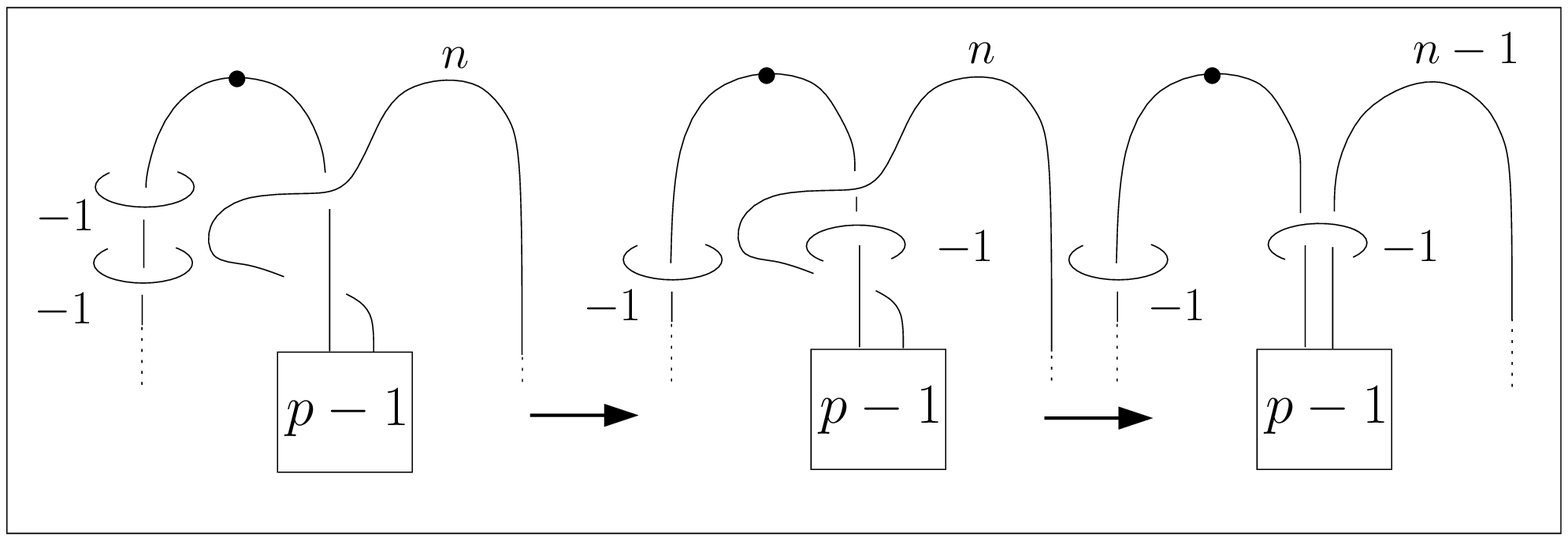}
\end{minipage}
\end{center}
\begin{center}
Figure \ref{rbdsec}.5
\end{center}

\pagebreak

\begin{center}
\begin{minipage}{8cm}
\includegraphics[width=6cm]{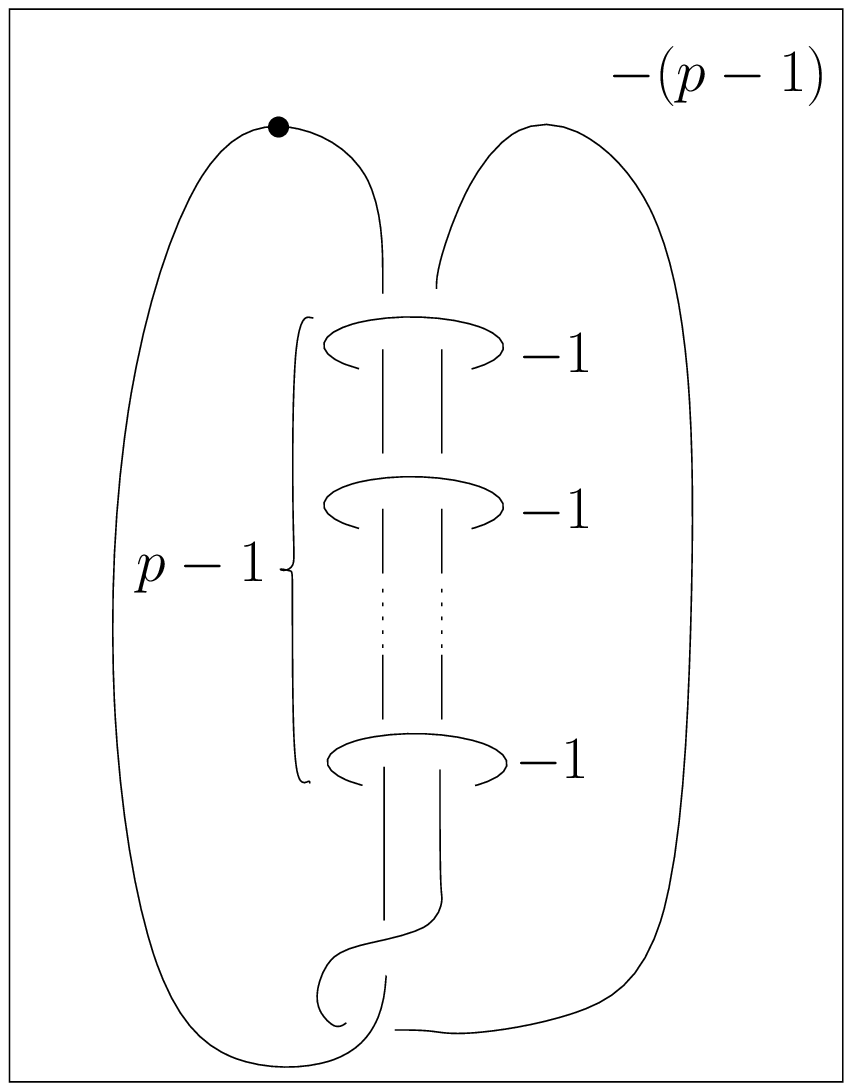}
\end{minipage}
\end{center}
\begin{center}
Figure \ref{rbdsec}.6
\end{center}

\begin{center}
\begin{minipage}{8cm}
\includegraphics[width=8cm]{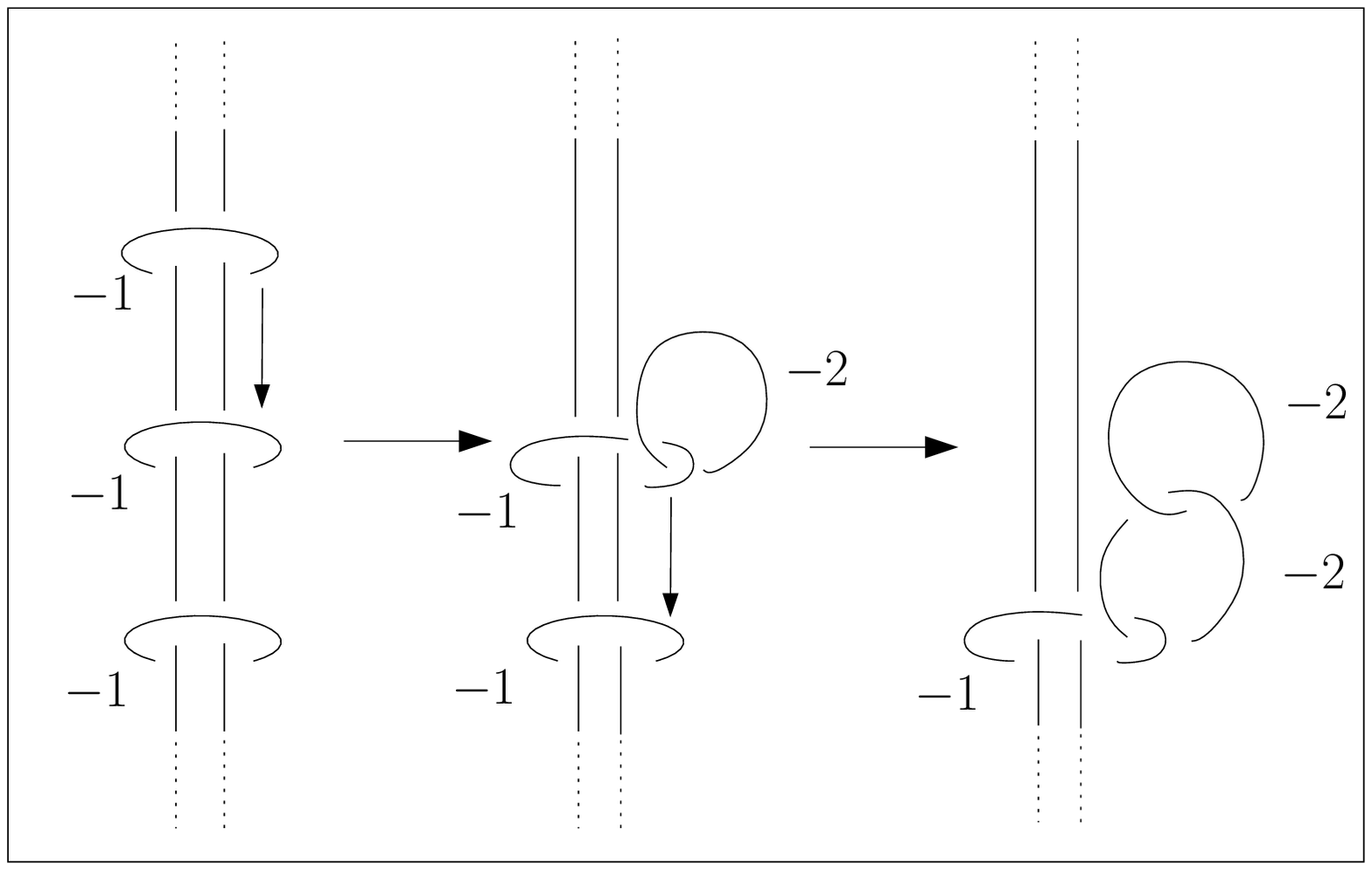}
\end{minipage}
\end{center}
\begin{center}
Figure \ref{rbdsec}.7
\end{center}

\begin{center}
\begin{minipage}{8cm}
\includegraphics[width=8cm]{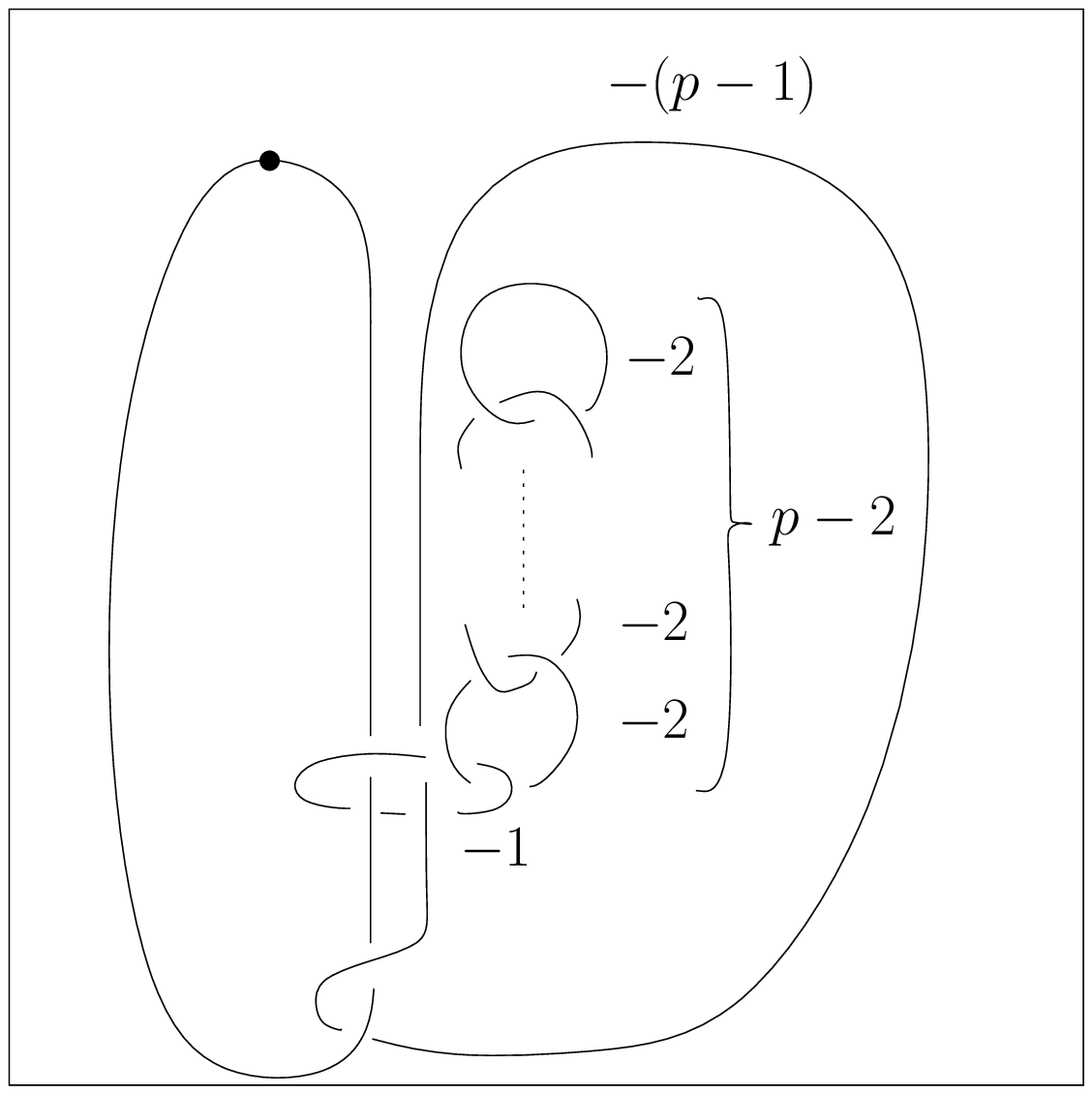}
\end{minipage}
\end{center}
\begin{center}
Figure \ref{rbdsec}.8
\end{center}

\begin{center}
\begin{minipage}{8cm}
\includegraphics[width=8cm]{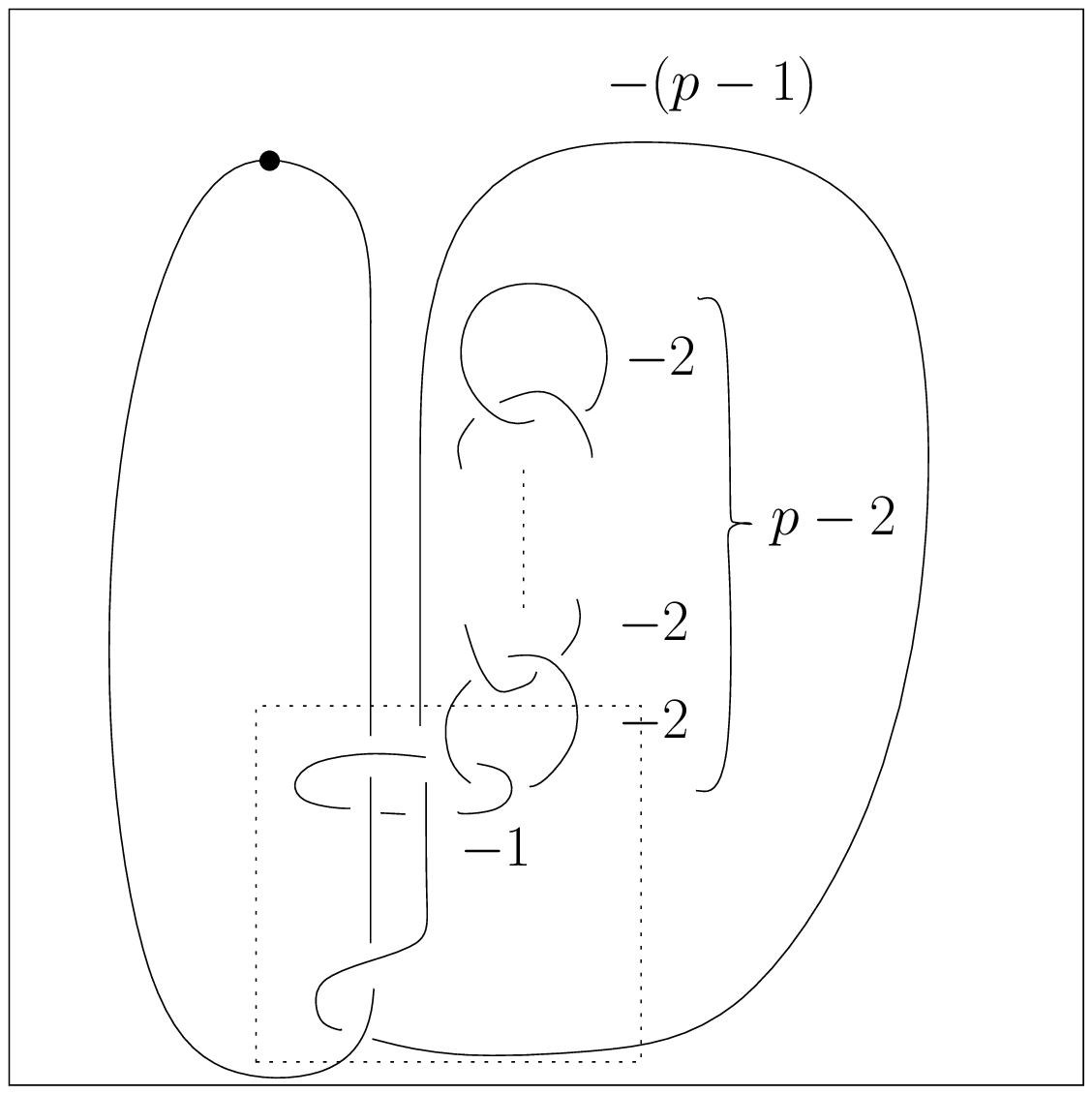}
\end{minipage}
\end{center}
\begin{center}
Figure \ref{rbdsec}.9
\end{center}

\begin{center}
\begin{minipage}{8cm}
\includegraphics[width=8cm]{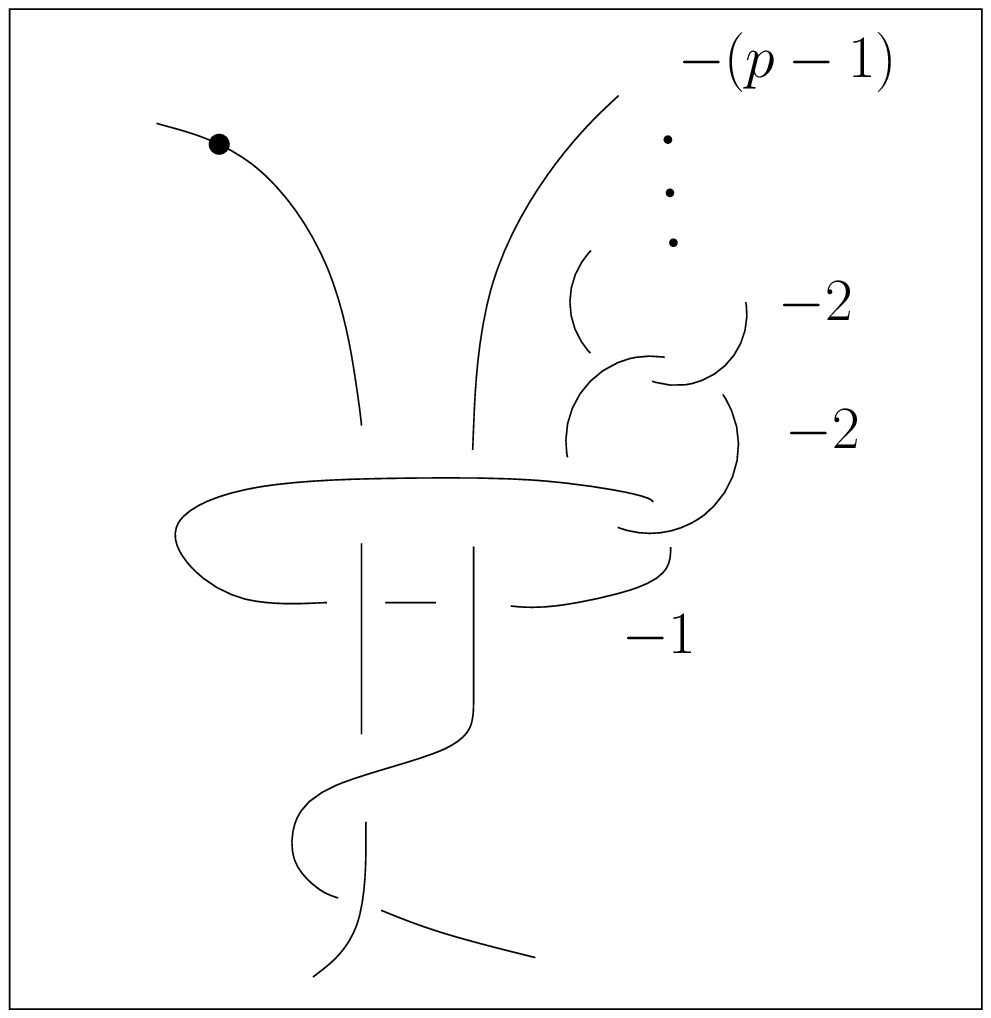}
\end{minipage}
\end{center}
\begin{center}
Figure \ref{rbdsec}.10
\end{center}

\begin{center}
\begin{minipage}{8cm}
\includegraphics[width=8cm]{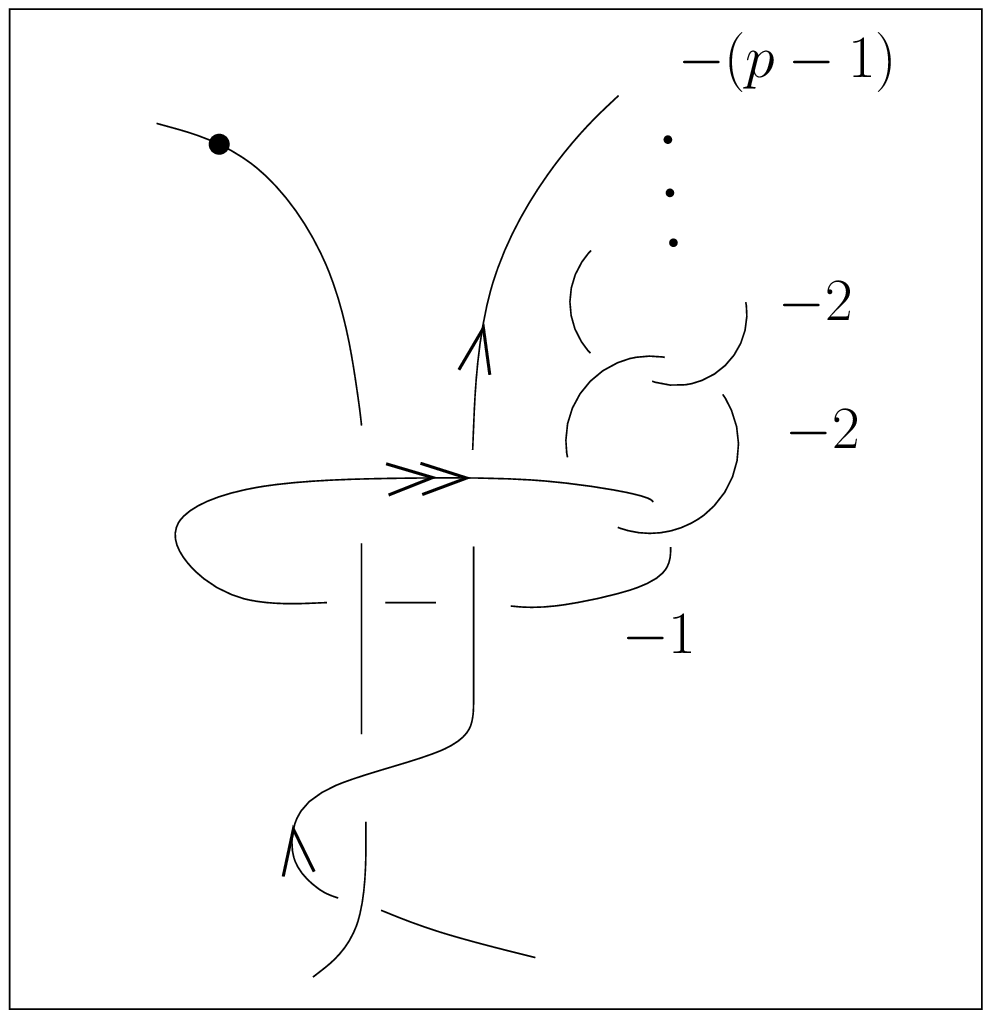}
\end{minipage}
\end{center}
\begin{center}
Figure \ref{rbdsec}.11
\end{center}

\begin{center}
\begin{minipage}{8cm}
\includegraphics[width=8cm]{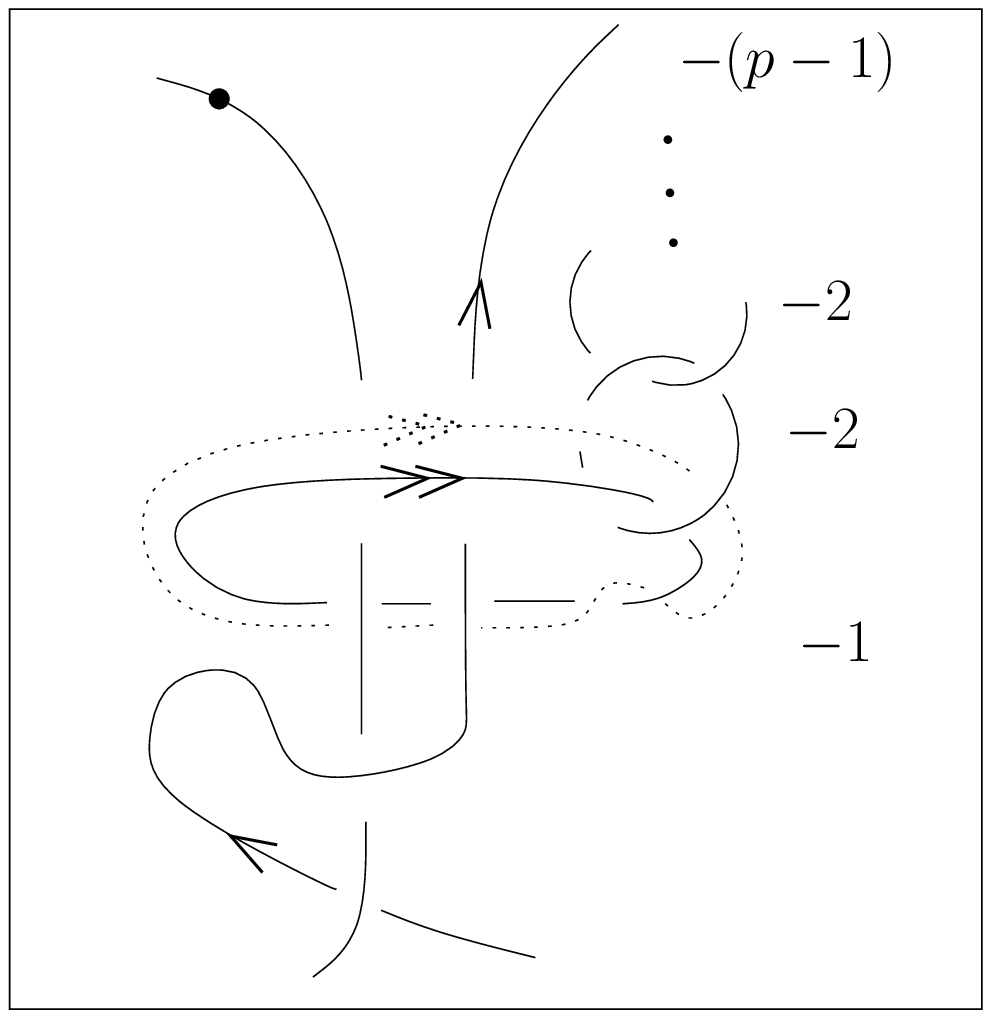}
\end{minipage}
\end{center}
\begin{center}
Figure \ref{rbdsec}.12
\end{center}

\begin{center}
\begin{minipage}{8cm}
\includegraphics[width=8cm]{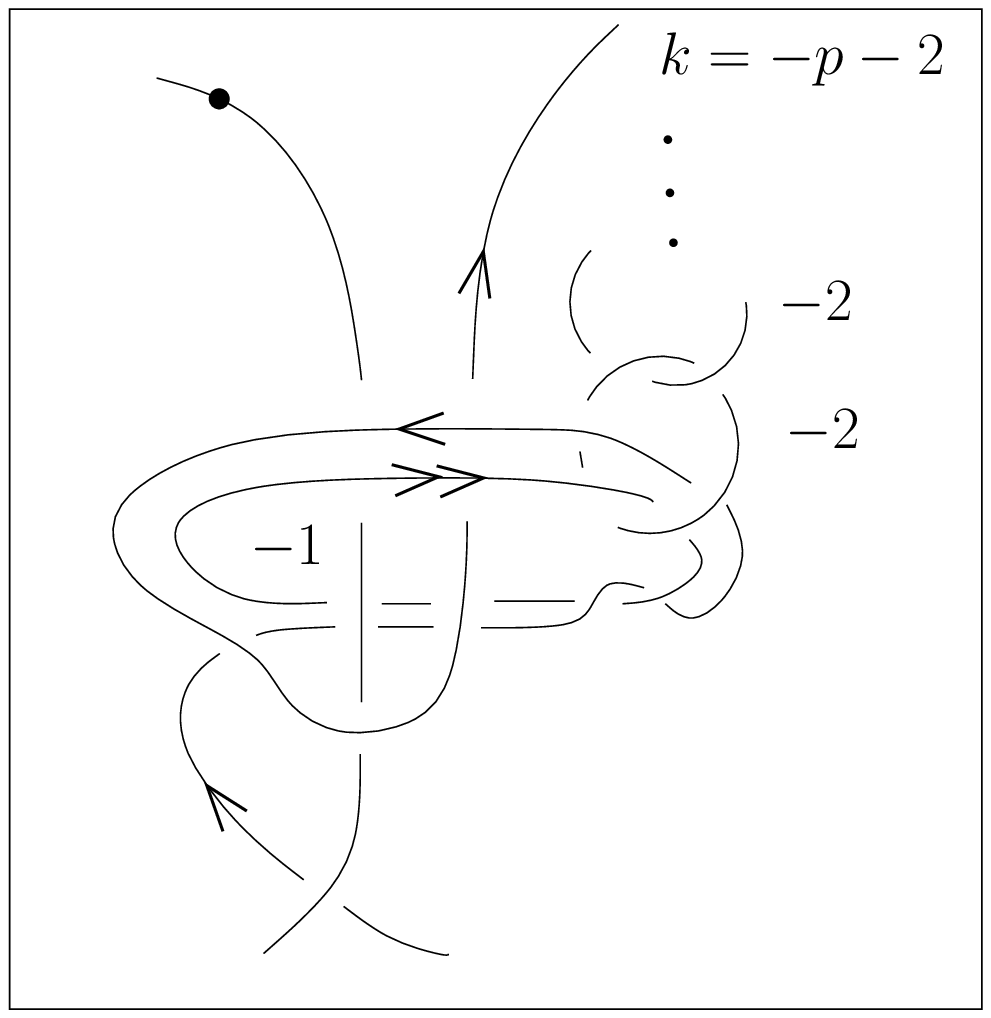}
\end{minipage}
\end{center}
\begin{center}
Figure \ref{rbdsec}.13
\end{center}

\begin{center}
\begin{minipage}{8cm}
\includegraphics[width=8cm]{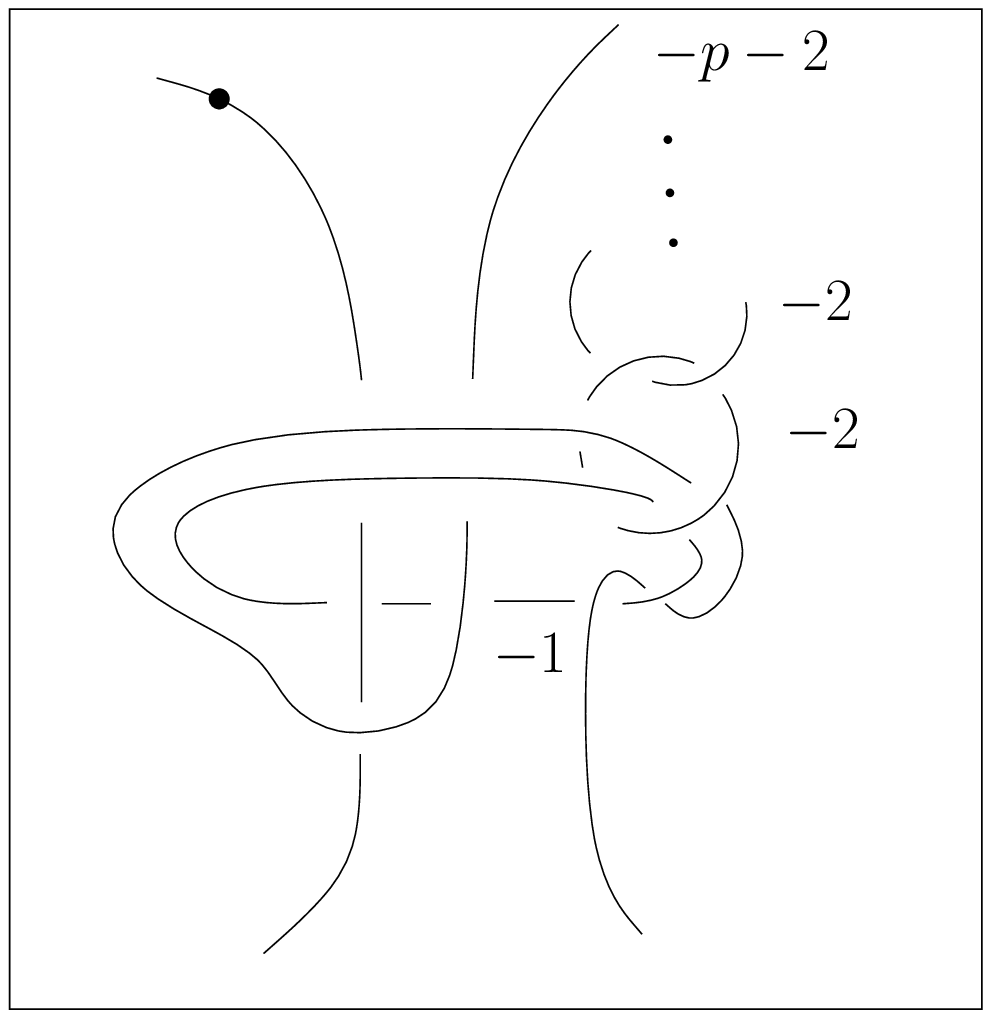}
\end{minipage}
\end{center}
\begin{center}
Figure \ref{rbdsec}.14
\end{center}

\begin{center}
\begin{minipage}{8cm}
\includegraphics[width=8cm]{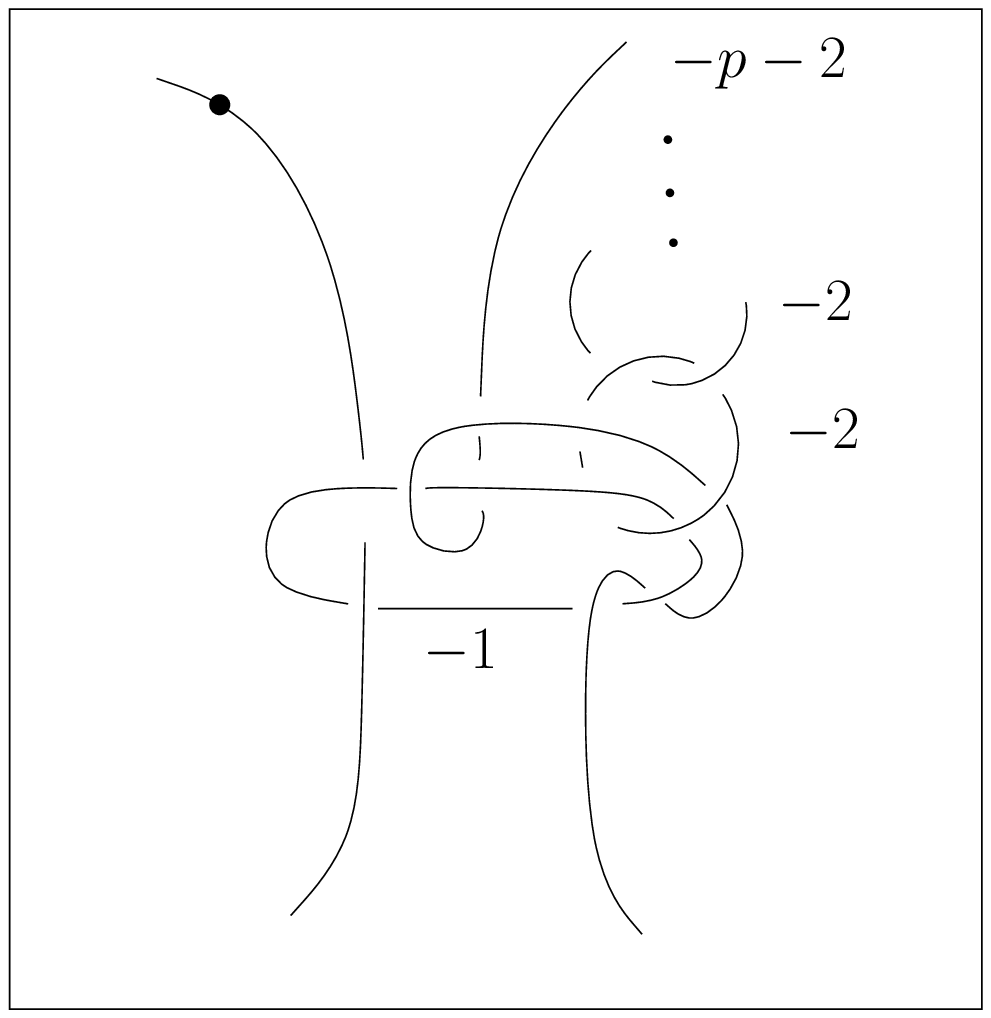}
\end{minipage}
\end{center}
\begin{center}
Figure \ref{rbdsec}.15
\end{center}

\begin{center}
\begin{minipage}{8cm}
\includegraphics[width=8cm]{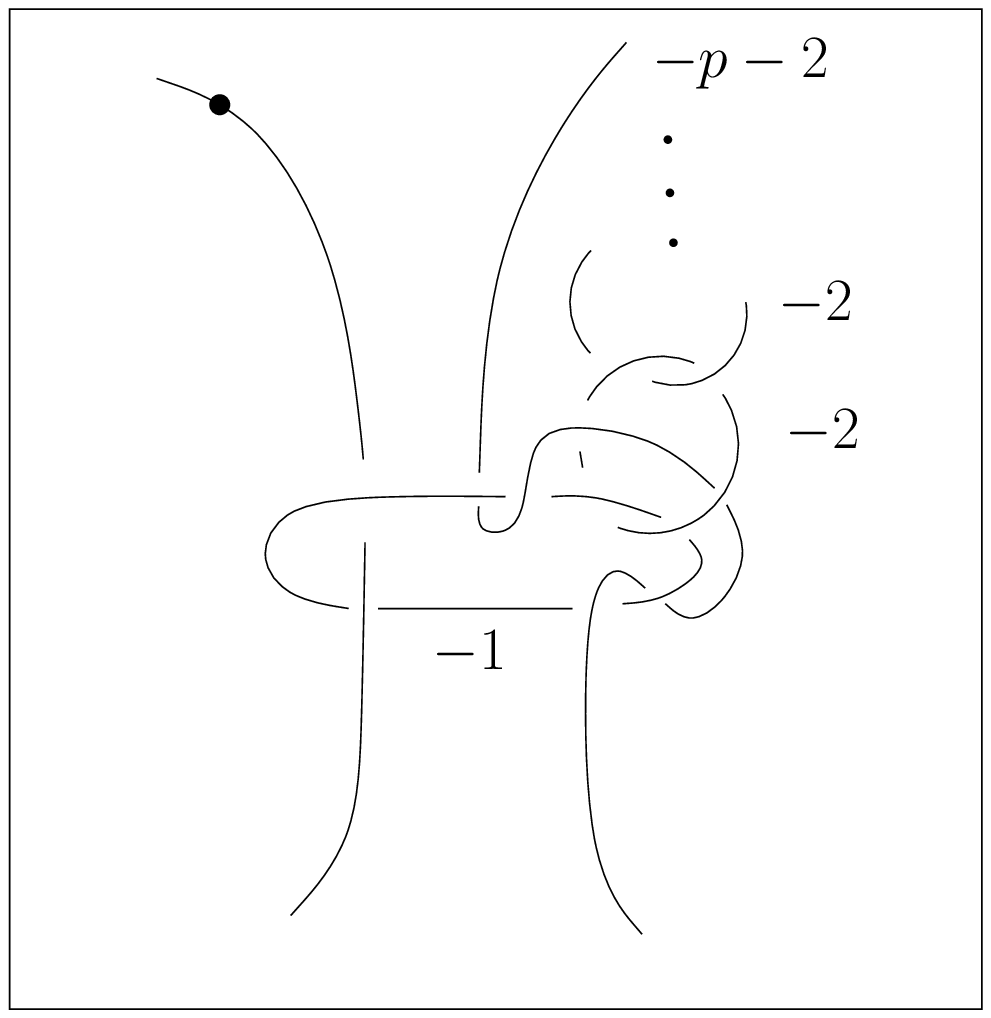}
\end{minipage}
\end{center}
\begin{center}
Figure \ref{rbdsec}.16
\end{center}

\begin{center}
\begin{minipage}{8cm}
\includegraphics[width=8cm]{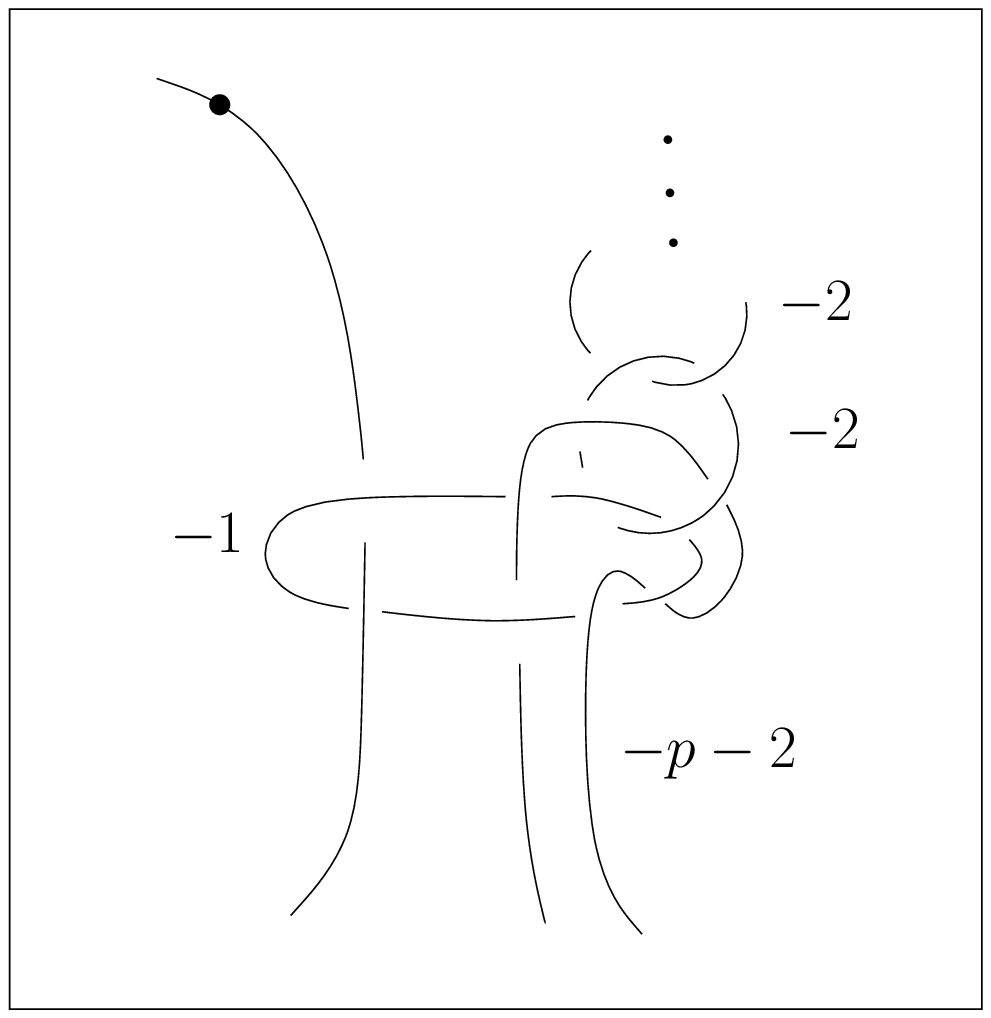}
\end{minipage}
\end{center}
\begin{center}
Figure \ref{rbdsec}.17
\end{center}

\begin{center}
\begin{minipage}{8cm}
\includegraphics[width=8cm]{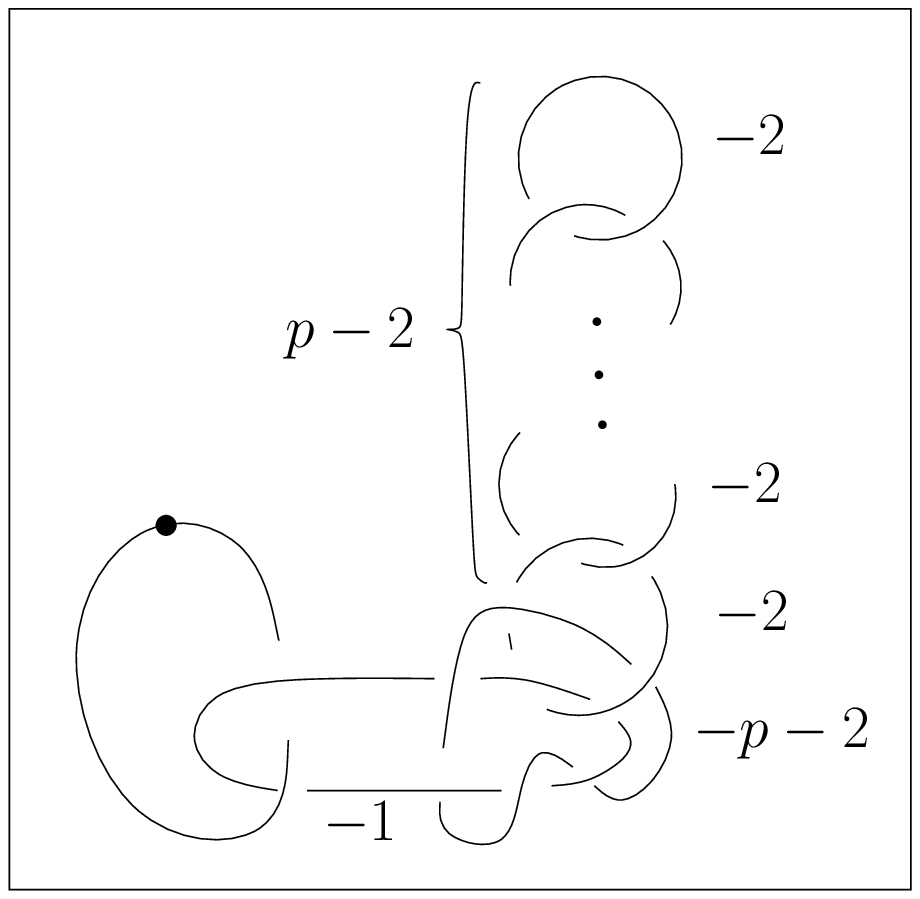}
\end{minipage}
\end{center}
\begin{center}
Figure \ref{rbdsec}.18
\end{center}

\begin{center}
\begin{minipage}{11cm}
\includegraphics[width=11cm]{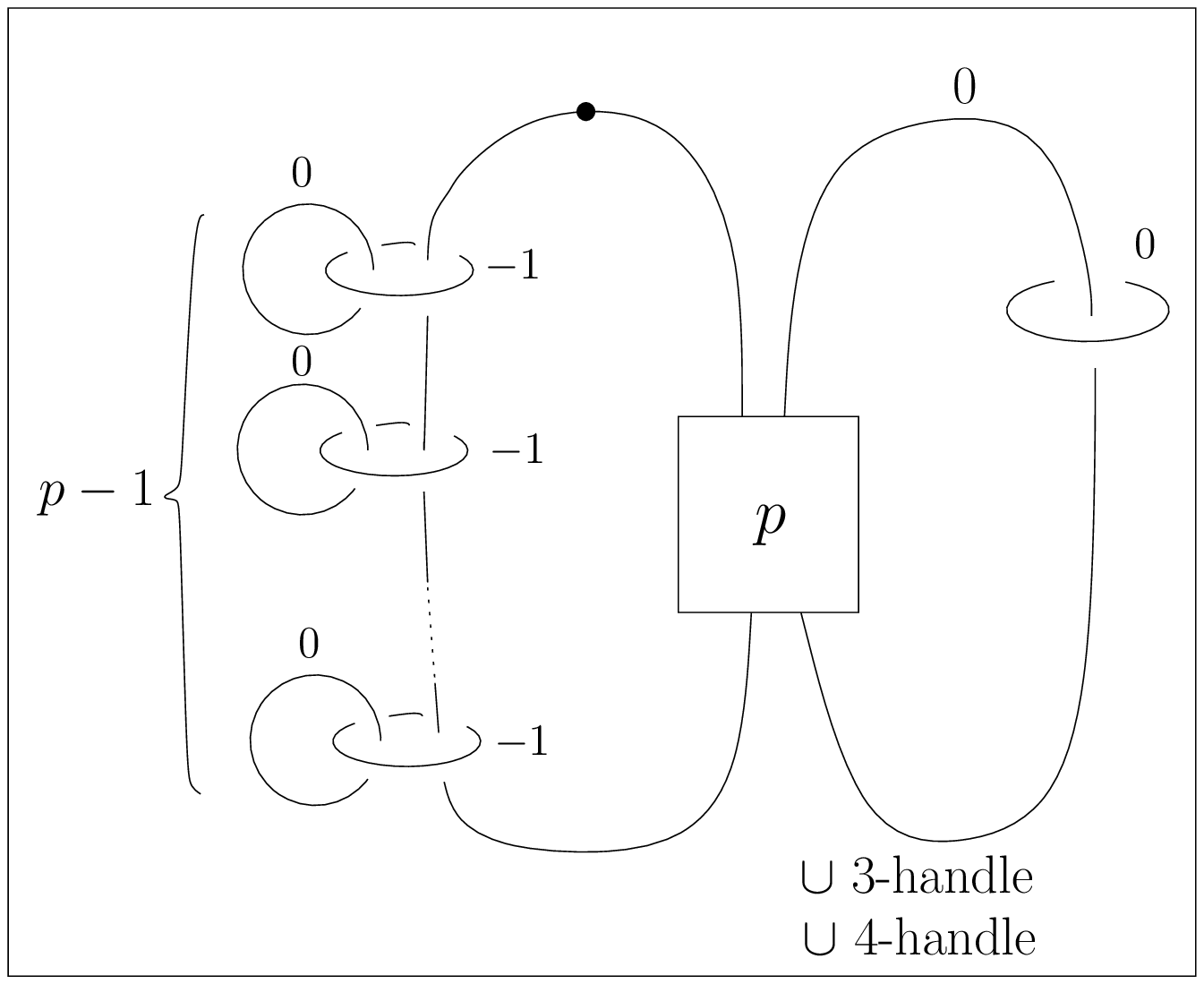}
\end{minipage}
\end{center}
\begin{center}
Figure \ref{rbdsec}.19
\end{center}

\begin{center}
\begin{minipage}{10cm}
\includegraphics[width=10cm]{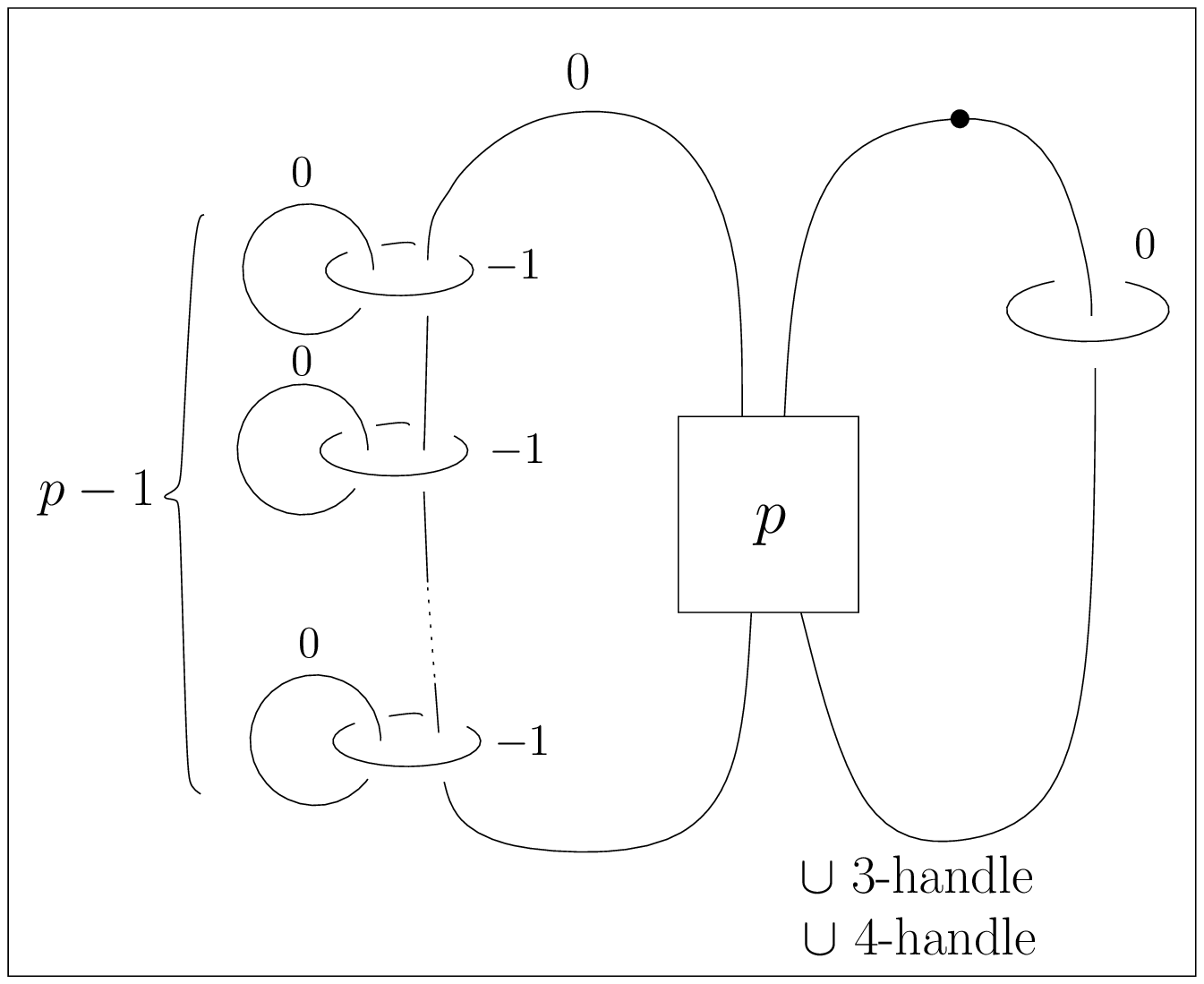}
\end{minipage}
\end{center}
\begin{center}
Figure \ref{rbdsec}.20
\end{center}

\begin{center}
\begin{minipage}{10cm}
\includegraphics[width=10cm]{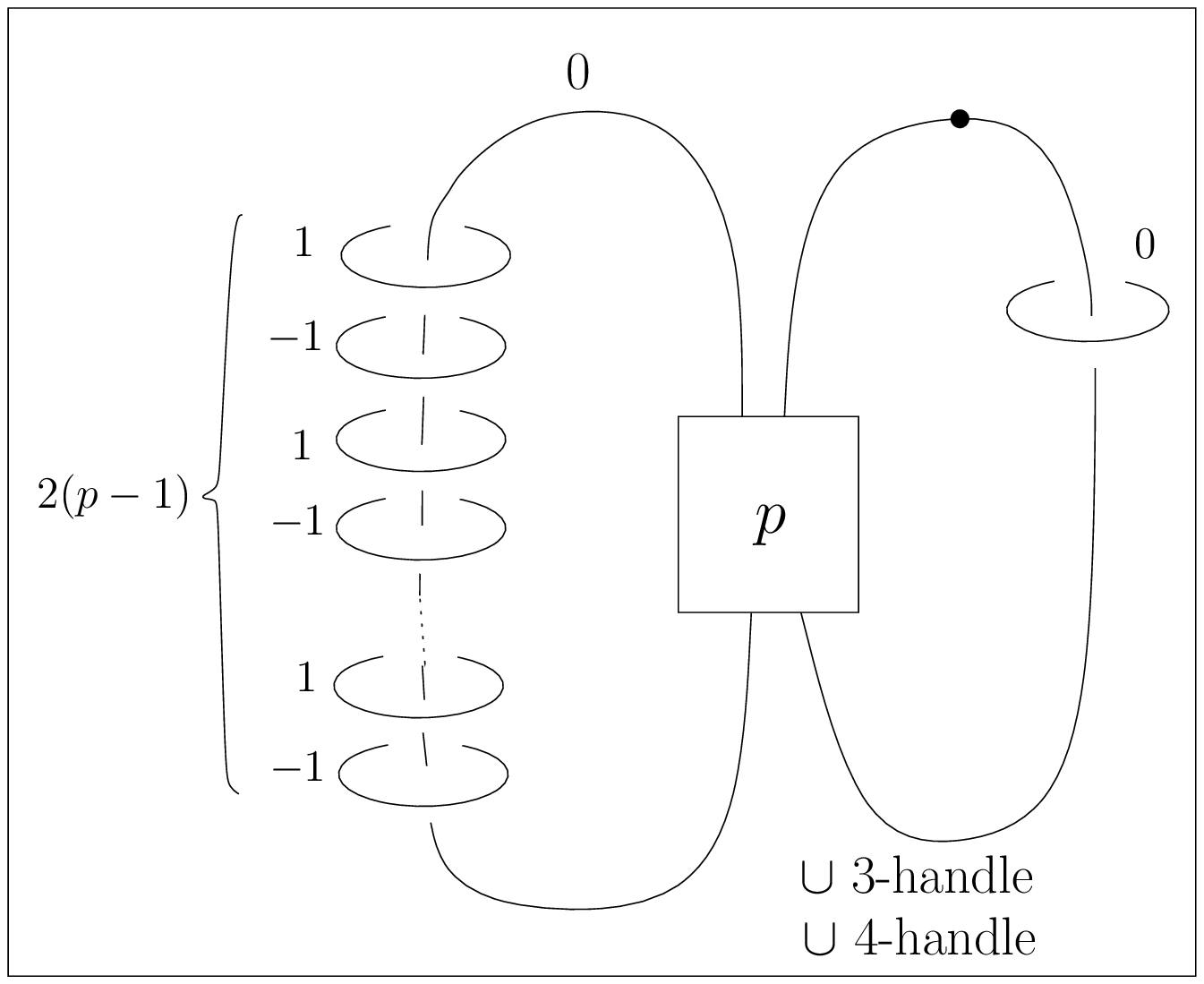}
\end{minipage}
\end{center}
\begin{center}
Figure \ref{rbdsec}.21
\end{center}

\begin{center}
\begin{minipage}{10cm}
\includegraphics[width=10cm]{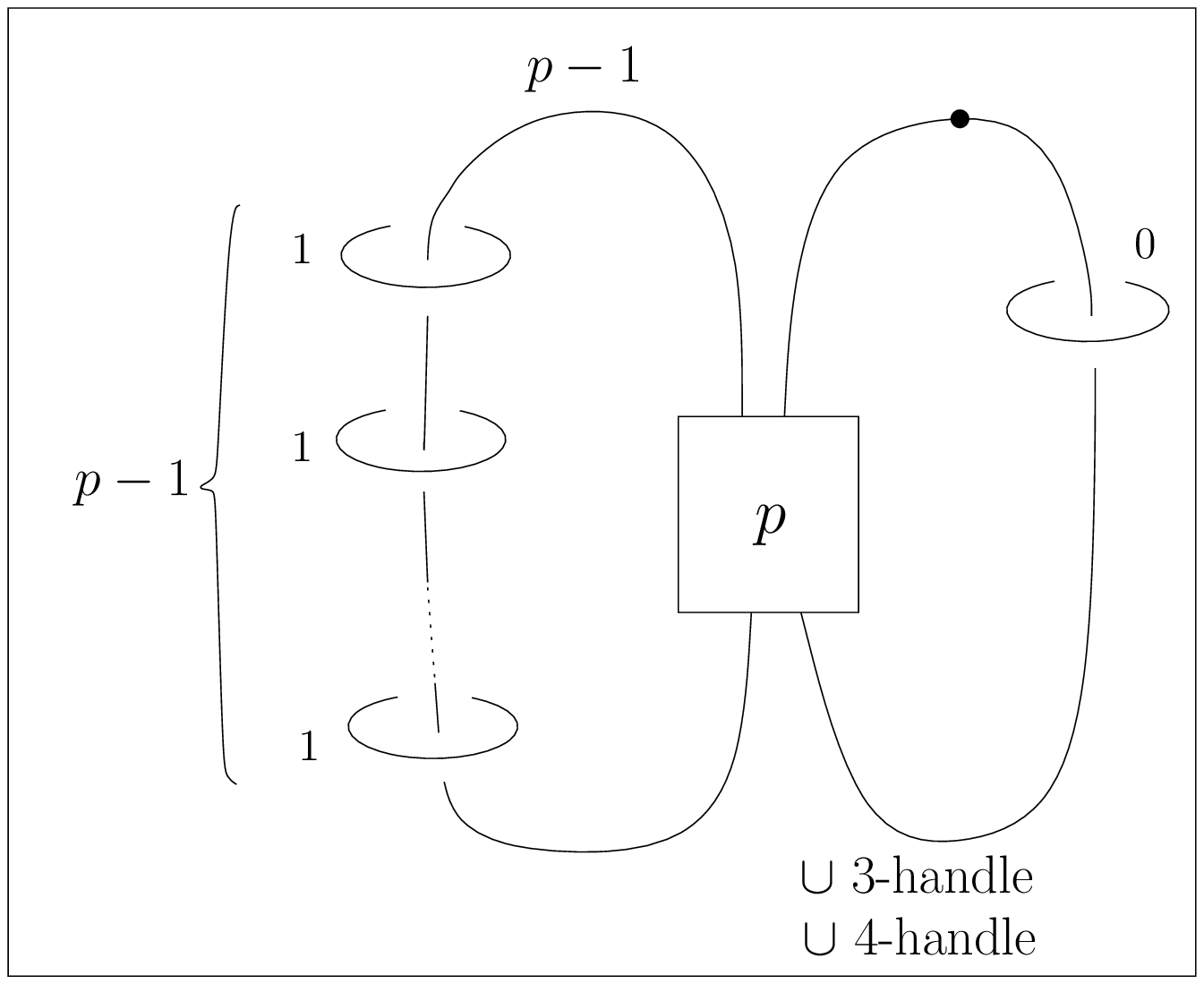}
\end{minipage}
\end{center}
\begin{center}
Figure \ref{rbdsec}.22
\end{center}

\begin{center}
\begin{minipage}{10cm}
\includegraphics[width=10cm]{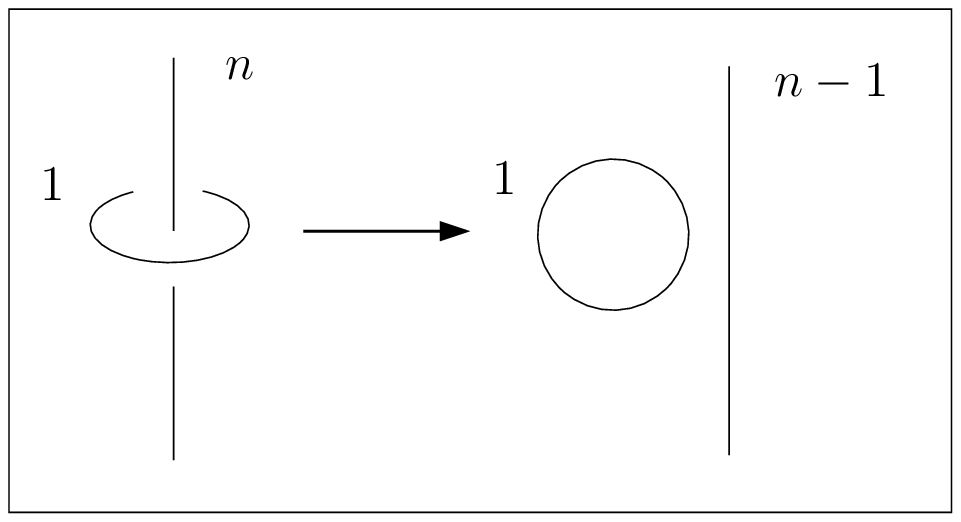}
\end{minipage}
\end{center}
\begin{center}
Figure \ref{rbdsec}.23
\end{center}

\begin{center}
\begin{minipage}{10cm}
\includegraphics[width=10cm]{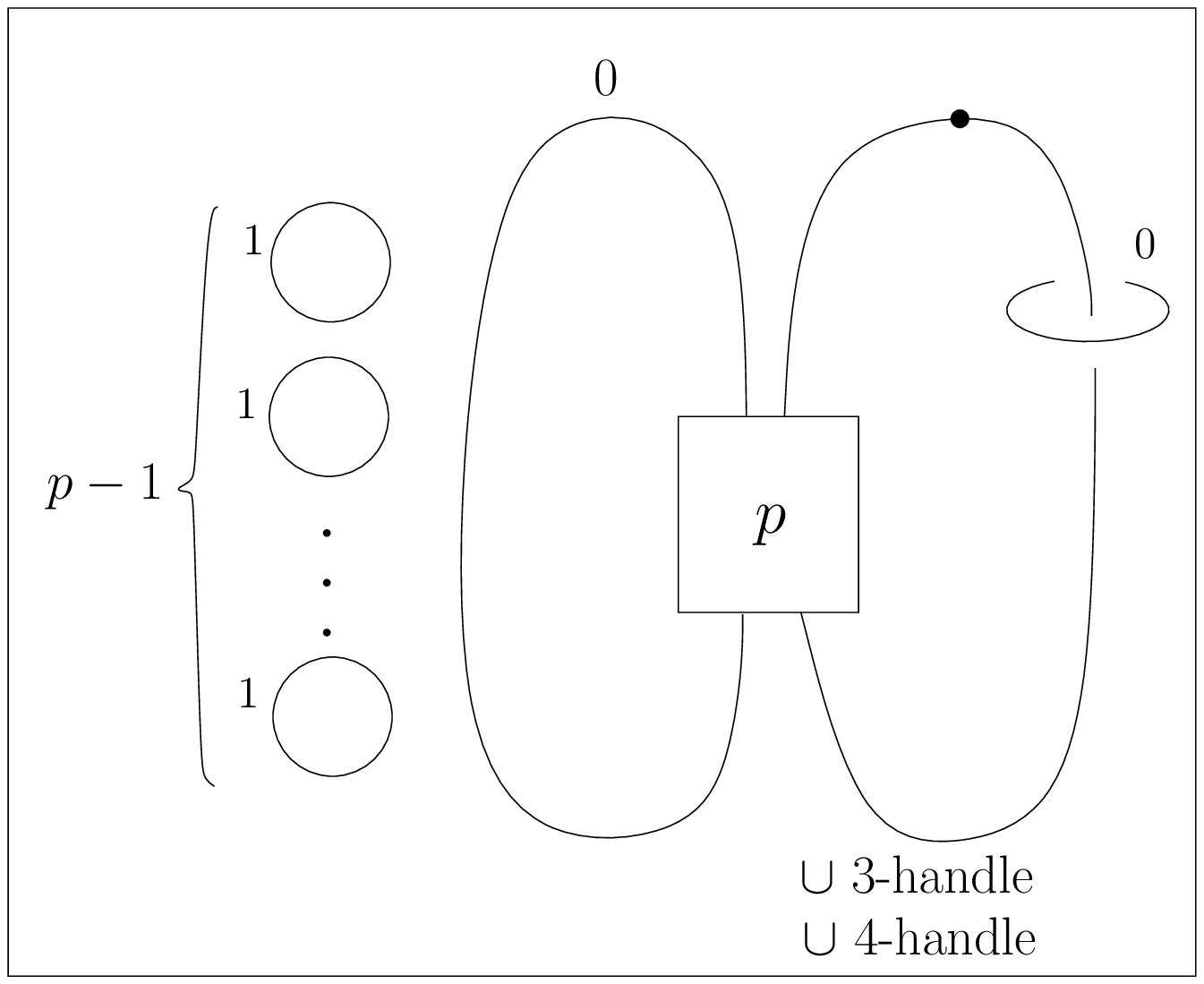}
\end{minipage}
\end{center}
\begin{center}
Figure \ref{rbdsec}.24
\end{center}

\begin{center}
\begin{minipage}{8cm}
\includegraphics[width=8cm]{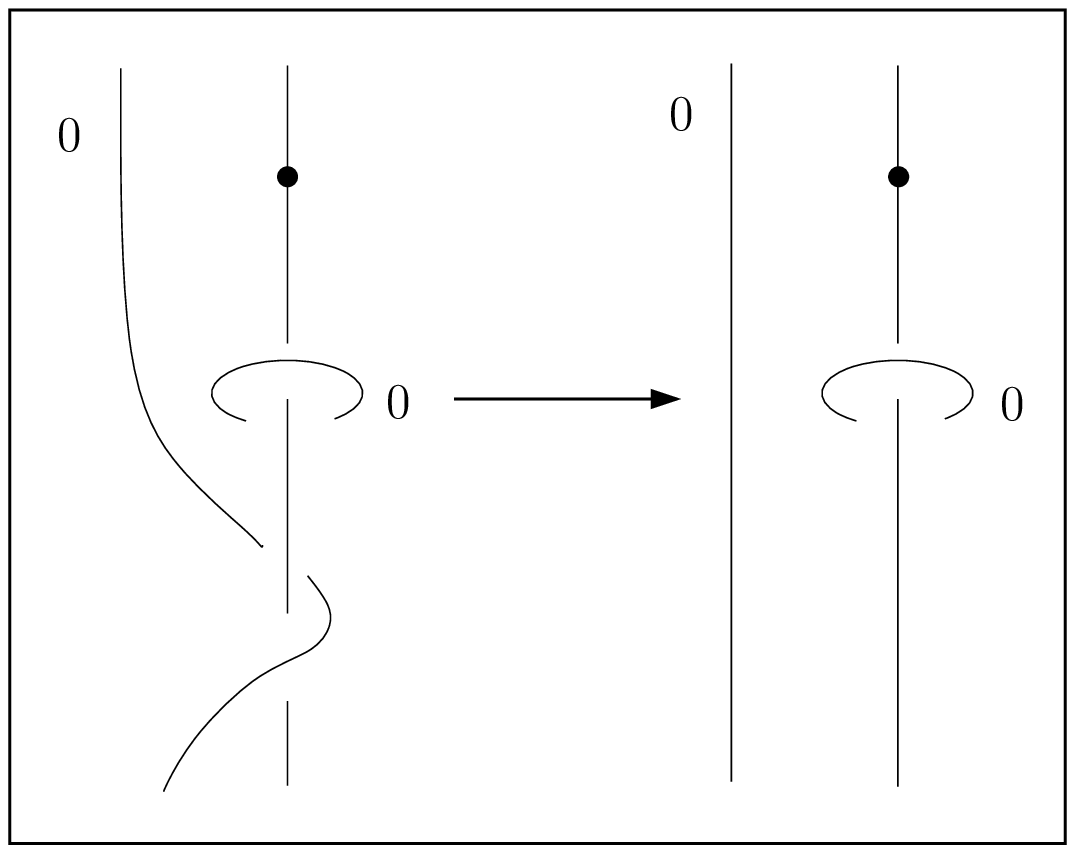}
\end{minipage}
\end{center}
\begin{center}
Figure \ref{rbdsec}.25
\end{center}

\begin{center}
\begin{minipage}{12cm}
\includegraphics[width=12cm]{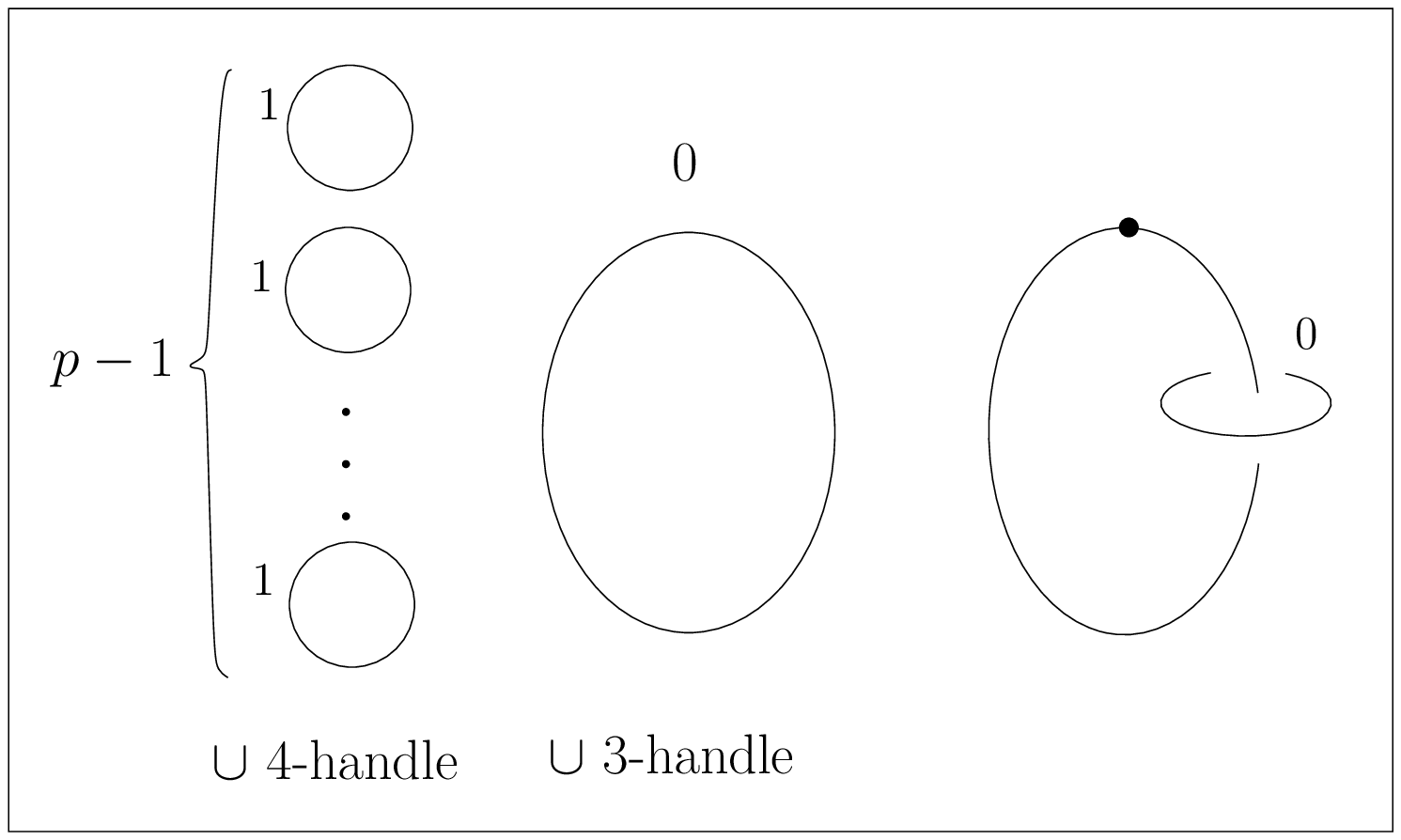}
\end{minipage}
\end{center}
\begin{center}
Figure \ref{rbdsec}.26
\end{center}

\begin{center}
\begin{minipage}{6cm}
\includegraphics[width=6cm]{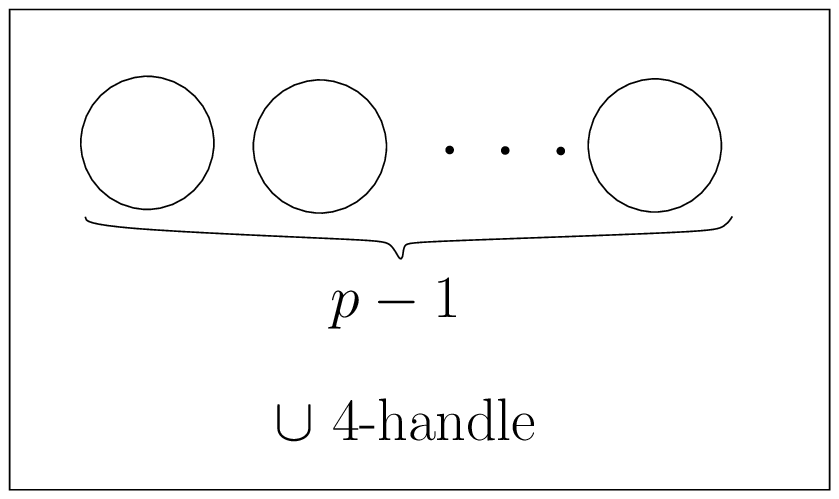}
\end{minipage}
\end{center}
\begin{center}
Figure \ref{rbdsec}.27
\end{center}

\begin{center}
\begin{minipage}{10cm}
\includegraphics[width=10cm]{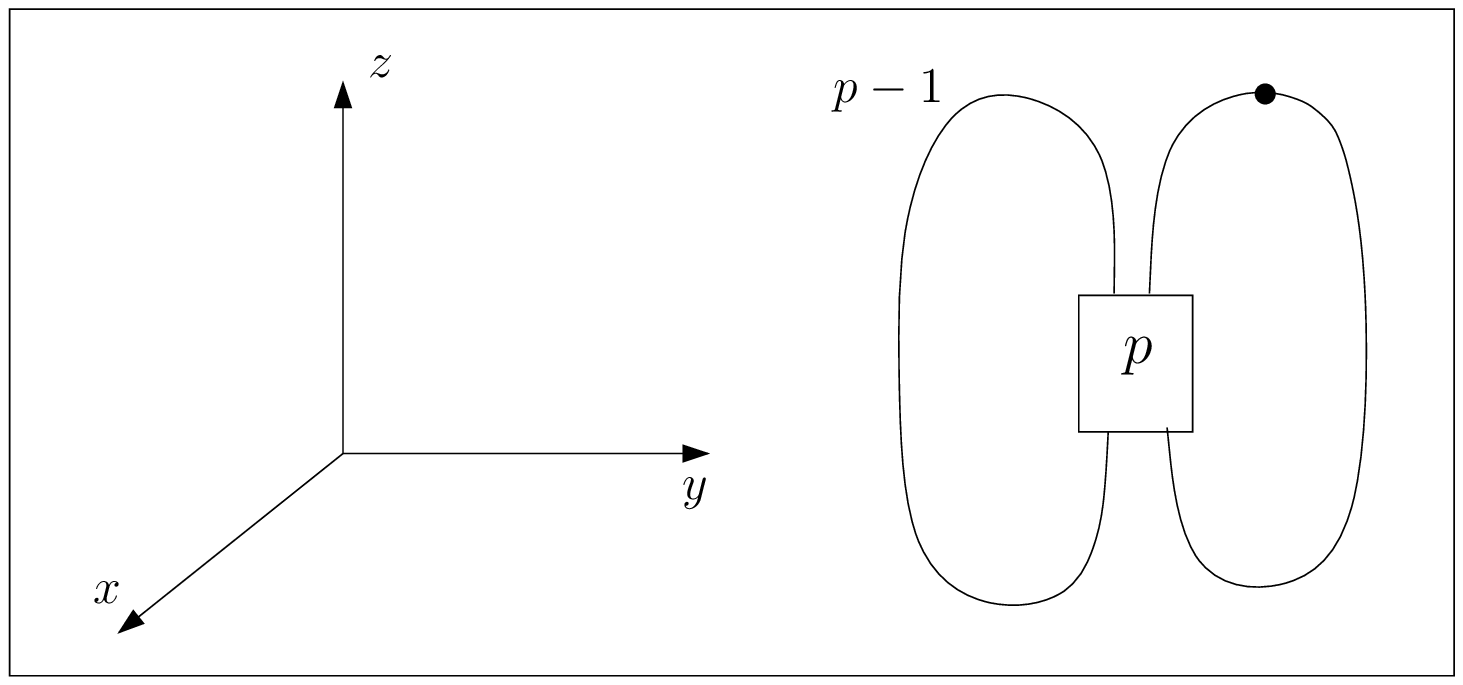}
\end{minipage}
\end{center}
\begin{center}
Figure \ref{rbdsec}.28
\end{center}

\pagebreak

\section{Results Concerning Rational Blowdowns} \label{x7proofsection}

In this section we prove that if $X_7$ is the rational blowdown of $E(1) \# 4 \cptbar$ along the configuration $C_7$, then $X_7$ is simply-connected.

We need the following theorem from \cite{OzSt} (see pages 37-43).
\newtheorem{ythm}{Theorem}[section]
\begin{ythm}
If $Y$ is a 3-manifold given by Dehn surgery along $(K_1, K_2, \dots, K_n) \subset S^3$ with surgery coefficients $\frac{p_i}{q_i}$ $(i = 1,2,\dots, n)$, then $H_1(Y; \mz)$ can be presented by the meridians $\mu_i$ as generators and the expressions 
\begin{equation*}
p_i \mu_i + q_i \sum_{j=1}^n \ell k(K_i,K_j) \mu_j = 0
\end{equation*}
as relations.
\end{ythm}
We first calculate $H_1(\partial C_7; \mz)$. We consider Figure \ref{x7proofsection}.1, which is a diagram of $C_7$ with the meridians $a_0, a_1, \dots, a_5$ drawn in. We have labelled $a_0$ as the meridian of the $-2$-sphere $S$ of square $-9$, and $a_1, \dots, a_5$ are the meridians of the $-2$-spheres $S_1, \dots, S_5$ in the $\tilde{E_6}$ singular fibre, respectively, where the spheres in Figure \ref{e6constrsec}.6 are given the labels as follows:
\begin{align*}
S_1:& \quad e_4-e_7\\
S_2:& \quad e_1-e_4\\
S_3:& \quad h-e_1-e_2-e_3\\
S_4:& \quad e_2-e_5\\
S_5:& \quad e_5-e_9\\
S_6:& \quad e_3-e_6\\
S_7:& \quad e_6-e_8 
\end{align*}

we can read off the relators (note that we use $a_i$ as a label for both the meridian and for the homology class the meridian represents):
\begin{align*}
r_0:& -9 a_0 + a_1 = 0 \\
r_1:& -2a_1 + a_0 + a_2 = 0 \\
r_2:& -2a_2 + a_1 + a_3 = 0 \\
r_3:& -2a_3 + a_2 + a_4 = 0 \\
r_4:& -2a_4 + a_3 + a_5 = 0 \\
r_5:& -2a_5 + a_4 = 0
\end{align*}
since integral surgery with coefficient $n$ is just rational surgery $\frac{p}{q}$ with $p=n$ and $q=1$, and so a presentation of the group is
\begin{equation*}
H_1(\partial C_7; \mz) \cong \; <a_0, a_1, a_2, a_3, a_4, a_5 | r_0, r_1, r_2, r_3, r_4, r_5>
\end{equation*}
where the relators are as above.

\newtheorem{a3thm}[ythm]{Lemma}
\begin{a3thm} \label{a3thmlabel}
$a_3$ is a generator of $\pi_1(\partial C_7)$. 
\end{a3thm}
Proof: \newline
First, we use the relators $r_i$ to write all the generators of $H_1(\partial C_7; \mz)$ in terms of $a_0$. We start with $r_0$:\newline
\begin{equation*}
r_0 \Rightarrow a_1 = 9a_0
\end{equation*}
Now, $r_1$ and $a_1 = 9a_0$ together give
\begin{align*}
a_2 &= 2a_1-a_0 \\
&= 18a_0 - a_0 \\
\Rightarrow a_2 &= 17a_0
\end{align*}
Similarly, the other relators become
\begin{align*}
a_3 &= 25a_0 \\
a_4 &= 33a_0 \\
a_5 &= 41a_0 \\
49a_0 &= 0
\end{align*}
So, $H_1(\partial C_7; \mz)$ can be presented as
\begin{equation} \label{preslabel1}
<a_0, a_1, a_2, a_3, a_4, a_5 | a_1 = 9a_0, a_2=17a_0, a_3=25a_0, a_4=33a_0, a_5=41a_0, 49a_0 =0>
\end{equation}
Since all the $a_i$ can be written in terms of $a_0$, this shows that $a_0$ is a generator. This presentation can be reduced to 
\begin{equation*}
<a_0 | 49a_0 =0> \cong \mz_{49}
\end{equation*}
We expected $H_1(\partial C_7; \mz) \cong \mz_{49}$ from Lemma \ref{cpzp2lemma}. We could also reduce $\eqref{preslabel1}$ to
\begin{equation*}
<a_0, a_3 | a_3 = 25a_0,  49a_0 =0>
\end{equation*}
which can now be reduced as follows:
\begin{align*}
&<a_0, a_3 | a_3 = 25a_0,  49a_0 =0> \\ 
\sim & <a_0, a_3 | a_3 = 25a_0, 49a_0 =0, 2a_3 = 50a_0> \\
\sim & <a_0, a_3 | a_3 = 25a_0, 2a_3 = a_0> \\
\sim & <a_0, a_3 | a_3 = 25a_0, 2a_3 = a_0, 50a_3 = 25a_0> \\
\sim & <a_0, a_3 | 2a_3=a_0, 50a_3 = a_3> \\
\sim & <a_0, a_3 | a_0 = 2a_3, 49a_3 = 0> \\
\sim & <a_3 | 49a_3 = 0>
\end{align*}
which shows that $a_3$ is a generator of $H_1(\partial C_7; \mz)$, as desired. \newline

Now, since $\partial C_7 \cong L(49,-6)$, we have $\pi_1(\partial C_7) \cong \pi_1(L(49,-6)) \cong \mz_{49}$. This shows that $\pi_1(\partial C_7)$ is abelian, which implies $\pi_1(\partial C_7) \cong H_1(\partial C_7; \mz)$. So, since $a_3$ is a generator of $H_1(\partial C_7; \mz)$, it is also a generator of $\pi_1(\partial C_7)$. $\Box$

\newtheorem{ythmrem}[ythm]{Remark}
\begin{ythmrem}
\upshape
It can be checked that actually every generator $a_i$ $(i=0, 1, \dots 5)$ is a generator of the group (although, we are only interested in the fact that $a_3$ is a generator). It should be noted that this does not always happen with every chain, and one reason it happens in this case is that $7$ is prime. For example, it can be checked that not every meridian in $C_6$ is a generator.
\end{ythmrem}

\newtheorem{x0thm}[ythm]{Proposition}
\begin{x0thm}
If we define $X = \cpt \# 13 \cptbar$, and write $X = X_0 \cup_{L(49,-6)} C_7$ so that $X_7 = X_0 \cup_{L(49,-6)} B_7$ is the rational blowdown along $C_7$, then $X_0$ is simply-connected.
\end{x0thm}

Proof: \newline
Firstly, we know $\pi_1(X) = 1$. Let us define $\Sigma =L(49,-6) \cong \partial C_7$. Then $X = X_0 \cup_{\Sigma} C_7$. We also have from Van Kampen's Theorem that if 
\begin{align*}
& j_0: \Sigma \hookrightarrow X_0 \\
& j_1: \Sigma \hookrightarrow C_7
\end{align*}
are the inclusion maps which induce the homomorphisms
\begin{align*}
& j_{0*}: \pi_1(\Sigma) \longrightarrow \pi_1(X_0) \\
& j_{1*}: \pi_1(\Sigma) \longrightarrow \pi_1(C_7)
\end{align*}

then $\pi_1(X) \cong (\pi_1(X_0) * \pi_1(C_7)) / N$, where $N$ is the normal subgroup generated by $j_{0*}(\omega) (j_{1*}(\omega))^{-1}$ for all $\omega \in \pi_1(\Sigma)$. \newline

Now, as Park observes in \cite{P1}, the generator $a_3$ of $\pi_1(\Sigma)$ intersects the $-2$-sphere labelled $S_6$ in the $\tilde{E_6}$-fibre. Note that $S_6$ is not in $C_7$, and so $S_6 \subset X_0$. In fact, $a_3$ intersects $S_6$ in such a way that it bounds a disk which is a hemisphere of $S_6$, and so $j_{0*}(a_3) = 1$. Therefore, by Lemma \ref{a3thmlabel}, $j_{0*}(\omega) = 1$ for all $\omega \in \pi_1(\partial C_7) \cong \pi_1(\Sigma)$. \newline

Therefore, the normal subgroup $N$ is just generated by all elements of the form $(j_{1*}(\omega))^{-1}$, or equivalently all elements of the form $j_{1*}(\omega) \in \pi_1(C_7)$, where $\omega \in \pi_1(\Sigma)$. \newline

Claim: $j_{1*}: \pi_1(\Sigma) \longrightarrow \pi_1(C_7)$ is a surjection. \newline

If this claim is true (and we shall prove it), then $N \cong \pi(C_7)$. Therefore,

\begin{equation*}
\pi_1(X) \cong (\pi_1(X_0) * \pi_1(C_7)) / N \cong \pi_1(X_0)
\end{equation*}

and so because $X$ is simply-connected, we have that $X_0$ is simply-connected. $\Box$ \newline

Proof of Claim (\cite{S}): \newline


Let $M$ be a 4-manifold that has only $0$-, $1$- and $2$-handles (no $3$-handles or $4$-handles). This manifold has non-empty boundary, since it does not have a 4-handle. So, $M$ has a handle decomposition consisting of a unique 0-handle, some 1-handles, some 2-handles, and a boundary $\partial M$. As for CW-complexes, 1-handles ``give'' generators for the $\pi_1(M)$ and 2-handles give relators for $\pi_1(M)$. \newline

We can look at $M$ ``upside-down'', by considering $k$-handles as $(4-k)$-handles. Then, $M$ has a handle decomposition with a unique 4-handle, some 3-handles, some 2-handles and a boundary $\partial M$. Since no generators of $\pi_1(M)$ come from $2$-, $3$- or $4$-handles, all the generators of $\pi_{1}(M)$ must come from its boundary $\partial M$. So, $\pi_1(M)$ has the same presentation as $\pi_1(\partial M)$, except it has extra relators coming from the $2$-handles of $M$. Therefore, there is a surjection $j_{*}: \pi_1(\partial M) \longrightarrow \pi_1(M)$. \newline

Since (for every $p$) $C_p$ is a manifold with no  $3$- or $4$-handles, the claim is proved for $C_7$. $\Box$ \newline

Finally, we need the following lemma:

\newtheorem{x7sclem}[ythm]{Lemma}
\begin{x7sclem}
If the configuration $C_{p}$ is embedded in $X$, and $X_p = (X \setminus C_p ) \cup_{L(p^2,p-1)} B_p$ is the rational blowdown of $X$ along $C_p$, and if $X$ and $X \setminus C_p$ are simply-connected, then $X_p$ is also simply connected.
\end{x7sclem}
Proof: \newline
Assume $X$ and $X \setminus C_p$ are both simply connected, and define $\Sigma = L(p^2,p-1)$. Then by Van Kampen's Theorem, 
\begin{align*}
\pi_1(X_p) &\cong (\pi_1(X \setminus C_p) * \pi_1(B_p) ) / N \\ 
&\cong \pi_1(B_p) ) / N
\end{align*} 
where $N$ is the normal subgroup generated by the elements $(j_{0*}(\omega))(j_{1*}(\omega))^{-1}$ for all $\omega \in \pi_1(\Sigma)$, where $j_{0*}$ and $j_{1*}$ are the homomorphisms induced from the inclusion maps
\begin{align*}
& j_0: \Sigma \hookrightarrow X \setminus C_p \\
& j_1: \Sigma \hookrightarrow B_p
\end{align*}
as before. \newline

However, since we assume $X \setminus C_p$ is simply connected, $j_{0*}$ is the trivial map, and so $N$ is generated by elements of the form $j_{1*}(\omega)$ for $\omega \in \pi_{1}(\Sigma)$, as before. By \cite{FS1}, this map $j_{1*}$ is surjective, and so $\pi_1(X_p)$ is trivial. $\Box$

\newtheorem{x7sclemr1}[ythm]{Remark}
\begin{x7sclemr1}
\upshape
In particular, $X_7$ is simply-connected.
\end{x7sclemr1}

\begin{center}
\begin{minipage}{10cm}
\includegraphics[width=10cm]{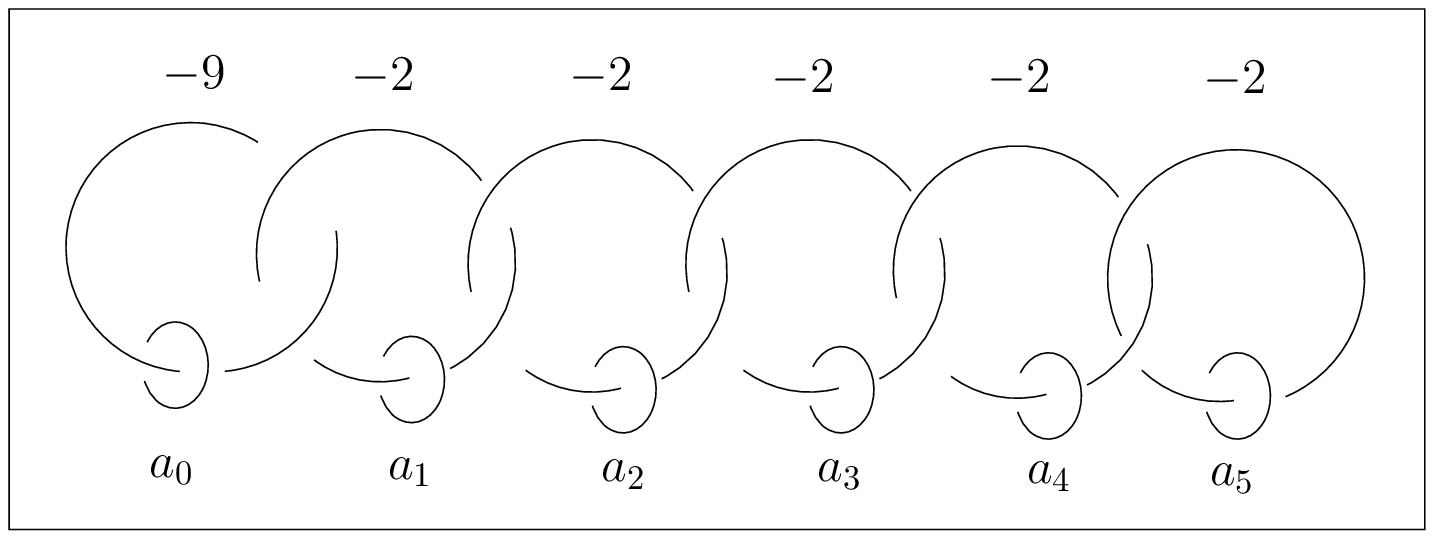}
\end{minipage}
\end{center}
\begin{center}
Figure \ref{x7proofsection}.1
\end{center}

\pagebreak

\section{Complex surfaces}

Recall that the complex projective line $H= \{[x:y:z] \in \cpt | \; x=0 \}$ defines a generator $h = [H] \in H_2(\cpt; \mz)$.

\newtheorem{cs41}{Proposition}[section]
\begin{cs41} \label{propdh}
The set $D = \{ [x:y:z] \in \cpt | \; x^{d} + y^{d} + z^{d} =0 \}$ is a smooth, connected submanifold of $\cpt$ representing $dh \in H_2(\cpt; \mz)$.
\end{cs41}

Since this result is so fundamental to many of our calculations later, we give a detailed proof, which comes from \cite{GS}. \newline

Proof: \newline
We first calculate how many points are in the intersection 
\begin{align*}
D \cap H &= \{[x:y:z] \in \cpt | x^{d} + y^{d} + z^{d} =0, x =0  \} \\
&=  \{[x:y:z] \in \cpt | y^{d} + z^{d} =0 \}
\end{align*}

We notice that if $y = 0$ then $z = 0$, and vice-versa, and then the point is $[0:0:0]$, which is not a point in $\cpt$. So, $y$ and $z$ are both non-zero. If we divide through by $z^d$, the equation $y^{d} + z^{d} =0$ becomes
\begin{equation*}
\alpha^d + 1 =0
\end{equation*}
where $\alpha = \frac{y}{z}$ is a non-zero complex number, and this equation has $d$ distinct solutions 
\begin{equation*}
\alpha_{k} = e^{i(\pi / d + 2k\pi / d)}; k = 0, 1, \dots, d-1
\end{equation*}
and so we get the $d$ points in $\cpt$
\begin{equation*}
D \cap H  = \{ [0: \alpha_k: 1] \in \cpt | \; \alpha_{k} = e^{i(\pi / d + 2k\pi / d)}, k = 0, 1, \dots, d-1 \}
\end{equation*}

If we define $g(x,y,z) = x^d + y^d + z^d$, then the Implicit Function Theorem says that 
\begin{equation*}
\tilde{D} = \{(x,y,z) \in \mathbb{C}^3 | \; g(x,y,z) = 0 \}
\end{equation*}
is a smooth manifold if $(g_{x}(x,y,z), g_y(x,y,z), g_z(x,y,z)) \neq (0,0,0)$ on every point in $\tilde{D}$ (where $g_x$ means $\frac{\partial g}{\partial x}$). \newline
\pagebreak

However, we are working with $D$, a subset of $\cpt$ and not $\tilde{D}$, which is a subset of $\mathbb{C}^3$. We therefore need to check the affine charts of $D$ (for example, the set $\{[x:y:1] \in \cpt\}) = \{(x,y) \in \mathbb{C}^2 \}$). \newline

We start by checking the chart $[x:y:1]$. Let $p(x,y) = g(x,y,1) = x^d + y^d + 1$. Then
\begin{align*}
p_x(x,y) &= dx^{d-1} \\
p_y(x,y) &= dy^{d-1}
\end{align*}
The only point $(x_0, y_0)$ that gives $(p_x(x_0,y_0), p_y(x_0,y_0)) = (0,0)$ is the point $(0,0)$. Fortunately, $p(0,0) = 1$, and so $(0,0)$ is not in the set $\{(x,y) \in \mathbb{C}^2 | p(x,y) = 0 \}$, and so $D$ is a smooth manifold in this chart. Similar calculations show that $D$ is also smooth in the charts $[x:1:z]$ and $[1:y:z]$, and so $D$ is a smooth submanifold of $\cpt$.\newline

Since $D$ and $H$ are both complex submanifolds of $\cpt$, each of their intersection points is a positive intersection, and so $Q_{\cpt}<[D],[H]> = d$, which shows that $D$ represents the homology class $dh \in H_2(\cpt; \mz)$. To see that $D$ is connected, notice that $g(x,y,z) = x^d + y^d + z^d$ is an irreducible polynomial, and so its zero set must be a connected component (see \cite{GH}). $\Box$ \newline

We now have the following interesting proposition, the proof of which can be found in \cite{GS}.

\newtheorem{s42}[cs41]{Proposition}
\begin{s42}
If $p_1$ and $p_2$ are two homogeneous polynomials with equal degree (and not powers of other polynomials) and the hypersurfaces $F_1 = \{P \in \mathbb{CP}^{n} |\; p_1(P) = 0 \}$ and $F_2 = \{P \in \mathbb{CP}^{n} |\; p_2(P) = 0 \}$ are smooth submanifolds of $\mathbb{CP}^{n}$, then $F_1$ is diffeomorphic to $F_2$.
\end{s42}

\newtheorem{s42r}[cs41]{Remark}
\begin{s42r}
\upshape
This proposition with Proposition \ref{propdh} above shows that any degree $d$ homogeneous polynomial is represented by $dh \in H_2(\cpt; \mz)$.
\end{s42r}

\pagebreak

\section{Resolving Singular Points} \label{resolvesection}

This section comes from Chapter 2 of \cite{GS}, and we follow the notation used there. \newline

Assume that we have two smooth surfaces $\Sigma_1$ and $\Sigma_2$ in $\cpt$ that are both closed and oriented, and assume that $\Sigma_1$ and $\Sigma_2$ intersect each other transversally in the single point $P \in \cpt$. Then $\Sigma = \Sigma_1 \cup \Sigma_2$ is not a smooth surface in $\cpt$, since at the point $P$ it fails to be a manifold, but it still defines a homology class $[\Sigma] = [\Sigma_1] + [\Sigma_2] \in H_2(\cpt; \mz)$. \newline

We now describe a process that changes $\Sigma$ into a smooth surface $\tilde{\Sigma}$. We consider a neighbourhood of the intersection point $P$ that is a 4-ball, which we call $D$. Inside $D$, a neighbourhood of the intersection point looks like
\begin{equation*}
F = \{(z_1,z_2) \in \mathbb{C}^2 \st z_1 z_2 = 0, |z_1|^2 + |z_2|^2 \leq 1 \}
\end{equation*}
which is a model for two 2-dimensional disks (note, $z_1, z_2 \in \mathbb{C}$) intersecting each other in a single point in the 4-ball $D$. To ``remove'' the singular intersection point at $P$, we cut out the pair $(D,F)$ and replace it with a pair $(D,R)$ that does not have a singular point at $P$, does not change the manifold $\cpt$, and also is such that the homology class of $[\tilde{\Sigma}] = [\Sigma]$. \newline

We choose $R$ to be the subset of $D$ that is obtained by perturbing the subset
\begin{equation*}
R'_{\epsilon} = \{(z_1,z_2) \in \mathbb{C}^2 | z_1 z_2 = \epsilon, |z_1|^2 + |z_2|^2 \leq 1\}
\end{equation*}
where $(0 < |\epsilon| \ll 1)$, so that $\partial R = \partial F \subset \partial D $. Since $R'_{\epsilon}$ is the graph of $z_2 = \frac{\epsilon}{z_1}$, it is topologically an annulus, and therefore $R$ is also topologically an annulus. \newline 

So, replacing the pair $(D,F)$ with the pair $(D,R)$ ``removes'' the singular point $P$, but since we are simply removing $D$ and gluing it back in, it does not change the $\cpt$ which contains $\Sigma_1$ and $\Sigma_2$. Furthermore, since the subsets $F$ and $R$ are homologous in $(D,\partial D)$, the homology class of $\tilde{\Sigma}$ (what $\Sigma$ becomes after this operation) is still $[\tilde{\Sigma}] = [\Sigma_1] + [\Sigma_2]$.

\pagebreak

\newtheorem{resrem1}{Remark}[section]
\begin{resrem1}
\upshape
The process above, as \cite{GS} say, ``removes the singular point $P$''. Since we worked locally around $P$, the method is valid for every 4-manifold $X$, and every pair of intersecting surfaces in $X$, even if the intersections are transverse self-intersections.
\end{resrem1}

\newtheorem{resrem2}[resrem1]{Remark}
\begin{resrem2}
\upshape
A nice way of looking at this operation is as follows: We have two surfaces $\Sigma_1$ and $\Sigma_2$ in some 4-manifold $X$ that intersect each other in a point $P$. We take a disk neighbourhood of $P$ in $\Sigma_1$, call it $D_1$, and a disk neighbourhood of $P$ in $\Sigma_2$, call it $D_2$. We remove the disks $D_1$ and $D_2$ and replace them with an annulus joining $\partial D_1$ to $\partial D_2$. This is exactly the operation of connect-summing the two surfaces together, from which it can be seen that the homology of the resulting surface $[\tilde{\Sigma}]$ is clearly $[\tilde{\Sigma}] = [\Sigma_1] + [\Sigma_2]$.
\end{resrem2}

\newtheorem{resrem3}[resrem1]{Remark}
\begin{resrem3} \label{resolvesympl}
\upshape
Although we have done this operation \emph{smoothly}, it is possible to resolve singular points \emph{symplectically}. The method is the same, except we use a function such as
\begin{equation*}
f(x) = 
\begin{cases}
\exp (-(\frac{1}{x})^2) & \text{if $x \neq 0$}, \\
0 & \text{if $x = 0$}
\end{cases}
\end{equation*}
to symplectically ``smooth corners''.
\end{resrem3}

\pagebreak
\section{Elliptic Fibrations} \label{ellfibsec}

This section is based on Chapter 3 of \cite{GS}. The definitions are from \cite{GS} and we follow the notation presented there. 

\newtheorem{esd}{Definition}[section]
\begin{esd}
\upshape
A complex surface $S$ is called an \emph{elliptic surface} if there is a holomorphic map $\pi: S \longrightarrow C$ from $S$ to a complex curve $C$ such that for generic $t \in C$ the inverse image $\pi^{-1}(t)$ is a smooth \emph{elliptic curve}. We call the map $\pi$ a \emph{(holomorphic) elliptic fibration}.
\end{esd}

\newtheorem{esdr1}[esd]{Remark}
\begin{esdr1}
\upshape
Recall that an elliptic curve is topologically a real 2-dimensional torus (\cite{GS}).
\end{esdr1}

\newtheorem{esd2}[esd]{Definition}
\begin{esd2} \label{correcteddefnlabel}
\upshape
Let $X$ be a smooth, closed, oriented 4-manifold, and let $C$ be a complex curve. A smooth map $\pi: X \longrightarrow C$ will be called a (\emph{$C^{\infty}$-}) \emph{elliptic fibration} if each fibre $\pi^{-1}(t)$ (which may be a singular fibre) has a neighbourhood $U \subset X$ and an orientation preserving diffeomorphism $\phi: U \longrightarrow \phi(U)$, where $\phi(U)$ is a subset of an elliptic surface $S$, such that $\phi$ commutes with the maps $\pi$.
\end{esd2}

We present three examples from \cite{GS} in order to illustrate what an elliptic fibration actually is.

\newtheorem{cp1eg}[esd]{Example}
\begin{cp1eg} \label{cp2n1eg}
\upshape
We first construct a $\mathbb{CP}^1$-fibration over $\mathbb{CP}^1$. We consider all the complex projective lines in $\cpt$ passing through the point $P = [0:0:1] \in \cpt$. To each line passing through $P$ we can associate a point $[t_0:t_1] \in \mathbb{CP}^1$ such that the line is $L_{[t_0:t_1]} = \{[x:y:z] \in \cpt \st t_0 x = t_1 y \}$. Essentially, we are picking the point in $\cpt$ where the line though $P$ crosses the projective line $\{ [x:y:z] \in \cpt \st z=0\} \cong \cpo$, and this association of a line through $P$ to a point $[t_0:t_1] \in \cpo$ parametrizes the set of such lines. \newline

It is easy to see that the family of lines $\{L_{[t_0:t_1]} \st [t_0:t_1] \in \cpo \}$ is a one-sheet cover of $\cpt \setminus \{P\}$, i.e. for every point $Q \in \cpt \setminus \{P\}$, there is a unique line in the family that passes through $Q$. Therefore, we can define a map $f: \cpt \setminus \{P\} \longrightarrow \cpo$ as follows: for $Q \in \cpt \setminus \{P\}$, there is a unique line $L_{[t_{0_Q}:t_{1_Q}]}$ in the above family that passes through $Q$, and we define $f(Q) = [t_{0_Q}:t_{1_Q}]$. \newline

We notice that all the lines $L_{[t_0:t_1]}$ intersect each other transversally in $P$, and so we cannot extend this map to all of $\cpt$. However, blowing up a point $P$ in a manifold essentially replaces the point $P$ with the set of lines going through that point, so we can extend $f$ to $\cpt \# \cptbar$, and therefore we have $\tilde{f}: \cpt \# \cptbar \longrightarrow \cpo$, a $\cpo$-fibration of $\cpt \# \cptbar$ over $\cpo$.
\end{cp1eg}

\newtheorem{cp2eg}[esd]{Example}
\begin{cp2eg} \label{cp2n2eg}
\upshape
We generalize the construction above. Above, the polynomials that define the lines passing through $P = [0:0:1]$ are linear polynomials (namely, $p_0(x,y,z) = x$ and $p_1(x,y,z) = y$). \newline

Instead of using linear polynomials, suppose we choose $p_0$ and $p_1$ to be quadratic (and homogeneous) polynomials in the variables $x,y,z$. Suppose further that we choose $p_0$ and $p_1$ ``generically enough'', so that their zero sets,
\begin{align}
V_{p_0} = \{[x:y:z] \in \cpt  \st p_0(x,y,z) = 0\}\\ 
V_{p_1} = \{[x:y:z] \in \cpt  \st p_1(x,y,z) = 0\}
\end{align}
which are the curves in $\cpt$ corresponding to the polynomials $p_0$ and $p_1$, are such that $V_{p_0}$ intersects $V_{p_1}$ in four points $P_1, P_2, P_3, P_4$. \newline

To see what is meant by choosing polynomials ``generically enough'', let us look at quadratic polynomials in two real variables. Figure \ref{ellfibsec}.1 shows how we could choose two quadratic polynomials that intersect in only two points, while Figure \ref{ellfibsec}.2 shows how we could choose two quadratic polynomials that intersect each other in four points. \newline

It should be remarked that we often identify a polynomial with the curve to which it corresponds (its zero set), in order to give meaning to the phrase of how ``polynomials intersect each other''. \newline

So, to recap, we have two quadratics which give curves $C_0$ and $C_1$ (the zero sets $V_{p_0}$ and $V_{p_1}$) which intersect each other in four points $P_1, P_2, P_3, P_4$. We now consider the family of polynomials
\begin{equation*}
\mathcal{Q} = \{t_0 p_0 + t_1 p_1 \st [t_0:t_1] \in \cpo \}
\end{equation*}
Such a family is called a \emph{pencil of curves}. \newline

Again, we blur the distinction of a polynomial and the curve to which a polynomial corresponds. So, we consider an element of $\mathcal{Q}$ to be both a curve and a polynomial, depending on context. \newline

As in Example \ref{cp2n1eg} above, this family $\mathcal{Q}$ gives a one-sheeted cover of $\cpt \setminus \{P_1, P_2, P_3, P_4 \}$, and any two curves in $\mathcal{Q}$ intersect each other transversally in every point of $\{P_1, P_2, P_3, P_4 \}$ (since for $i=1,2,3,4$, $t_0 p_0 (P_i) + t_1 p_1 (P_i) = 0+0 = 0$). Therefore, we can define the map
\begin{equation*}
f: \cpt \setminus \{P_1, P_2, P_3, P_4 \} \longrightarrow \cpo
\end{equation*}
and although we cannot extend this map to the whole of $\cpt$, we can blow up at the points $P_1, P_2, P_3, P_4$ to define a map
\begin{equation*}
\tilde{f}: \cpt \# 4 \cptbar \longrightarrow \cpo
\end{equation*}
\end{cp2eg}

\newtheorem{cp2egr1}[esd]{Remark}
\begin{cp2egr1} 
\upshape
Although $\tilde{f}$ in Example \ref{cp2n1eg} was a bundle map, in Example \ref{cp2n2eg} it is not. First of all, we do have generic fibres (a generic quadric curve is irreducible, and so since ``a generic quadric curve in $\cpt$ is a copy of $\cpo$'' (\cite{GS}), a generic fibre is $\cpo$). \newline

However, there are singular fibres as well, which correspond to those ``non-generic'' quadratic polynomials that are reducible, and so we get fibres that are the union of two lines (e.g. $x^2 + y^2 = 0$ gives two lines, $x = iy$ and $x=-iy$), which is not simply a copy of $\cpo$. \newline

The reason we did not encounter this problem in Example \ref{cp2n1eg} is that all linear polynomials are irreducible, and so, as \cite{GS} puts it ``there are no singular linear subspaces of $\cpt$.'' \newline

So, $\tilde{f}: \cpt \# 4 \cptbar \longrightarrow \cpo$ is not a fibre bundle, since there is the possibility that two different fibres may be non-diffeomorphic. However, we still call $\tilde{f}: \cpt \# 4 \cptbar \longrightarrow \cpo$ a (singular) fibration.
\end{cp2egr1}

\newtheorem{cp2egr2}[esd]{Remark}
\begin{cp2egr2} 
\upshape
Recall from Example \ref{cp2n2eg} that every polynomial $p_{[t_0:t_1]}$ in the pencil $\mathcal{Q} = \{t_0 p_0 + t_1 p_1 \st [t_0:t_1] \in \cpo \}$ has the property that $p_{[t_0:t_1]}(P) = 0$ for each $P \in \{P_1, P_2, P_3, P_4 \}$. Since each fibre $\tilde{f}^{-1}([t_0:t_1])$ is simply a curve $p_{[t_0:t_1]}$ corresponding to the point $[t_0:t_1] \in \cpo$, and each point $P_i$ was blown-up to become the exceptional sphere $E_{i}$ ($i = 1,2,3,4$), we observe that each exceptional sphere intersects each fibre in a unique point. Therefore, the exceptional spheres of the blow-ups are sections of the fibration $\tilde{f}: \cpt \# 4 \cptbar \longrightarrow \cpo$.
\end{cp2egr2}

\newtheorem{cp3eg}[esd]{Example}
\begin{cp3eg} 
\upshape
We now start with a generic pair of cubics, $p_0$ and $p_1$ that intersect each other in 9 points $\{P_1, \dots, P_9 \}$ (for an example of a generic pair of cubics in two real variables see Figure \ref{ellfibsec}.3, for a non-generic pair of cubics, see Figure \ref{ellfibsec}.4). \newline

We consider the pencil of curves $\mathcal{Q} = \{t_0 p_0 + t_1 p_1 \st [t_0:t_1] \in \cpo \}$ and we define a map
\begin{equation*}
f: \cpt \setminus \{P_1, \dots, P_9 \} \longrightarrow \cpo
\end{equation*}
as before, i.e. for a point $Q \in \cpt \setminus \{P_1, \dots, P_9 \}$ there is a unique cubic of the form $p_{[t_0:t_1]} = t_0 p_0 + t_1 p_1$ that passes through $Q$, and we define $f(Q) = [t_0:t_1] \in \cpo$. \newline

We blow up at each of the 9 points $P_1, \dots, P_9$ and extend $f$ to a fibration 
\begin{equation*}
\pi: \cpt \# 9 \cptbar \longrightarrow \cpo
\end{equation*}
whose fibres are cubic curves, and so generic fibre a smooth elliptic curve, which is topologically a torus. Therefore, we have just shown that there is a holomorphic elliptic fibration on $\cpt \# 9 \cptbar$.
\end{cp3eg}

\newtheorem{cp3egr1}[esd]{Remark}
\begin{cp3egr1}
\upshape
If a fibration $\pi: X \longrightarrow \cpo$ has only generic fibres (fibres that are tori), then since we know the Euler characteristic of a torus is $\chi(T^2) = 0$, this would imply that $\chi(X) =0$. \newline

However, we know $H_{2}(\cpt; \mz) \cong \mz$ and $\cpt$ is simply connected, so by the argument at the beginning of section \ref{class4msec} and the definition of the Euler characteristic we have 
\begin{equation*}
\chi(\cpt) = \sum_{i=0}^{4}H_{i}(\cpt; \mz) = 1+0+1+0+1= 3 
\end{equation*}
So $\chi(\cpt \# 9 \cptbar) = 3 + 9 = 12 \neq 0$, and so there must be fibres that are not diffeomorphic to the torus. These fibres are the \emph{singular fibres} we have mentioned, and the possible singular fibres will be discussed in the next section. Note that our choice of polynomials $p_0$ and $p_1$ decides which types of singular fibres will occur in the elliptic fibration.
\end{cp3egr1}

\newtheorem{cp3egr2}[esd]{Remark}
\begin{cp3egr2}
\upshape
When we consider $\cpt \# 9 \cptbar$ as being equipped with an elliptic fibration, we denote it by $E(1)$.
\end{cp3egr2}

We note that this process of constructing (singular) fibrations can be extended for higher-order polynomials:
\newtheorem{cp4thm}[esd]{Lemma}
\begin{cp4thm}
\upshape \cite{GS} \itshape The manifold $\cpt \# d^2 \cptbar$ admits a (singular) fibration $\cpt \# d^2 \cptbar \longrightarrow \cpo$, where the generic fibre is a complex curve of genus $\frac{1}{2}(d-1)(d-2)$.
\end{cp4thm}

\subsection*{The general case}
The following is taken from \cite{SSS}. \newline 

Suppose that $p_0$ and $p_1$ are two homogeneous polynomials of degree $n$ in the variables $x,y,z$. Suppose furthermore that $\mathrm{gcd}(p_0,p_1) = 1$, which means that the curves $C_0, C_1 \subset \cpt$ defined by $p_0$ and $p_1$ do not share a common component, and intersect each other only in finitely many points. We call the family of polynomials 
\begin{equation*}
p_t = \{ t_0 p_0 + t_1 p_1 \st t = [t_0:t_1] \in \cpo \} 
\end{equation*}
the \emph{pencil} generated by $p_0$ and $p_1$. The zero-set of $p_t$,
\begin{equation*}
C_t = \{[x:y:z]\in \cpt \st p_t(x,y,z) =0 \} \subset \cpt
\end{equation*}
is a complex curve. We define the set
\begin{equation*}
B = \{P \in \cpt \st p_0(P) = p_1(P) = 0 \}
\end{equation*}
and $B$ is called the set of \emph{base points} of the pencil. By the assumption above, this set must have only finitely many points. \newline

As in the examples above, we notice that the map
\begin{equation*}
P \mapsto [p_0(P):p_1(P)]
\end{equation*}
defines a map from $\cpt \setminus B$ to $\cpo$, and after blowing up at the base points of the pencil, this map extends to a well-defined holomorphic map from $\cpt \# k \cptbar$ to $\cpo$, where $k$ is the number of points in $B$.


\pagebreak

\begin{center}
\begin{minipage}{6cm}
\includegraphics[width=6cm]{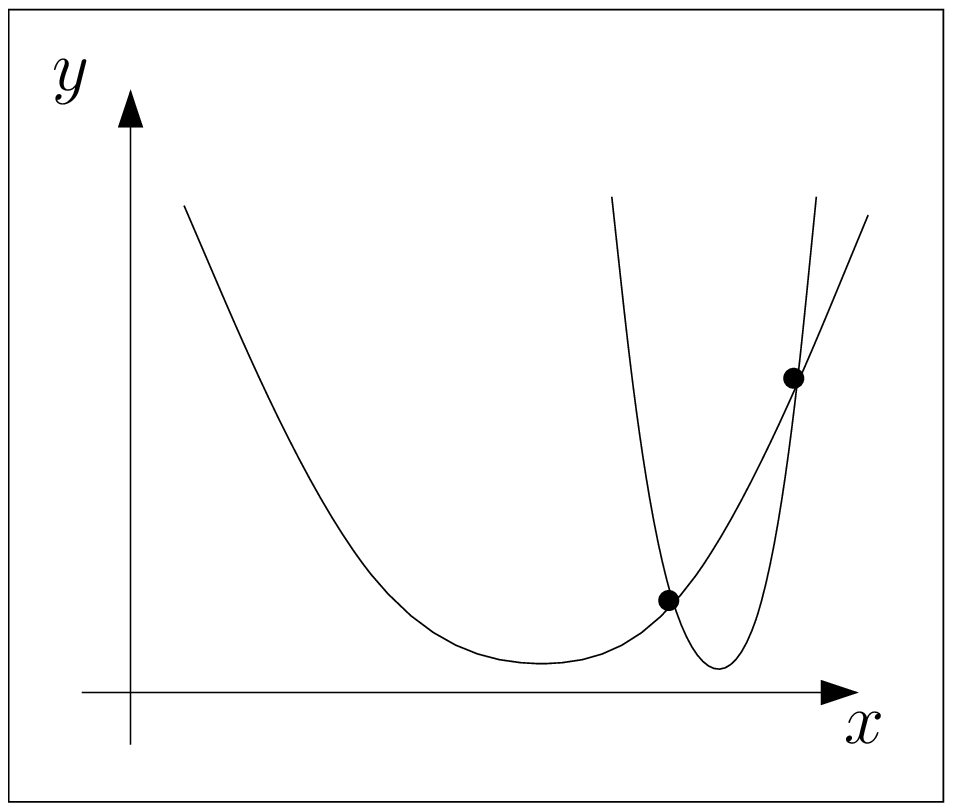}
\end{minipage}
\end{center}
\begin{center}
Figure \ref{ellfibsec}.1
\end{center}

\begin{center}
\begin{minipage}{6cm}
\includegraphics[width=6cm]{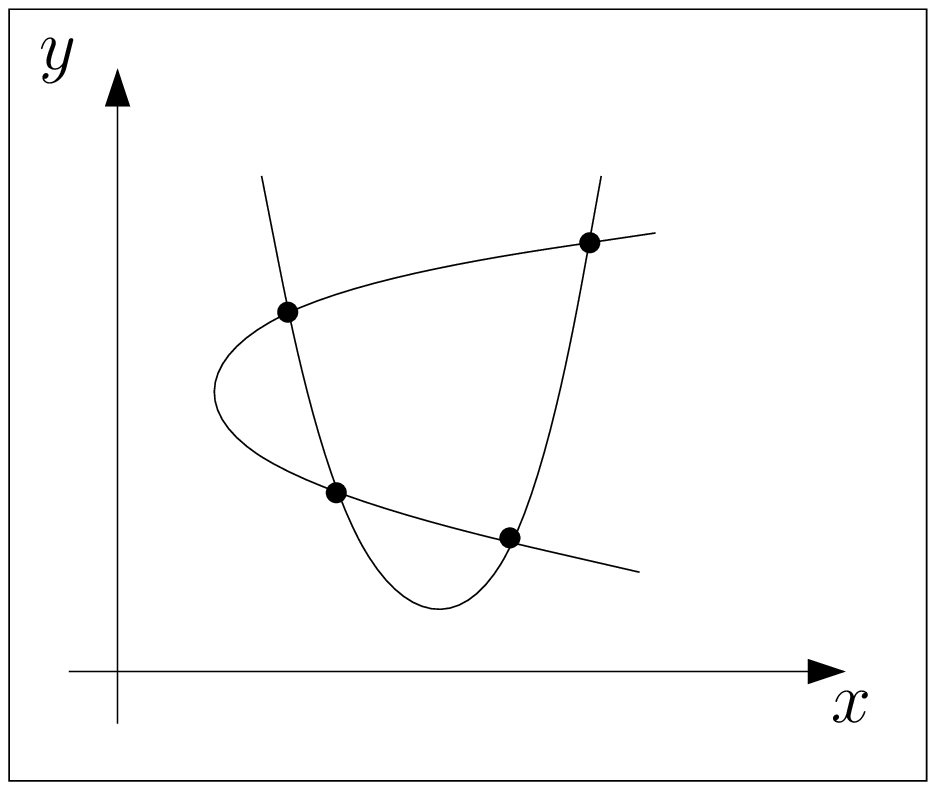}
\end{minipage}
\end{center}
\begin{center}
Figure \ref{ellfibsec}.2
\end{center}

\begin{center}
\begin{minipage}{6cm}
\includegraphics[width=6cm]{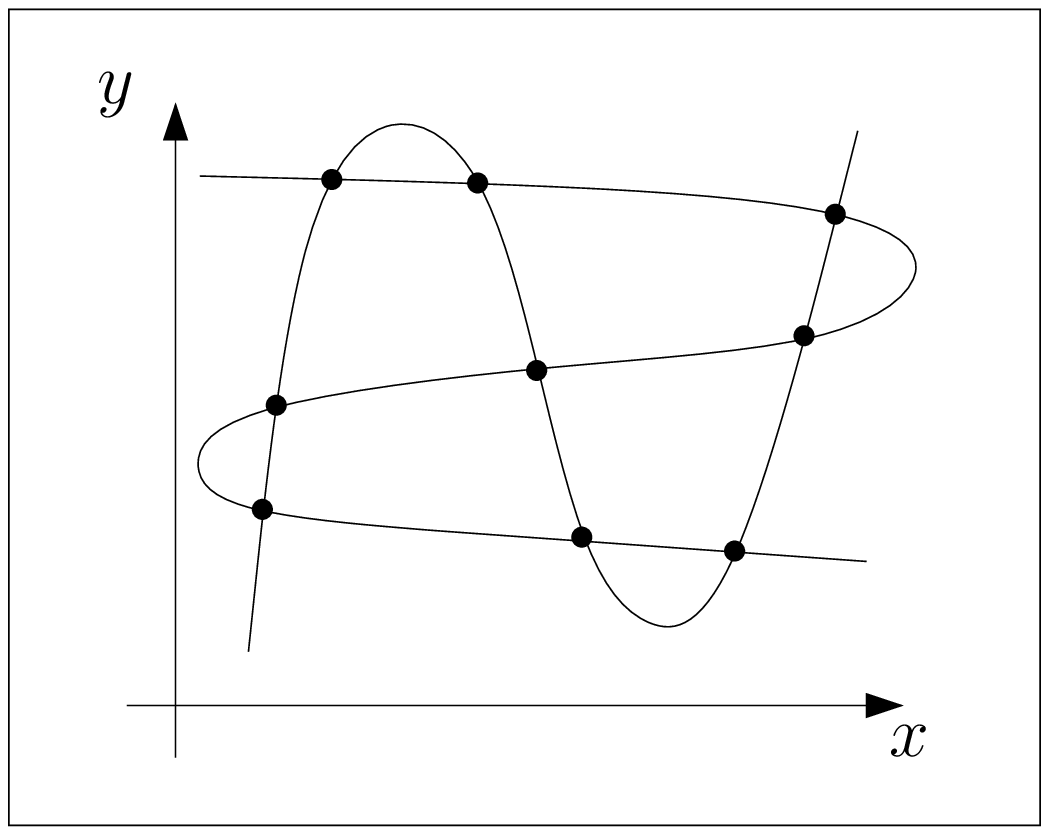}
\end{minipage}
\end{center}
\begin{center}
Figure \ref{ellfibsec}.3
\end{center}

\begin{center}
\begin{minipage}{6cm}
\includegraphics[width=6cm]{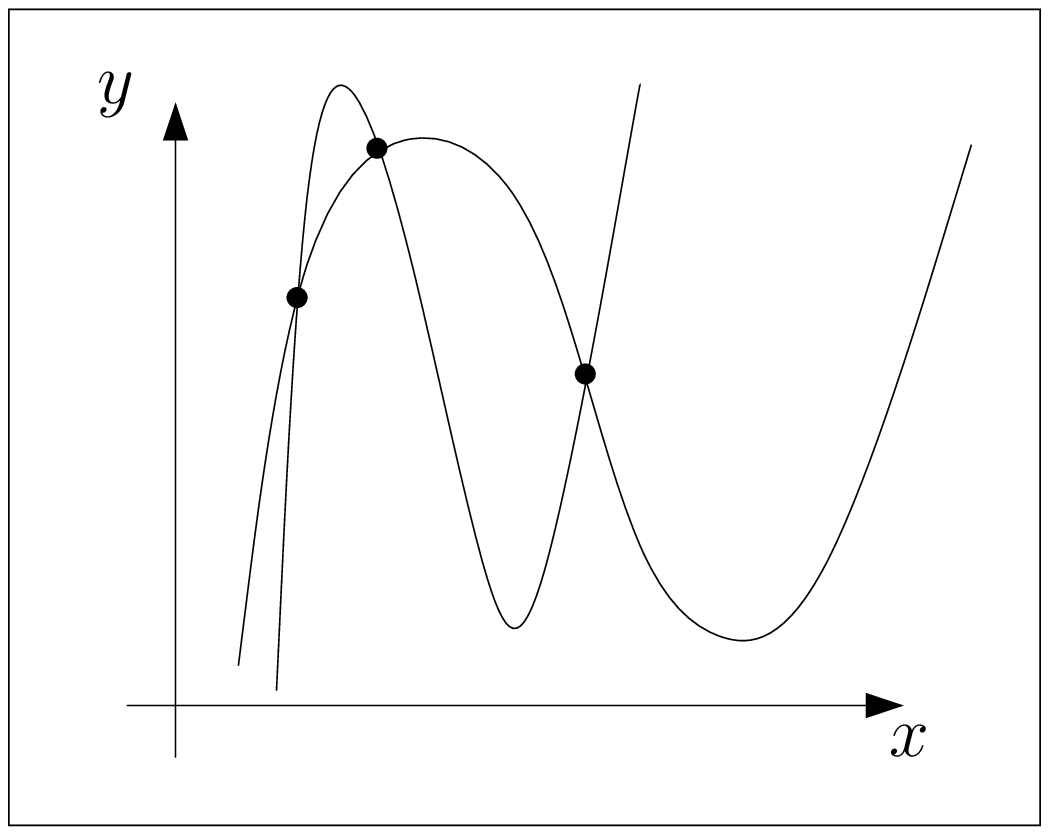}
\end{minipage}
\end{center}
\begin{center}
Figure \ref{ellfibsec}.4
\end{center}

\pagebreak
\section{Singular Fibres in Elliptic Fibrations} \label{kodairasec}

K. Kodaira, in \cite{Ko2}, classified the possible singular fibres of a locally holomorphic elliptic fibration. We follow \cite{SSS} mainly, although the original source is \cite{Ko2}, and other good sources are \cite{HKK} and \cite{BPV}. Two other sources are \cite{Pe} and \cite{M3}, which are quoted by \cite{SSS}. There is also a good, short review in \cite{Sc}. The following theorem is quoted from \cite{SSS}.

\newtheorem{kthm}{Theorem}[section]
\begin{kthm}
\upshape \cite{Ko2} \itshape A singular fibre of a locally holomorphic elliptic fibration without multiple fibres is either of type $I_n$ ($n \geq 1$), of type $II$, $III$, $IV$, or of type $I_n^*$ ($n \geq 0$), or an $\tilde{E_6}$-, $\tilde{E_7}$- or $\tilde{E_8}$-fibre.
\end{kthm}

\newtheorem{kthmr1}[kthm]{Remark}
\begin{kthmr1}
\upshape
We shall not discuss \emph{multiple fibres} here, or how they occur in elliptic fibrations. 
\end{kthmr1}

We now describe the topological properties of these singular fibres. See Table \ref{kodairasec}.1 (from \cite{Sc} and \cite{SSS}).

\newtheorem{typeiiirem1}[kthm]{Remark}
\begin{typeiiirem1}
\upshape
We call a 2-sphere of self-intersection $-2$ a \emph{$-2$-sphere}.
\end{typeiiirem1}
\subsection*{Type $I_n$ fibres $(n \geq 1)$}
The $I_1$-fibre is also known as a \emph{fishtail} fibre or as a \emph{nodal} fibre. It is an immersed sphere of homological self-intersection zero with one positive double point. Hence, its Euler characteristic is $\chi(I_1) = 1$. \newline

The $I_n$-fibre for $n \geq 2$ is a plumbing of $n$ $-2$-spheres plumbed along a circle. Such a fibre is sometimes called a \emph{necklace} fibre (\cite{Sc}). The Euler characteristic is $\chi(I_n) = n$.

\subsection*{Type $II$ fibre}
A type $II$ fibre is also known as a \emph{cusp} fibre, since it is topologically a 2-sphere with a cusp singularity, where the singularity is a cone on the trefoil knot (see \cite{GS}). Its Euler characteristic is $2$.

\subsection*{Type $III$ fibre}
A type $III$ fibre is topologically the union of two $-2$-spheres intersecting each other (not transversally) in a unique point, with multiplicity 2. Therefore, its Euler characteristic is 3.

\subsection*{Type $IV$ fibre}
The type $IV$ fibre is topologically the union of three $-2$-spheres intersecting each other transversally in a unique point. Therefore its Euler characteristic is $4$.

\newtheorem{eulcharrem1}[kthm]{Remark}
\begin{eulcharrem1}
\upshape
The Euler characteristic of the a singular fibre is computed quite easily. Each sphere in the fibre contributes 2 to the Euler characteristic, and then each time a pair of points is identified, the Euler characteristic decreases by 1. For example, the type $III$ singular fibre consists two spheres with two points identified, so the Euler characteristic is $2 +2 -1 = 3$. The type $IV$-fibre consists of three spheres, with three points identified, so the Euler characteristic is $2+2+2-1-1 = 4$ (think of it as replacing three points with a single point).
\end{eulcharrem1}

\subsection*{Type $I_n^*$ fibres $(n \geq 0)$}
The type $I_n^*$ fibre is described by the plumbing given in Table \ref{kodairasec}.1, where all the spheres are $-2$-spheres, and the multiplicities are indicated on the tree. There are $n+1$ spheres of multiplicity $2$, and so the there are $n+5$ spheres in total. From the plumbing diagram, it is easy to calculate that the Euler characteristic if the $I_n^*$ fibre is $n+6$.

\subsection*{The $\tilde{E_8}$ fibre}
The type $\tilde{E_8}$ fibre is described by the plumbing given in Table \ref{kodairasec}.1, where all nine of the vertices are $-2$-spheres. The numbers next to the vertices are the multiplicities the spheres have as homology classes in the fibre. Again, we can use the plumbing diagram to calculate that the Euler characteristic of this fibre is 10.

\subsection*{The $\tilde{E_7}$ fibre}
The type $\tilde{E_7}$ fibre is described by the plumbing given in Table \ref{kodairasec}.1, where all eight of the vertices are $-2$-spheres, and again the multiplicities are indicated. The Euler characteristic is 9.
 
\subsection*{The $\tilde{E_6}$ fibre}
The type $\tilde{E_6}$ fibre is described by the plumbing diagram given in Table \ref{kodairasec}.1, where all seven of the vertices are $-2$-spheres and the multiplicities are again given. The Euler characteristic is 8.

\subsection*{More remarks}

\newtheorem{i0rem}[kthm]{Remark}
\begin{i0rem}
\upshape
We define the $I_0$ fibre to be simply the generic fibre (the torus), and hence it is not a singular fibre. 
\end{i0rem}

\newtheorem{e6e7e8rem}[kthm]{Remark}
\begin{e6e7e8rem}
\upshape
The $\tilde{E_8}$, $\tilde{E_7}$ and $\tilde{E_6}$ fibres are also known as the type $II^*$, type $III^*$ and type $IV^*$ fibres, respectively. One reason for this is the Euler characteristics:  
\begin{align*}
\chi(II) + \chi(II*) = 2 + 10 = 12 \\
\chi(III) + \chi(III*) = 3 + 9 = 12 \\
\chi(IV) + \chi(IV*) = 4 + 8 = 12
\end{align*}
and so the fibres are ``dual to each other''. 
\end{e6e7e8rem}

\newtheorem{e6e7e8rem2}[kthm]{Remark}
\begin{e6e7e8rem2}
\upshape
A section of an elliptic fibration can be thought of as a curve that intersects each fibre in a unique point. Consequently, a section can only intersect those $-2$-spheres in an elliptic fibration that have multiplicity 1. Therefore, a section of the an elliptic fibration that contains an $\tilde{E_8}$ fibre can only intersect the one $-2$-sphere that has multiplicity 1. Similarly, an $\tilde{E_7}$ fibre has two spheres (of multiplicity 1) that a section can interesect, while an $\tilde{E_6}$ has three such spheres.
\end{e6e7e8rem2}

\newtheorem{obsrem1}[kthm]{Remark}
\begin{obsrem1}
\upshape
The reason for calculating the Euler characteristic is that it provides us with an obstruction to the existence of an elliptic fibration with a certain collection of singular fibres (which we shall call a \emph{configuration} of fibres). In other words, since $\chi (\cpt \# \cptbar) = 12$, the sum of the Euler characteristics of the singular fibres in a particular elliptic fibration must equal 12. \newline
\pagebreak

So, while it may be possible for an elliptic fibration to have 12 $I_1$ fibres, or an $\tilde{E_8}$ fibre and 2 $I_1$ fibres (since the sum of the Euler characteristics is 12), it is not possible to have an elliptic fibration with an $\tilde{E_6}$ fibre, a type $III$ fibre and two $I_1$ fibres, since $\chi(\tilde{E_6}) + \chi(III) + 2\chi(I_1) = 8+3+2(1) = 13$. \newline

Note however, that just because a certain configuration of fibres has (collectively) an Euler characteristic of 12, there is no guarantee that there can exist an elliptic fibration with that collection of singular fibres. \cite{SSS} discusses this issue in detail, and also uses the signature of the singular fibres to show when certain configurations cannot exist. \newline
\end{obsrem1}


\pagebreak

\begin{center}
\begin{minipage}{13cm}
\includegraphics[width=13cm]{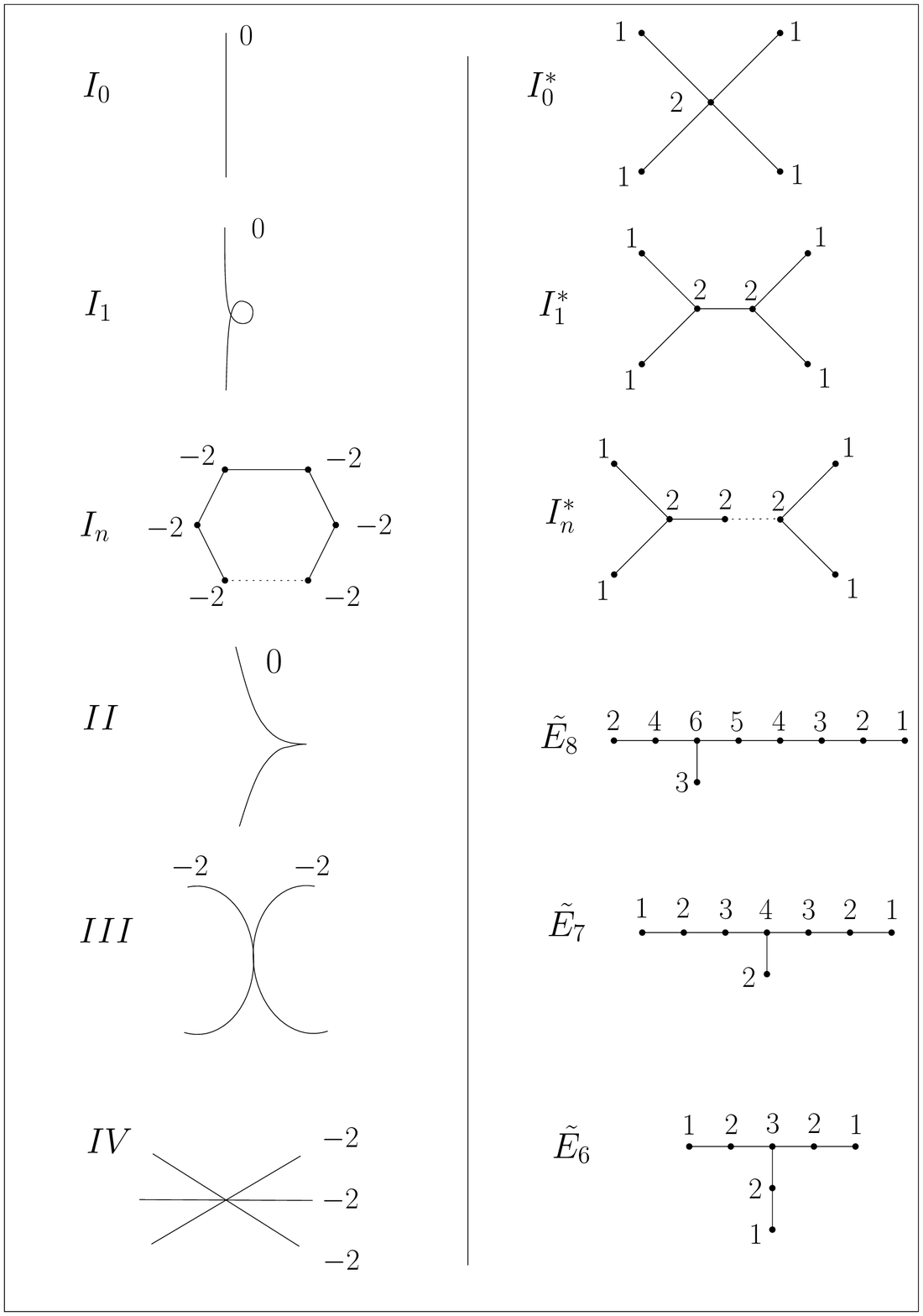}
\end{minipage}
\end{center}
\begin{center}
Table \ref{kodairasec}.1
\end{center}

\pagebreak

\section{The Fishtail Fibre and the Cusp Fibre} \label{fishtailsection}

This section comes from Section 2.3 in \cite{GS} and we follow their notation. We take a closer look at the fishtail and cusp fibres.

\newtheorem{ff}{Definition}[section]
\begin{ff}
\upshape
Consider the following singular curve
\begin{equation*}
C_1 = \{[x:y:z] \in \cpt \st zy^2 = x^3 + zx^2 \}
\end{equation*}
We call any fibre in $E(1)$, that comes from blowing up a curve ambiently isotopic to $C_1$, a \emph{fishtail fibre}.
\end{ff}

\newtheorem{ffr1}[ff]{Remark}
\begin{ffr1}
\upshape
Figure \ref{fishtailsection}.1 shows an example of a curve $C_1$.
\end{ffr1}

\newtheorem{fff}[ff]{Proposition}
\begin{fff} \label{fishtailprop}
The curve $C_1$ is smooth except at the point $P = [0:0:1] \in \cpt$ and is homeomorphic to a sphere with two points identified.
\end{fff}
Proof: \newline
To prove that $C_1$ is smooth except at the point $P=[0:0:1]$, we use the Implicit Function Theorem. Let $p_1(x,y,z) = zy^2 - zx^2 - x^3$. In the chart $[x:y:1]$, we need to find all the points $(x,y) \in \mathbb{C}^2$ that satisfy the following three equations
\begin{align*}
p_1(x,y) &= y^2 - x^2 - x^3 = 0 \\
\frac{\partial p_1}{\partial x}(x,y) &= - 2x - 3x^2 = 0 \\
\frac{\partial p_1}{\partial y}(x,y) &= 2y = 0
\end{align*}
and the only point in $[x:y:1]$ that satisfies all three equations is $[0:0:1]$. It can be checked that there are no points in the charts $[1:y:z]$ and $[x:1:z]$ that satisfy the relevant three equations, and therefore $C_1$ is a smooth curve except at the point $P=[0:0:1]$. \newline

Recall that $\cpo$ is homeomorphic to the 2-sphere. We now want to show that there is a map $\cpo \longrightarrow C_1$ that is one-to-one except that two points get mapped to $P$. We define this map as follows: consider all the projective lines that pass through $P$. This space can be parametrized by $[t_0:t_1] \in \cpo$ as
\begin{equation*}
L_{[t_0:t_1]} = \{[x:y:z] \in \cpt | t_0x = t_1y \}
\end{equation*}
(every projective line through $P$ is of the form $ax + by + cz = 0$ for some $a,b,c \in \mathbb{C}$ not all zero. Since each line goes through $P$, we must have $a(0)+b(0)+c(1) = 0$, which implies that $c=0$, and so $ax = -by$, where $a,b \in \mathbb{C}$ are both not zero. Therefore, we can choose $a = t_0$, $b = -t_1$, where $[t_0:t_1] \in \cpo$). \newline

Let us calculate the number of intersection points of $L_{[t_0:t_1]}$ and $C_1$ (we define $\alpha = \frac{t_0}{t_1}$ and assume $t_1 \neq 0$):
\begin{align*}
& zy^2 - zx^2 - x^3 = 0 \; \; \mathrm{and} \; \; t_0 x = t_1 y \\
\Rightarrow & zy^2 - zx^2 - x^3 = 0 \; \; \mathrm{and} \; \; y = \frac{t_0}{t_1} x = \alpha x \\
\Rightarrow & z (\alpha x)^2 - zx^2 - x^3 = 0 \\
\Rightarrow & \alpha ^2 zx^2 - zx^2 - x^3 = 0 \\
\Rightarrow & (\alpha ^2-1) zx^2 - x^3 = 0 \\
\Rightarrow & x^2 ((\alpha ^2-1) z - x) = 0 \\
\Rightarrow & (\alpha ^2-1) z = x \qquad(x \neq 0) \\
\Rightarrow & z = \frac{1}{\alpha ^2-1} x
\end{align*}

So the line $L_{[t_0:t_1]}$ intersects $C_1$ in $P$ and the point
\begin{equation*}
[x: \alpha x : \frac{1}{\alpha^2 - 1} x]
\end{equation*}
which is the same as the point 
\begin{equation*}
[1: \frac{t_0}{t_1} : \frac{1}{(\frac{t_0}{t_1}) ^2 - 1}]
\end{equation*}

However, when $(\frac{t_0}{t_1}) ^2 = 1$, $\frac{1}{(\frac{t_0}{t_1})^2- 1}$ is undefined. If $(\frac{t_0}{t_1}) ^2 = 1$, then $[t_0:t_1] = [\pm 1: 1]$, and then

\begin{align*}
& zy^2 - zx^2 - x^3 = 0 \; \; \mathrm{and} \; \; x = \pm y \\
\Rightarrow & z(\pm x)^2 - zx^2 - x^3 = 0 \\
\Rightarrow & zx^2 - zx^2 - x^3 = 0 \\
\Rightarrow & x^3 = 0 \\
\Rightarrow & x = 0 \\
\Rightarrow & y = 0 \\
\Rightarrow & [x:y:z] = [0:0:1]
\end{align*}
and so $L_{[t_0:t_1]} \cap C_1 = \{P, Q_{[t_0:t_1]}\}$ except for the cases $[t_0:t_1] = [\pm 1:1]$, and then we define $Q_{[\pm1:1]} = P$, so $L_{[\pm 1: 1]} \cap C_1 = \{P\}$. We define the map $\psi :\cpo \longrightarrow C_1$ by $\psi([t_0:t_1]) = Q_{[t_0:t_1]}$, which is one-to-one except that $\psi([1:1]) = \psi([-1:1]) = P$ (see Figures \ref{fishtailsection}.2 and \ref{fishtailsection}.3). So, $C_1$ is homeomorphic to $\cpo$ with $[1:1]$ and $[-1:1]$ identified, and so is homeomorphic to the 2-sphere with two points identified. $\Box$

\newtheorem{cfd}[ff]{Definition}
\begin{cfd}
\upshape
We call a fibre in $E(1)$ that comes from blowing up the curve ambiently isotopic to
\begin{equation*}
C_2 = \{[x:y:z] \in \cpt | zy^2 = x^3 \}
\end{equation*}
a \emph{cusp fibre}.
\end{cfd}

\newtheorem{cfdr1}[ff]{Remark}
\begin{cfdr1}
\upshape
See Figure \ref{fishtailsection}.4 for an example of a curve $C_2$.
\end{cfdr1}

\newtheorem{cf}[ff]{Proposition}
\begin{cf}
The curve $C_2$ is smooth except at the point $P = [0:0:1] \in \cpt$ and is homeomorphic to a sphere.
\end{cf}
Proof: \newline
The calculations are almost identical to the proof for Proposition \ref{fishtailprop}. The only difference is that $C_2 \cap L_{[t_0:t_1]} = \{P, Q_{[t_0:t_1]} \}$, and $Q_{[0:1]} = P$, in which case we have $C_2 \cap L_{[0:1]} = \{P \}$. Therefore the map $\psi : \cpo \longrightarrow C_2$ defined by  $\psi([t_0:t_1]) = Q_{[t_0:t_1]}$ gives a homeomorphism between $\cpo$ and $C_2$. $\Box$


\pagebreak

\begin{center}
\begin{minipage}{12cm}
\includegraphics[width=12cm]{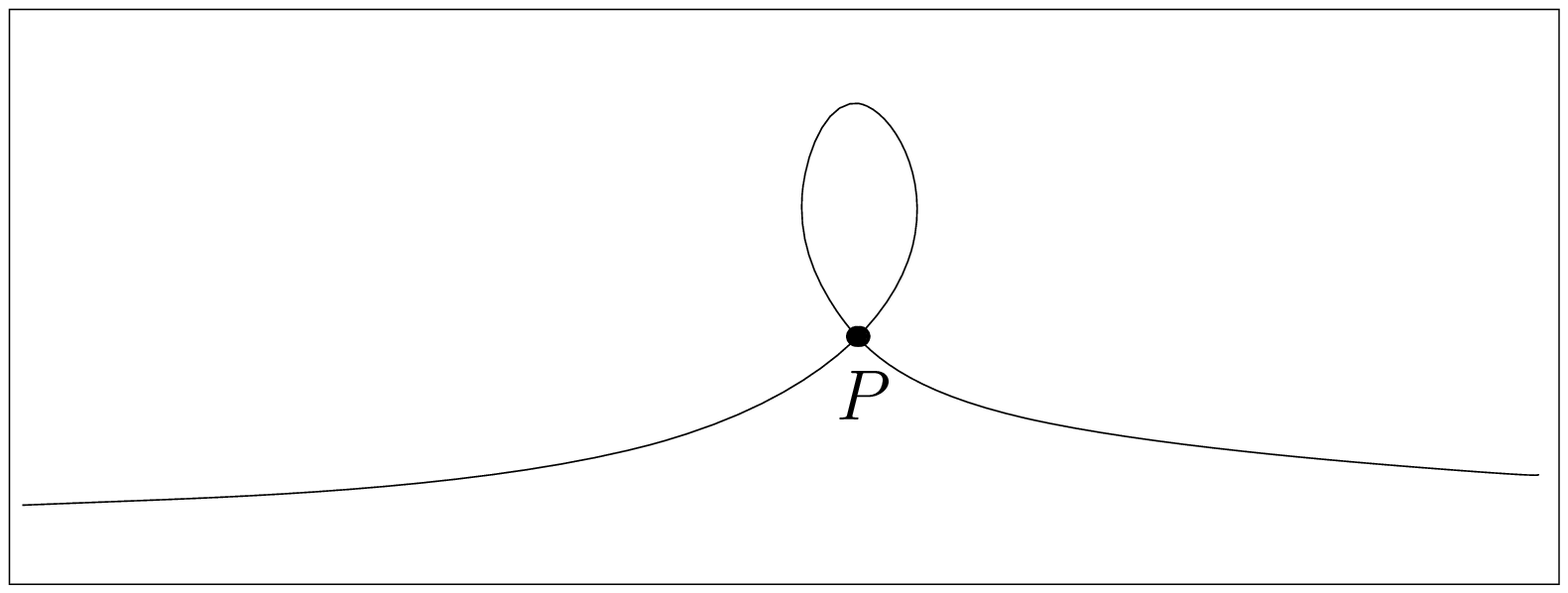}
\end{minipage}
\end{center}
\begin{center}
Figure \ref{fishtailsection}.1
\end{center}

\begin{center}
\begin{minipage}{12cm}
\includegraphics[width=12cm]{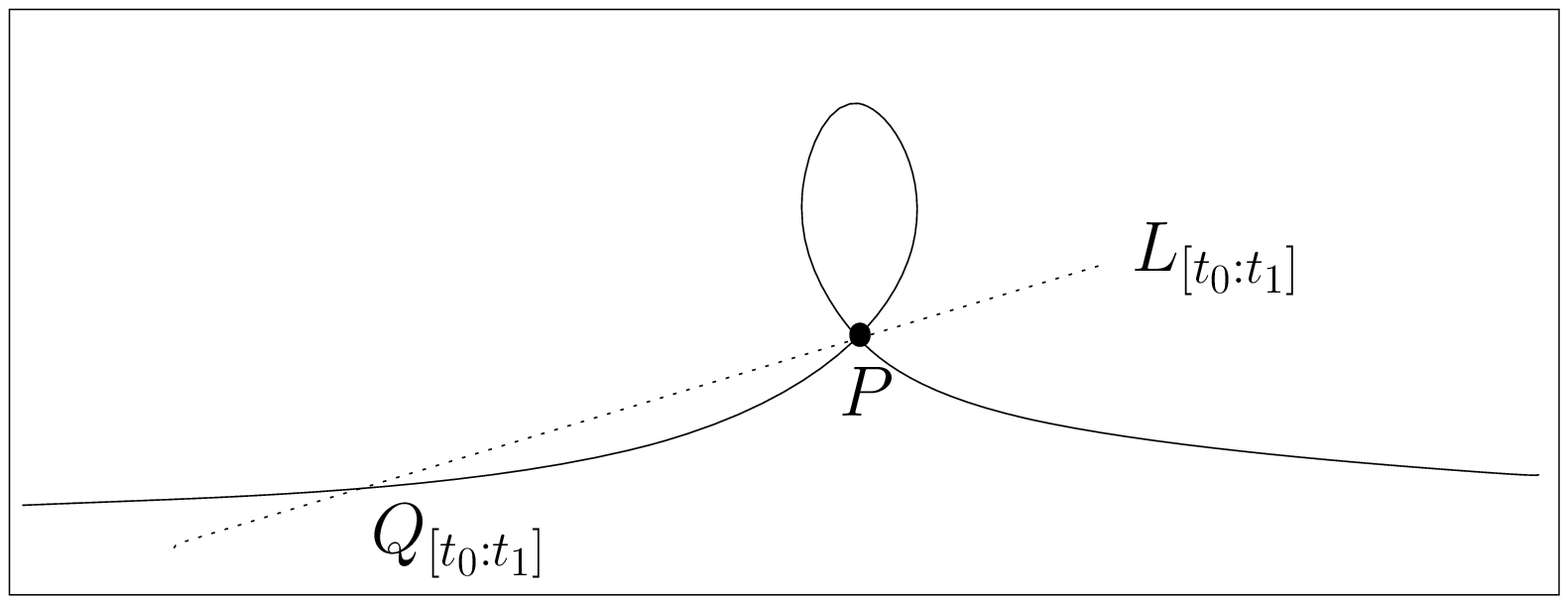}
\end{minipage}
\end{center}
\begin{center}
Figure \ref{fishtailsection}.2
\end{center}

\begin{center}
\begin{minipage}{12cm}
\includegraphics[width=12cm]{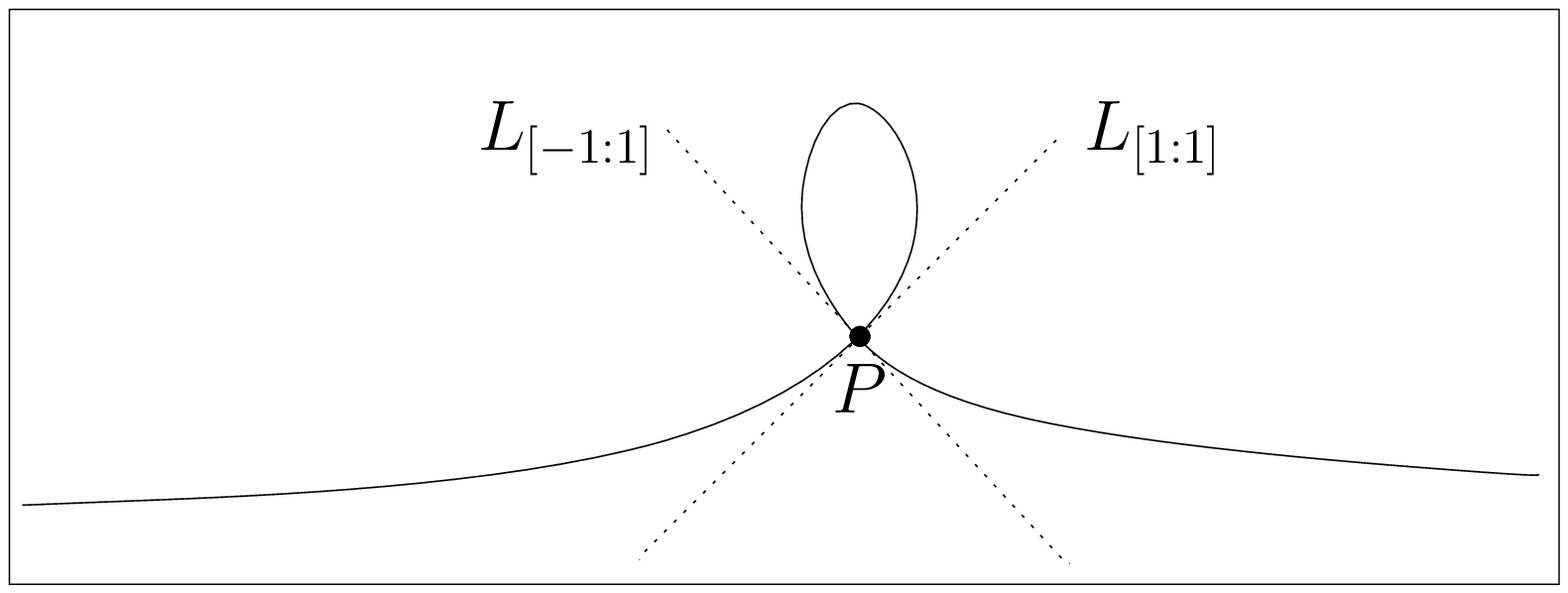}
\end{minipage}
\end{center}
\begin{center}
Figure \ref{fishtailsection}.3
\end{center}

\begin{center}
\begin{minipage}{5cm}
\includegraphics[width=5cm]{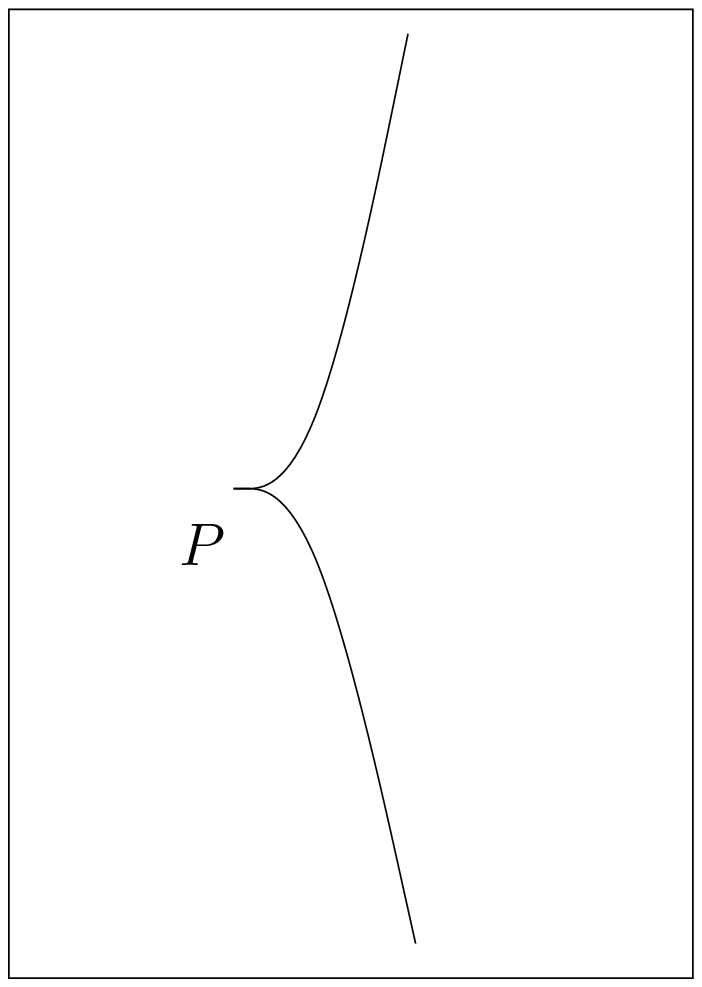}
\end{minipage}
\end{center}
\begin{center}
Figure \ref{fishtailsection}.4
\end{center}

\pagebreak
\section{Blowing Up To Create an $\tilde{E_8}$ Fibre} \label{e8sec1}

This example follows the outline given in section \ref{ellfibsec} (source: \cite{S}). \newline

We start with with a cubic polynomial $p_1$ and a linear polynomial $l$, which are chosen so that the corresponding curves $C_1$ and $L$ intersect in only one point $P$. This point will therefore be a point of multiplicity 3. If we let the homology class of $L$ be $h \in H_2(\cpt; \mz)$, then the homology class of $C_1$ is $3h \in H_2(\cpt; \mz)$, and both homology classes have multiplicity 1. \newline

We then define the cubic $p_0$ as $p_0 = l^3$, and then we have two cubics, $p_0$ and $p_1$ such that their corresponding curves $C_0$ and $C_1$ intersect in only one point $P \in \cpt$. We can think of $C_0$ as three copies of $L$ lying on top of each other, each intersecting $C_1$ three times in the single point $P$; therefore $P$ has multiplicity 9. The homology class of $C_0$ is $h$, but it has multiplicity 3. Therefore, the homology class of $C_0$ in $H_2(\cpt; \mz)$ is $h(3)=3h$. Notice that the homology class of $C_1$ in $H_2(\cpt; \mz)$ is $3h(1) = 3h$. Any element of the pencil on $\cpt$ must represent the same homology class (in this case, $3h$) $-$ this fact (mentioned in \cite{SSS}) will be used repeatedly to calculate the multiplicity of the exceptional curve $e_i$ in the $i$th blow-up at $P$. There will be 9 blow-ups at $P$, since it has multiplicity 9. \newline

The starting point is illustrated in Figure \ref{e8sec1}.1, and the homology classes are indicated next to their respective curves, with their multiplicities in parentheses. We draw $C_0$ as being tangent to $C_1$ at $P$, and also indicate the multiplicity of the point $P$ in parentheses (this multiplicity will decrease by 1 with each blow-up at $P$, until it is 0). \newline

We show that the pencil of curves generated by the polynomials $p_0$ and $p_1$ provides a fibration with an $\tilde{E_8}$ fibre. We shall need to blow up nine times at the point $P$. \newline 

The first blow-up introduces the exceptional curve $e_1$. The proper transform of $C_1$ represents the homology class $3h-e_1 \in H_2(\cpt \# \cptbar; \mz)$, with multiplicity 1 (if a curve has multiplicity $m$, its proper transform will also have multiplicity $m$). The proper transform of $C_0$ represents the homology class  $h - e_1 \in H_2(\cpt \# \cptbar; \mz)$, with multiplicity 3 (so actually, it represents $3h- 3e_1 \in H_2(\cpt \# \cptbar; \mz)$). We consider $C_0$ and the exceptional curve(s) to be one element of the pencil. Therefore, we need $e_1$ to have multiplicity $m$ so that $3h-3e_1 + me_1 = 3h-e_1$, since every element in the pencil on $\cpt \# \cptbar$ must represent the same homology class (in this case, $3h - e_1$). Therefore, $e_1$ must have have multiplicity $2$. Since $L$ intersected $C_1$ in $P$ with multiplicity $3$, it would take three blow-ups at $P$ before the proper transforms of $L$ and $C_1$ would not intersect each other anymore. In fact, the calculation
\begin{align*}
[3h -e_1] \cdot [h - e_1] & = [3h \cdot h] + [3h \cdot (-e_1)] + [-e_1 \cdot h] + [(-e_1) \cdot (-e_1)]  \\
&=  (3) + (0) + (0) + (-1) \\
&= 2 
\end{align*}
shows this (in fact, $h-e_1$ has multiplicity 3, so the curves actually intersect with multiplicity 6). Therefore, we still need to draw the proper tranform of $C_0$ so that it intersects the proper transform of $C_1$. To indicate a change, instead of drawing the curves tangent to each other, we draw the proper transform of $C_0$ intersecting the proper transform of $C_1$ transversally. The proper transforms of $C_0$ and $C_1$ still intersect each other in the point $P$, but now only with multiplicity 8. See Figure \ref{e8sec1}.2.  \newline

From now on, we shall label a curve with its homology class (and only indicate its multiplicity in the diagrams). Therefore, we shall say that $3h - e_1$ is a curve with multiplicity 1, and $h - e_1$ is a curve with multiplicity 3. The diagrams will make this terminology clear. \newline

We blow-up a second time at $P$ (repeated blow-ups at a point are called \emph{infinitely close blow-ups}; see \cite{PSS}). This blow-up introduces the exceptional curve $e_2$, which separates $e_1$ and $h-e_1$ since their proper tranforms, which have homology classes $e_1 - e_2$ and $h - e_1 - e_2$, no longer intersect, the following calculation shows:
\begin{align*}
[e_1 - e_2] \cdot [h - e_1 - e_2] & = [e_1 \cdot h] + [e_1 \cdot (-e_1)] + [e_1 \cdot (-e_2)]  \\
& \; \; - [e_2 \cdot h] - [e_2 \cdot (-e_1)] - [e_2 \cdot (-e_2)] \\
&= (0) + (-(-1)) + (0) - (0) - (0) - (-(-1)) \\
&= 1 - 1\\
&= 0 
\end{align*}
Note that the multiplicities of $e_1 - e_2$ and $h-e_1-e_2$ makes no difference to the calculation. The proper transform of $3h - e_1$ is $3h - e_1 - e_2$. This is only the second blow-up at $P$, and so $h - e_1 - e_2$ still passes through $P$ transversally (we could again check algebraically that $3h-e_1-e_2$ and $h-e_1-e_2$ intersect). We calculate the multiplicity of $e_2$: since 2 is the multiplicity of $e_1-e_2$, 3 is the multiplicity of $h-e_1-e_2$ and 1 is the multiplicity of $3h-e_1-e_2$, if we let $m$ be the multiplicity of $e_2$, we have 
\begin{align*}
& m (e_2) + 2(-e_2) + 3(-e_2) = 1(-e_2) \\
\Rightarrow & m = 4 
\end{align*}
which shows that the multiplicity of $e_2$ is 4. See Figure \ref{e8sec1}.3. \newline

The third blow-up at $P$ introduces the exceptional curve $e_3$. The proper transforms of  $3h-e_1 - e_2$ and $h - e_1 - e_2$ are $3h-e_1 - e_2-e_3$ and $h - e_1 - e_2-e_3$, respectively, and it can be calculated that they no longer intersect. $e_3$ has multiplicity 6, and also separates $e_2$ and $h - e_1 - e_2$ (again, their proper tranforms $e_2-e_3$ and $h-e_1-e_2-e_3$ no longer intersect). See Figure \ref{e8sec1}.4. \newline

The fourth blow-up at $P$ introduces the exceptional curve $e_4$. The proper transform of $3h-e_1-e_2-e_3$ is $3h-e_1-e_2-e_3-e_4$ and the proper tranform of $e_3$ is $e_3-e_4$, and the other curves remain the same. $e_4$ has multiplicity 5. See Figure \ref{e8sec1}.5. \newline

The fifth, sixth, seventh and eighth blow-ups at $P$ follow the same pattern as the fourth blow-up, and are illustrated in Figures \ref{e8sec1}.6,\ref{e8sec1}.7,\ref{e8sec1}.8 and \ref{e8sec1}.9, respectively. It can be seen in the diagrams that the multiplicities of $e_5$, $e_6$, $e_7$ and $e_8$ are 4,3,2 and 1, respectively. \newline

Finally, the ninth blow-up at $P$ introduces the exceptional curve $e_9$. The proper tranform of $e_8$ is $e_8-e_9$ (with multiplicity 1) and the proper transform of $3h-e_1-\dots - e_8$ is $3h-e_1-\dots - e_8-e_9$ (still with multiplicity 1). This means that $e_9$ has multiplicity 0, and note that $P$ also has multiplicity 0, which indicates that it is no longer a singular point (after nine blow-ups at $P$). $e_9$ is therefore a section between the two fibres, $3h-e_1-\dots-e_9 $ (which has self-intersection $0$ and is a fishtail fibre) and the rest of the Figure \ref{e8sec1}.10, which is an $\tilde{E_8}$ fibre. \newline

To see that it is an $\tilde{E_8}$ fibre, observe that every curve $e_j - e_{j+1}$ ($j=1,\dots,8$) and $h-e_1-e_2-e_3$ has self-intersection $-2$, and so are $-2$-spheres. Furthermore, these $-2$-spheres intersect each other and have multiplicities as in Figure \ref{e8sec1}.11, which is the diagram of an $\tilde{E_8}$ as in Table \ref{kodairasec}.1, where the vertices represent the $-2$-spheres (which are labelled above with their multiplicities, and below with the curve to which the sphere corresponds), and the lines joining the vertices show that two $-2$-spheres intersect one another. \newline

So, we have given one construction of an $\tilde{E_8}$ fibre. In section \ref{e8algeomsection}, we shall give an explicit construction using two well-chosen polynomials.

\newtheorem{e8secrem1}{Remark}[section]
\begin{e8secrem1} \label{e8secrem1rem}
\upshape
Note that this $\tilde{E_8}$ fibre has one section, the exceptional sphere $e_9$, connecting the fishtail fibre and the $\tilde{E_8}$ fibre. One definition of a section (of an elliptic fibration) is that it is a curve (the image of the section map) that intersects each fibre exactly once. Therefore, a section only can intersect a singular fibre in one of its spheres that has multiplicity 1. $\tilde{E_8}$ has only one such $-2$-sphere, the sphere $e_8-e_9$, and therefore this is the only sphere through which a section can pass, and we say that $\tilde{E_8}$ admits only one section. In the case of $\tilde{E_6}$, which has three $-2$-spheres of multiplicity 1, there are three spheres that a section can possibly intersect, and we say that $\tilde{E_6}$ has 3 sections.
\end{e8secrem1}

There is a useful lemma from \cite{SSS} (which also contains a proof of this lemma):
\newtheorem{fibseclemsss}[e8secrem1]{Lemma}
\begin{fibseclemsss} \label{fibseclemssslab}
Suppose that two cubic polynomials $p_1, p_2$ define a pencil of elliptic curves in $\cpt$ with $k$ base points. Suppose furthermore that the pencil contains at least one smooth cubic curve. Then the fibration has $k$ sections.
\end{fibseclemsss}
In the next section we construct an elliptic fibration over $E(1)$ with an $\tilde{E_8}$ fibre.

\pagebreak

\begin{center}
\begin{minipage}{8cm}
\includegraphics[width=8cm]{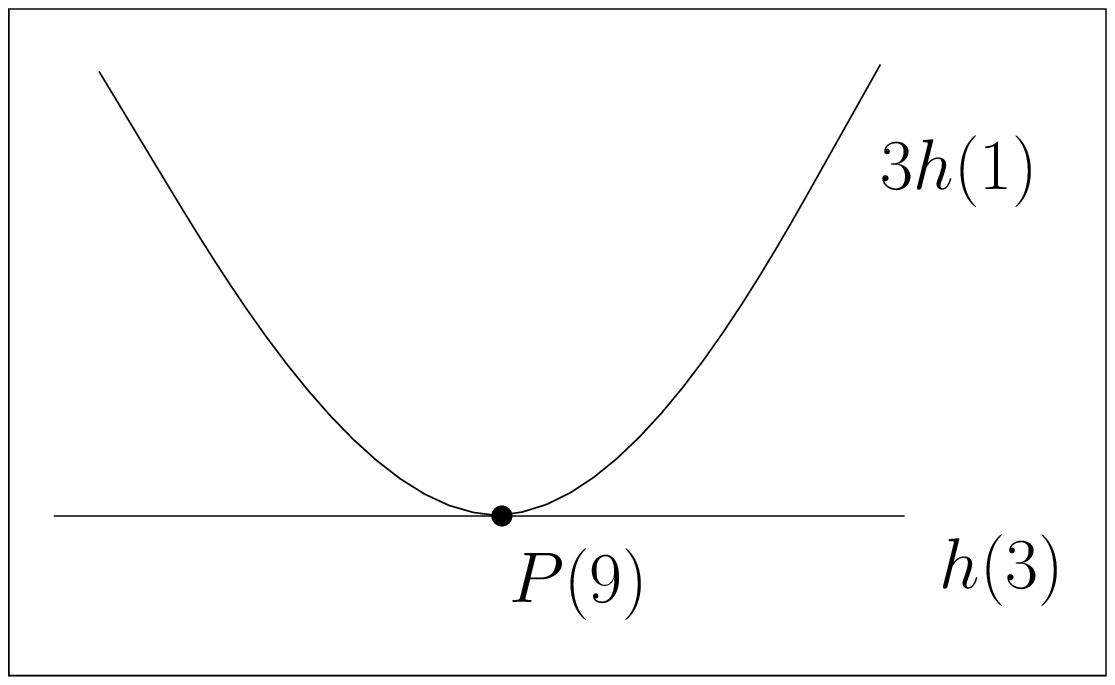}
\end{minipage}
\end{center}
\begin{center}
Figure \ref{e8sec1}.1
\end{center}

\begin{center}
\begin{minipage}{8cm}
\includegraphics[width=8cm]{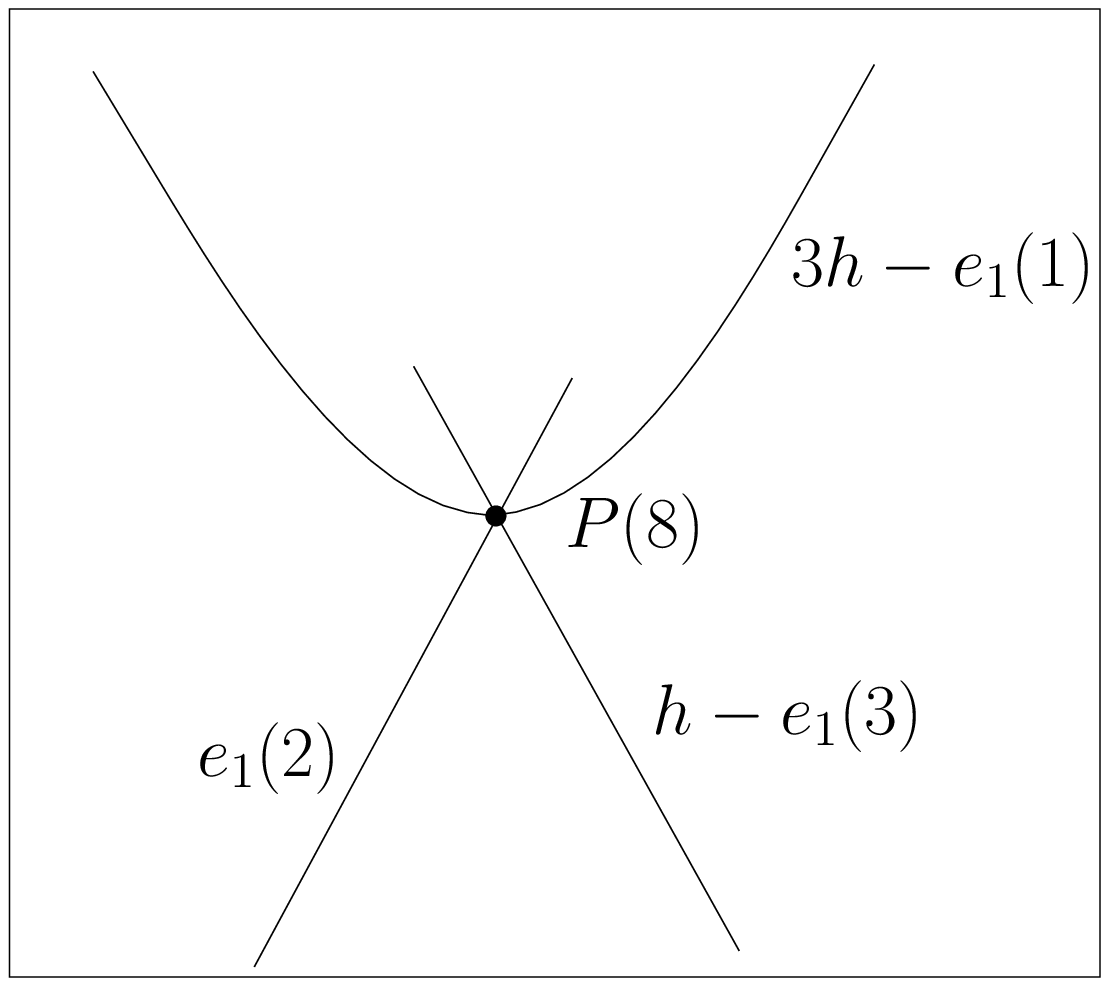}
\end{minipage}
\end{center}
\begin{center}
Figure \ref{e8sec1}.2
\end{center}

\begin{center}
\begin{minipage}{10cm}
\includegraphics[width=10cm]{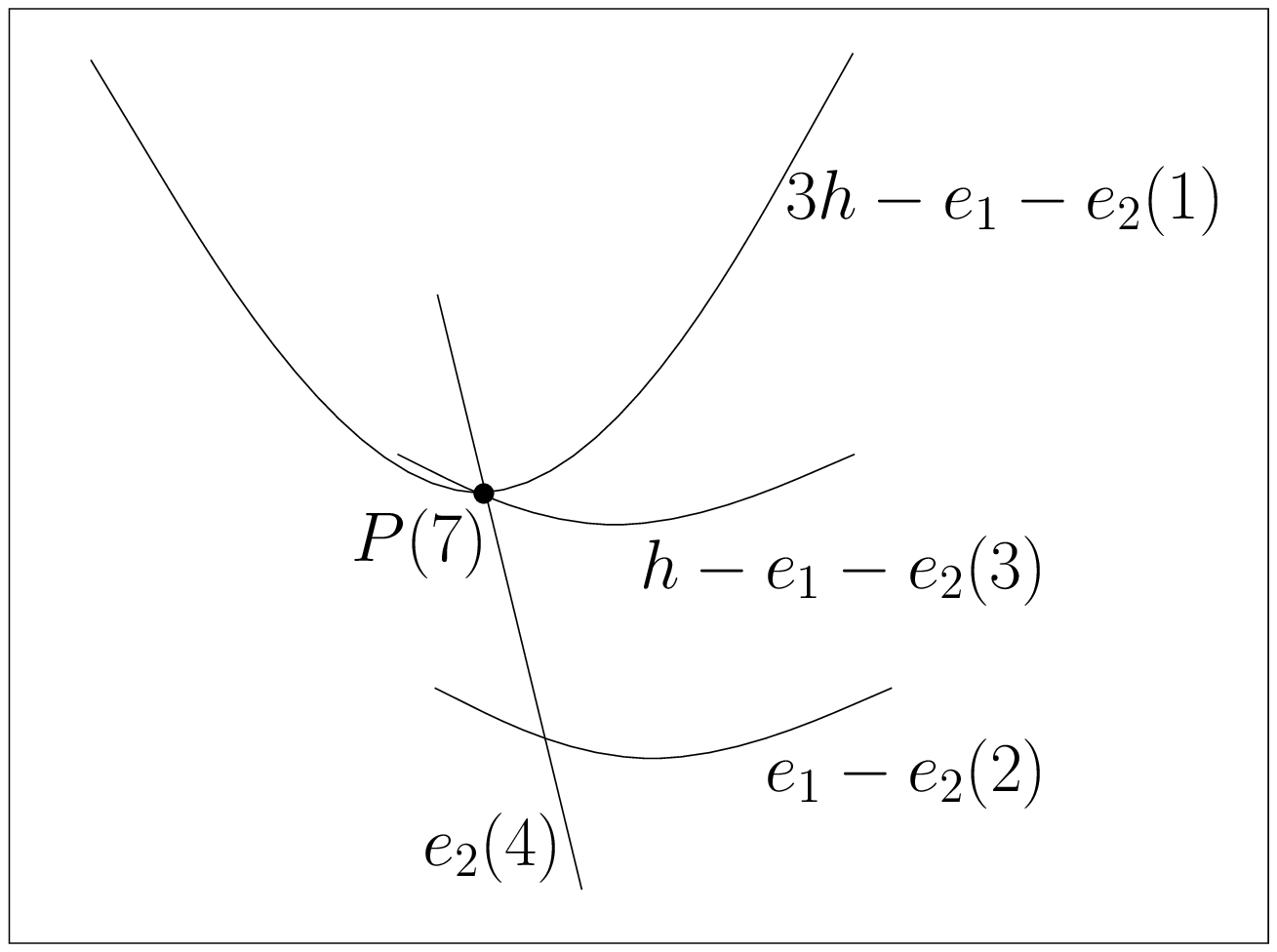}
\end{minipage}
\end{center}
\begin{center}
Figure \ref{e8sec1}.3
\end{center}

\begin{center}
\begin{minipage}{10cm}
\includegraphics[width=10cm]{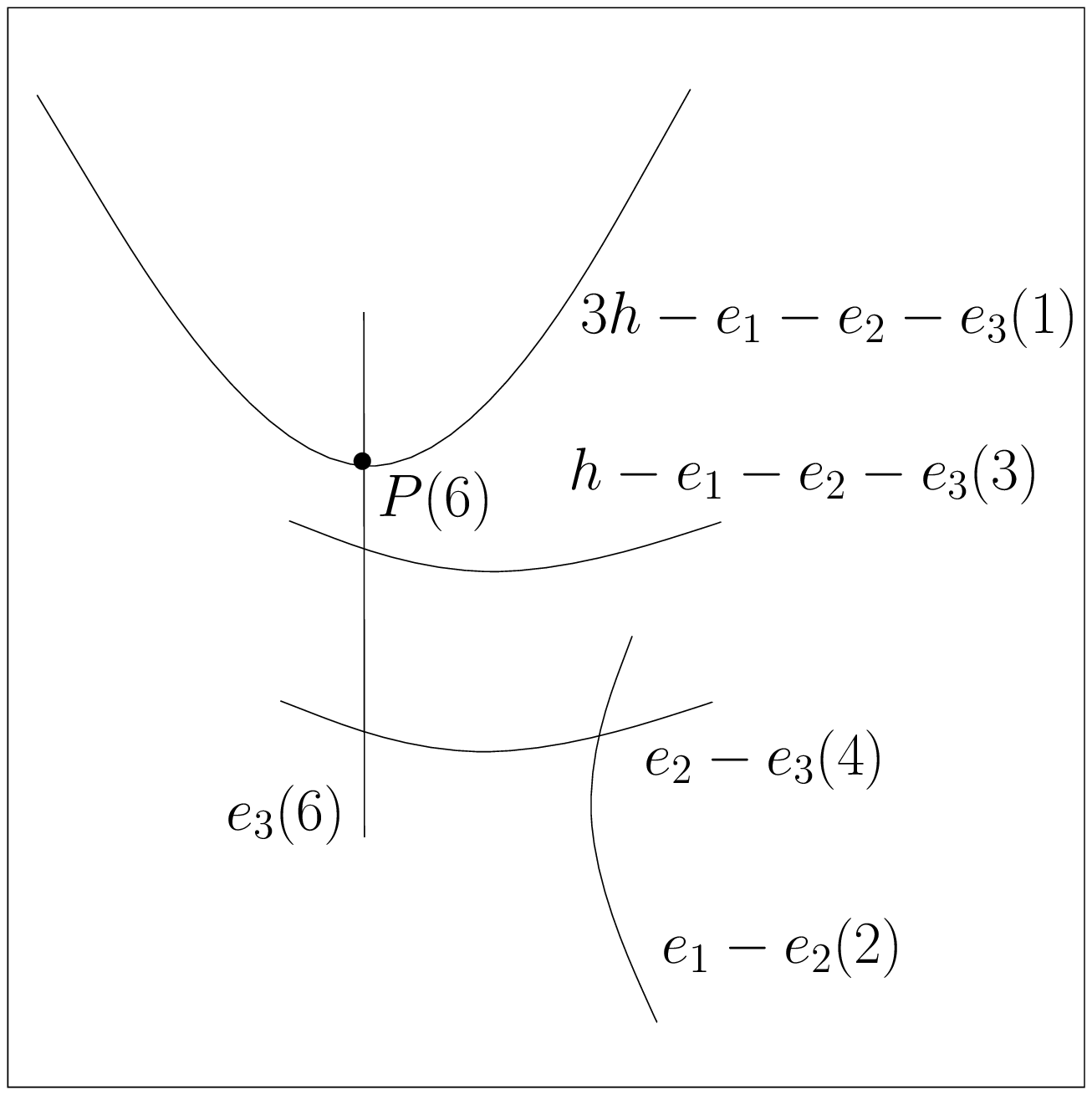}
\end{minipage}
\end{center}
\begin{center}
Figure \ref{e8sec1}.4
\end{center}

\begin{center}
\begin{minipage}{10cm}
\includegraphics[width=10cm]{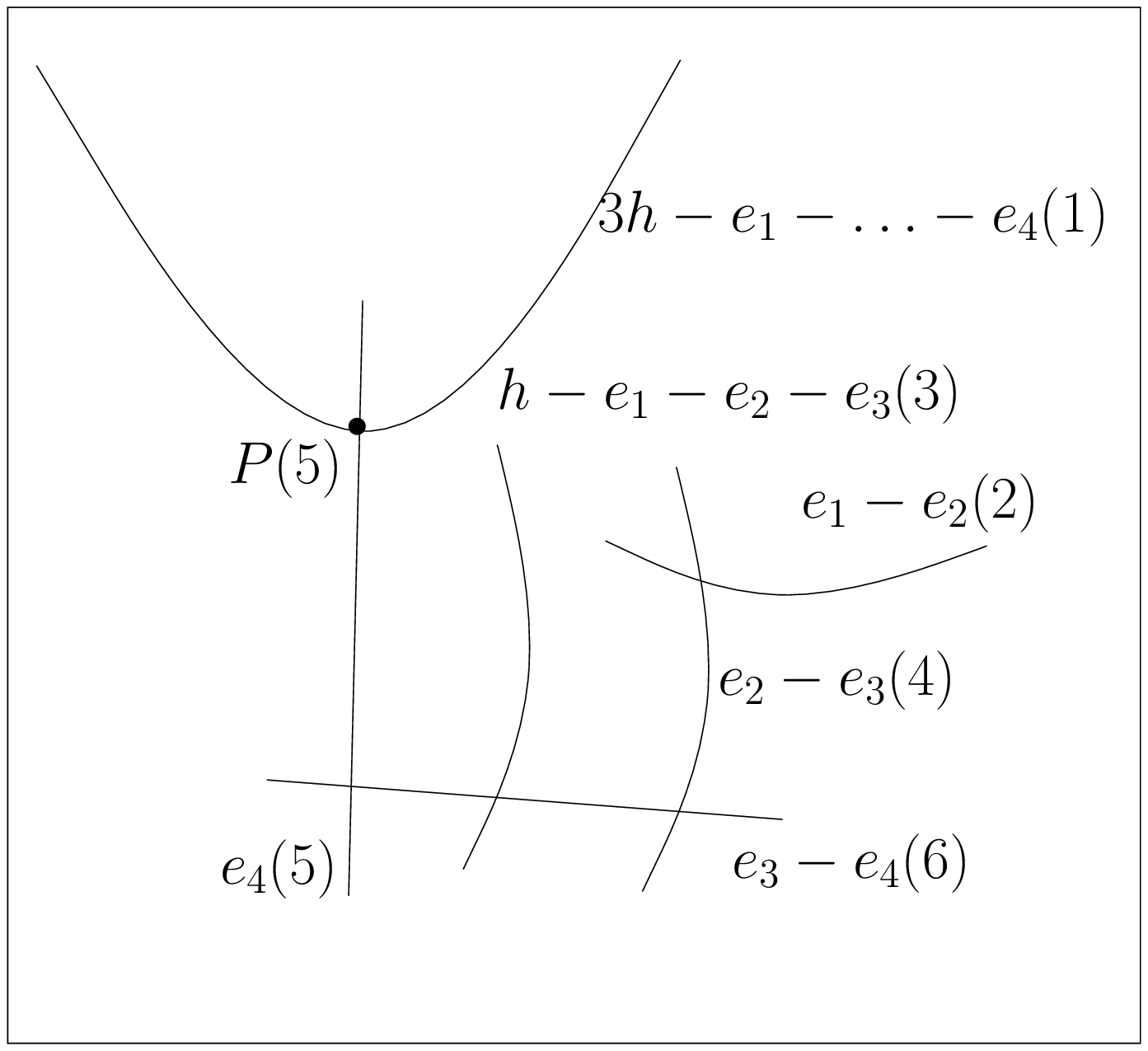}
\end{minipage}
\end{center}
\begin{center}
Figure \ref{e8sec1}.5
\end{center}

\begin{center}
\begin{minipage}{10cm}
\includegraphics[width=10cm]{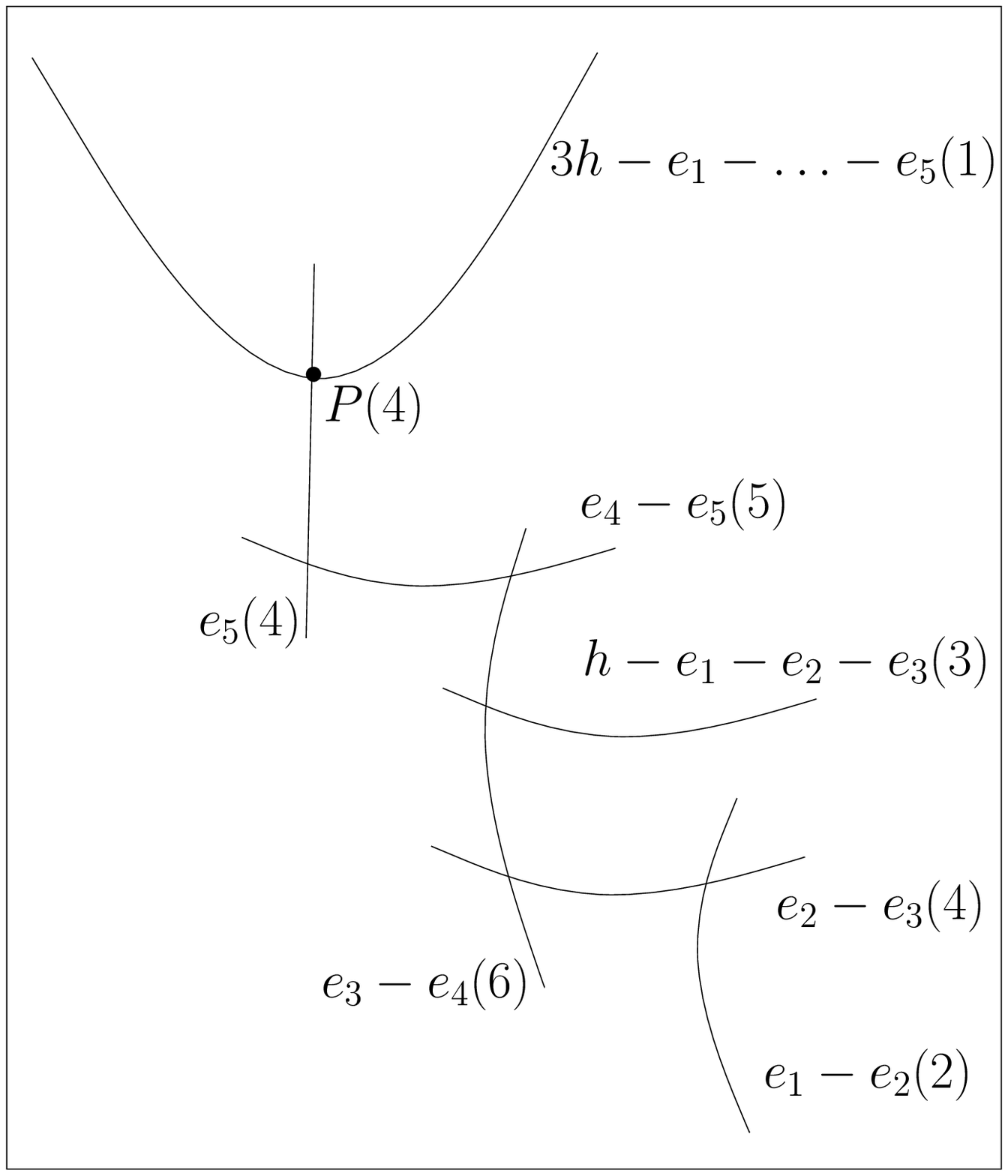}
\end{minipage}
\end{center}
\begin{center}
Figure \ref{e8sec1}.6
\end{center}

\begin{center}
\begin{minipage}{12cm}
\includegraphics[width=12cm]{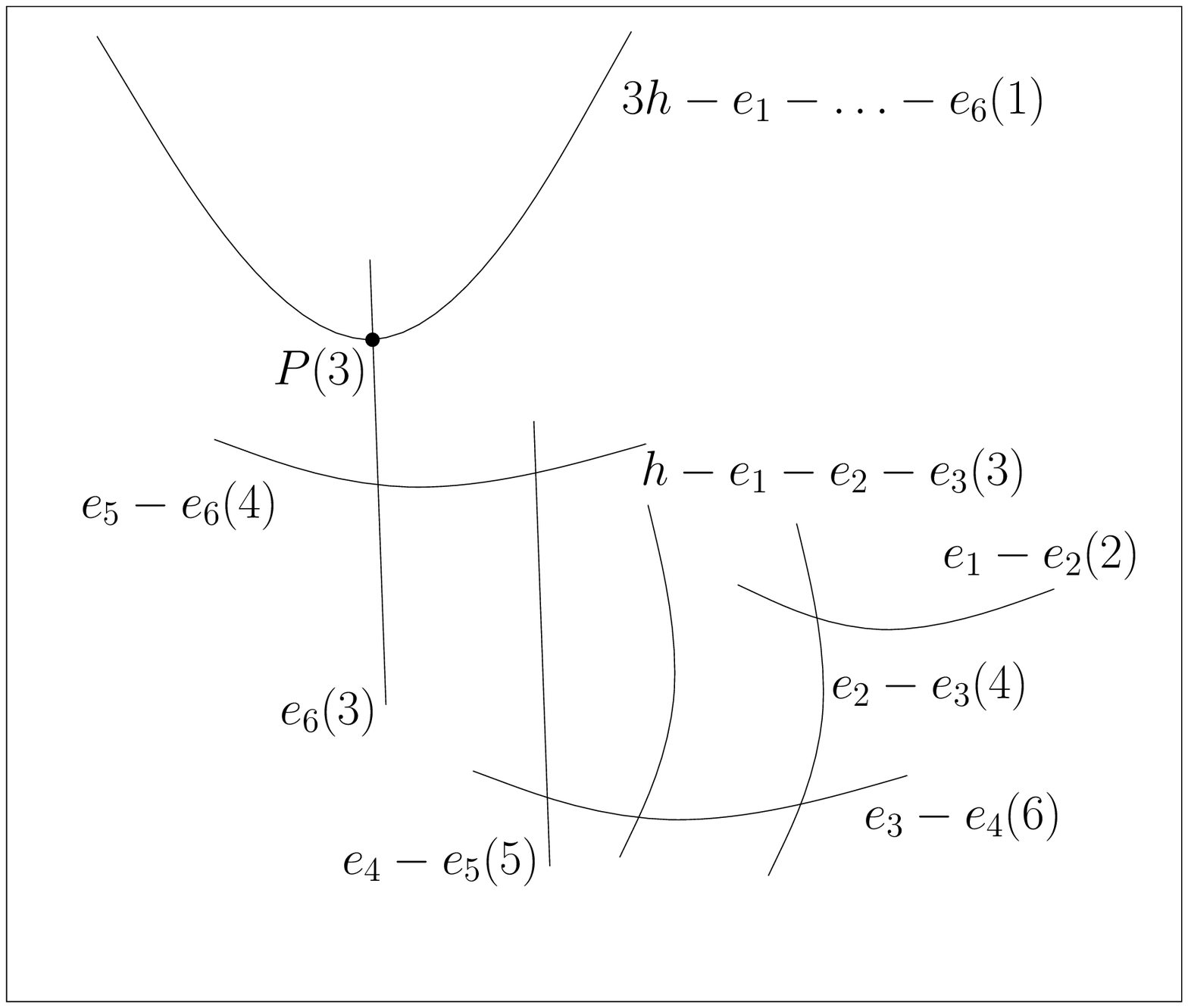}
\end{minipage}
\end{center}
\begin{center}
Figure \ref{e8sec1}.7
\end{center}

\begin{center}
\begin{minipage}{12cm}
\includegraphics[width=12cm]{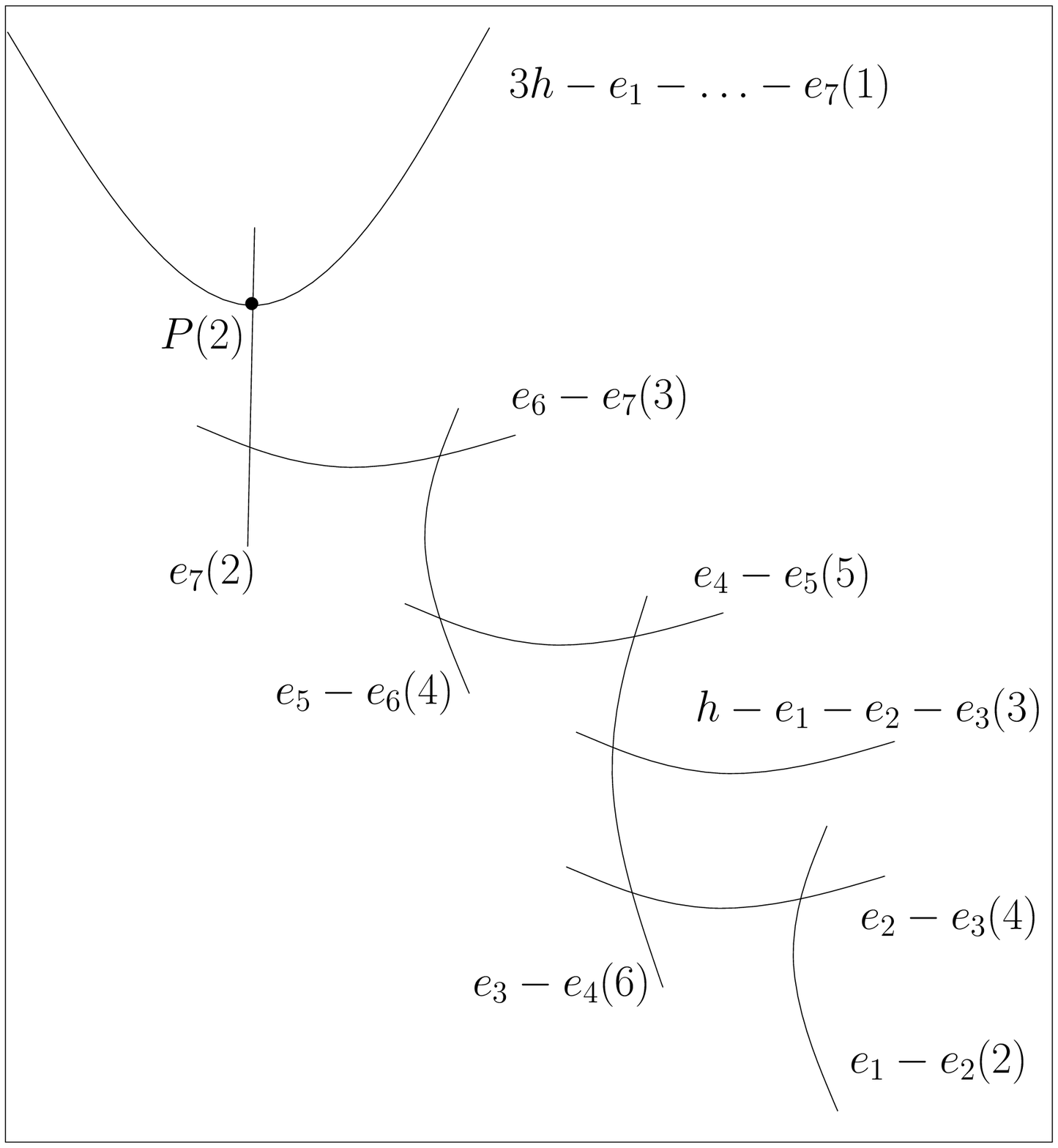}
\end{minipage}
\end{center}
\begin{center}
Figure \ref{e8sec1}.8
\end{center}

\begin{center}
\begin{minipage}{13cm}
\includegraphics[width=13cm]{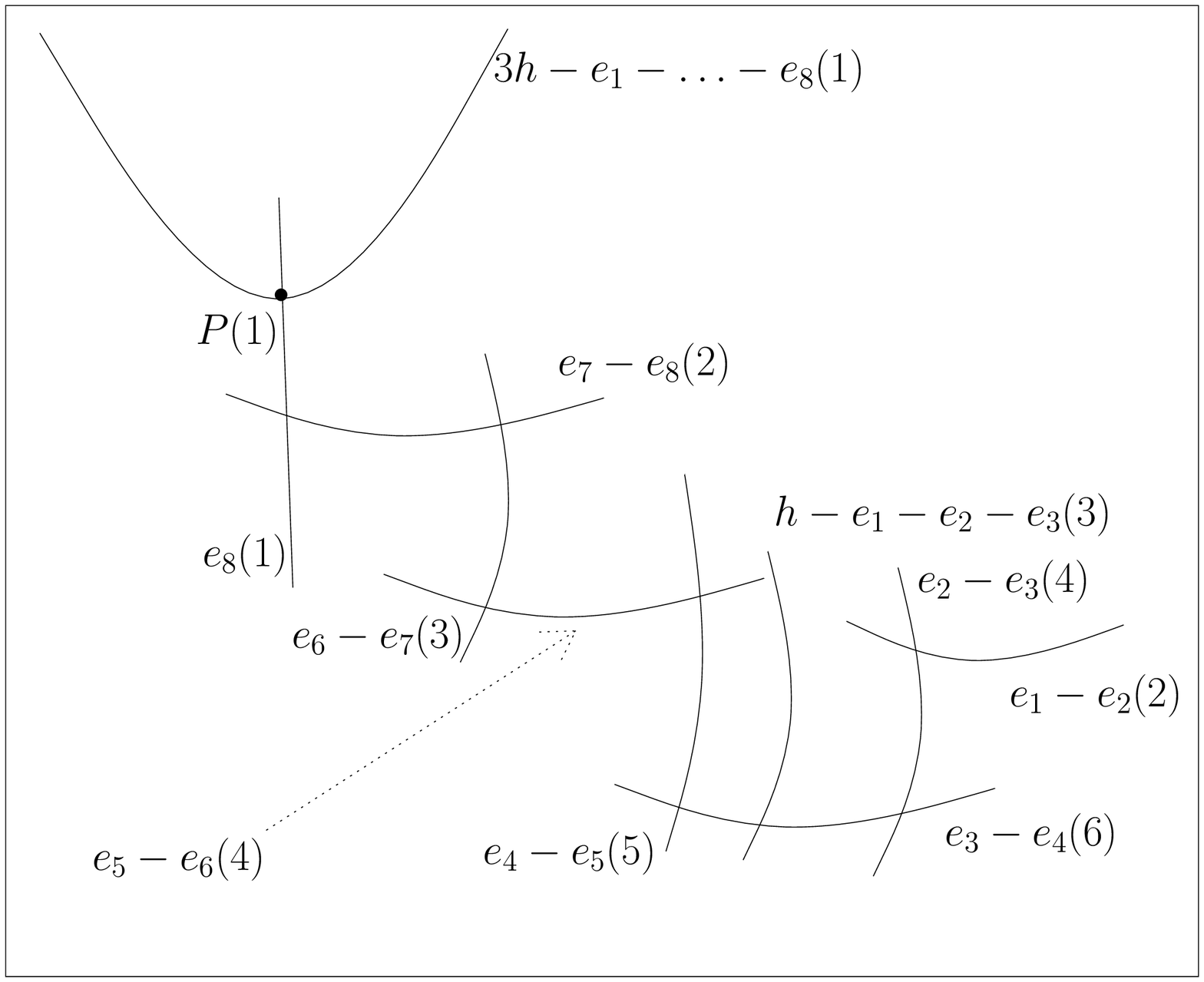}
\end{minipage}
\end{center}
\begin{center}
Figure \ref{e8sec1}.9
\end{center}

\begin{center}
\begin{minipage}{13cm}
\includegraphics[width=13cm]{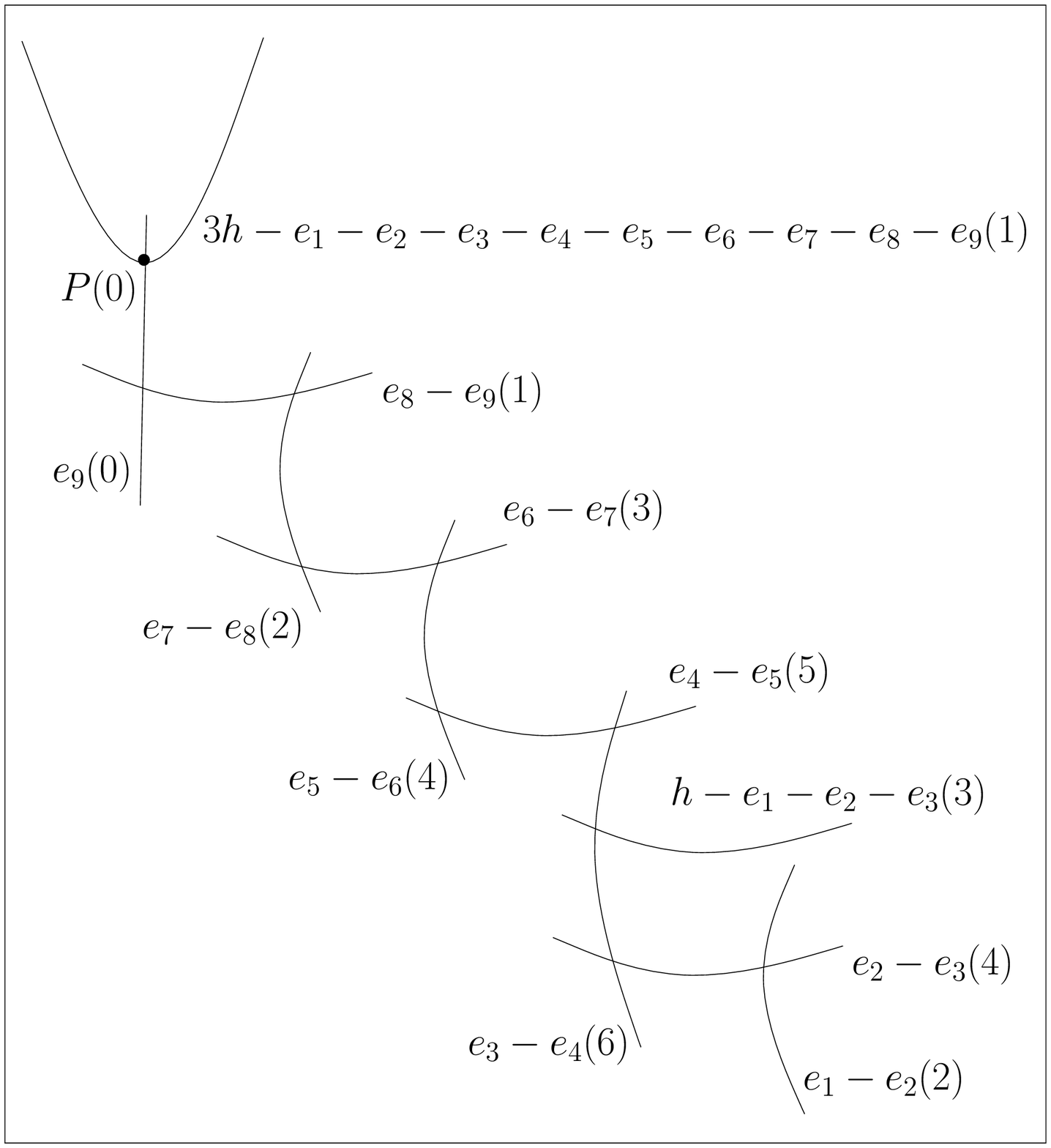}
\end{minipage}
\end{center}
\begin{center}
Figure \ref{e8sec1}.10
\end{center}

\begin{center}
\begin{minipage}{12cm}
\includegraphics[width=12cm]{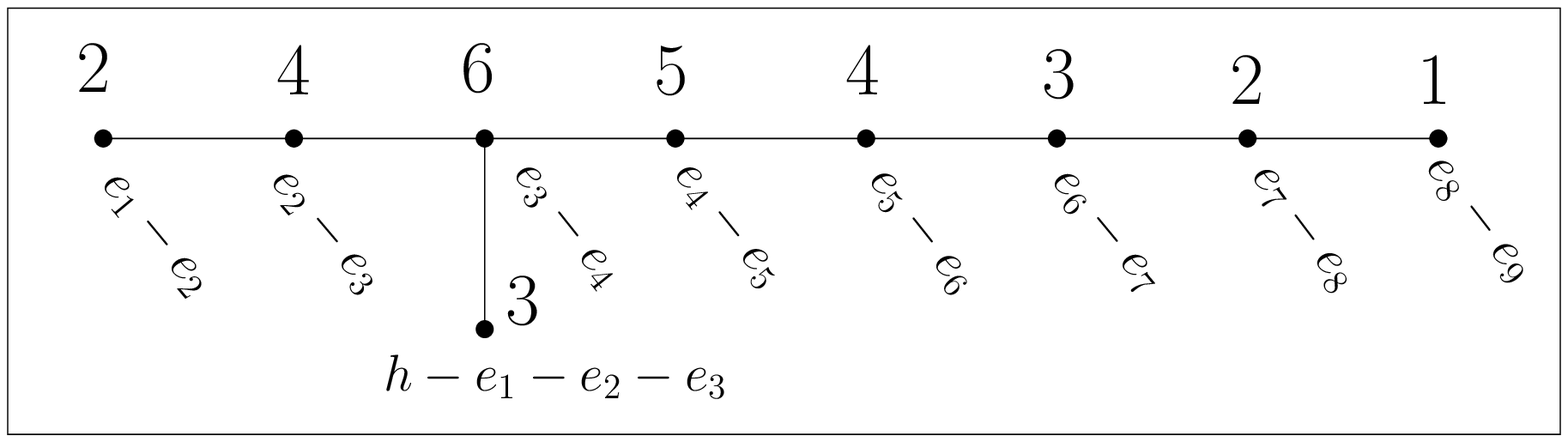}
\end{minipage}
\end{center}
\begin{center}
Figure \ref{e8sec1}.11
\end{center}

\pagebreak

\section{An Elliptic Fibration on $E(1)$ with an $\tilde{E_8}$ Fibre} \label{e8algeomsection}

The method used in this section can be found in \cite{SS}. \newline

We explicitly construct an elliptic fibration on $E(1)$ with one $\tilde{E_8}$ fibre and two fishtail fibres. We do the calculations in detail in order in order to show how such constructions are done. Again, in this section there are times when we shall not distinguish between a polynomial and the curve to which it corresponds.\newline

The construction in section \ref{e8sec1} will guide our choice of cubics. We choose one cubic to be $p_1(x,y,z) = zy^2 - zx^2 - x^3$, the ``prototype'' of the fishtail fibre. We now need to choose another cubic $p_0(x,y,z)$ that intersects $p_1(x,y,z)$ in exactly one point, and for ease of calculation, we would like this cubic to be fairly simple. Furthermore, this cubic must be the perfect cube of a homogeneous polynomial of degree one (since it must represent a line with multiplicity 3) that intersects the curve $C_1$ exactly once (because of Lemma \ref{fibseclemssslab} and the fact that $\tilde{E_8}$ has only one section). We try $p_0(x,y,z) = z^3$. We define the curves
\begin{align}
C_0 &= \{[x:y:z] \in \cpt | p_0(x,y,z) = z^3 = 0 \} \\
C_1 &= \{[x:y:z] \in \cpt | p_1(x,y,z) = zy^2 - zx^2 - x^3 = 0 \}
\end{align}
Let us calculate the intersection point(s) of $C_0$ and $C_1$:
\begin{align}
& z^3 = 0 \nonumber \\
\Rightarrow & z = 0 \nonumber \\
& zy^2 - zx^2 - x^3 = 0 \qquad (\mathrm{and} \;\; z = 0) \nonumber \\
\Rightarrow & -x^3 = 0 \nonumber \\
\Rightarrow & x = 0  \nonumber \\
& x= 0 \; \; \mathrm{and} \;\; z =0 \nonumber \\
\Rightarrow & P = [0:1:0] \qquad(\mathrm{since} \; \; [x:y:z]\in \cpt)
\end{align}
and so the only intersection point of $z^3=0$ and $zy^2-zx^2-x^3=0$ is the point $P=[0:1:0]$. This point will have multiplicity 9. The pencil 
\begin{equation*}
C_{[t_0:t_1]} = \{(t_0 p_0 + t_1 p_1)^{-1}(0) \; | \; [t_0:t_1] \in \cpo \}
\end{equation*}
of elliptic curves defined by $C_0$ and $C_1$ provides a map from $\cpt$ to $\cpo$ well-defined away from the point $P = [0:1:0]$. We perform nine (infinitely close) blow-ups at $P$ (as in section \ref{e8sec1}) in order to get the desired elliptic fibration on $E(1)$ with one $\tilde{E_8}$ and one fishtail fibre. Since the Euler characteristic of an $\tilde{E_8}$ fibre is 10 and  that of a fishtail is 1, and since the Euler characteristic of $E(1)$ is 12, we know (from section \ref{kodairasec}), that there must be one other fishtail fibre in this elliptic fibration. We shall now explicitly ``find'' this fibre. \newline

We know that our singular fibres are curves $C_{[t_0:t_1]}$ which correspond to polynomials, where $[t_0:t_1] \in \cpo$. We define
\begin{equation*}
p_t(x,y,z) = t_0 p_0(x,y,z) + t_1 p_1(x,y,z) \qquad (t = [t_0:t_1], [x:y:z] \in \cpt)
\end{equation*} 
Since we are working in $\cpt$, in order to find singular points using the Implicit Function Theorem, we need to work in the charts $[1:y:z]$, $[x:1:z]$ and $[x:y:1]$ (which we shall call charts 1, 2 and 3, respectively), which are all isomorphic to $\mathbb{C}^2$. For example, working in the chart $z=1$, where each point is of the form $[x:y:1]$, our polynomial $p_t$ becomes
\begin{equation*}
p_t(x,y) = t_0 p_0(x,y) + t_1 p_1(x,y) \qquad ((x,y) \in \mathbb{C}^2)
\end{equation*}
and then a singular point $(x_0,y_0) \in \mathbb{C}^2$ is a point that satisfies the following three equations:
\begin{align*}
p_t(x_0,y_0) = 0 \\
\frac{\partial p_t}{ \partial x}(x_0,y_0) = 0 \\
\frac{\partial p_t}{ \partial y}(x_0,y_0) = 0
\end{align*}
and then we say $[x_0:y_0:1] \in \cpt$ is a singular point of $p_t$. We see that the polynomial for $[t_0:t_1] = [0:1]$, $p_{[0:1]}(x,y,z) = zy^2 - zx^2 - x^3$, corresponds to the fishtail fibre (see section \ref{fishtailsection}), and the polynomial for $[t_0:t_1] = [1:0]$, $p_{[1:0]}(x,y,z) = z^3$, corresponds to the $\tilde{E_8}$ fibre (and every point $[x:y:0]$ is a singular point). \newline

Since we have covered the cases when either $t_0 = 0$ or $t_1 = 0$, we now look for all the polynomials $p_t$ that have singular points and $t_0 \neq 0$ and $t_1 \neq 0$. Let us return to our example
\begin{align}
p_0(x,y,z) &=z^3 \\
p_1(x,y,z) &= zy^2 - zx^2 - x^3 \\
p_t(x,y,z) &= t_0z^3 + t_1(zy^2 - zx^2 - x^3) \qquad (t = [t_0:t_1] \in \cpo) \label{e8polyform}
\end{align}
from now on, we shall just write $p_t$ as $p$.

\subsection*{Chart 1} 
Let us start\footnote{The reader could skip ahead to chart 3.} by looking for singular points in chart 1, $[1:y:z]$. We need to find points $(y,z) \in \mathbb{C}^2$ that satisfy the following three equations
\begin{align}
p(y,z) &= t_0z^3 + t_1(zy^2 - z - 1)  = 0 \label{e8peq}\\
\frac{\partial p}{\partial y}(y,z) &= t_1 2zy  = 0 \label{e8pyeq}\\ 
\frac{\partial p}{\partial z}(y,z) &= t_0 3z^2 + t_1(y^2 -1) =0 \label{e8pzeq}
\end{align}
In order for equation $\eqref{e8pyeq}$ to hold, we need
\begin{itemize}
\item[(i)] $z = 0$, or
\item[(ii)] $y = 0$
\end{itemize}
since $t_1 \neq 0$. In case (i), if $z=0$, then equation $\eqref{e8peq}$ implies
\begin{equation*}
p(y,0) = t_1(- 1) = 0
\end{equation*}
which is a contradiction, since $t_1 \neq 0$. \newline

In case (ii), if $y=0$, then equation $\eqref{e8pzeq}$ implies
\begin{align}
t_0 3z^2 - t_1 &= 0 \nonumber \\
\Rightarrow  -t_1 &= -t_0 3z^2 \nonumber \\
\Rightarrow  \frac{t_1}{t_0} &= 3z^2 \label{e8t0t11}
\end{align}
and with $y=0$ equation $\eqref{e8peq}$ becomes
\begin{align*}
& t_0z^3 + t_1(- z - 1)  = 0 \\
\Rightarrow & z^3 + \frac{t_1}{t_0}(- z - 1)  = 0 \\
\end{align*}
and substituting in equation $\eqref{e8t0t11}$, this becomes
\begin{align*}
& z^3 + 3z^2(- z - 1)  = 0 \\
\Rightarrow & z^3 - 3z^3 - 3z^2 = 0 \\
\Rightarrow & -2z^3 - 3z^3 = 0 \\
\Rightarrow & -2z^2 \Big(z + \frac{3}{2} \Big) = 0 \\
\Rightarrow & z=0 \; \; \mathrm{or} \; \; z = -\frac{3}{2}
\end{align*}
The case $z=0$ was shown above to lead to a contradiction. However, substituting $z = -\frac{3}{2}$ back into $\eqref{e8t0t11}$ gives
\begin{align*}
& \frac{t_1}{t_0} = 3(-\frac{3}{2})^2 \\
\Rightarrow & \frac{t_1}{t_0} = \frac{27}{4} \\
\Rightarrow & \frac{t_0}{t_1} = \frac{4}{27} \\
\Rightarrow & [t_0:t_1] = [\frac{4}{27}:1]
\end{align*}
and it can be checked the point $[1:0:-\frac{3}{2}]$ is indeed a singular point of $\frac{4}{27}z^3 + zy^2 - zx^2 - x^3 = 0$. \newline 

Although we were only expecting to find one more singular curve with a single singular point (a fishtail fibre), we check the other charts for completeness sake.

\subsection*{Chart 2}
In chart 2, $y=1$ and every point is of the form $[x:1:z]$. We need to find points $(x,z)\in \mathbb{C}^2$ that satisfy the following three equations
\begin{align}
p(x,z) &= t_0z^3 + t_1(z - zx^2 - x^3) = 0 \label{e8peq2}\\
\frac{\partial p}{\partial x}(x,z) &= t_1 (-2zx - 3x^2)  = 0 \label{e8pxeq2}\\ 
\frac{\partial p}{\partial z}(x,z) &= t_0 3z^2 + t_1(1-x^2) =0 \label{e8pzeq2}
\end{align}
From equation $\eqref{e8pxeq2}$ we get
\begin{align*}
& -2zx - 3x^2 =0\\
\Rightarrow & -2x \Big( z+\frac{3}{2}x  \Big) 
\end{align*}
which gives us two cases
\begin{itemize}
\item[(i)] $x = 0$, or
\item[(ii)] $z =-\frac{3}{2}x$
\end{itemize}
In case (i), if $x=0$, then equation $\eqref{e8pzeq2}$ becomes
\begin{align}
& t_0 3z^2 + t_1 =0 \nonumber \\
\Rightarrow & t_1 = -t_0 3z^2 \nonumber \\
\Rightarrow & \frac{t_1}{t_0} = -3z^2 \label{e8alphaeq2}
\end{align}
and equation $\eqref{e8peq2}$ becomes
\begin{align*} 
& t_0z^3 + t_1z = 0 \\
\Rightarrow & z^3 + \frac{t_1}{t_0}z = 0 \\
\Rightarrow & z^3 + (- 3z^2) z = 0 \qquad \mathrm{using} \;\; \eqref{e8alphaeq2} \\
\Rightarrow & -2z^3 = 0 \\
\Rightarrow & z =0
\end{align*}
and $z=0$ in equation  \eqref{e8alphaeq2} implies $t_1 = 0$, which contradicts our choice of $t_1 \neq 0$. \newline

In case (ii), $z =-\frac{3}{2}x$, equation $\eqref{e8pzeq2}$ becomes
\begin{align}
& t_0 3 \Big( -\frac{3}{2}x \Big) ^2 + t_1(1-x^2) =0 \nonumber \\
\Rightarrow & \frac{t_0}{t_1} \frac{27}{4}x^2 + 1-x^2 =0 \nonumber \\
\Rightarrow & \Big( \frac{t_0}{t_1} \frac{27}{4} - 1 \Big)x^2 = -1 \nonumber \\
\Rightarrow & \frac{t_0}{t_1} \frac{27}{4} - 1 = -\frac{1}{x^2} \nonumber \\
\Rightarrow & \frac{t_0}{t_1} \frac{27}{4} = -\frac{1}{x^2} + 1 \label{e8alphaeq3}
\end{align}
and equation $\eqref{e8peq2}$ becomes
\begin{align*}
& t_0 \Big(-\frac{3}{2}x \Big)^3 + t_1\Bigg( \Big(-\frac{3}{2}x \Big) - \Big(-\frac{3}{2}x \Big) x^2 - x^3 \Bigg) = 0 \\
\Rightarrow & t_0 \Big(-\frac{27}{8} \Big) x^3 + t_1\Bigg( -\frac{3}{2}x  +  \frac{3}{2}x^3 - x^3 \Bigg) = 0 \\
\Rightarrow & -\frac{t_0}{t_1} \frac{27}{8} x^3 -\frac{3}{2}x  +  \frac{1}{2}x^3 = 0 \\
\Rightarrow & x \Bigg( -\frac{t_0}{t_1} \frac{27}{8} x^2 -\frac{3}{2}  +  \frac{1}{2}x^2 \Bigg) = 0 \\
\Rightarrow & -\frac{t_0}{t_1} \frac{27}{8} x^2 -\frac{3}{2}  +  \frac{1}{2}x^2 = 0 \qquad (x \neq 0) \\
\Rightarrow & x^2 \Bigg(-\frac{t_0}{t_1} \frac{27}{8} + \frac{1}{2} \Bigg) = \frac{3}{2} \\
\Rightarrow & x^2 \Bigg(-\frac{1}{2}\frac{t_0}{t_1} \frac{27}{4} + \frac{1}{2} \Bigg) = \frac{3}{2}\\
\Rightarrow & x^2 \Bigg(-\frac{1}{2} \Big(-\frac{1}{x^2} + 1 \Big) + \frac{1}{2} \Bigg) = \frac{3}{2} \qquad (\mathrm{using} \; \; \eqref{e8alphaeq3})
\end{align*}
\begin{align*}
\Rightarrow & x^2 \Bigg(\frac{1}{2}\frac{1}{x^2} + \Big(-\frac{1}{2}\Big) + \frac{1}{2} \Bigg) = \frac{3}{2} \\
\Rightarrow & x^2 \Bigg(\frac{1}{2}\frac{1}{x^2} \Bigg) = \frac{3}{2} \\
\Rightarrow & \frac{1}{2} = \frac{3}{2}
\end{align*}
which is a contradiction. Therefore, there are no singular points in chart 2.

\subsection*{Chart 3}
In chart 3, $z=1$ and every point is of the form $[x:y:1]$. We need to find points $(x,y)\in \mathbb{C}^2$ that satisfy the following three equations
\begin{align}
p(x,y) &= t_0 + t_1(y^2 - x^2 - x^3) = 0 \label{e8peq3}\\
\frac{\partial p}{\partial x}(x,y) &= t_1 (-2x - 3x^2)  = 0 \label{e8pxeq3}\\ 
\frac{\partial p}{\partial y}(x,y) &= t_1(2y) =0 \label{e8pyeq3}
\end{align}
Equation $\eqref{e8pyeq3}$ implies that $y=0$, since $t_1 \neq 0$. Equation $\eqref{e8pxeq3}$ gives
\begin{align*}
& -2x - 3x^2  = 0 \\
\Rightarrow & -x(2 + 3x) = 0 \\
\Rightarrow & x = 0 \; \; \mathrm{or} \;\; x = -\frac{2}{3}
\end{align*}
If $x=0$, since $y=0$ equation $\eqref{e8peq3}$ becomes
\begin{equation*}
t_0 = 0
\end{equation*}

which contradicts our choice of $t_0 \neq 0$. If $x = -\frac{2}{3}$, then equation $\eqref{e8peq3}$ becomes

\begin{align}
& t_0 + t_1\Bigg(- \Big(-\frac{2}{3} \Big)^2 - \Big(-\frac{2}{3} \Big)^3 \Bigg) = 0 \\
\Rightarrow & t_0 + t_1\Bigg(- \frac{4}{9}  + \frac{8}{27} \Bigg) = 0 \\
\Rightarrow & t_0 - \frac{4}{27}t_1 = 0 \\
\Rightarrow & t_0 = \frac{4}{27}t_1
\end{align}
which is the point $[\frac{4}{27}:1] \in \cpo$. We expected this, since $[-\frac{2}{3}:0:1] \sim [1:0:-\frac{3}{2}] \in \cpt$ (which we calculated is a singular point in chart 1). Therefore, if both $t_0 \neq 0$ and $t_1 \neq 0$, there is only one polynomial of the form as in $\eqref{e8polyform}$, namely
\begin{equation}
p(x,y,z) = \frac{4}{27} z^3 + zy^2 - zx^2 - x^3
\end{equation}
which corresponds to the singular curve
\begin{equation}
C_3 = \{[x:y:z] \in \cpt | \frac{4}{27} z^3 + zy^2 - zx^2 - x^3 = 0\} \label{thirdfflab}
\end{equation}
which has a singular point at $[-\frac{2}{3}:0:1]\in \cpt$.

\subsection*{$C_3$ is a fishtail fibre}
We proceed as in the proof of Proposition \ref{fishtailprop}. We first note that the calculations above show that $C_3$ has $P=[-\frac{2}{3}:0:1]$ as its only singular point. The space of projective lines through the point $P=[-\frac{2}{3}:0:1]$ consists of line of the form 
\begin{equation*}
ax + by + cz = 0 \qquad [a:b:c] \in \cpt
\end{equation*}
that satisfy 
\begin{align*}
& a \Big(-\frac{2}{3} \Big) + b(0) + c(1)= 0 \\
\Rightarrow & c = \frac{2}{3} a
\end{align*}
and so we can parametrize this space of lines by $[u_0:u_1] \in \cpo$ by 
\begin{equation*}
L_{[u_0:u_1]} = \{[x:y:z] \in \cpt | u_0x + u_1y + \frac{2}{3}u_0 z = 0 \}
\end{equation*}
Defining
\begin{equation}
\alpha = \frac{u_0}{u_1} \label{alphau0u1lab}
\end{equation}
we have
\begin{align}
u_0x &+ u_1y + \frac{2}{3}u_0 z = 0 \nonumber \\
\Rightarrow u_1 y &= -u_0x - u_0\frac{2}{3}z \nonumber \\
\Rightarrow y &= -\frac{u_0}{u_1} \Big(x + \frac{2}{3}z \Big) \nonumber \\
\Rightarrow y &= -\alpha \Big(x + \frac{2}{3}z \Big) \label{yeqne8}
\end{align}
and then substituting this into our polynomial $\frac{4}{27} z^3 + zy^2 - zx^2 - x^3 = 0$, we get
\begin{align}
& \frac{4}{27} z^3 + zy^2 - zx^2 - x^3 = 0 \nonumber \\
\Rightarrow & \frac{4}{27} z^3 + z\Bigg( -\alpha \Big(x + \frac{2}{3}z \Big) \Bigg)^2 - zx^2 - x^3 = 0 \nonumber \\
\Rightarrow & \frac{4}{27} z^3 + \alpha ^2 z \Big(x + \frac{2}{3}z \Big)^2 - zx^2 - x^3 = 0 \label{zxfact1}
\end{align}
We shall now try and find a factor of $\Big( x+\frac{2}{3}z \Big)^2$ in the above line:
\begin{align*}
\frac{4}{27} z^3 - zx^2 - x^3 &= \frac{4}{27} z^3 - x \Big( x^2 + xz \Big) \\
&= \frac{4}{27} z^3 - x \Big( x^2 + \frac{4}{3}xz + \frac{4}{9}z^2 \Big) + \frac{1}{3}x^2 z + \frac{4}{9}xz^2 \\
&= \frac{1}{3}x^2 z + \frac{4}{9}xz^2 + \frac{4}{27} z^3 - x \Big(x + \frac{2}{3}z \Big)^2 \\
&= \frac{1}{3}z \Big( x^2 + \frac{4}{3}xz + \frac{4}{9} z^2 \Big) - x \Big(x + \frac{2}{3}z \Big)^2 \\
&= \frac{1}{3}z \Big( x + \frac{2}{3} z \Big)^2 - x \Big(x + \frac{2}{3}z \Big)^2 \\
&= \Big( x + \frac{2}{3} z \Big)^2 \Big(\frac{1}{3}z  - x \Big)
\end{align*}
and therefore $\eqref{zxfact1}$ becomes
\begin{align*}
& \Big( x+\frac{2}{3} z \Big)^2 \Big(\frac{1}{3}z  - x \Big) + \alpha ^2 z \Big(x + \frac{2}{3}z \Big)^2 = 0 \\
\Rightarrow & \Big( x+\frac{2}{3} z \Big)^2 \Big(\frac{1}{3}z  - x + \alpha ^2 z \Big) = 0 \\
\Rightarrow & x = -\frac{2}{3}z \; \; \mathrm{or} \;\; x= \Big( \frac{1}{3} + \alpha ^2 \Big) z 
\end{align*}
$x = -\frac{2}{3}z$ corresponds to the point $P = [-\frac{2}{3}:0:1]$, but for the second value of $x$, we have from equation $\eqref{yeqne8}$ we have
\begin{align*}
y &= -\alpha \Bigg( \Big( \frac{1}{3} + \alpha ^2 \Big) z + \frac{2}{3}z \Bigg) \\
\Rightarrow y &= -\alpha ( 1 + \alpha ^2)z
\end{align*}
and so we have the point (where we recall that $\alpha = \frac{u_0}{u_1}$ for $[u_0:u_1] \in \cpo$)
\begin{equation*}
Q_{[u_0:u_1]} = [\frac{1}{3} + \Big(\frac{u_0}{u_1}\Big) ^2 : -\frac{u_0}{u_1} ( 1 + \Big(\frac{u_0}{u_1}\Big) ^2) : 1 ]
\end{equation*}
It is perhaps worthwhile to recap what we have just done and try to interpret the results of our calculations. We have just calculated the intersection points of $C_3$ and a line $L_{[u_0:u_1]}$. Our solution $x = -\frac{2}{3}z$ shows that each line $L_{[u_0:u_1]}$ passes through the point $P=[-\frac{2}{3}:0:1]$. We expected this solution, since we are considering the pencil of all curves passing through $P$. The point $Q_{[u_0:u_1]}$ is the other intersection point of $C_3$ and $L_{[u_0:u_1]}$. \newline

Now, we notice that if $\Big( \frac{u_0}{u_1} \Big)^2 = -1$ then $Q_{[u_0:u_1]} = P = [-\frac{2}{3}:0:1]$. Therefore, for the points $[i:1]$ and $[-i:1]$, $Q_{[i:1]} = Q_{[-i:1]} = P$. So, $L_{[u_0:u_1]} \cap C_3 = \{P, Q_{[u_0:u_1]} \}$, except if $[u_0:u_1] = [\pm i:1]$, in which case $L_{[i:1]} \cap C_3 = L_{[-i:1]} \cap C_3 = \{P\}$. We define the map $\psi : \cpo \longrightarrow C_3$ by $\psi ([u_0:u_1]) = Q_{[u_0:u_1]}$, and we define $Q_{[\pm i: 1]} = P$. We have already shown that $C_3$ has only one singular point, and this calculation shows that $C_3$ is homeomorphic to $\cpo$ with two points identified. We can therefore conclude that $C_3$ is a fishtail fibre. $\Box$

\pagebreak

\section{A Construction of an $\tilde{E_6}$ Fibre} \label{e6constrsec}

In this section, we use the techniques explained in section \ref{e8sec1} to give a construction of an $\tilde{E_6}$ fibre. \newline

By Remark \ref{e8secrem1rem} and Lemma \ref{fibseclemssslab}, we need to choose two cubics $p_0$ and $p_1$ that intersect transversally in exactly three points. We choose a cubic $p_1$ so that the curve corresponding to $p_1$ will have homology class $3h$, with multiplicity 1. As in section \ref{e8sec1}, we shall label the curves with their homology classes, so we shall call the curve corresponding to $p_1$ the curve $3h$. Next, we choose a projective line $L$ that intersects $3h$ in exactly three points in $\cpt$. The line $L$ corresponds to a linear polynomial $l$ which has homology class $h$. We then choose the cubic $p_0$ to be $p_0 = l^3$, and therefore the curve corresponding to $p_0$ has homology class $h$ with multiplicity 3 (and we shall also call this curve $h$). We label the intersection points $P_1$, $P_2$ and $P_3$. Each intersection point has multiplicity 3. See Figure \ref{e6constrsec}.1. \newline

We start by blowing-up at $P_1$. This introduces the exceptional curve $e_1$, and the proper transform of $h$ is $h-e_1$ and the proper transform of $3h$ is $3h-e_1$. Therefore, $e_1$ has multiplicity 2. See Figure \ref{e6constrsec}.2. \newline

The second and third blow-ups (at $P_2$ and $P_3$, respectively) are similar to the first, each introducing the exceptional curves $e_2$ and $e_3$, respectively, both with multiplicity 2, and the proper transforms of the curves corresponding to $p_0$ and $p_1$ are now $h-e_1-e_2-e_3$ (with multiplicity 3) and $3h-e_1-e_2-e_3$ (with multiplicity 1), respectively. It can be calculated that these two curves no longer intersect each other. This gives us Figure \ref{e6constrsec}.3. \newline

We perform the fourth blow-up at $P_1$, the fifth blow-up at $P_2$ and the sixth blow-up at $P_3$. These blow-ups introduce the exceptional curves $e_4$, $e_5$ and $e_6$, which each have multiplicity 1. The proper tranforms of $e_1$, $e_2$ and $e_3$ are $e_1-e_4$, $e_2-e_5$ and $e_3-e_6$, respectively, and the proper transform of $3h-e_1-e_2-e_3$ is $3h-e_1 - \dots - e_6$. See Figure \ref{e6constrsec}.4. \newline

\pagebreak

We perform the seventh blow-up at $P_1$, introducing the exceptional curve $e_7$, the eighth blow-up at $P_3$ introducing $e_8$ and the ninth blow-up at $P_2$, introducing $e_8$. This unexpected ordering of blow-ups is done so that the homology classes of the resulting $-2$-spheres will be the same as in \cite{P1}. The proper transforms of $e_4$, $e_5$ and $e_6$ are $e_4-e_7$, $e_5-e_9$ and $e_6-e_8$, respectively. These are therefore $-2$-spheres each with multiplicity 1. The proper transform of $3h-e_1-\dots-e_6$ is $3h-e_1-\dots-e_6-e_7-e_8-e_9$ (still with multiplicity 1), and therefore the three exceptional curves $e_7$, $e_8$ and $e_9$ each have multiplicity 0, and are therefore sections of the $\tilde{E_6}$ fibre. See Figures \ref{e6constrsec}.5 and \ref{e6constrsec}.6.

\pagebreak

\begin{center}
\begin{minipage}{12cm}
\includegraphics[width=12cm]{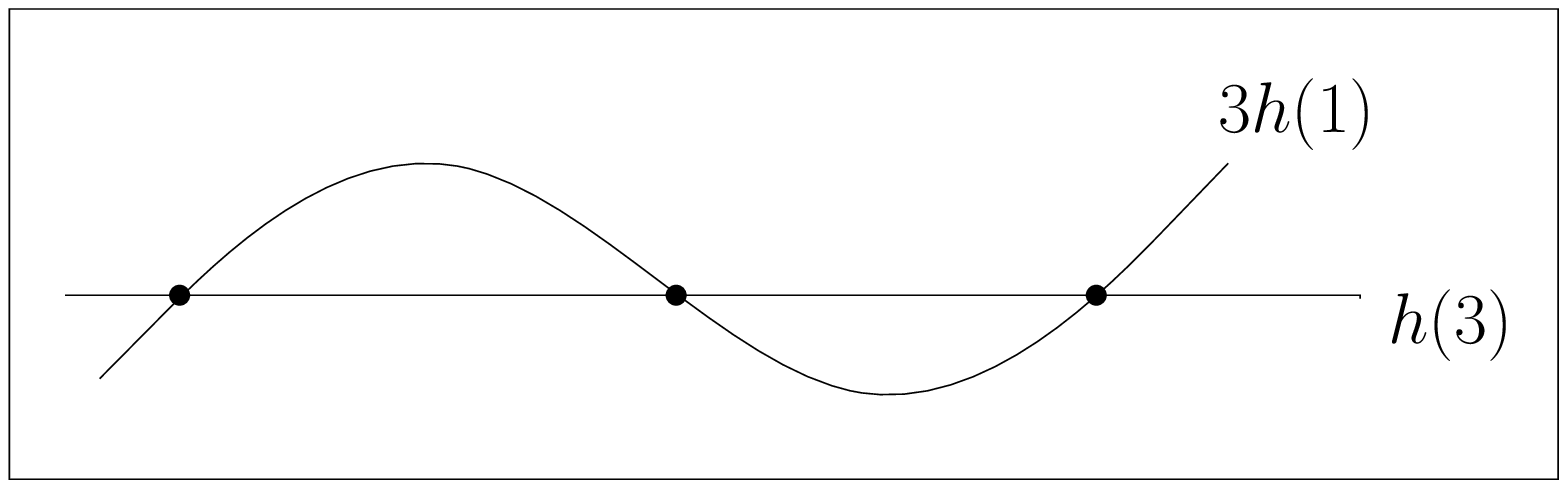}
\end{minipage}
\end{center}
\begin{center}
Figure \ref{e6constrsec}.1
\end{center}

\begin{center}
\begin{minipage}{12cm}
\includegraphics[width=12cm]{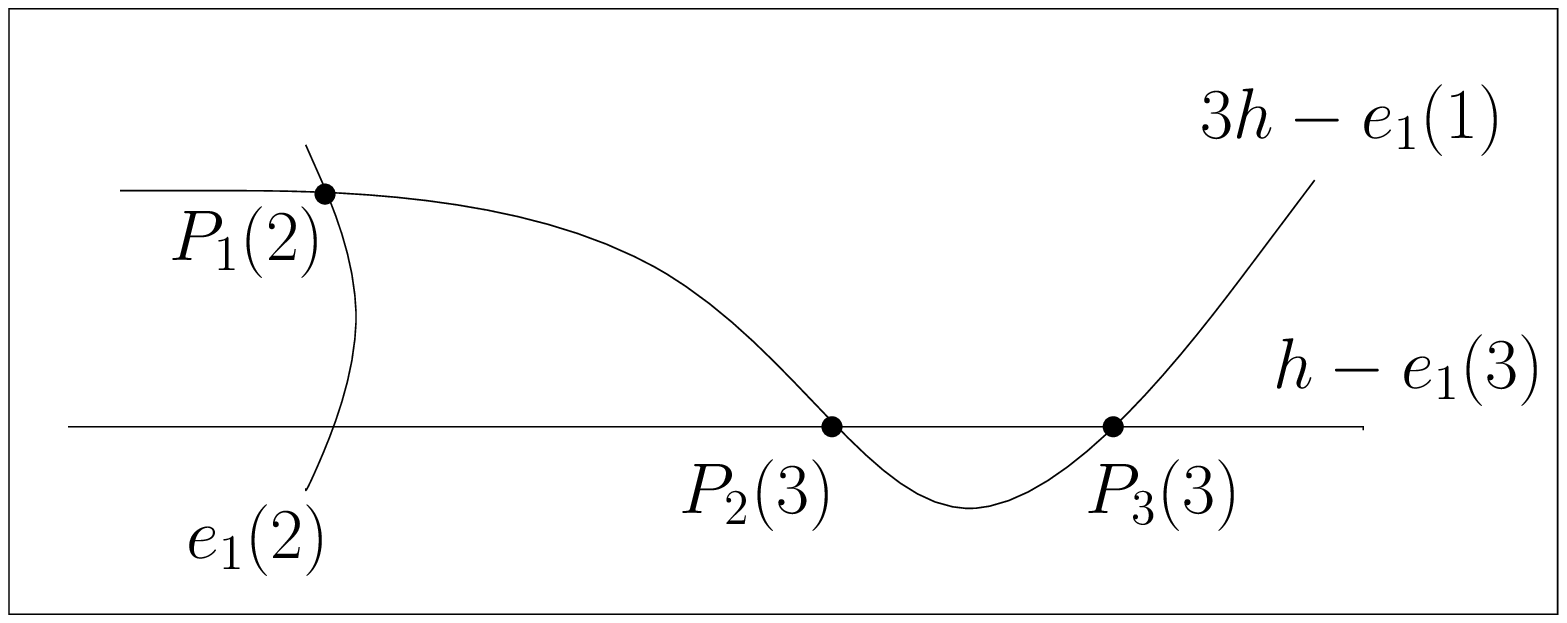}
\end{minipage}
\end{center}
\begin{center}
Figure \ref{e6constrsec}.2
\end{center}

\begin{center}
\begin{minipage}{12cm}
\includegraphics[width=12cm]{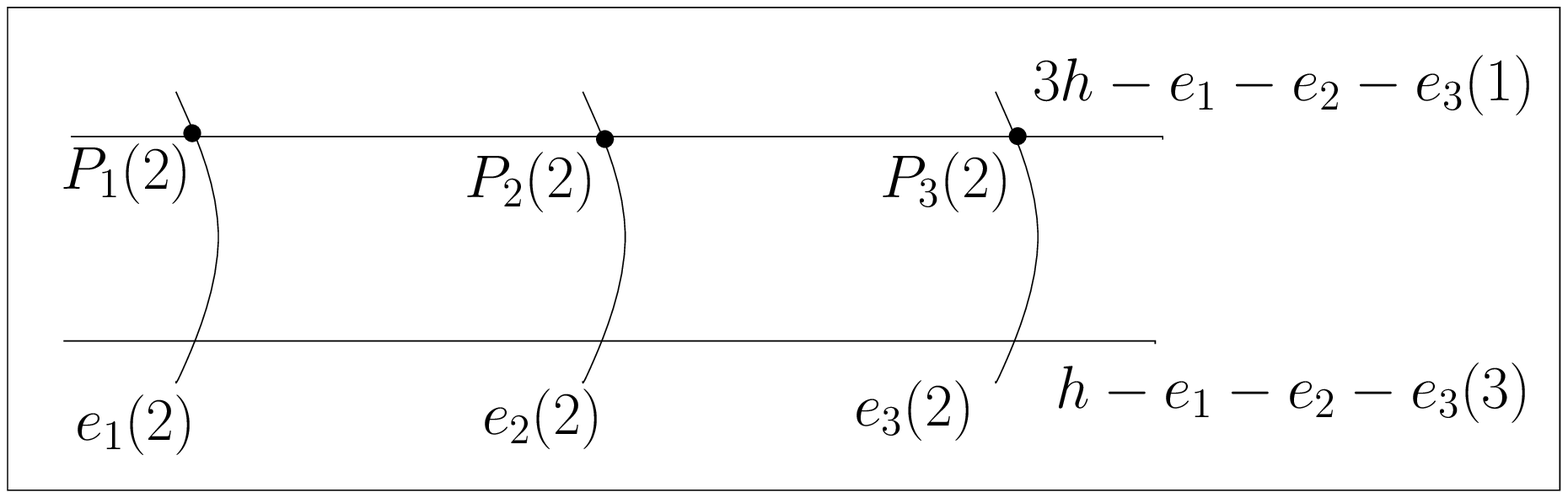}
\end{minipage}
\end{center}
\begin{center}
Figure \ref{e6constrsec}.3
\end{center}

\begin{center}
\begin{minipage}{12cm}
\includegraphics[width=12cm]{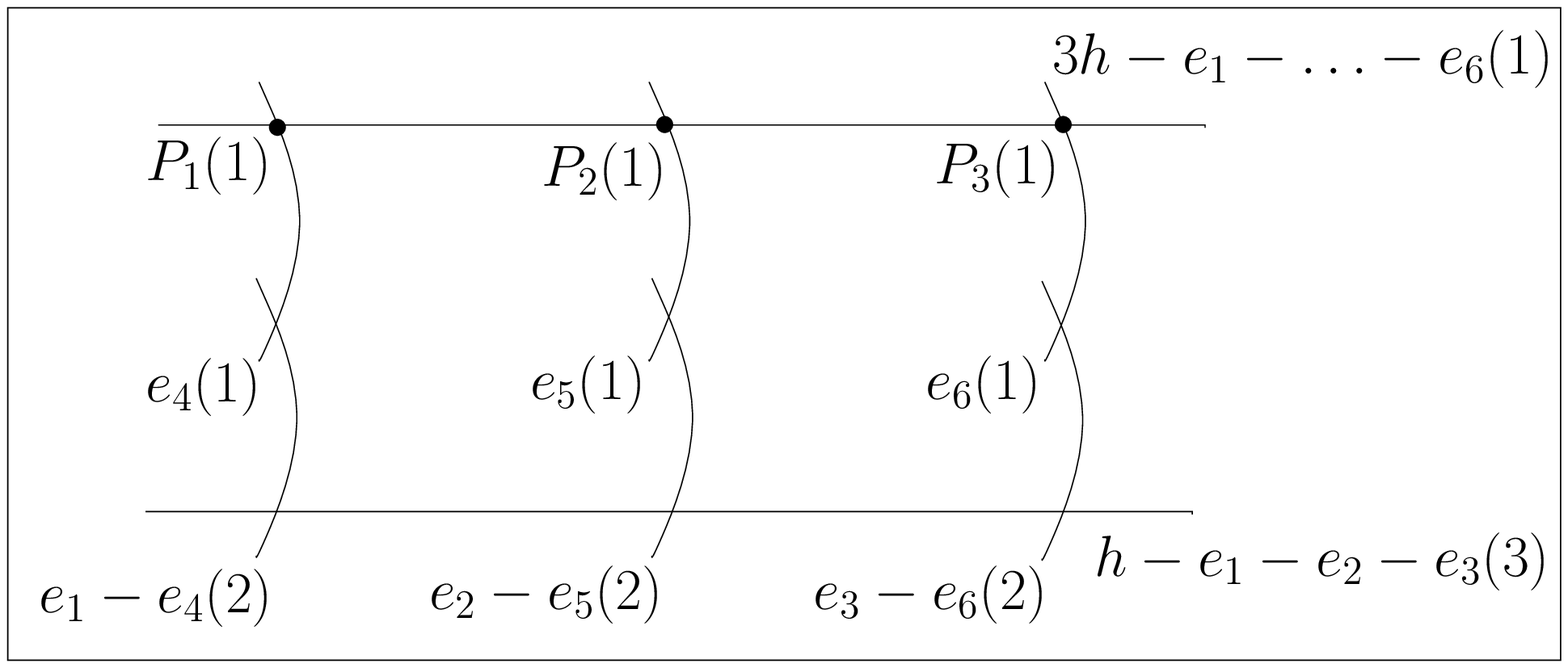}
\end{minipage}
\end{center}
\begin{center}
Figure \ref{e6constrsec}.4
\end{center}

\begin{center}
\begin{minipage}{12cm}
\includegraphics[width=12cm]{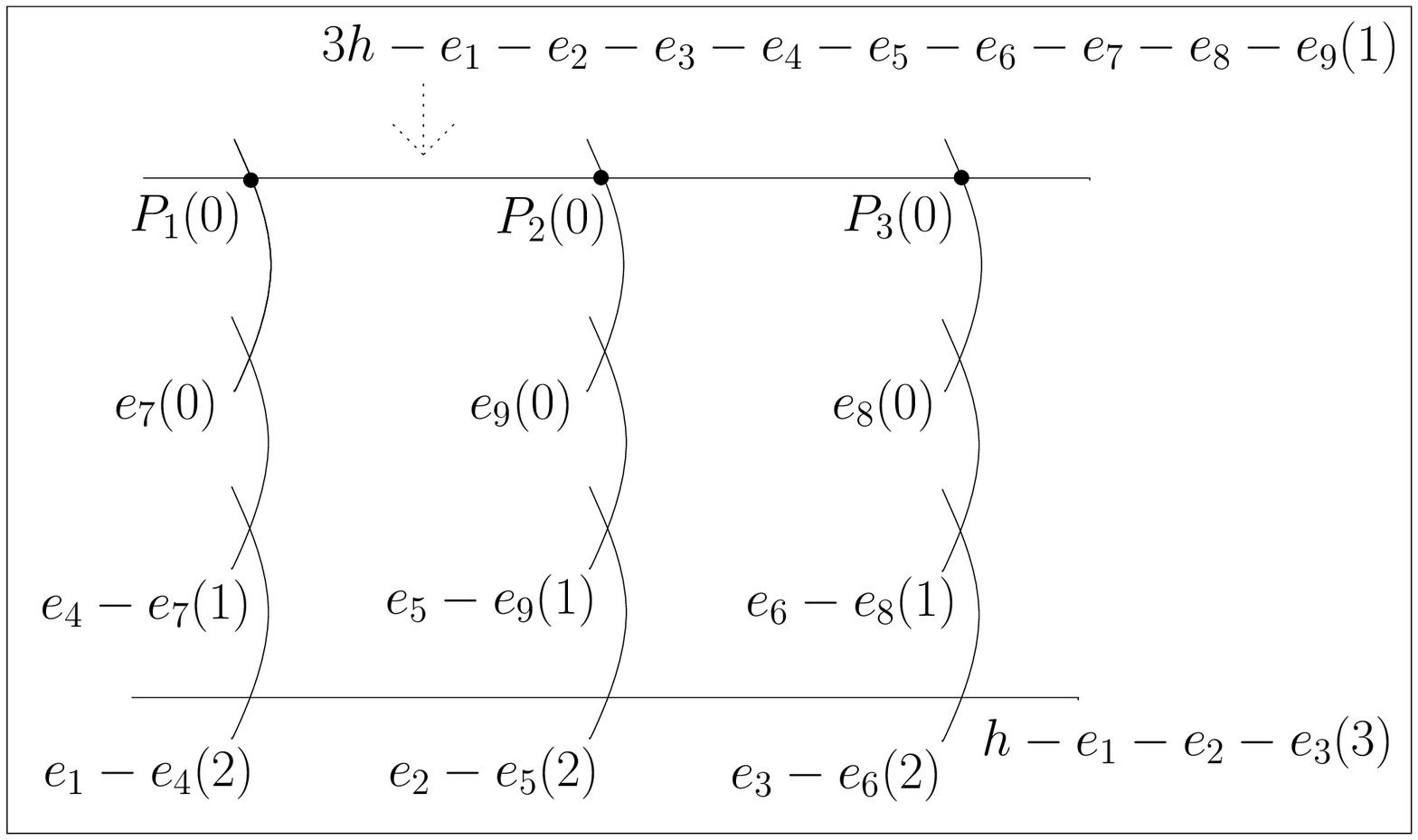}
\end{minipage}
\end{center}
\begin{center}
Figure \ref{e6constrsec}.5
\end{center}

\begin{center}
\begin{minipage}{12cm}
\includegraphics[width=12cm]{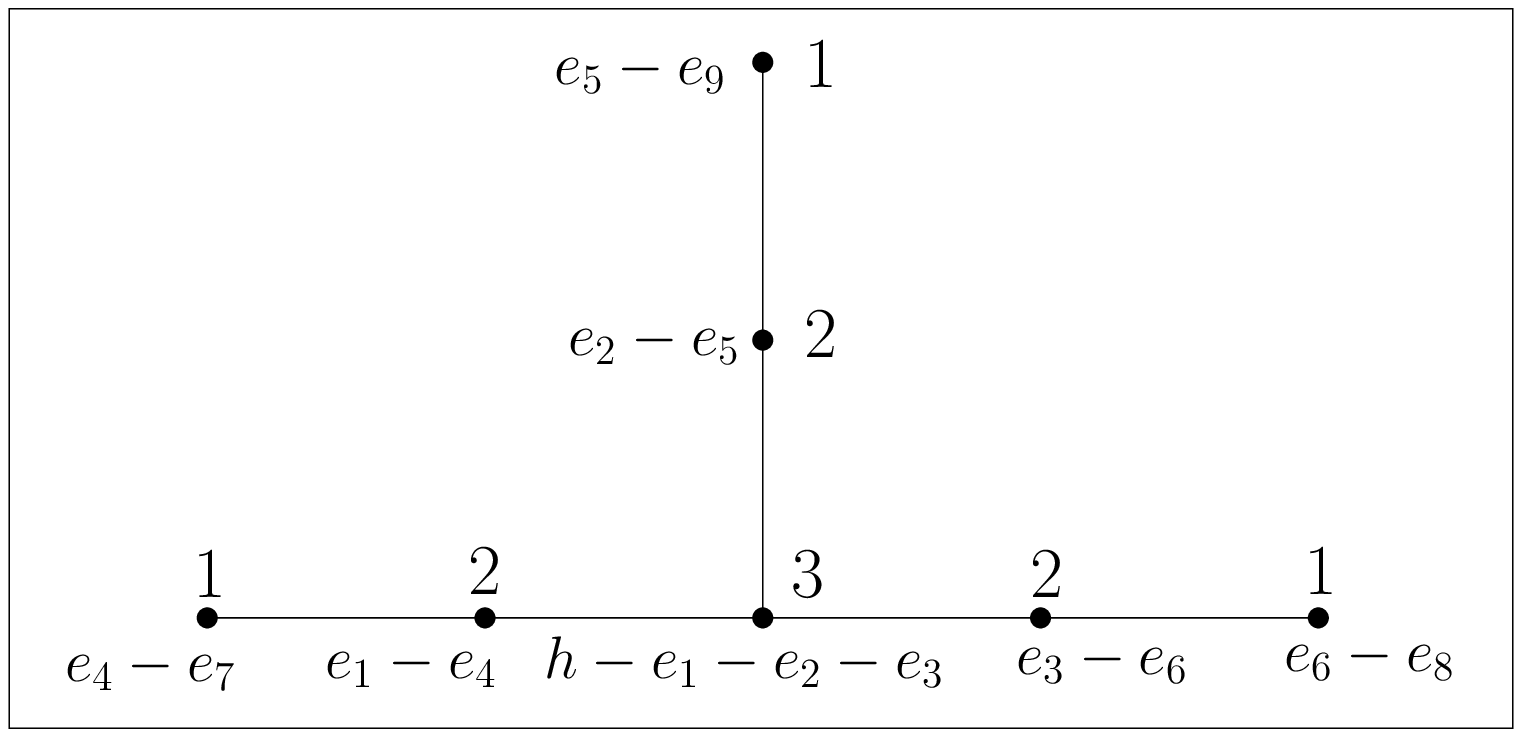}
\end{minipage}
\end{center}
\begin{center}
Figure \ref{e6constrsec}.6
\end{center}

\pagebreak


\section{An Elliptic Fibration on $E(1)$ with an $\tilde{E_6}$ Fibre} \label{explicitE6sec}

We explicitly construct an elliptic fibration on $E(1)$ with one $\tilde{E_6}$ fibre and four fishtail fibres.\newline

We follow the construction in section \ref{e6constrsec}, where we started with two cubics that intersected each other transversally in three points, and then blew up three times at each of those three points. \newline

We know from \cite{GS} that a fishtail fibre is ambiently isotopic to
\begin{equation}
C_1 = \{ [x:y:z] \in \cpt | zy^2 = x^3 + zx^2 \} \label{fishtaillab}
\end{equation}
so, we shall take one cubic to be 
\begin{equation}
p_1(x,y,z) = zy^2 - zx^2- x^3 \label{fish1lab}
\end{equation}

We now need to choose a cubic $p_0$ such that $p_0$ and $p_1$ intersect transversally in three points, which means that the set 
\begin{equation*}
\{P \in \cpt | p_0(P) = 0 \; and \; p_1(P) = 0 \}
\end{equation*}
contains exactly three distinct points, none of them tangent points. \newline

Furthermore, we need $p_0$ to be a perfect cube, i.e. a polynomial of the form
\begin{equation*}
p_0(x,y,z) = (ax + by + cz)^3
\end{equation*}
for some constants $a,b,c \in \mathbb{C}$. \newline

We shall choose $p_0(x,y,z) = (y + \frac{1}{2}z)^3$ to be $p_0$. See Remark \ref{choiceremlab}. So, we have
\begin{align}
p_0(x,y,z) &= (y + \frac{1}{2}z)^3 \label{p0lab} \\
p_1(x,y,z) &= zy^2 - zx^2- x^3 \label{p1lab}
\end{align}

\subsection*{The intersection points}

Firstly,

\begin{align}
& p_0(x,y,z) = 0 \nonumber \\
\Rightarrow & (y + \frac{1}{2}z)^3 = 0 \nonumber \\
\Rightarrow & y + \frac{1}{2}z = 0 \nonumber  \\ 
\Rightarrow & y = -\frac{1}{2}z \nonumber  \\
\Rightarrow & y^2 = \frac{1}{4}z^2 \label{p0ylab}
\end{align}

and substituting $\eqref{p0ylab}$ into $\eqref{p1lab}$, we get
\begin{align}
& p_1(x,y,z) = 0 \nonumber \\
\Rightarrow & zy^2 - zx^2- x^3 = 0 \nonumber \\
\Rightarrow & z \Big( \frac{1}{4}z^2 \Big) - zx^2 - x^3 = 0 \nonumber \\
\Rightarrow & z^3 - 4zx^2 - 4x^3 = 0 \label{p1yzlab}
\end{align}

Now, if $x=0$, we get from $\eqref{p1yzlab}$
\begin{align*}
& z^3 - 4zx^2 - 4x^3 = 0 \;\; \mathrm{and} \;\; x=0 \\
\Rightarrow & z^3 = 0 \\
\Rightarrow & z = 0
\end{align*}
and, from $\eqref{p0ylab}$, if $z=0$, then $y=0$. Since $[0:0:0]$ is not a point in $\cpt$, we can assume that $x \neq 0$. Dividing $\eqref{p1yzlab}$ through by $x^3$ and setting $t = \frac{z}{x}$, we get
\begin{align}
& z^3 - 4zx^2 - 4x^3 = 0 \nonumber \\
\Rightarrow & \Big( \frac{z}{x} \Big) ^3 - 4 \Big( \frac{z}{x} \Big) - 4 = 0 \nonumber\\
\Rightarrow & t^3 - 4t - 4 = 0 \nonumber\\ 
\Rightarrow & t^3 - 4t = 4 \label{tcublab}
\end{align}

From \cite{We}, we can explicitly calculate the roots of a cubic equation of the form
\begin{equation*}
x^3 + px = q \qquad p,q \in \mathbb{R}
\end{equation*}

We have the following intermediate variables:
\begin{align}
Q &= \frac{1}{3}p \label{peqlab} \\
R &= \frac{1}{2}q \label{qeqlab} \\
D &= \frac{p^3}{27} + \frac{q}{4} \\
S &= \sqrt[3]{R + \sqrt{D}} \\
T &= \sqrt[3]{R - \sqrt{D}} \\
A &= S+T \\
B &= S-T \label{Beqlab}
\end{align}

Then the three roots $t_1, t_2, t_3$ are:
\begin{align}
t_1 &= B \label{root1lab} \\
t_2 &= -\frac{1}{2} B + i \frac{\sqrt{3}}{2} A \label{root2lab}\\
t_3 &= -\frac{1}{2} B - i \frac{\sqrt{3}}{2} A \label{root3lab}
\end{align} 

For our particular cubic in $\eqref{tcublab}$, $t^3 - 4t = 4$, we have $p = -4$ in $\eqref{peqlab}$ and $q=4$ in $\eqref{qeqlab}$, and so from equations $\eqref{peqlab}$ to $\eqref{Beqlab}$ we get 
\begin{align}
B &= \sqrt[3]{2 + \sqrt{ \frac{44}{27} } } + \sqrt[3]{2 - \sqrt{ \frac{44}{27} } } \\
A &= \sqrt[3]{2 + \sqrt{ \frac{44}{27} } } - \sqrt[3]{2 - \sqrt{ \frac{44}{27} } }
\end{align}

Substituting these values for $B$ and $A$ back into equations $\eqref{root1lab}$, $\eqref{root2lab}$ and $\eqref{root3lab}$, it is simple algebra to verify that these are indeed solutions to our polynomial in $\eqref{tcublab}$. \newline

Or, we could check numerically that 
\begin{align}
t_1 &\approx 2.3830 \\
t_2 &\approx -1.1915 + 0.5089i \\
t_3 &\approx -1.1915 - 0.5089i
\end{align} 
are solutions to $\eqref{tcublab}$. Recalling that $t = \frac{z}{x}$, $y = -\frac{1}{2}z$, $x \neq 0$ and that $[x:y:z] \in \cpt$, we can let $x =1$ and then we get the three intersection points as

\begin{align}
P_1 &= [1: -\frac{1}{2}t_1 :t_1] \\
P_2 &= [1: -\frac{1}{2}t_2 :t_2] \\
P_3 &= [1: -\frac{1}{2}t_3 :t_3]
\end{align} 

Notice that these are the three distinct intersection points between the line $y + \frac{1}{2}z$ and the cubic $zy^2 - zx^2- x^3$, and so must be transverse intersection points. Therefore, the cubic $p_0(x,y,z) = (y + \frac{1}{2}z)^3$ and the cubic $p_1(x,y,z) = zy^2 - zx^2- x^3$ intersection transversely in three points, and so blowing up the intersection points of the two cubic curves corresponding to $p_0$ and $p_1$ will give rise to an elliptic fibration with an $\tilde{E_6}$ fibre (as in section \ref{e6constrsec}), and at least one fishtail fibre. The next calculation shows that there are in fact four fishtail fibres.

\subsection*{The singular points}

We now look at every polynomial of the form
\begin{equation}
p_{[t_0:t_1]} = t_0 p_0 + t_1 p_1
\end{equation}
where $[t_0:t_1] \in \cpo$ and $p_0$, $p_1$ are as in $\eqref{p0lab}$ and $\eqref{p1lab}$. Every polynomial of that form is a cubic which passes through the intersection points $P_1,P_2$ and $P_3$. We wish to find the values of $[t_0:t_1] \in \cpo$ such that $p_{[t_0:t_1]}$ is a singular curve. Again, we use the Implicit Function Theorem. We start with the polynomial
\begin{equation}
p(x,y,z) = t_0 \Big(y+\frac{1}{2}z \Big)^3  + t_1(zy^2-zx^2-x^3) \label{polylab}
\end{equation}

The point $[0:1]$ corresponds to a fishtail fibre (see equation $\eqref{fishtaillab}$ and \cite{GS}) and the point $[1:0]$ corresponds to the polynomial $p_0(x,y,z) = (y + \frac{1}{2}z)^3$ which is singular at every point $[x:-\frac{1}{2}z:z] \in \cpt$ and corresponds to the $\tilde{E_6}$ fibre. We therefore look for points (polynomials) $[t_0:t_1]$ such that $t_0 \neq 0$ and $t_1 \neq 0$, and so defining $\alpha = \frac{t_0}{t_1}$, instead of $\eqref{polylab}$ we can consider the polynomial
\begin{equation}
p(x,y,z) = \alpha \Big(y+\frac{1}{2}z\Big)^3  + zy^2-zx^2-x^3 \label{polylab2}
\end{equation}

Since we are working in $\cpt$, we need to work in charts in order to use the Implicit Function Theorem. We first look in the chart $z=1$, i.e. $[x:y:1]$. Then, the polynomial in $\eqref{polylab2}$ becomes
\begin{equation}
p(x,y,1) = \alpha \Big(y+\frac{1}{2} \Big)^3  + y^2-x^2-x^3
\end{equation}
and the partial derivatives are
\begin{align}
\frac{\partial p}{\partial x}(x,y,1) &= -2x- 3x^2 \label{partialxlab}\\
\frac{\partial p}{\partial y}(x,y,1) &= 3\alpha \Big(y+\frac{1}{2} \Big)^2  + 2y \label{partialylab}
\end{align}

Recall that a singular point is a point $P$ such that
\begin{equation*}
p(P) = \frac{\partial p}{\partial x}(P) = \frac{\partial p}{\partial y}(P) = 0
\end{equation*}

One can check that there are only three values of $\alpha$, and therefore only three points $[t_0:t_1]$, that give polynomials $p_{[t_0:t_1]}$ that each have a single singular point $P$. These values are:
\begin{align}
[t_0:t_1] &= [-\frac{8}{27}:1]  &&P = [0:1:1] \label{singcurve1}\\
[t_0:t_1] &= [-\frac{32}{121}:1] &&P = [-\frac{2}{3}: \frac{4}{3} :1] \label{singcurve2} \\
[t_0:t_1] &= [1:\frac{1}{8}] &&P = [-\frac{2}{3}: -\frac{1}{3} :1] \label{singcurve3}
\end{align}

\pagebreak

\newtheorem{choicerem}{Remark}[section]
\begin{choicerem} \label{choiceremlab}
\upshape
It should be remarked that the polynomial $p_0(x,y,z) = \Big(y + \frac{1}{2} \Big)^3$ was chosen specifically so that the singular points would be rational, to make the calculations easier. The trade-off is that the intersection points are irrational, but since these points are not used in any further calculations (after we have shown that there are three intersection points, we don't use them anymore), this is not a problem.
\end{choicerem}

We now need to determine the type of singular fibre(s) to which these three polynomials correspond. We show that the polynomial given by $\eqref{singcurve1}$ is a fishtail fibre. \newline

Consider the polynomial
\begin{equation}
p(x,y,z) = -\frac{8}{27} \Big(y+\frac{1}{2}z \Big)^3 + zy^2 - zx^2 - x^3 \label{singpolylab}
\end{equation}
and let $C$ be the curve
\begin{equation}
C = \{[x:y:z] \in \cpt | -\frac{8}{27} \Big(y+\frac{1}{2}z \Big)^3 + zy^2 - zx^2 - x^3 = 0 \}
\end{equation}
We know this curve has a singular point $P = [0:1:1]$, and we look at the space of projective lines that pass through $P$. A line of the form
\begin{equation*}
ax+by+cz = 0 \qquad a,b,c \in \mathbb{C}
\end{equation*}
must satisfy
\begin{align*}
& a(0)+b(1)+c(1) = 0 \\
\Rightarrow & c = -b
\end{align*}
and so we can parametrize these lines by $\cpo$ by setting $a=u_0$, $b=u_1$ for $[u_0:u_1]\in \cpo$, and then these lines are of the form
\begin{equation*}
L_{[u_0:u_1]} = \{[x:y:z] \in \cpt | u_0x + u_1 y - u_1 z = 0 \}
\end{equation*}
which shows that for every line of this form,
\begin{equation*}
z = y + \frac{u_0}{u_1}x 
\end{equation*}
and letting 
\begin{equation}
\beta = \frac{u_0}{u_1} \label{betau0lab}
\end{equation}
we have
\begin{equation}
z = y + \beta x \label{zbetalabel}
\end{equation}
Substituting this back into $\eqref{singpolylab}$, we get
\begin{align*}
& p(x,y,z) = 0\\
\Rightarrow & -\frac{8}{27} \Big(y+\frac{1}{2}z \Big)^3 + zy^2 - zx^2 - x^3 =0 \\
\Rightarrow & -\frac{8}{27} \Big(y+\frac{1}{2}(y + \beta x) \Big)^3 + (y + \beta x)y^2 - (y + \beta x)x^2 - x^3 =0 \\
\Rightarrow & -\frac{8}{27} \Big(\frac{3}{2}y + \frac{1}{2}\beta x \Big)^3 + y^3 + \beta xy^2 - x^2 y - \beta x^3 - x^3 =0 \\
\Rightarrow & -\frac{8}{27} \Bigg(\frac{27}{8}y^3 + 3 \Big( \frac{9}{4} \Big) \Big( \frac{1}{2} \Big) \beta xy^2 + 3 \Big( \frac{3}{2} \Big) \Big( \frac{1}{4} \Big) \beta ^2 x^2y + \frac{1}{8} \beta ^3 x^3 \Bigg) \\
& \qquad + y^3 + \beta xy^2 - x^2 y - (\beta + 1) x^3 =0 \\
\Rightarrow & -y^3 - \beta xy^2 - \frac{1}{3} \beta ^2 x^2y - \frac{1}{27} \beta ^3 x^3 + y^3 + \beta xy^2 - x^2 y - (\beta +1) x^3=0 \\
\Rightarrow & -y^3 + y^3 - \beta xy^2 + \beta xy^2 - \frac{1}{3} \beta ^2 x^2y - \frac{1}{27} \beta ^3 x^3 - x^2 y - (\beta +1) x^3=0 \\
\Rightarrow & - \frac{1}{3} \beta ^2 x^2y - \frac{1}{27}\beta ^3 x^3 - x^2 y - (\beta +1) x^3=0 \\
\Rightarrow & \Big(- \frac{1}{3} \beta ^2 - 1 \Big)  x^2y - (\beta +1 +\frac{1}{27}\beta ^3 ) x^3=0
\end{align*}
and dividing through by $y^3$ and letting $t= \frac{x}{y}$, we get
\begin{align*}
& \Big(- \frac{1}{3} \beta ^2-1\Big) \Big(\frac{x}{y}\Big)^2 - (\beta +1 +\frac{1}{27}\beta ^3 ) \Big( \frac{x}{y} \Big)^3=0 \\
\Rightarrow & \Big(- \frac{1}{3} \beta ^2 -1 \Big) t^2 - (\beta +1 +\frac{1}{27}\beta ^3 ) t ^3=0
\end{align*}
and so if 
\begin{equation*}
- \frac{1}{3} \beta ^2 -1 = 0
\end{equation*}
(which, after a short calculation, implies that $(\beta +1 +\frac{1}{27}\beta ^3 ) \neq 0$) \newline
we then must have $t = 0$ with multiplicity 3, which implies
\begin{align}
& t = 0 \nonumber  \\ 
\Rightarrow & \frac{x}{y} = 0 \nonumber  \\
\Rightarrow & x= 0 \nonumber \\
\Rightarrow & z= y + \beta (0) \label{zeqylab} \\
\Rightarrow & z = y \nonumber 
\end{align}
which is the point $[0:1:1]$ (the singular point), and where we have used $\eqref{zbetalabel}$ in line $\eqref{zeqylab}$. But
\begin{align*}
- \frac{1}{3} & \beta ^2 -1  = 0 \\
\Rightarrow & \beta ^2 = -3 \\
\Rightarrow & \beta  = \pm \sqrt{3}i \\
\Rightarrow & \frac{u_0}{u_1}  = \pm \sqrt{3}i
\end{align*}
where the last line follows from $\eqref{betau0lab}$. This means that the two lines corresponding to $[\sqrt{3}i:1]$ and $[-\sqrt{3}i:1]$ each intersect the curve $C$ in the point $[0:1:1]$ (and in no other points), which shows that the polynomial in $\eqref{singpolylab}$ corresponds to a fishtail fibre.

\subsection*{The other singular fibres}
So far, we have shown that the two polynomials
\begin{align*}
p_0(x,y,z) &= \Big(y+\frac{1}{2}z \Big)^3 \\
p_1(x,y,z) &= (zy^2-zx^2-x^3)
\end{align*}
give rise to an elliptic fibration with an $\tilde{E_6}$ fibre and two fishtail fibres. By Kodaira's classification of singular fibres in section \ref{kodairasec}, the Euler characteristic of an $\tilde{E_6}$ fibre is 8, and the Euler characteristic of a fishtail fibre is 1, so the sum of the Euler characteristics of these three fibres is 10. \newline

There are two more singular curves given by lines $\eqref{singcurve2}$ and $\eqref{singcurve3}$,
\begin{equation*}
C_1 = \{[x:y:z] \in \cpt | -\frac{32}{121}\Big(y + \frac{1}{2}z \Big)^3 + zy^2 - zx^2 - x^3 =0  \}
\end{equation*}
which has a singular point at $[-\frac{2}{3}: \frac{4}{3}: 1]$, and 
\begin{equation*}
C_2 = \{[x:y:z] \in \cpt | \Big(y + \frac{1}{2}z \Big)^3 + \frac{1}{8}\Big(zy^2 - zx^2 - x^3\Big) =0  \}
\end{equation*}
which has a singular point at $[-\frac{2}{3}: -\frac{1}{3}: 1]$. \newline

Both of these singular curves must have Euler characteristic at least 1. However, the Euler characteristic of $E(1) = \cpt \# 9 \cptbar$ is 12, and so if one of these singular fibres has Euler characteristic greater than 1, then the sum of the Euler characteristics of the five fibres is greater than 12, which is a contradiction. Therefore, both $C_1$ and $C_2$ must be fishtail fibres, and thus we have found an elliptic fibration on $E(1)$ with an $\tilde{E_6}$ fibre and four fishtail fibres.

\pagebreak


\section{An Exotic $\mathbb{CP}^2 \# 7 \overline{\mathbb{CP}^2}$} \label{mainsectionlabel}

In this section, we give an exposition of J. Park's construction in \cite{P1} of an exotic $\cpt \# 7 \cptbar$. \newline

Let $C_{7}$ be the plumbing given by Figure \ref{x7proofsection}.1, without the meridians $a_i$ $(i=0,1, \dots, 5)$.  \newline

It is proved in \cite{FS1} that if the intersection form of $H_2(C_7; \mz)$ with respect to the basis $\{u_{1}, u_{2}, \hdots , u_{6}\}$ is given by 

\begin{displaymath}
P = 
\left( \begin{array} {cccccc}
-2 &  1 &  0 &  0 &  0 &  0 \\
 1 & -2 &  1 &  0 &  0 &  0 \\
 0 &  1 & -2 &  1 &  0 &  0 \\
 0 &  0 &  1 & -2 &  1 &  0 \\
 0 &  0 &  0 &  1 & -2 &  1 \\
 0 &  0 &  0 &  0 &  1 & -9
\end{array} \right)
\end{displaymath}

then if $\{\gamma_{1}, \gamma_{2}, \hdots, \gamma_{6} \}$ is the basis of $H_{2}(C_{7}, \partial C_{7} ; \mathbb{Z}) \cong H^{2}(C_{7} ; \mathbb{Z})$ which is dual to the basis $\{u_{1}, u_{2}, \hdots , u_{6}\}$ of $ H_{2}(C_{7} ; \mathbb{Z})$, i.e. $<\gamma_{i}, u_{j}> = \delta_{ij}$, then the intersection form of $H^{2}(C_{7} ; \mathbb{Q})$ with respect to this basis is 

\begin{displaymath}
T = P^{-1} = \frac{-1}{49} 
\left( \begin{array} {cccccc}
41 & 33 & 25 & 17 &  9 &  1 \\
33 & 66 & 50 & 34 & 18 &  2 \\
25 & 50 & 75 & 51 & 27 &  3 \\
17 & 34 & 51 & 68 & 36 &  4 \\
 9 & 18 & 27 & 36 & 45 &  5 \\
 1 &  2 &  3 &  4 &  5 &  6
\end{array} \right)
\end{displaymath}

\newtheorem{remp1}{Remark}[section]
\begin{remp1} \label{sec21rem1}
\upshape
It will turn out that we only need to know the values in the last row/column.
\end{remp1}

For the proof of the main theorem, we shall need the fact that the manifold on which we perform the rational blowdown has a symplectic structure. In \cite{Sy}, M. Symington proved the following theorem (stated in \cite{P1}, as below) which gives conditions when a manifold that has been rationally blown-down admits a symplectic structure. \pagebreak

\newtheorem{sy1}[remp1]{Theorem}
\begin{sy1}
Suppose that $(X, \omega)$ is a symplectic 4-manifold that contains a configuration $C_{p}$. Suppose further that all the 2-spheres in $C_{p}$ are symplectically embedded in $X$ and intersect each other orthogonally. Then the manifold which is the rational blow-down of $X$ along $C_{p}$, denoted by $X_{p} = X_{0} \cup_{L(p^{2}, 1-p)} B_{p}$, admits a symplectic 2-form $\omega_{p}$ such that $(X_{0}, \omega_{p}|_{X_{0}})$ is symplectomorphic to $(X_{0}, \omega|_{X_{0}})$.
\label{sy1ref}
\end{sy1}

As noted in \cite{Sy}, the fact that the 2-spheres are symplectic means that the orthogonal intersections are positive. \newline

Now suppose we have a 4-manifold $X$ that satisifes the conditions for Theorem \ref{sy1ref} above; $(X, \omega)$ is a symplectic 4-manifold, $X = X_{0} \cup_{L(p^{2}, 1-p)} C_{p}$, and $C_{p}$ is symplectically embedded in $X$. So, $X_{p} = X_{0} \cup_{L(p^{2}, 1-p)} B_{p}$ is a symplectic 4-manifold with symplectic 2-form $\omega_{p}$, such that there is a symplectomorphism $\psi_{p} : (X_{0}, \omega_{p}|_{X_{0}}) \longrightarrow (X_{0}, \omega|_{X_{0}})$. Let $K$ be the canonical class on $X$ induced by $\omega$, and let $K_{p}$ be the canonical class on $X_{p}$ induced by $\omega_{p}$ (see Section \ref{canclasssection}). \newline

We now claim that $H^{1}(L(p^{2}, 1-p); \mathbb{Q})$ and $H^{2}(L(p^{2}, 1-p); \mathbb{Q})$ are both trivial. 
\newtheorem{claim1}[remp1]{Claim}
\begin{claim1} \label{claim1ref}
$H^{1}(L(p^{2}, 1-p); \mathbb{Q})$ is trivial
\end{claim1}
Proof: \newline
For ease of reading the calculations, we shall denote the lens space $L(p^2, 1-p)$ simply by $L$. Firstly, by the Universal Coefficient Theorem,
\begin{equation*} 
H^{1}(L; \mathbb{Q}) \cong \mathrm{Hom}(H_{1}(L; \mathbb{Z}); \mathbb{Q}) \oplus \mathrm{Ext}(H_{0}(L; \mathbb{Z}); \mathbb{Q})
\end{equation*}
We know $\pi_1(L) \cong \mathbb{Z}_{p^2}$, and since $H_1(L; \mathbb{Z})$ is the abelianization of $\pi_1(L)$ (which is already abelian, in this case), we have $H_1(L; \mathbb{Z}) \cong \mathbb{Z}_{p^2}$. \newline 

Recall that if $F$ is a finite abelian group and $G$ is an abelian group of infinite order, then $\mathrm{Hom}(F,G)$ must be trivial, since if $\phi: F \longrightarrow G$ is a homomorphism, then for all $x \in F$, $n\phi(x) = 0_{G}$, where $n$ is the order of $F$, and the only way this can happen is if $\phi(x) = 0_{G}$ for each $x \in F$. \newline

In particular $\mathrm{Hom}(\mathbb{Z}_{p^2}; \mathbb{Q})$ is trivial, and so $\mathrm{Hom}(H_{1}(L; \mathbb{Z}); \mathbb{Q})$ is trivial. \newline 

Recall from \cite{H1} that if $F$ is a free group, then $\mathrm{Ext}(F; \mathbb{Q})$ is trivial. Since $H_{0}(L; \mathbb{Z}) \cong \mathbb{Z}$ (see \cite{H1}), $\mathrm{Ext}(H_{0}(L; \mathbb{Z}); \mathbb{Q})$ is trivial, and so finally we have that $H^{1}(L; \mathbb{Q})$ is trivial. This proves the claim. $\Box$ \newline

\newtheorem{claim2}[remp1]{Claim}
\begin{claim2} \label{claim2ref}
$H^{2}(L(p^{2}, 1-p); \mathbb{Q})$ is trivial
\end{claim2}
Proof: \newline
By Poincar\'e Duality, since $L(p^{2}, 1-p)$ is three-dimensional
\begin{equation*}
H^{2}(L(p^{2}, 1-p); \mathbb{Q}) \cong H_{1}(L(p^{2}, 1-p); \mathbb{Q})
\end{equation*}
$H_{1}(L(p^{2}, 1-p); \mathbb{Z}) \cong \mathbb{Z}_{p^2}$, and changing from integral to rational coefficients shows us that $H_{1}(L(p^{2}, 1-p); \mathbb{Q})$ is trivial. $\Box$ \newline

We now claim that the triviality of these cohomology groups allows us to decompose $K$ and $\omega$ as $K = K|_{X_{0}} + K|_{C_{p}}$ and $[\omega] = [\omega|_{X_{0}}] + [\omega|_{C_{p}}]$ (where $K|_{X_{0}}, [\omega|_{X_{0}}] \in H^{2}(X_{0}, \mathbb{Q})$, $K|_{C_{p}}, [\omega|_{C_{p}}] \in H^{2}(C_{p}; \mathbb{Q})$ and $K, [\omega] \in H^{2}(X; \mathbb{Q})$). \newline 

If we choose $A$ to be a small open neighbourhood\footnote{I say ``small open neighbourhood'' to mean $A$ deformation retracts onto $X_{0}$, so the homology groups of $A$ and $X_{0}$ are isomorphic. The same goes for $B$ and $A \cap B$.} of $X_0$ and $B$ to be a small open neighbourhood of $C_{p}$, so that $A \cap B$ is a small open neighbourhood of $L(p^2, 1-p)$ and $X = \mathrm{int}(A) \cup \mathrm{int}(B)$, then the Mayer-Vietoris sequence (see \cite{H1})
\begin{equation*}
\dots \longrightarrow H_{2}(A \cap B) \longrightarrow H_{2}(A) \oplus H_{2}(B) \longrightarrow H_{2}(X) \longrightarrow H_{1}(A \cap B) \longrightarrow \dots
\end{equation*} 
implies, by Poincar\'e duality and the claims above, that
\begin{align*}
 H_{2}(L(p^{2}, 1-p)) &\longrightarrow H_{2}(X_{0}) \oplus H_{2}(C_{p}) \longrightarrow H_{2}(X) \longrightarrow H_{1}(L(p^{2}, 1-p)) \\
\Rightarrow  0 &\longrightarrow H_{2}(X_{0}) \oplus H_{2}(C_{p}) \longrightarrow H_{2}(X) \longrightarrow 0
\end{align*}
which implies $H_{2}(X_{0}) \oplus H_{2}(C_{p}) \cong H_{2}(X)$ (and the omitted coefficients are the rationals). This proves the claim. $\Box$\newline

So, we have the following decompositions:
\begin{equation*}
K = K|_{X_{0}} + K|_{C_{p}}
\end{equation*}
where $K \in H^{2}(X; \mathbb{Q})$, $K|_{X_{0}} \in H^{2}(X_{0}, \mathbb{Q})$ and $K|_{C_{p}} \in H^{2}(C_{p}; \mathbb{Q})$, and 
\begin{equation*}
[\omega] = [\omega|_{X_{0}}] + [\omega|_{C_{p}}]
\end{equation*}
where $[\omega] \in H^{2}(X; \mathbb{Q})$, $[\omega|_{X_{0}}] \in H^{2}(X_{0}, \mathbb{Q})$ and $[\omega|_{C_{p}}] \in H^{2}(C_{p}; \mathbb{Q})$. \newline


Similarly, $K_{p}$ and $[\omega_{p}]$ decompose as:
\begin{align*}
&K_{p} = K_{p}|_{X_{0}} + K_{p}|_{B_{p}} \\
&[\omega_{p}] = [\omega_{p}|_{X_{0}}] + [\omega_{p}|_{B_{p}}] 
\end{align*}
where $K_{p}, [\omega_{p}] \in H^{2}(X_{p}; \mathbb{Q})$, $K_{p}|_{X_{0}}, [\omega_{p}|_{X_{0}}] \in H^{2}(X_{0}, \mathbb{Q})$ and $K_{p}|_{B_{p}}, [\omega_{p}|_{B_{p}}] \in H^{2}(B_{p}; \mathbb{Q})$. \newline

(Basically, all the cohomology classes are where we expect them to be.) \newline

These decompositions lead to the following lemma (Lemma 2.1 in \cite{P1}).

\newtheorem{L21}[remp1]{Lemma}
\begin{L21} \label{kpomegaplemma}
Under the same hypothesis on $(X, K, \omega)$ and $(X_{p}, K_{p}, \omega_{p})$ as above, we have
\begin{equation*}
K_{p} \cdot [\omega_{p}] = K \cdot [\omega] - K|_{C_{p}} \cdot [\omega|_{C_{p}}] 
\end{equation*}
\end{L21}

Proof:\newline
We have the short exact sequence
\begin{equation*}
H_{2}(\partial B_{p}; \mathbb{Q}) \longrightarrow H_{2}(B_{p}; \mathbb{Q}) \longrightarrow H_{2}(B_{p}, \partial B_{p}; \mathbb{Q}) \longrightarrow H_{1} (\partial B_{p}; \mathbb{Q})
\end{equation*}

Since $\partial B_{p} \cong L(p^{2}, 1-p)$, from Claims \ref{claim1ref} and \ref{claim2ref} above, both $H_{2}(\partial B_{p}; \mathbb{Q})$ and $H_{1}(\partial B_{p}; \mathbb{Q})$ are trivial, and so the short exact sequence becomes
\begin{equation*}
0 \longrightarrow H_{2}(B_{p}; \mathbb{Q}) \longrightarrow H_{2}(B_{p}, \partial B_{p}; \mathbb{Q}) \longrightarrow 0
\end{equation*}
proving $H_{2}(B_{p}; \mathbb{Q}) \cong H_{2}(B_{p}, \partial B_{p}; \mathbb{Q})$. \newline

By Poincar\'e duality, $H_{2}(B_{p}, \partial B_{p}; \mathbb{Q}) \cong H^{2}(B_{p}; \mathbb{Q})$, so $H^{2}(B_{p}; \mathbb{Q}) \cong H_{2}(B_{p}; \mathbb{Q})$. Since $B_{p}$ is a rational ball, $H_{2}(B_{p}; \mathbb{Q})$ is trivial, and so $H^{2}(B_{p}; \mathbb{Q})$ is also trivial. \newline

Therefore, $K_{p}|_{B_{p}}$ and $[\omega|_{B_{p}}]$ are zero elements in $H^{2}(B_{p}; \mathbb{Q})$, and consequently zero elements in $H^{2}(X_{p}; \mathbb{Q})$. \newline

Now let us explicitly calculate $K_{p} \cdot [\omega_{p}]$:
\begin{align*}
K_{p} \cdot [\omega_{p}] &= (K_{p}|_{X_{0}} + K_{p}|_{B_{p}}) \cdot ([\omega_{p}|_{X_{0}}] + [\omega_{p}|_{B_{p}}] ) \\
&= K_{p}|_{X_{0}} \cdot [\omega_{p}|_{X_{0}}] +   K_{p}|_{X_{0}} \cdot [\omega_{p}|_{B_{p}}] + K_{p}|_{B_{p}} \cdot [\omega_{p}|_{X_{0}}] +  K_{p}|_{B_{p}} \cdot [\omega_{p}|_{B_{p}}] \\
&= K_{p}|_{X_{0}} \cdot [\omega_{p}|_{X_{0}}]
\end{align*}
since the last three terms of the second line are all zero, since each contains $K_{p}|_{B_{p}}$ or $[\omega_{p}|_{B_{p}}]$. \newline

By the definition of the symplectomorphism $\psi_{p} : (X_{0}, \omega_{p}|_{X_{0}}) \longrightarrow (X_{0}, \omega|_{X_{0}})$, we must have $\psi_{p}^{*} ([\omega|_{X_{0}}]) = [\omega_{p}|_{X_{0}}]$ and $\psi_{p}^{*} (K|_{X_{0}}) = K_{p}|_{X_{0}}$, so we have

\begin{equation} \label{psieq1}
K_{p} \cdot [\omega_{p}] = K_{p}|_{X_{0}} \cdot [\omega_{p}|_{X_{0}}] = \psi_{p}^{*}(K|_{X_{0}}) \cdot \psi_{p}^{*}([\omega|_{X_{0}}])
\end{equation}


Furthermore, since $\psi_{p}^{*}$ is a homomorphism,

\begin{equation} \label{psieq2}
\psi_{p}^{*}(K|_{X_{0}}) \cdot \psi_{p}^{*}([\omega|_{X_{0}}]) = \psi_{p}^{*}(K|_{X_{0}} \cdot [\omega|_{X_{0}}])
\end{equation}

Since $\psi_{p}$ is a symplectomorphism, in particular, it is an orientation-preserving diffeomorphism, and so $\psi_{p}^{*}$ is an isomorphism between the homology and cohomology classes of $(X_0, \omega_{p}|_{X_0})$ and $(X_0, \omega|_{X_0})$. \newline

Since $H^4(X_0; \mz) \cong \mz$, $\psi_{p}^{*}$ maps $1 \in H^4((X_0, \omega_{p}|_{X_0}); \mz)$ to $1 \in H^4((X_0, \omega|_{X_0}); \mz)$ (because it is orientation-preserving, $1$ does not get mapped to $-1$). Therefore,

\begin{equation} \label{psieq3}
\psi_{p}^{*}(K|_{X_{0}} \cdot [\omega|_{X_{0}}]) = K|_{X_{0}} \cdot [\omega|_{X_{0}}]
\end{equation}

and using equations $\eqref{psieq2}$ and $\eqref{psieq3}$, we have

\begin{equation} \label{psieq4}
\psi_{p}^{*}(K|_{X_{0}}) \cdot \psi_{p}^{*}([\omega|_{X_{0}}]) = K|_{X_{0}} \cdot [\omega|_{X_{0}}]
\end{equation}


Finally, using equations $\eqref{psieq1}$ and $\eqref{psieq4}$, we have
\begin{equation} \label{kdo1ref}
K_{p} \cdot [\omega_{p}] = K|_{X_{0}} \cdot [\omega|_{X_{0}}]
\end{equation}

We now have $K_{p} \cdot [\omega_{p}] = K|_{X_{0}} \cdot [\omega|_{X_{0}}]$. Calculating
\begin{align*}
K \cdot [\omega] &= (K|_{X_{0}} + K|_{C_{p}}) \cdot ([\omega|_{X_{0}}] + [\omega_{p}|_{C_{p}}]) \\
&= K|_{X_{0}} \cdot [\omega|_{X_{0}}] + K|_{X_{0}} \cdot [\omega_{p}|_{C_{p}}] + K|_{C_{p}} \cdot [\omega|_{X_{0}}] + K|_{C_{p}} \cdot [\omega_{p}|_{C_{p}}]
\end{align*}

Now, $K|_{X_{0}}$ and $[\omega_{p}|_{C_{p}}]$ are forms restricted to different submanifolds, and so $K|_{X_{0}} \cdot [\omega_{p}|_{C_{p}}] = 0$. Algebraically, this is because each term belongs to a different summand of $H_{2}(X_{0}; \mathbb{Q}) \oplus H_{2}(C_{p}; \mathbb{Q}) \cong H_{2}(X; \mathbb{Q})$. Similarly, $K|_{C_{p}} \cdot [\omega|_{X_{0}}] = 0$ and so we have
\begin{equation} \label{kdo2ref}
K \cdot [\omega] = K|_{X_{0}} \cdot [\omega|_{X_{0}}] + K_{C_{p}} \cdot [\omega_{p}|_{C_{p}}]
\end{equation}

Putting $\eqref{kdo1ref}$ and $\eqref{kdo2ref}$ together, we have
\begin{equation*}
K_{p} \cdot [\omega_{p}] = K|_{X_{0}} \cdot [\omega|_{X_{0}}] = K \cdot [\omega] - K_{C_{p}} \cdot [\omega_{p}|_{C_{p}}]
\end{equation*}

which proves the lemma. $\Box$ \newline

Next, we shall need the fact that $E(1) = \mathbb{CP}^{2} \# 9 \overline{\mathbb{CP}^{2}}$ can be described to be an elliptic fibration over $\mathbb{CP}^{1}$ with one $\tilde{E_{6}}$ singular fibre and four fishtail fibres. The existence of such a fibration was shown in section \ref{explicitE6sec}. \newline

The following is Lemma 3.1 in \cite{P1}, credited to D. Auroux and R. Fintushel.

\newtheorem{L31}[remp1]{Lemma}
\begin{L31} \label{aurouxfintlemma}
The second (co)homology classes $[S_{i}]$ ($1 \leq i \leq 7$) of the 2-spheres $S_{i}$ embedded in $\tilde{E_{6}}$ can be represented by: \newline
$[S_{1}] = e_{4} - e_{7}$, $[S_{2}] = e_{1} - e_{4}$, $[S_{3}] = h - e_{1} - e_{2} - e_{3}$, \newline
$[S_{4}] = e_{2} - e_{5}$, $[S_{5}] = e_{5} - e_{9}$, $[S_{6}] = e_{3} - e_{6}$, $[S_{7}] = e_{6} - e_{8}$; \newline
where $h$ denotes a generator of $H_{2}(\mathbb{CP}^{2}; \mathbb{Z})$ and each $e_{i}$ denotes the (co)homology class represented by the i-th exceptional curve in $\overline{\mathbb{CP}^{2}} \subset E(1) = \mathbb{CP}^{2} \# 9 \overline{\mathbb{CP}^{2}}$.
\end{L31}

Section \ref{e6constrsec} gives a proof of this Lemma, and we can take the labelling of the spheres $S_1, \dots, S_7$ to be as given above. \newline

Below are two theorems that we shall need in the proof of the lemma below. The first is Corollary 1.4 in \cite{MS}, credited to A. Liu and H. Ohta and K. Ono. Recall that a 4-manifold $X$ is said to be \emph{minimal} if it contains no exceptional spheres, and that a 4-manifold $X$ is said to be \emph{rational} if it is the blow-up of $\mathbb{CP}^{2}$ or $S^{2} \times S^{2}$. For example, $E(1)$ is rational but not minimal. We now quote a series of important results. \pagebreak

\newtheorem{cor14ms}[remp1]{Lemma}
\begin{cor14ms} \label{cor14mslabel}
Let $X$ be a symplectic 4-manifold. Then the following are equivalent:
\begin{itemize}
\item[(i)] $X$ admits a metric of positive scalar curvature.
\item[(ii)] $X$ admits a symplectic structure $\omega$ with $K \cdot \omega < 0$
\item[(iii)] $X$ is either rational or ruled
\end{itemize}
\end{cor14ms}

First, note that in the statement of the lemma in \cite{MS}, $X$ is a \emph{minimal} symplectic 4-manifold. However, it is then noted in \cite{MS} that the lemma extends to the case when $X$ is not minimal (which is how we have stated it above). Before we can quote another lemma, we need:

\newtheorem{D33LL}[remp1]{Definition}
\begin{D33LL}
\upshape
For a non-minimal rational manifold with a standard decomposition $\mathbb{CP}^{2} \# n \overline{\mathbb{CP}^{2}}$ and a standard basis $\{h, e_{1}, e_{2},\dots, e_{n} \}$, a class $\xi = ah - b_{1}e_{1}-b_{2}e_{2}- \dots - b_{n}e_{n}$ is called \emph{reduced} if 
\begin{align*}
& b_{1} \geq b_{2} \geq \dots \geq b_{n} \geq 0, and\\
& a \geq b_{1} + b_{2} + b_{3}
\end{align*}
\end{D33LL}

Note that the second condition (with the first condition) implies $a \geq b_{i}$ for $1 \leq i \leq n$. With this definition in mind, the following lemma is the first part of Lemma 3.4 (which is actually stronger and has other implications) in \cite{LL1}.

\newtheorem{L34LL}[remp1]{Lemma}
\begin{L34LL} \label{reducedsquarelemma}
Let $M$ be a non-minimal rational manifold with a standard decomposition and a standard basis. Then any class of non-negative square is equivalent to a reduced class under the action of orientation-preserving diffeomorphisms. (Furthermore, we can find such a diffeomorphism by a simple algorithm.)
\end{L34LL}

The canonical class of $E(1)$,  $K_{E(1)} \in H^{2}(E(1); \mathbb{Z})$, is represented by $K_{E(1)} = -3h + e_{1}+ \dots + e_{9} = -[f]$, following the notation in \cite{P1} (although in \cite{P1} the notation seems to change from `$[f]$' to just `$f$'). Using the lemma above, we later get an important relation between this canonical class and a compatible symplectic 2-form on a non-minimal rational surface, which we shall need in the proof of our main result. \newline

In $\cite{LiLiu2}$ it is proved that $\mathbb{CP}^{2} \# k \overline{\mathbb{CP}^{2}}$ has a `unique' symplectic structure for certain $k \geq 2$:

\newtheorem{L33P}[remp1]{Theorem}
\begin{L33P} \label{usscctref}
There is a unique symplectic structure on $\mathbb{CP}^{2} \# k \overline{\mathbb{CP}^{2}}$ for $ 2 \leq k \leq 9$ up to diffeomorphisms and deformation. For $k \geq 10$, the symplectic structure is still unique for the standard canonical class.
\end{L33P}

\newtheorem{L33PCor}[remp1]{Corollary}
\begin{L33PCor} \label{komegacorr}
$\cpt \# k \cptbar$, for $2 \leq k \leq 9$, with canonical class $K$ does not admit a symplectic 2-form $\omega$ for which $K \cdot \omega > 0$.
\end{L33PCor}

\newtheorem{L33PCorRem}[remp1]{Remark}
\begin{L33PCorRem}
\upshape
This result is of great importance to us, as it will be used to show that the 4-manifold that we contruct, although homeomorphic to $\cpt \# 7 \cptbar$, is not diffeomorphic to it.
\end{L33PCorRem}

We now use these results to prove the following lemma.

\newtheorem{L32P}[remp1]{Lemma}
\begin{L32P} \label{redclasslemma}
For each integer $k \geq 1$, $E(1)\#k \overline{\mathbb{CP}^{2}}$ admits a symplectic 2-form $\omega$ which is compatible with the standard canonical class $K_{E(1)\#k \overline{\mathbb{CP}^{2}}} = -3h + e_{1} + e_{2} +\dots+e_{9+k}$, such that its cohomology class $[\omega]$ can be represented by $ah- b_{1}e_{1}-b_{2}e_{2}-\dots - b_{9+k}e_{9+k}$, where $a, b_{1}, b_{2}, \dots, b_{9+k}$ are some rational numbers satisfying 
\begin{itemize}
\item [(i)] $a \geq b_{1} \geq \dots \geq b_{9+k} \geq 0$, and 
\item [(ii)] $3a > b_{1}+b_{2}+\dots +b_{9+k}$.
\end{itemize}
\end{L32P}

Proof: \newline
By the equivalence between $(iii)$ and $(ii)$ in Lemma $\ref{cor14mslabel}$ above, since $E(1)\#k \overline{\mathbb{CP}^{2}}$ is a rational surface it admits symplectic 2-form $\omega$, which is compatible with the standard canonical class $K_{E(1)\#k \overline{\mathbb{CP}^{2}}} = -3h + e_{1} + e_{2} +\dots+e_{9+k}$ and satisfies the inequality $K_{E(1)\#k \overline{\mathbb{CP}^{2}}} \cdot \omega < 0$. \newline

In \cite{GS} we are reminded that a symplectic form is non-degenerate and so $\omega \wedge \omega > 0$. This then implies that $[\omega]$ has non-negative square. Since $[\omega]$ has non-negative square, Lemma \ref{reducedsquarelemma} above implies that $[\omega]$ can be represented by $ah- b_{1}e_{1}-b_{2}e_{2}-\dots - b_{9+k}e_{9+k}$ for some rational numbers satisfying $a \geq b_{1} \geq \dots \geq b_{9+k} \geq 0$ (part (i) of the lemma). \newline

Since $h \cdot h = 1$, $h \cdot e_{i} = 0$ for all $e_{i}$, and $e_{i} \cdot e_{j} = -\delta_{ij}$, 
\begin{align*}
K_{E(1)\#k \overline{\mathbb{CP}^{2}}} \cdot \omega &= (-3h + e_{1} +\dots+e_{9+k}) \cdot (ah - b_{1}e_{1}- \dots - b_{9+k}e_{9+k}) \\
&= -3a h \cdot h - b_{1}e_{1}\cdot e_{1} - \dots - b_{9+k} e_{9+k} \cdot e_{9+k} \\
&= -3a +b_{1} + b_{2} + \dots + b_{9+k}
\end{align*}

and this together with $K_{E(1)\#k \overline{\mathbb{CP}^{2}}} \cdot \omega < 0$ then implies

\begin{equation*}
3a > b_{1} + b_{2} + \dots + b_{9+k}
\end{equation*}

which is part (ii) of the lemma. $\Box$ \newline

Next comes an important proposition, concerning the existence of a specific configuation $C_{p}$ in $E(1) \# k \overline{\mathbb{CP}^{2}}$. We also need that the 2-spheres in the configuration are symplectically embedded, in order to use the theorem proved in \cite{Sy} and stated above.

\newtheorem{P31}[remp1]{Proposition}
\begin{P31} \label{c7sympprop}
There exists a configuration $C_{7}$ in the rational surface $E(1) \# 4 \overline{\mathbb{CP}^{2}}$ such that all the 2-spheres $u_{i}$ lying in $C_{7}$ are symplectically embedded.
\end{P31}

Proof: \newline
Recall that $E(1)$ can be viewed as an elliptic fibration with an $\tilde{E_{6}}$-singular fibre and 4 singular fishtail fibres (section \ref{explicitE6sec}). \newline

Recall that the homology class $[f] = 3h - e_{1} - \dots - e_{9}$ of the elliptic fibre $f$ in $E(1)$ can be represented by an immersed 2-sphere with one positive double point, which is equivalent to a fishtail fibre (sections \ref{e6constrsec} and  \ref{explicitE6sec}). \newline

$E(1)$ contains at least 4 such immersed 2-spheres, since it contains 4 singular fishtail fibres. If we blow up at each of these 4 double points, there exist embedded 2-spheres $f-2e_{10}$, $f-2e_{11}$, $f-2e_{12}$, $f-2e_{13}$ in $E(1) \# 4 \overline{\mathbb{CP}^{2}}$. \newline

The reason it is $f-2e_{k}$ and not simply $f-e_{k}$ (for $10 \leq k \leq 13$) is that $f$ has a positive double point, and when we blow-up at this double point, the exceptional divisor will intersect $f$ with multiplicity $2$, and so the proper transform of $f$ is $f-2e_{k}$. \newline

Since $[f] \cdot e_{9} = (3h - e_{1} - \dots - e_{9})\cdot e_9 = -e_{9} \cdot e_{9} = -(-1) = 1$, and $e_{k} \cdot e_{9} = 0$ for $10 \leq k \leq 13$, we have $(f-2e_{k}) \cdot e_{9} = 1$ for $10 \leq k \leq 13$. So, each of the embedded 2-spheres $f-2e_{10}$, $f-2e_{11}$, $f-2e_{12}$, $f-2e_{13}$ intersects a section $e_{9}$ of $E(1)$ positively. Let us call these intersection points $p_{10}, p_{11}, p_{12}, p_{13}$, respectively. \newline

We resolve these 4 intersection points (see section \ref{resolvesection}) to get a sphere $S$. The homology class representing $S$ is 
\begin{equation*}
[S] = (f-2e_{10}) + (f-2e_{11}) + (f-2e_{12}) + (f-2e_{13}) + e_{9} = 4f + e_{9} -2e_{10} - 2e_{11}-2e_{12}-2e_{13}
\end{equation*}
since when we resolve the points, we add the homology classes of the surfaces together (see section \ref{resolvesection}). Recall that this ``resolving'' can be done symplectically (see Remark \ref{resolvesympl}). \newline

And recalling that $f \cdot f = 0$, $f \cdot e_{9} = 1$, $f \cdot e_{k} =0 $ for $10 \leq k \leq 13$, and $e_{i} \cdot e_{j} = -\delta_{ij}$, we can compute the square of $[S]$ easily:
\begin{align*}
[S] \cdot [S] &= (4f + e_{9} -2(e_{10} - \dots - e_{13})) \cdot (4f + e_{9} -2(e_{10} - \dots - e_{13})) \\
&= 16f \cdot f + 8f \cdot e_{9} + e_{9}\cdot e_{9} + 4(e_{10} \cdot e_{10} + e_{10} + \cdots + e_{13} \cdot e_{13}) \\
&= 16 \cdot 0 + 8 \cdot 1 + (-1) + 4\cdot ((-1) + (-1) + (-1) + (-1)) \\
&= -9
\end{align*}

Therefore, we have found a symplectically embedded 2-sphere $S$ with square $-9$ in $E(1) \# 4 \overline{\mathbb{CP}^2}$. \newline

We now recall constructing an elliptic fibration on $E(1)$ with an $\tilde{E_{6}}$ fibre. Lemma \ref{aurouxfintlemma} showed that we can consider the $\tilde{E_6}$ fibre as consisting of the spheres $S_1,\dots, S_7$ where
$[S_{1}] = e_{4} - e_{7}$, $[S_{2}] = e_{1} - e_{4}$, $[S_{3}] = h - e_{1} - e_{2} - e_{3}$, $[S_{4}] = e_{2} - e_{5}$ and $[S_{5}] = e_{5} - e_{9}$. Note that

\begin{equation*}
[S_{1}]\cdot[S_{2}] = [S_{2}]\cdot[S_{3}] = [S_{3}]\cdot[S_{4}] = [S_{4}]\cdot[S_{5}] = 1
\end{equation*}
and 
\begin{align*}
[S]\cdot [S_{5}] &= (4f + e_{9} -2(e_{10} - \dots - e_{13})) \cdot (e_{5}-e_{9}) \\
&= 4f \cdot e_{5} - 4f\cdot e_{9} - e_{9} \cdot e_{9} \\
&= 4(-e_{5}) \cdot e_{5} - 4(-e_{9}) \cdot e_{9} - e_{9} \cdot e_{9}\\
&= 4 - 4 + 1 \\
&= 1
\end{align*}

So, if we set $u_{1}= S_{1}$, $u_{2}= S_{2}$, \dots,$u_{5}= S_{5}$ and $u_{6}= S$, we obtain a configuration $C_{7}$ lying in $E(1) \# 4 \overline{\mathbb{CP}^{2}}$. See Figures \ref{e6constrsec}.6 and \ref{x7proofsection}.1. \newline

Note that all the 2-spheres $u_{i}$ lying in $C_{7}$ are symplectically embedded. (since $S$ is symplectic and $S_1, \dots, S_5$ were constructed using algebraic techniques, and so are just algebraic curves, and are therefore symplectic). This proves the proposition. $\Box$ \newline

We are finally ready to prove the main result in \cite{P1}.

\newtheorem{pmain}[remp1]{Theorem}
\begin{pmain}
The exists a simply connected 4-manifold with which is homeomorphic, but not diffeomorphic, to $\mathbb{CP}^{2} \# 7 \overline{\mathbb{CP}^{2}}$. 
\end{pmain}

Proof: \newline
We shall construct this exotic manifold, which we denote by $X_{7}$. \newline

Let us denote the manifold $E(1) \# 4 \overline{\mathbb{CP}^{2}}$ by $X$. By Proposition \ref{c7sympprop}, there exists a symplectically embedded configuration $C_{7}$ in $X$. If we blow down along this configuration $C_{7}$ in $X = X_{0} \cup _{L(49,-6)} C_{7}$, we get a new smooth 4-manifold, which we shall denote by $X_{7} = X_{0} \cup _{L(49,-6)}B_{7}$. \newline

By the theorem proved in \cite{Sy} and stated above, since $C_{7}$ was symplectically embedded in $X$, there exists a symplectic structure on $X_{7}$. \newline

We proved in Section \ref{x7proofsection} that $X_7$ is simply-connected. \newline

Claim: $X_{7}$ is homeomorphic to $\mathbb{CP}^{2} \# 7 \overline{\mathbb{CP}^{2}}$ \newline
Let us look at $X = \cpt \# 7 \cptbar$ first. It has Betti numbers $b_2^{+}(X) = 1$ and $b_2^{-}(X) = 7$, so the rank of $X$ is $rk(X) =8$ and the signature is $\sigma(X) = 6$. Since $X$ is a smooth manifold (actually, it is symplectic) and since $\sigma(X)$ is not divisible by 8, we have by Lemma \ref{L210} that the intersection form of $X$ must be odd. \newline

Now, let us look at $X_7 = X_{0} \cup _{L(49,-6)}B_{7}$. Since this was constructed from $E(1) \# 4 \cptbar$ by blowing down a configuration $C_7$ (and so, removing 6 spheres of negative self-intersection), we know that

\begin{align*}
b_2^{+}(X_7) &= b_2^{+}(E(1) \# 4 \cptbar) = 1 \\
b_2^{-}(X_7) &= b_2^{-}(E(1) \# 4 \cptbar) - 6 = 13 - 6 = 7
\end{align*}
So, the rank of $X_7$ is $rk(X_7) = 8$ and the signature is $\sigma(X_7) = 6$. Since $X_7$ is a smooth manifold (in fact, also symplectic) we again have by Lemma \ref{L210} that its intersection form must be odd. \newline

By Corollary \ref{frcor1} (Freedman's Theorem) since $X$ and $X_7$ are both smooth manifolds with the same rank, signature and parity, they must be homeomorphic. \newline

We now need to show that $X_{7}$ is not diffeomorphic to $\mathbb{CP}^{2} \# 7 \overline{\mathbb{CP}^{2}}$. \newline

Let $K_{7}$ be the canonical class on $X_{7}$, and let $\omega_{7}$ be the corresponding symplectic 2-form on $X_{7}$. We claim that $K_{7} \cdot \omega_{7} > 0$ (and we shall prove this claim in the lemma following this theorem). \newline

Therefore, by Corollary \ref{komegacorr}, $X_7$ is not diffeomorphic to $X$. So, we have shown that $X_{7}$ is homeomorphic, but not diffeomorphic, to $\mathbb{CP}^{2} \# 7 \overline{\mathbb{CP}^{2}}$. $\Box$

\newtheorem{LKw}[remp1]{Lemma}
\begin{LKw}
As defined above, $K_{7} \cdot [\omega_{7}] > 0$.
\end{LKw}

Proof: \newline
We shall denoted the homology class of a sphere $S$ by $[S]$. \newline

The canonical class $K_{E(1)}$ of $E(1)$ is represented by $-[f] = -3h + (e_{1} + \dots + e_{9})$, and the canonical class $K$ of $X = E(1) \# 4 \overline{\mathbb{CP}^{2}}$ is represented by $K = -3h + (e_{1} + \dots + e_{13}) = -[f] + (e_{10} + \dots + e_{13})$. \newline

Using Lemma \ref{redclasslemma}, we may assume that the cohomology class $[\omega]$ of the symplectic 2-form $\omega$ on $X$, that is compatible with the canonical class $K$, can be represented by $ah - (b_{1}e_{1} + \dots + b_{13}e_{13})$ for some rational numbers $a, b_{1}, \dots, b_{13}$ satisfying $a \geq b_{1} \geq \dots \geq b_{13} \geq 0$ and $3a \geq b_{1} + \dots b_{13}$. \linebreak Recall from Proposition \ref{c7sympprop} above that 
\begin{align*} 
[u_{6}] &= [S] \\
&= 4f + e_{9} - 2(e_{10} + \dots + e_{13})\\
\Rightarrow [u_6]&= 12h + e_{9} -4(e_{1} + \dots + e_{9}) - 2(e_{10} + \dots + e_{13}) 
\end{align*}
 
Let us now recall that we defined $u_{i} = S_{i}$ for $1 \leq i \leq 5$ in Proposition \ref{c7sympprop} above, and the so the homology classes are: 
\begin{align*}
[u_{1}] &= e_{4} - e_{7}, [u_{2}] = e_{1} - e_{4}, [u_{3}] = h - e_{1} - e_{2} - e_{3}, \\
[u_{4}] &= e_{2} - e_{5}, [u_{5}] = e_{5} - e_{9}
\end{align*}

and recall $K = -3h + (e_{1} + \dots + e_{13})$. \newline

From these definitions we can see that $[u_{i}] \cdot K = 0$ (for $1 \leq i \leq 5$). \newline

In order to calculate $K_{7} \cdot [\omega_{7}]$, we shall first calculate $K \cdot [\omega]$ and $K|_{C_{7}} \cdot [\omega|_{C_{7}}]$, and then use the result we proved in Lemma \ref{kpomegaplemma} above:
\begin{equation} \label{lemma215rel}
K_{7} \cdot [\omega_{7}] = K \cdot [\omega] - K|_{C_{7}} \cdot [\omega|_{C_{7}}] 
\end{equation}

to calculate $K_{7} \cdot [\omega_{7}]$. \newline

Firstly, if we define $[\omega] = ah - (b_{1}e_{1} + \dots + b_{13}e_{13})$,
\begin{align}
K \cdot [\omega] &= (-3h + (e_{1} + \dots + e_{13})) \cdot (ah - (b_{1}e_{1} + \dots + b_{13}e_{13})) \nonumber \\
\Rightarrow K \cdot [\omega] &= -3a + (b_{1} + \dots + b_{13}) \label{kcdotomegalab}
\end{align}
We now express the two cohomology classes $K|_{C_{7}}$ and $[\omega|_{C_{7}}]$ using the dual basis $ \{ \gamma_{i} : 1 \leq i \leq 6 \}$ (such that $<\gamma_{i}, u_{j} > = \delta_{ij}$) for $H^{2}(C_{7}; \mathbb{Q})$. We have,
\begin{align*}
K|_{C_{7}} &= (K \cdot [u_{1}]) [\gamma_{1}] + (K \cdot [u_{2}])[\gamma_{2}] +  \dots + (K \cdot [u_{6}])[\gamma_{6}] \\
&= 7 [\gamma_{6}]
\end{align*} 
since $K \cdot [u_{i}] = 0$ for $1 \leq i \leq 5$, and 
\begin{align*}
K \cdot [u_{6}] &= (-3h + (e_{1} + \dots + e_{13})) \cdot (12h + e_{9} -4(e_{1} + \dots + e_{9}) \\
& - 2(e_{10} + \dots + e_{13})) \\
&=  (-3)(12)h \cdot h + e_{9} \cdot e_{9} -4 (e_{1}\cdot e_{1} + \dots + e_{9} \cdot e_{9}) \\
&- 2(e_{10} \cdot e_{10} + \dots + e_{13} \cdot e_{13}) \\
&= -36 + (-1) - 4(-9) - 2(-4)\\
&= -36 -1 +36 + 8 \\
&= 7
\end{align*}

Similarly, 
\begin{align*}
[\omega|_{C_{7}}] &= ([\omega] \cdot [u_{1}])\gamma_{1} + ([\omega] \cdot [u_{2}])\gamma_{2} + \dots + ([\omega] \cdot [u_{6}])\gamma_{6} 
\end{align*}

and using 
\begin{align*}
[\omega] \cdot [u_{1}] &= (ah - (b_{1}e_{1} + \dots + b_{13}e_{13})) \cdot (e_{4} - e_{7}) = b_{4} - b_{7} \\
[\omega] \cdot [u_{2}] &= (ah - (b_{1}e_{1} + \dots + b_{13}e_{13})) \cdot (e_{1} - e_{4}) = b_{1} - b_{4} \\
[\omega] \cdot [u_{3}] &= (ah - (b_{1}e_{1} + \dots + b_{13}e_{13})) \cdot (h-e_{1}-e_{2}-e_{3}) = a-b_{1}-b_{2}-b_{3} \\
[\omega] \cdot [u_{4}] &= (ah - (b_{1}e_{1} + \dots + b_{13}e_{13})) \cdot (e_{2} - e_{5}) = b_{2} - b_{5} \\
[\omega] \cdot [u_{5}] &= (ah - (b_{1}e_{1} + \dots + b_{13}e_{13})) \cdot (e_{5} - e_{9}) = b_{5} - b_{9} \\
[\omega] \cdot [u_{6}] &= (ah - (b_{1}e_{1} + \dots + b_{13}e_{13})) \cdot [u_{6}] \\
&= 12a + b_{9} - 4(b_{1} + \dots + b_{9}) - 2(b_{10} + \dots + b_{13})
\end{align*}

we finally get
\begin{align*}
[\omega|_{C_{7}}] &= (b_{4} - b_{7}) [\gamma_{1}] + (b_{1} - b_{4}) [\gamma_{2}] + (a-b_{1}-b_{2}-b_{3}) [\gamma_{3}] + (b_{2} - b_{5}) [\gamma_{4}] \\
& + (b_{5} - b_{9}) [\gamma_{5}] + (12a + b_{9} - 4(b_{1} + \dots + b_{9}) - 2(b_{10} + \dots + b_{13})) [\gamma_{6}]
\end{align*}

Then, using the intersection form of $H^{2}(C_{7}; \mathbb{Q})$ given above Remark \ref{sec21rem1}, we have $[\gamma_{6}] \cdot [\gamma_{k}] = <\gamma_6 , \gamma_k> = k(\frac{-1}{49})$ for $1 \leq k \leq 6$, and so
\begin{align*}
K|_{C_{7}} \cdot [\omega|_{C_{7}}] &= 7[\gamma_6] \cdot [\omega|_{C_{7}}] \\
&= (7) \Big(\frac{-1}{49} \Big) \Big( 1(b_{4} - b_{7}) + 2(b_{1} - b_{4}) + 3(a-b_{1}-b_{2}-b_{3}) \\
& + 4(b_{2} - b_{5}) + 5(b_{5} - b_{9}) \\ 
& + 6 \big(12a + b_{9} - 4(b_{1} + \dots + b_{9}) - 2(b_{10} + \dots + b_{13}) \big) \Big)
\end{align*}

Therefore
\begin{align}
K|_{C_{7}} \cdot [\omega|_{C_{7}}] &= \Big( \frac{-1}{7} \Big) \Big( 75a -25 b_{1} - 23 b_{2} - 27 b_{3} - 25 b_{4} - 23 b_{5} - 24 b_{6} \nonumber\\
& - 25 b_{7} - 24 b_{8} - 23 b_{9} - 12 (b_{10} + b_{11} + b_{12} + b_{13}) \Big) \label{Kc7cdotomegac7}
\end{align}

And using the relation in $\eqref{lemma215rel}$ above, along with $\eqref{kcdotomegalab}$ and $\eqref{Kc7cdotomegac7}$, we have
\begin{align*}
K_{7} \cdot [\omega_{7}] &= K \cdot [\omega] - K|_{C_{7}} \cdot [\omega|_{C_{7}}] \\
&= (-3a + (b_{1} + \dots + b_{13})) - K|_{C_{7}} \cdot [\omega|_{C_{7}}] \\
&= \Big( \frac{1}{7} \Big) \Big(-21a + 7(b_{1} + \dots + b_{13}) \Big) - K|_{C_{7}} \cdot [\omega|_{C_{7}}] \\
&= \Big( \frac{1}{7} \Big) \Big(-21a + 7(b_{1} + \dots + b_{13}) + 75a -25 b_{1} - 23 b_{2} \\
& \qquad \qquad - 27 b_{3} - 25 b_{4} - 23 b_{5} - 24 b_{6} - 25 b_{7} - 24 b_{8} - 23 b_{9} \\
& \qquad \qquad - 12 (b_{10} + b_{11} + b_{12} + b_{13}) \Big) \\
&= \Big( \frac{1}{7} \Big) \Big( 54a -18 b_{1} -16 b_{2} - 20 b_{3} - 18 b_{4} - 16 b_{5} - 17 b_{6} \\
& - 18 b_{7} - 17 b_{8} - 16 b_{9} - 5 (b_{10} + b_{11} + b_{12} + b_{13}) \Big)
\end{align*}

and using the inequality

\begin{equation*}
3a > b_{1} + \dots + b_{13}
\end{equation*}

which implies
\begin{equation*}
54a > 18 b_{1} + \dots + 18 b_{13}
\end{equation*}

we get
\begin{equation*}
K_{7} \cdot [\omega_{7}] > \Big( \frac{1}{7} \Big) \Big( 2 b_{2} - 2 b_{3} + 2 b_{5} + b_{6} + b_{8} + 2 b_{9} + 13 (b_{10} + b_{11} + b_{12} + b_{13}) \Big)
\end{equation*} 

and since $b_{2}> b_{3}$, we have $2 b_{2} - 2 b_{3} > 0$, which in turn implies
\begin{equation*}
K_{7} \cdot [\omega_{7}] > \Big( \frac{1}{7} \Big) \Big( 2 b_{5} + b_{6} + b_{8} + 2 b_{9} + 13 (b_{10} + b_{11} + b_{12} + b_{13}) \Big)
\end{equation*} 

and since $b_{i} \geq 0$ for $1 \leq i \leq 13$, we finally have
\begin{equation*}
K_{7} \cdot [\omega_{7}] > 0
\end{equation*} 
This proves the lemma. $\Box$ \newline

\pagebreak
\section{The Next Constructions}

The first example of a manifold homeomorphic but not diffeomorphic to $\cpt \# 6 \cptbar$ was constructed by A. Stipsicz and Z. Szab\'o in \cite{SS}. They started with an elliptic fibration on $\cpt \# 9 \cptbar$ which had an $\tilde{E_7}$ fibre and three fishtail fibres, and (after a few blow-ups) then performed a generalized version of a rational blowdown to obtain a manifold homeomorphic to $\cpt \# 6 \cptbar$. Then, computation of the Seiberg-Witten invariants of this manifold showed that it was not diffeomorphic to $\cpt \# 6 \cptbar$. \newline

Up until this point, only finitely many non-diffeomorphic exotic smooth 4-manifolds had been found for $\cpt \# n \cptbar$ (for $n = 6,7,8$). R. Fintushel and R. Stern introduced a new technique in \cite{FS3} which can be used to construct infinite families of non-diffeomorphic exotic $\cpt \# n \cptbar$'s, for $n=6,7,8$. \newline

A few days after a preprint of \cite{FS3} was posted on the arXiv, J. Park, A. Stipsicz and Z. Szab\'o posted a preprint of \cite{PSS} on the arXiv, in which they used this technique to construct an infinite family of non-diffeomorphic exotic $\cpt \# 5 \cptbar$'s. \newline

We give an outline of the technique presented in \cite{FS3} below, as well as an outline of how it was applied to construct exotic $\cpt \# n \cptbar$'s (for $n=5,6$). First, however, we need to review R. Fintushel and R. Stern's \emph{knot surgery} technique presented in \cite{FS2}, as well as the concept of a \emph{mapping class group} and the \emph{monodromy} of a singular fibre.

\pagebreak


\section{Fibre Sums and Embedded Tori}

For fibre sums, one can refer to \cite{GS}, and for knot surgery the original source is \cite{FS2}, but a good presentation of both can be found in \cite{Sc}.

\subsection*{Fibre sums}
Suppose we have two $C^{\infty}$-elliptic fibrations (see section \ref{ellfibsec})
\begin{align*}
& \pi_1 : S_1 \longrightarrow C_1 \\
& \pi_2 : S_2 \longrightarrow C_2
\end{align*}
Choose $t_1 \in C_1$ and $t_2 \in C_2$ such that the fibres $F_1 = \pi_1 ^{-1}(t_1)$ and $F_2 = \pi_2 ^{-1}(t_2)$ are tori (generic fibres). \newline

A tubular neighbourhood $\nu F_i$ of each $F_i$ in $S_i$ is a copy of $D^2 \times T^2$ in each $S_i$ ($i=1,2$). Then $\partial(S_i \setminus \nu F_i) \cong T^3$ (for $i=1,2$), and we choose a fibre-preserving, orientation-reversing diffeomorphism of $T^3$
\begin{equation*}
\phi : \partial(S_1 \setminus \nu F_1) \longrightarrow \partial(S_2 \setminus \nu F_2)
\end{equation*}
Then, the \emph{fibre sum} $S_1 \# _{f} S_2$ is defined as the manifold $(S_1 \setminus \nu F_1) \cup _{\phi} (S_2 \setminus  \nu F_2)$. \newline

Note that $S_1 \# _{f} S_2$ will admit a $C^{\infty}$-elliptic fibration $\pi : S_1 \# _{f} S_2 \longrightarrow C_1 \# C_2$.

\newtheorem{fibresumrem1}{Remark}[section]
\begin{fibresumrem1}
\upshape
Although the diffeomorphism type of $S_1 \#_{f} S_2$ might depend on the choice of the diffeomorphism $\phi$, if either elliptic fibration $\pi_i : S_i \longrightarrow C_i$ ($i=1,2$) contains a cusp fibre, then for any choice of $\phi$ the manifolds $(S_1 \setminus \nu F_1) \cup _{\phi} (S_2 \setminus  \nu F_2)$ will be diffeomorphic, and then $S_1 \#_{f} S_2$ is a well-defined 4-manifold. See Chapter 8 in \cite{GS}.
\end{fibresumrem1}

\newtheorem{fibresumrem2}[fibresumrem1]{Remark}
\begin{fibresumrem2}
\upshape
It should be remarked that a tubular neighbourhood is often called a \emph{regular} neighbourhood.
\end{fibresumrem2}

\subsection*{Near-cusp embedded tori}

Before we get to the knot surgery technique itself, we first need a definition from \cite{Sc}.
\newtheorem{ncetor}[fibresumrem1]{Definition}
\begin{ncetor}
\upshape
Let $X$ be a simply-connected 4-manifold. Let $T$ be a torus embedded in $X$ that is homologically nontrivial and has zero self-intersection. Such a torus $T$ is called \emph{near-cusp embedded} if and only if a neighbourhood of $T$ in $X$ is diffeomorphic to a neighbourhood $U$ of a generic torus fibre inside some elliptic fibration, so that $U$ contains a cusp fibre and so that $T$ corresponds to a regular (generic) fibre.
\end{ncetor}

\pagebreak


\section{Knot Surgery} \label{knotsurgerysec}
We now describe the knot surgery technique described in \cite{FS2}. Another good explanation can be found in \cite{Sc}. \newline

Let us start with a closed simply-connected 4-manifold $X$, which contains a near-cusp embedded torus $T$. Let $\nu T$ be a tubular neighbourhood of $T$ (i.e. a copy of $T \times D^2$), and consider the manifold $X \setminus \nu T$. Then $\partial(X \setminus \nu T) \cong \partial (\nu T) \cong S^1 \times S^1 \times S^1 \cong T^3$. \newline

Now, let $K$ be a knot in $S^3$. Let $\nu K$ be a tubular neighbourhood of $K$, and consider the manifold $S^3 \setminus \nu K$ (the \emph{knot complement}). Homologically, $S^3 \setminus \nu K$ is indistinguishable from a solid torus $S^1 \times D^2$, and it therefore has boundary $\partial (S^3 \setminus \nu K) \cong S^1 \times S^1$. Therefore, the boundary of $S^1 \times (S^3 \setminus \nu K)$ is also the 3-torus $S^1 \times S^1 \times S^1 \cong T^3$. \newline

We shall glue $X \setminus \nu T$ and $S^1 \times (S^3 \setminus \nu K)$ together along their boundaries. However, there are several choices we can make in how we glue their boundaries to each other. \newline

Firstly, $\partial(\nu K) \cong \partial (S^3 \setminus \nu K)  \cong T^2$, a torus. Let $\lambda$ be a longitude of $\partial(\nu K)$ (and therefore a meridian of $\partial (S^3 \setminus \nu K)$). Consider the homology class $[pt \times \lambda] \in H_1(S^1 \times (S^3 \setminus \nu K); \mz)$. \newline

Secondly, $\partial(\nu T) \cong \partial(T \times D^2)$, and we consider the homology class $[pt \times \partial D^2] \in H_1( T \times D^2;\mz)$. Since $\partial(\nu T) \cong \partial(X \setminus \nu T)$, we also consider $[pt \times \partial D^2]$ as a homology class in $H_1( X \setminus \nu T;\mz)$ \newline

We define the manifold
\begin{equation*}
X_K = (X \setminus \nu T) \cup _{\phi} (S^1 \times (S^3 \setminus \nu K))
\end{equation*}
where $\phi: \partial(X \setminus \nu T) \longrightarrow \partial(S^1 \times (S^3 \setminus \nu K))$ is an orientation-reversing diffeomorphism such that $[pt \times \partial D^2]$ is identified with $[pt \times \lambda]$.

\newtheorem{ksrem1}{Remark}[section]
\begin{ksrem1}
\upshape
We did not specify where the homology class $[pt \times \mu] \in H_1(S^1 \times (S^3 \setminus \nu K); \mz) $ is mapped to, where $\mu$ is the meridian of $\partial(\nu K)$. In fact, $[pt \times \mu]$ can be mapped to any generator of $H_1(T; \mz)$ in $T \times \partial D^2$. 

Since we assumed that $T$ is near-cusp embedded, the diffeomorphism type of $X_K$ is completely determined by our choice of mapping $[pt \times \lambda]$ to $[pt \times D^2]$, since (because $T$ is near-cusp embedded) the Seiberg-Witten invariant of $X_K$ is completely determined by the Seiberg-Witten invariant of $X$ and the Alexander polynomial of $K$. \newline

This means that all such constructed $X_K$'s (where $[pt \times \mu]$ can be mapped to any generator of $H_1(T; \mz)$) have the same Seiberg-Witten invariant. It is not known whether these (different) $X_K$'s are all diffeomorphic (\cite{FS2}); recall that if $Y$ and $Z$ are two smooth 4-manifolds that have different Seiberg-Witten invariants, then they are definitely non-diffeomorphic, but if $Y$ and $Z$ have the same Seiberg-Witten invariant, then they could be diffeomorphic or non-diffeomorphic.
\end{ksrem1}

Let us review the knot surgery construction. We started with a closed, simply-connected smooth 4-manifold $X$. We removed $\nu T \cong T \times D^2$, and then ``glued back'' a homological copy of $T^2 \times D^2$ (our $S^1 \times (S^3 \setminus \nu K)$), and called this new manifold $X_K$. Therefore, the homology of $X_K$ is the same as the homology of $X$, and by the corollary of Freedman's Classification Theorem, Corollary \ref{frcor1}, $X_K$ is homeomorphic to $X$. \newline

However, the Seiberg-Witten invariants of $X_K$ and $X$ will be different for most choices of $K$, and therefore in general $X_K$ will not be diffeomorphic to $X$.

\newtheorem{ksrem2}[ksrem1]{Remark}
\begin{ksrem2}
\upshape
If we also assume (or choose our torus $T$) so that $X \setminus \nu T$ is simply-connected, then the fact that $X$ and $X \setminus \nu T$ are simply-connected implies that $X_K$ is also simply-connected (stated in \cite{FS3}).
\end{ksrem2}

\pagebreak


\section{Fishtail and Cusp Fibres Revisited} \label{fibresrevisitedsec}

We now look at fishtail and cusp fibres from a viewpoint different to that of section \ref{fishtailsection}. \newline

By definition, a generic fibre in an elliptic fibration is a torus. In section \ref{fishtailsection}, it was shown that a fishtail fibre is a sphere with a point of transverse self-intersection. The fishtail fibre ``appears'' in an elliptic fibration by collapsing a homologically nontrivial circle, in a nearby generic torus fibre, to a point (\cite{Sc}). Such a circle is called a \emph{vanishing cycle}. It is a circle that bounds a disk of self-intersection $-1$ in the fibre's complement (for a computation of this fact, see \cite{GS} pages 292-293). \newline

Explicitly, suppose that $\pi : S \longrightarrow C$ is an elliptic fibration and that $t_1 \in C$ is a point such that $F_1 = \pi^{-1}(t_1)$ is a fishtail fibre. Then, in a neighbourhood of $t_1$ in $C$, there is a point $t$ such that $F = \pi^{-1}(t)$ is a generic torus fibre with a vanishing cycle $v$. As $t \rightarrow t_1$, $v$ collapses to a point. See Figure \ref{fibresrevisitedsec}.1 and \ref{fibresrevisitedsec}.2. \newline

In section \ref{fishtailsection}, it is shown that a cusp fibre is also a sphere, but with one singular point. Suppose again that $\pi : S \longrightarrow C$ is an elliptic fibration, and suppose that $t_2 \in C$ is a point such that $F_2 = \pi^{-1}(t_2)$ is a cusp fibre. Then a nearby (generic) torus fibre has two vanishing cycles $v_1$ and $v_2$ which collapse to a single point (we can think of these circles, which are generators of $H_1(T^2; \mz)$, as the meridan and longitude of the torus). See Figure \ref{fibresrevisitedsec}.3. The singular point of a cusp fibre has a neighbourhood which looks like a cone over the trefoil knot (see \cite{Sc}, \cite{GS}). \newline

For Kirby diagrams of a fishtail fibre and of a cusp fibre, see \cite{GS}, page 299.

\pagebreak

\begin{center}
\begin{minipage}{12cm}
\includegraphics[width=12cm]{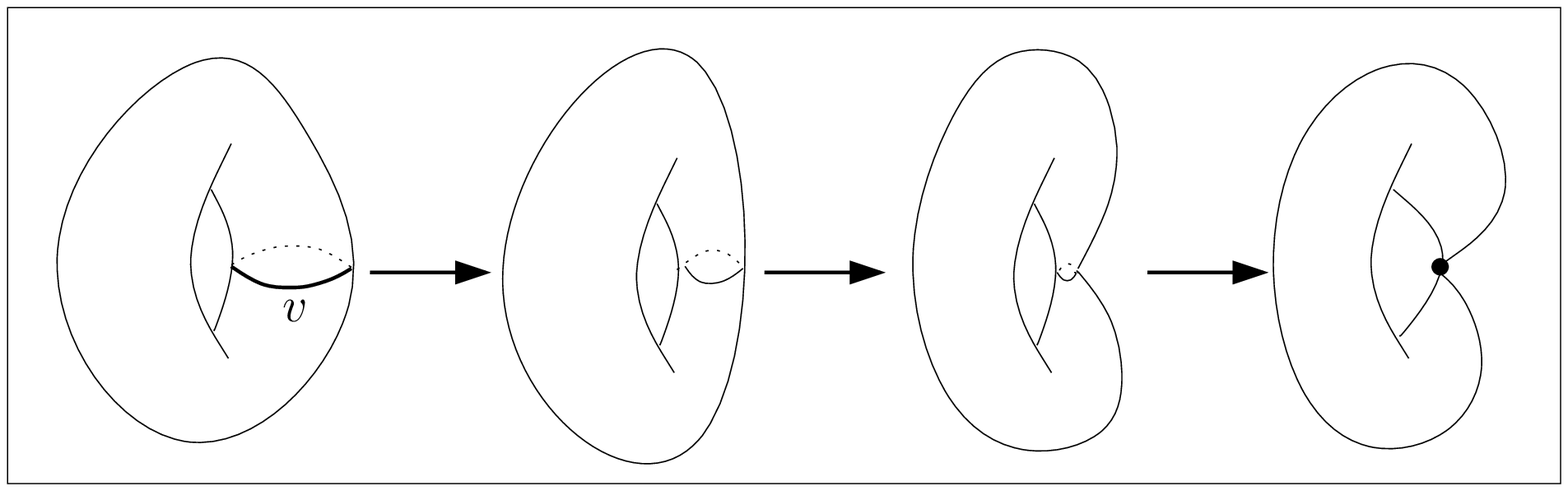}
\end{minipage}
\end{center}
\begin{center}
Figure \ref{fibresrevisitedsec}.1
\end{center}

\begin{center}
\begin{minipage}{12cm}
\includegraphics[width=12cm]{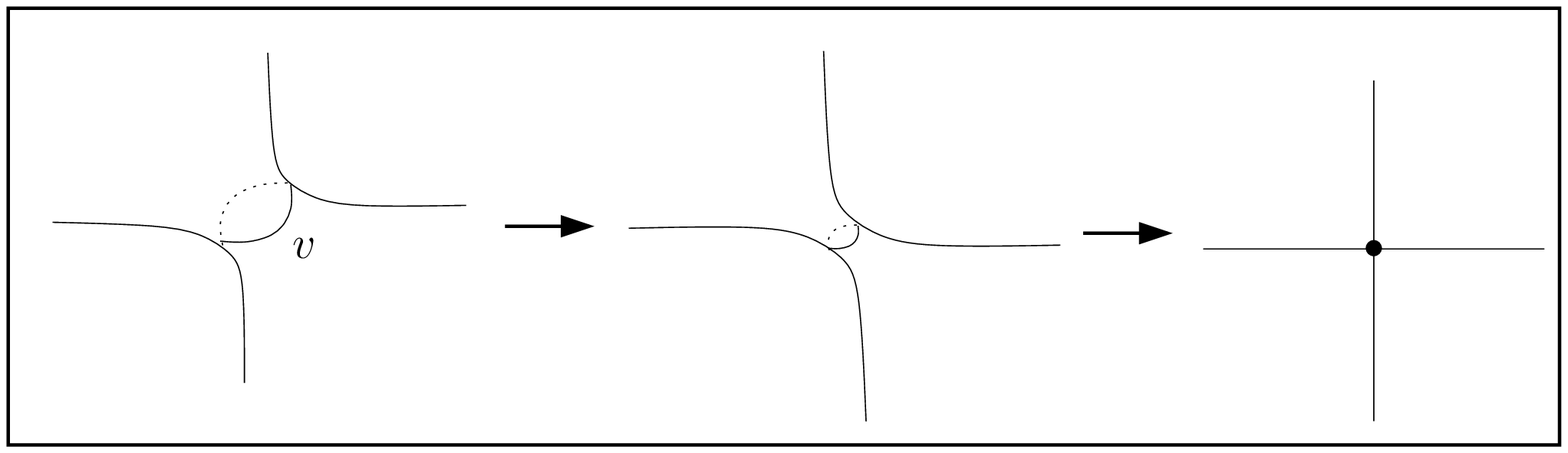}
\end{minipage}
\end{center}
\begin{center}
Figure \ref{fibresrevisitedsec}.2
\end{center}

\begin{center}
\begin{minipage}{12cm}
\includegraphics[width=12cm]{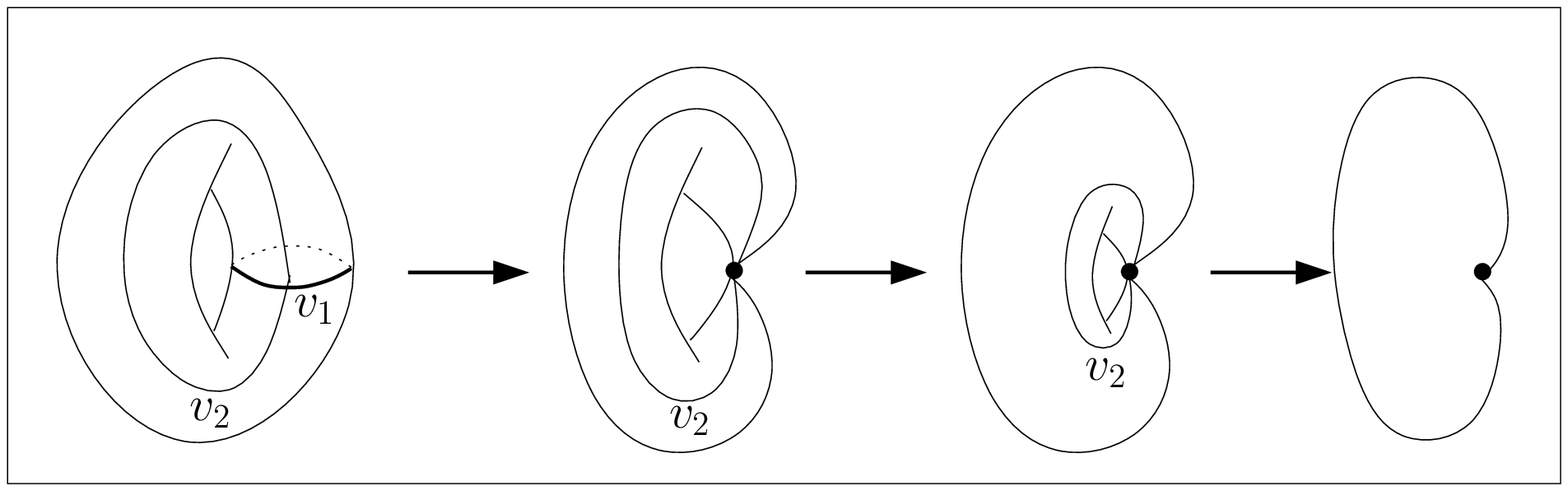}
\end{minipage}
\end{center}
\begin{center}
Figure \ref{fibresrevisitedsec}.3
\end{center}

\pagebreak

\section{The Mapping Class Group of the Torus} \label{mpgtorsec}

\newtheorem{mcgtdef}{Definition}[section]
\begin{mcgtdef}
\upshape
The \emph{mapping class group} of the torus, labelled $\Gamma_1$, is defined to be the set of isotopy classes of orientation-preserving diffeomorphisms of $T^2$.
\end{mcgtdef}

More formally, if we let $\mathrm{Diff}^{+}(T^2)$ be the set of orientation-preserving diffeomorphisms of $T^2$, and we define $\mathrm{Diff}_{0}(T^2)$ to be the set of diffeomorphisms of $T^2$ that are isotopic to the identity map, then
\begin{equation*}
\Gamma_1 \cong \mathrm{Diff}^{+}(T^2) / \mathrm{Diff}_0(T^2)
\end{equation*}

It was first proved by M. Dehn in \cite{De} that the mapping class group $\Gamma_g$ of a genus-$g$ surface is generated by finitely many \emph{twist homeomorphisms} (now called \emph{Dehn twists}). In fact, only $3g-1$ Dehn twists are needed to generate $\Gamma_g$ (\cite{L2}, see also \cite{L1} and \cite{L3}). The following definition is from \cite{GS}.

\newtheorem{dtfm}{Definition}[section]
\begin{dtfm}
\upshape
Let $S$ be a surface and let $C$ be a circle in $S$. A right-handed \emph{Dehn twist} $\psi:S \longrightarrow S $ is a diffeomorphism obtained by ``cutting'' $S$ along $C$, twisting a neighbourhood of one of the boundary components $360^\circ$ to the right, and then ``regluing'' the cut-out component ``back in''.
\end{dtfm}

We, of course, also have a more formal definition (also from \cite{GS}).
\newtheorem{dtfm2}[dtfm]{Definition}
\begin{dtfm2}
\upshape
Let $S$ be a surface and let $C$ be a circle in $S$. Identify $\nu C$, a neighbourhood of $C$, with $S^1 \times I$. Define the map $\psi$ such that on $\nu C$ 
\begin{equation} \label{psilab1}
\psi(\theta, t) = (\theta + 2 \pi t, t) 
\end{equation}
and $\psi$ is the identity map on $S \setminus \nu C$ (and $\psi$ goes smoothly from $\mathrm{id}|_{S \setminus \nu C}$ to $\eqref{psilab1}$). Then $\psi$ is a right-handed \emph{Dehn twist}.
\end{dtfm2}

For an example of a Dehn twist on a cylinder, see Figure \ref{mpgtorsec}.1 below (from \cite{GS}). \newline

According to \cite{SSS}, $\Gamma_1$ admits a presentation (see also chapter 7 in \cite{CM})
\begin{equation*}
\Gamma_1 = <a,b \; | \; aba=bab, (ab)^6 =1 >
\end{equation*}

It is shown in \cite{FM} that $\Gamma_1$ is isomorphic to $SL(2;\mz)$, the group of $2 \times 2$ matrices with integer entries and determinant 1. In fact, the map $\phi: \Gamma_1 \longrightarrow SL(2;\mz)$, defined in \cite{SSS} by
\begin{align*}
\phi(a) &=
\left( \begin{array}{cc}
1 & 1 \\
0 & 1
\end{array} \right) \\
\phi(b) &=
\left( \begin{array}{cc}
\phantom{1} 1 & 0 \\
-1 & 1
\end{array} \right) \\
\end{align*}
provides us with an isomorphism between $\Gamma_1$ and $SL(2;\mz)$. It is easily checked that $\left( \begin{array}{cc}
1 & 1 \\
0 & 1
\end{array} \right)$ and 
$\left( \begin{array}{cc}
\phantom{1}1 & 0 \\
-1 & 1
\end{array} \right)$ generate $SL(2;\mz)$ and that 
\begin{align*}
\phi(aba) &=
\left( \begin{array}{cc}
\phantom{1} 0 & 1 \\
-1 & 0
\end{array} \right) = \phi(bab) \\
\phi((ab)^6) &=
\left( \begin{array}{cc}
 1 & 0 \\
 0 & 1
\end{array} \right) \\
\end{align*}
For more about Dehn twists and mapping class groups, see \cite{FM}.

\begin{center}
\begin{minipage}{10cm}
\includegraphics[width=10cm]{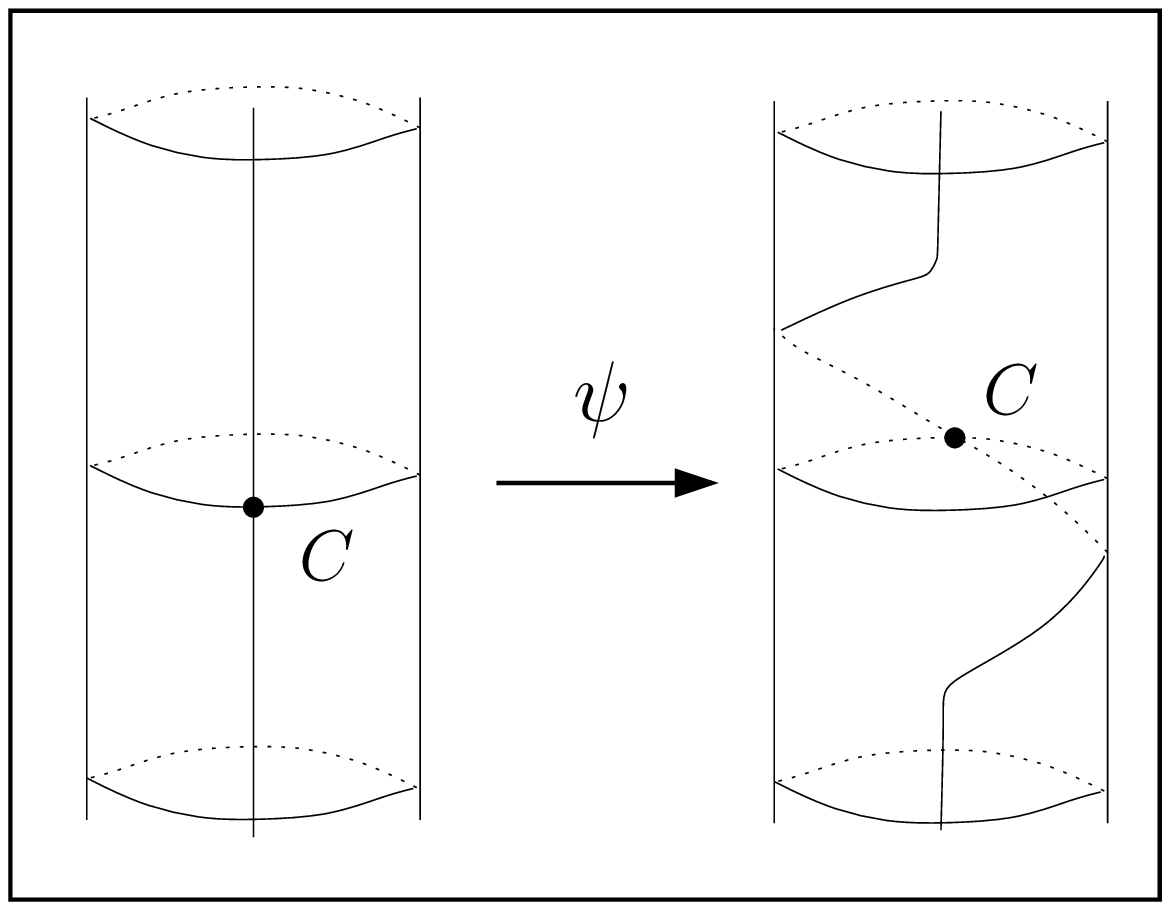}
\end{minipage}
\end{center}
\begin{center}
Figure \ref{mpgtorsec}.1
\end{center}
\pagebreak


\section{Monodromy and the Existence of Certain Elliptic Fibrations} \label{monodromysec}

We consider an elliptic fibration 
\begin{equation} \label{sfib1lab}
\pi : \cpt \# 9 \cptbar \longrightarrow \cpo
\end{equation}

There are only finitely-many singular points $\{s_1, s_2, \dots, s_n \} \subset \cpo$ (i.e. points $s$ such that $\pi ^{-1}(s)$ is a singular fibre of the type listed in Kodaira's table of singular fibres). In fact, there are at most 12 singular points, since each singular point corresponds to a singular fibre which has Euler characteristic at least 1, and $\chi(\cpt \# 9 \cptbar) = 12$. \newline

Since there are only finitely-many singular points $\{s_1, s_2, \dots, s_n \} \subset \cpo$, we can find disjoint disks $\{D_1, D_2, \dots D_n \}$ such that $s_i \in \mathrm{int}(D_i)$ and $D_i \cap D_j = \emptyset$ if $i \neq j$ ($i,j = 1,2,\dots n$). \newline

Now consider one of these singular points $s$, and its disk neighbourhood $D$, and let $C = \partial D$. Since every $t\in C \subset \cpo$ is such that $\pi^{-1}(t)$ is a (generic) $T^2$ fibre, if we restrict the fibration in $\eqref{sfib1lab}$ to $C$, by traversing along $C$ we get a diffeomorphism $\psi$ of $T^2$ (\cite{SSS}), and this diffeomorphism is defined up to isotopy and conjugation (\cite{SSS}). The corresponding element (defined only up to conjugation) of $\psi$ in $\Gamma_1$, the mapping class group of the torus, is called the \emph{monodromy} of the singular fibre $F = \pi^{-1}(s)$. \newline

The following argument is stated in \cite{SSS}: suppose that $w$ is a word $\Gamma_1$ which is composed of 12 right-handed Dehn twists $w_1, w_2, \dots, w_{12}$ satisfying
\begin{equation*}
w = w_1 \dots w_{12} = (w_1 \dots w_{i_1})(w_{i_1+1}\dots w_{i_2})\dots(w_{i_k+1}\dots w_{12}) = 1
\end{equation*}
in $\Gamma_1$. Then for the singular fibres $F_j$ $(j=1, \dots , k)$ with monodromies conjugate to $(w_{i_j+1}\dots w_{i_{j+1}})$ there is an elliptic fibration $\pi : \cpt \# 9 \cptbar \longrightarrow \cpo$ with singular fibres $F_1, \dots, F_k$. \newline

Table \ref{monodromysec}.1 gives a list of the monodromies (in $\Gamma_1$ and $SL(2;\mz)$) of the singular fibres in Kodaira's list of singular fibres. This list comes from \cite{SSS} and \cite{HKK}.
\pagebreak
\begin{align*}
& \mathrm{Singular \; \; fibre} && \Gamma_1 && SL(2;\mz) \\
& I_1 && a &&  \left( \begin{array}{cc} \ 1 & 1 \\ 0 & 1 \end{array} \right) \\
& I_n && a^n &&  \left( \begin{array}{cc} 1 & n \\ 0 & 1 \end{array} \right) \\
& II && ba &&  \left( \begin{array}{cc} \phantom{1} 1 & 1 \\ -1 & 0 \end{array} \right)\\
& III && aba=bab && \left( \begin{array}{cc} \phantom{1} 0 & 1 \\ -1 & 0 \end{array} \right) \\
& IV && (ba)^2 && \left( \begin{array}{cc} \phantom{1} 0 & 1 \\ -1 & -1 \end{array} \right)\\
& I_n^{*} && (ab)^3a^n && \left( \begin{array}{cc} -1 & -n \\ \phantom{1}0 & -1 \end{array} \right)\\
& \tilde{E_8} && (ba)^5 && \left( \begin{array}{cc} 0 & -1 \\ 1 & \phantom{1} 1 \end{array} \right) \\
& \tilde{E_7} && (ba)^4 b &&\left( \begin{array}{cc} 0 & -1 \\ 1 & \phantom{1} 0 \end{array} \right) \\
& \tilde{E_6} && (ba)^4 &&\left( \begin{array}{cc} -1 & -1 \\ \phantom{1}1 & \phantom{1} 0 \end{array} \right)
\end{align*}
\begin{center}
Table \ref{monodromysec}.1
\end{center}

\pagebreak

\section{Double Node Neighbourhoods and Exotic $\cpt \# 5 \cptbar$'s} \label{dnnsec}

We give an outline of the double node neighbourhood technique presented in \cite{FS3}. \newline

A \emph{double node neighbourhood} $D$ is a (fibred) neighbourhood of an elliptic fibration that contains two fishtail fibres with the same monodromy\footnote{This is equivalent to saying that the neightbourhood contains an $I_2$-fibre.}. This means that there is a smooth torus fibre of $D$ that has two vanishing cycles $C_1$ and $C_2$ that collapse to a point over the points $p_1$ and $p_2$ in $D$, respectively. \newline

We now consider the effect of performing knot surgery along a (generic) torus fibre $F$ in a double node neighbourhood using a twist knots $K = T(n)$, pictured in Figure \ref{dnnsec}.1 (from \cite{FS3}). Recall that although we are forced to send the homology class of the longitude of our knot $K$ to a certain homology class, we are free to send the meridian of our knot $K$ to any homology class (that is a generator) of our choice. We choose the gluing in the knot surgery construction to be such that we send the homology class of a meridian of $K$ to the class $a \times pt$, the class of the vanishing cycle. \newline

It is proved in \cite{FS3} that while $D$ had section which was a disk $D_1$, the effect of knot surgery is that $D_1$ has a smaller disk $D_2$ removed from it and a punctured torus glued onto the boundary of $D_2$ in $D_1$ (essentially, connect-summing with a torus). Furthermore, this punctured torus contains a loop which bounds a disk $U$ self-intersection $-1$. We can then perform surgery on an annular neighbourhood of $U$ that will result in the torus becoming an immersed sphere $S$ of self-intersection $-1$. Note that this step is nontrivial and is the key argument in \cite{FS3}.  \newline

Now, R. Fintushel and R. Stern proceed in \cite{FS3} to construct an infinite family of $\cpt \# 6 \cptbar$'s as follows: we recall that $E(1)$ has an elliptic fibration with one $\tilde{E_6}$ fibre and four fishtail fibres. Since $(ab)^6 = 1$ in $\Gamma_1$, the factorisation $(ab)^6 = (ab)^4 a^2 (a^{-1}ba) b$ shows the existence of such a fibration, since 
\begin{itemize}
\item[(i)] $b(ab)^4 b^{-1} \sim (ba)^4$ is the monodromy of the $\tilde{E_6}$ fibre, 
\item[(ii)] $a$ is the monodromy  of a fishtail fibre,
\item[(iii)] a conjugate of $b$ in $\Gamma_1$ is $(ba)b(ba)^{-1} = (bab)a^{-1}b^{-1} = (aba)a^{-1}b^{-1} = a$ (where we used the braid relation $aba=bab$ in the second step), and so $b$ is also a fishtail fibre
\item[(vi)] $a^{-1}ba$ is a conjugate of $b$, and therefore of $a$, in $\Gamma_1$ and so is also a fishtail fibre.
\end{itemize}
and, furthermore, two of the fishtail fibres have the same monodromy $a$. \newline

We can therefore find a double node neighbourhood $D \subset E(1)$ containing two fishtail fibres with the same monodromy, and $E(1)\setminus D$ contains an $\tilde{E_6}$ fibre and the two remaining fishtail fibres $F_1$ and $F_2$. Choosing the knot $K$ to be the twist knot $T(n)$, we perform the knot surgery on $D$ as above, and therefore $D_K$ contains an immersed sphere $S$ of self-intersection $-1$. We then glue $D_K$ back into $E(1) \setminus D$ to obtain the manifold $Y_{n} = (E(1) \setminus D) \cup D_K$, and then the two fishtail fibres $F_1$ and $F_2$ intersect $S$ transversely in a single point. \newline

We blow-up at the double points of $S$, $F_1$ and $F_2$ to obtain spheres of self-intersection $-5$, $-4$ and $-4$, respectively, which intersect in a pair of points. These intersection points are singular points which can be smoothed to obtain a sphere $R$ of self-intersection $-9$. Furthermore, $R$ intersects the $S_5$ sphere in $\tilde{E_6}$ (one of the spheres of multiplicity 1) in a single positive point (i.e. $R \cdot S_5 = +1$). We are therefore in the same position to perform a rational blowdown along a configuration of spheres $C_7$, as in Park's construction, except that we have only used three blow-ups of $E(1)$ (actually, $Y_n$, which is homeomorphic to $E(1)$) instead of the four done in \cite{P1} and in section \ref{mainsectionlabel} above. \newline

Therefore, the rational blowdown along $C_7$ will produce a manifold $X_n$ that is homeomorphic to $\cpt \# 6 \cptbar$. These manifolds $X_n$ will not be diffeomorphic to $\cpt \# 6 \cptbar$, and in fact because of the knot surgery, will not even be diffeomorphic to each other. Therefore, we have constructed an infinite family of 4-manifolds homeomorphic to $\cpt \# 6 \cptbar$ but not diffeomorphic to it.  

Several constructions are presented in \cite{PSS} to construct infinite families of manifolds homeomorphic but not diffeomorphic to $\cpt \# 5 \cptbar$. One particularly nice construction is started by first showing that there is an elliptic fibration on $E(1)$ with an $I_6$ fibre and six fishtail fibres, one with monodromy $m_1$, one with monodromy $m_2$, two with monodromy $m_3$ and the last two with monodromy $m_4$, where the $m_i$ ($i=1,2,3,4$) are conjugates of $a$ in $\Gamma _1$ (explicitly, this is done by using the braid relation to turn a conjugate of $(ab)^6$ into $(a^3 b)^3$ which factorises as $a^6 (a^{-3}b a^3)(bab^{-1})^2 b^2 (b^{-1}ab)$). \newline 

Two pairs of fishtail fibres with the same monodromy allow us to perform the double node neighbourhood knot surgery twice, and we get to the stage of having a sphere of self-intersection $-9$ in a manifold $V_n \# 2 \cptbar$, where $V_n$ is some manifold (the result of the knot surgery with the twist knots $T(1)$ and $T(n)$) that is homeomorphic to $E(1)$, and we are again ready to perform a rational blowdown along a configuration $C_7$. This time, the resulting 4-manifolds $Q_n$ will be homeomorphic (but not diffeomorphic) to $\cpt \# 5 \cptbar$, and the manifolds $Q_n$ will all be non-diffeomorphic. \newline

This short account of a few results in \cite{FS3} and \cite{PSS} do not do justice to these papers, which contain far more detail.

\begin{center}
\begin{minipage}{10cm}
\includegraphics[width=10cm]{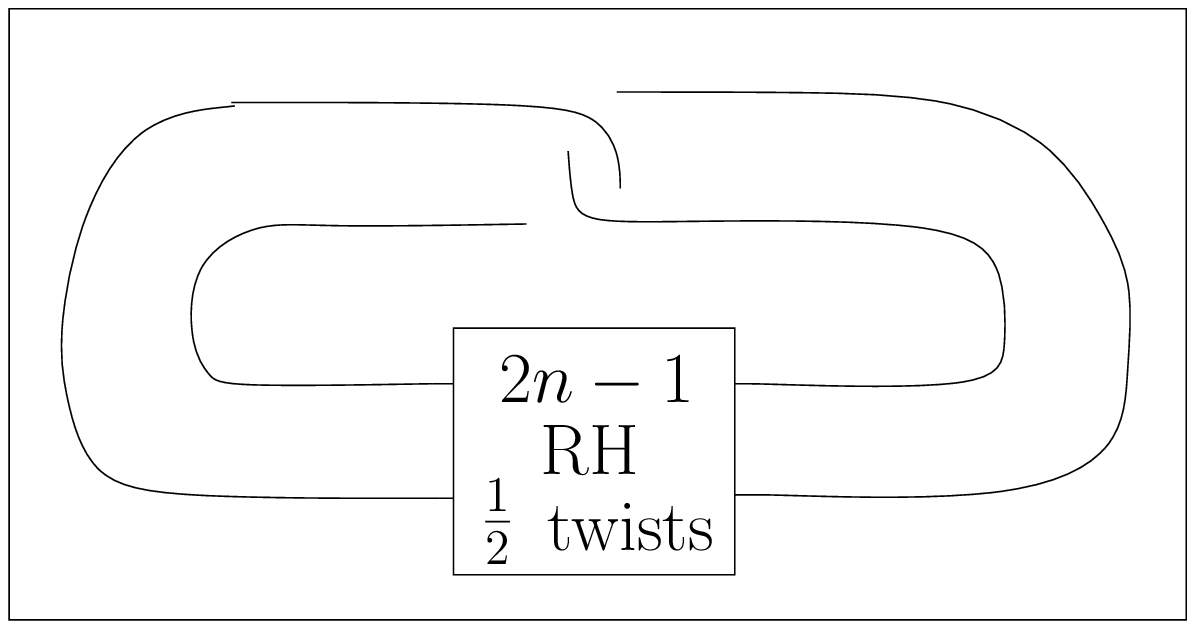}
\end{minipage}
\end{center}
\begin{center}
Figure \ref{dnnsec}.1
\end{center}

\pagebreak

\section{Epilogue}
In early 2007, the first example of an exotic $\cpt \# 3 \cptbar$ was given in \cite{AP} by A. Akhmedov and B. Doug Park. The reader is also referred to the subsequent papers \cite{BK1} and \cite{ABP}.  \newline

K. Yasui recently posted a paper (\cite{Ya}) on the arXiv giving constructions of $\cpt \# n \cptbar$ ($n = 5,6,7,8,9$) using rational blowdowns and without using elliptic fibrations. \newline

In conclusion, over the last few years mathematicians, in a way, have got closer to finding an exotic $\cpt$. However, it is still unknown whether any of the following 4-manifolds
\begin{equation*}
\cpt, \qquad \cpt \# \cptbar, \qquad \cpt \# 2 \cptbar, \qquad S^2 \times S^2, \qquad S^4 
\end{equation*}
admits an exotic smooth structure.

\pagebreak
\section{A Note About Sources, Diagrams and Corrections}

Before discussing the books and papers, I must give credit where it is due. My advisor, David Gay, helped me immensely with all the background material and sections 2$-$14 and 21$-$24, and answered many, many questions on the original thesis overall and topology in general. Andr\'as Stipsicz gave me several excellent lectures on singular fibres in elliptic fibrations and double-node neighbourhoods, and patiently answered my multitude of questions, which has resulted in sections 11, 15$-$20 and 25$-$28. I cannot overstate their contribution. \newline

Although the bibliography is quite large, most references are only used for a single result. Only a few books or articles have been used extensively. \newline

\cite{GS} is the key background reference, containing chapters on the classification of topological 4-manifolds, blowing up and blowing down, rational blowdowns, Kirby calculus, elliptic fibrations, $\dots$ the list goes on. It would be helpful to have already read \cite{GP}, or to have a copy at hand, before reading this article. \newline

\cite{P1} and \cite{PSS} are the papers this article is ``built around'', and the main goal of writing this article was to make these papers more accessible. \newline

The next most important are the papers \cite{FS1}, \cite{FS2} and \cite{FS3}, the original sources for rational blowdowns, knot surgery and double node neighbourhoods. \newline

\cite{SSS} is my favourite source for singular fibres, although other good sources are \cite{HKK} and \cite{BPV}. \newline

\cite{Sc} is an excellent book that covers a lot of the material presented in \cite{GS}, although is often not as detailed. It also has an excellent section on knot surgery. \newline

\cite{MS2} is an excellent reference for anything to do with symplectic topology. \newline

\cite{Ro} is my favourite source for the theory of knots and links. \cite{CF} also has a good section on group presentations not found in \cite{Ro}. \newline

Although these references listed above provided over 95\% of the source material for this article, the remaining 5\% provided by the rest of the references is no less important.\newline

Finally, I should mention that the diagrams included in this article were made using the wonderful program Ipe (version 6.0 preview 23).

\subsection*{Corrections}

There are several versions of this article available on the arXiv, of which this is version 4. \newline

Version 1 was the first draft, which I had hoped was error-free. \newline

In version 2 the order of attribution, regarding the construction of exotic $\cpt \# 3 \cptbar$s, at the beginning of the Epilogue was corrected. \newline

In version 3 the abstract was modified to make clear the purpose of this article.\newline

In version 4 Definition \ref{correcteddefnlabel} was modified, and the accompanying commutative diagram was removed. 
The details for the references for [R] and [U], which were neglected in the earlier versions, were finally filled in. Section 20 in version 3 was moved to become section 11, and was renamed ``Results Concerning Rational Blowdowns''. Consequently, a few sections have been shifted accordingly. Various typographical errors were corrected; however, I fear that there are still some out there.

\pagebreak



\end{flushleft}

\begin{thebibliography}{99}



\bibitem[ABP]{ABP}
\newblock A. Akhmedov, R. I. Baykur, B. Doug Park,
\newblock \emph{Constructing infinitely many smooth structures on small 4-manifolds},
\newblock Journal of Topology (2008) 1(2) 409-428

\bibitem[AP]{AP}
\newblock A. Akhmedov, B. Doug Park,
\newblock \emph{Exotic smooth structures on small 4-manifolds},
\newblock Invent. Math (2008) $\mathbf{173}$ 209-223


\bibitem[BK]{BK1}
\newblock S. Baldridge and P. Kirk,
\newblock \emph{A symplectic manifold homeomorphic but not diffeomorphic to $\cpt \# 3 \cptbar$}
\newblock Geometry and Topology $\mathbf{12}$ (2008), 919-940


\bibitem[BPV]{BPV}
\newblock W. Barth, C. Peters and A. Van de Ven,
\newblock \emph{Compact complex surfaces}, 
\newblock Ergebnisse der Mathematik und ihrer Grenzgebiete (3), vol. 4, Springer, Berlin, 1984

\bibitem[Bo]{Bo}
\newblock F. Bonahon,
\newblock \emph{Diff\'eotopes des espaces lenticulaires},
\newblock Topology $\mathbf{22}$ (1983), 305-314

\bibitem[CM]{CM}
\newblock H. S. M. Coxeter and W. O. J. Moser,
\newblock \emph{Generators and relations for discrete groups},
\newblock Springer-Verlag Ergebnisse der Mathematik (1957)



\bibitem[CF]{CF}
\newblock R. H. Crowell and R. H. Fox,
\newblock \emph{Introduction to knot theory},
\newblock Springer-Verlag, 1977

\bibitem[De]{De}
\newblock M. Dehn,
\newblock \emph{Die Gruppe der Abbildungsklassen}, 
\newblock Acta Math. $\mathbf{69}$ (1938), 135-206

\bibitem[Do]{Do}
\newblock S. K. Donaldson,
\newblock Irrationality and the $h$-cobordism conjecture,
\newblock J. Differential Geom. $\mathbf{26}$ (1987), 141-168

\bibitem[FM]{FM}
\newblock B. Farb and D. Margalit,
\newblock \emph{A primer on mapping class groups},
\newblock http://www.math.utah.edu/~margalit/primer/, version 2.95, 2007

\bibitem[FrM]{FrM}
\newblock R. Friedman and J. Morgan, 
\newblock \emph{Smooth 4-manifolds and complex surfaces},
\newblock Ergeb. Math. Grenzgeb., vol. 27, Springer-Verlag, 1994

\bibitem[FS1]{FS1}
\newblock R. Fintushel and R. Stern,
\newblock \emph{Rational blowdowns of smooth 4-manifolds},
\newblock Jour. Diff. Geom. $\mathbf{46}$ (1997), 181-235

\bibitem[FS2]{FS2}
\newblock R. Fintushel and R. Stern,
\newblock \emph{Knots, links and 4-manifolds},
\newblock Invent. Math. $\mathbf{134}$ (1998), 363-400

\bibitem[FS3]{FS3}
\newblock R. Fintushel and R. Stern,
\newblock \emph{Double node neighborhoods amd families of simply connected 4-manifolds with $b^+ = 1$},
\newblock Journal of the American Mathematical Society $\mathbf{19}$(1) (2006), 171-180

\bibitem[F]{F}
\newblock M. Freedman,
\newblock \emph{The topology of four-dimensional manifolds},
\newblock J. Diff. Geom. $\mathbf{17}$ (1982), 357-453

\bibitem[GH]{GH}
\newblock P. Griffiths and J. Harris,
\newblock \emph{Principles of algebraic geometry}, 
\newblock Wiley, New York, 1978





\bibitem[GP]{GP}
\newblock V. Guillemin and A. Pollack,
\newblock \emph{Differential Topology},
\newblock Prentice-Hall, 1974

\bibitem[GS]{GS}
\newblock R. Gompf and A. Stipsicz, 
\newblock \emph{4-manifolds and Kirby Calculus},
\newblock AMS Grad. Studies in Math. $\mathbf{20}$, 1999

\bibitem[HKK]{HKK}
\newblock J. Harer, A. Kas and R. Kirby, 
\newblock \emph{Handlebody decompositions of complex surfaces},
\newblock Memoirs of the AMS $\mathbf{62}$ no. 350, 1986

\bibitem[H1]{H1}
\newblock A. Hatcher, 
\newblock \emph{Algebraic Topology},
\newblock Cambridge University Press, Cambridge, 2002. http://www.math.cornell.edu/~hatcher/AT/ATpage.html

\bibitem[H2]{hatcher1}
\newblock A. Hatcher,
\newblock \emph{Vector Bundles and K-Theory},
\newblock http://www.math.cornell.edu/~hatcher/VBKT/VBpage.html, 2003

\bibitem[He]{He}
\newblock I. N. Herstein,
\newblock \emph{Topics in Algebra},
\newblock Blaisdell Publishing Company, 1964

\bibitem[Is]{I1}
\newblock C. Isham,
\newblock \emph{Modern Differential Geometry for Physicists},
\newblock World Scientific, 1999



\bibitem[Ko]{Ko2}
\newblock K. Kodaira,
\newblock \emph{On compact analytic surfaces: II}, 
\newblock Ann. Math, $\mathbf{77}$ (1963), 563-626

\bibitem[Ko2]{Ko3}
\newblock D. Kotschick,
\newblock \emph{On 4-manifolds homeomorphic to $\cpt \# 8 \cptbar$},
\newblock Invent. Math. $\mathbf{95}$ (1989), 591-600


\bibitem[LL]{LL1}
\newblock B. H. Li and T. J. Li,
\newblock \emph{Symplectic genus, minimal genus and diffeomorphisms},
\newblock preprint, http://front.math.ucdavis.edu/math.GT/0108227v1

\bibitem[LiLiu1]{LiLiu2}
\newblock T. J. Li and A. Liu,
\newblock \emph{Symplectic structures on ruled surfaces and a generalized adjunction formula},
\newblock Math. Res. Letters $\mathbf{2}$ (1995), 453-471

\bibitem[LiLiu2]{LiLiu1}
\newblock T. J. Li and A. Liu,
\newblock \emph{Uniqueness of the symplectic canonical class, surface cone and symplectic cone of 4-manifolds with $b^{+} = 1$},
\newblock Jour. Diff. Geom. $\mathbf{58}$ (2001), 331-370

\bibitem[LM]{LM}
\newblock P. Libermann and C.-M. Marle,
\newblock \emph{Symplectic Geometry and Analytical Mechanics},
\newblock D. Reidel Publishing Company, Dordrecht, Holland, 1987

\bibitem[L1]{L1}
\newblock W. B. R. Lickorish,
\newblock \emph{A representation of orientable, combinatorial 3-manifolds},
\newblock Annals of Math. $\mathbf{76}$ (1962), 531-540

\bibitem[L2]{L2}
\newblock W. B. R. Lickorish,
\newblock \emph{A finite set of generators for the homeotopy group of a 2-manifold},
\newblock Proc. Camb. Phil. Soc. $\mathbf{60}$ (1964), 769-778

\bibitem[L3]{L3}
\newblock W. B. R. Lickorish,
\newblock \emph{A finite set of generators for the homotopy group of a 2-manifold} corrigendum,
\newblock Proc. Camb. Phil. Soc. $\mathbf{62}$ (1966), 679-681



\bibitem[McS1]{MS2}
\newblock D. McDuff and D. Salamon,
\newblock \emph{Introduction to symplectic topology},
\newblock Oxford University Press, 1995

\bibitem[McS2]{MS}
\newblock D. McDuff and D. Salamon, 
\newblock \emph{A survey of symplectic 4-manifolds with $b_{2}^{+}=1$},
\newblock Turkish. Jour. Math. $\mathbf{20}$ (1996), 47-61




\bibitem[MH]{MH}
\newblock J. Milnor and D. Husemoller,
\newblock \emph{Symmetric bilinear forms},
\newblock Springer-Verlag, 1973

\bibitem[MiSt]{MiSt}
\newblock J. Milnor and J. Stasheff,
\newblock \emph{Characteristic Classes},
\newblock Princeton University Press, 1974

\bibitem[M]{M3}
\newblock R. Miranda
\newblock \emph{Persson's list of singular fibers for a rational elliptic surface}
\newblock Math. Zeit. Vol 205, no. 1, 1990, pages 191-211


\bibitem[Mo]{M}
\newblock B. Moishezon
\newblock \emph{Complex surfaces and connected sums of complex projective planes}
\newblock Lecture notes in Mathematics $\mathbf{603}$, Springer-Verlag, 1977

\bibitem[Mu]{Mu}
\newblock J. Munkres,
\newblock \emph{Elementary differential topology},
\newblock Ann. Math. Studies $\mathbf{54}$, Princeton University Press, 1966

\bibitem[O]{O1}
\newblock H. Osborn,
\newblock \emph{Vector Bundles, Volume 1, Foundations and Stiefel-Whitney Classes},
\newblock Academic Press, 1982

\bibitem[OzSt]{OzSt}
\newblock B. Ozbagci and A. I. Stipsicz,
\newblock \emph{Surgery on contact 3-manifolds and Stein surfaces}
\newblock Springer-Verlag, 2004


\bibitem[OS]{OS}
\newblock P. Ozsvath and Z. Szabo,  
\newblock \emph{On Park's exotic smooth four-manifolds},
\newblock `Geometry and Topology of manifolds', 253-260, Fields Inst. Commun., $\mathbf{47}$, Amer. Math. Soc., Providence, RI, 2005

\bibitem[P]{P1}
\newblock Jongil Park,
\newblock \emph{Simply connected 4-manifolds with $b_{2}^{+}=1$ and $c_{1}^{2}=2$.} 
\newblock Invent. Math. 159 (2005), no. $\mathbf{3}$, 657-667

\bibitem[PSS]{PSS}
\newblock Jongil Park, A. Stipsicz and Z. Szabo, 
\newblock \emph{Exotic smooth structures on $\mathbb{CP}^{2} \#  5 \overline{\mathbb{CP}^{2}}$}.
\newblock Math. Res. Lett. $\mathbf{12}$ (2005), no. 5-6, 701-712

\bibitem[Pe]{Pe}
\newblock U. Persson
\newblock \emph{Configurations of Kodaira fibers on rational elliptic surfaces}
\newblock Math. Zeit. Vol 205, no. 1, 1990, pages 1-47

\bibitem[R]{R2}
\newblock V. Rohlin,
\newblock \emph{New results in the theory of four dimensional manifolds},
\newblock Dokl. Akad. Nauk. USSR $\mathbf{84}$ (1952), 221-224; French transl. in [$\mathbf{RTP}$], pp. 14-16, 22-23

\bibitem[Ro]{Ro}
\newblock D. Rolfsen,
\newblock \emph{Knots and Links},
\newblock Publish or Perish, Berkeley, 1976

\bibitem[Sc]{Sc}
\newblock A. Scorpan,
\newblock \emph{The Wild World of 4-Manifolds},
\newblock American Mathematical Society, 2005

\bibitem[Se]{Se}
\newblock J.-P. Serre,
\newblock \emph{Formes bilin\'eaires sym\'etriques enti\`eres \`a discriminant $\pm 1$},
\newblock S\'eminaire Henri Cartan, tome 14 (1961-1962) Exp. 14-15, pgs 1-16, Secr\'etariat math\'ematique, Paris, 1961-1962

\bibitem[St]{St}
\newblock N. Steenrod,
\newblock \emph{The Topology of Fibre Bundles},
\newblock Princeton University Press, 1951

\bibitem[S]{S}
\newblock A. I. Stipsicz,
\newblock Personal communication

\bibitem[SS]{SS}
\newblock A. Stipsicz and Z. Szab\'o,
\newblock \emph{An exotic smooth structure on $\cpt \# 6 \cptbar$},
\newblock Geom. Topol. $\mathbf{9}$ (2005), 813 - 832


\bibitem[SSS]{SSS}
\newblock A. I. Stipsicz, Z. Szab\'o and A. Szil\'ard,
\newblock \emph{Singular fibers in elliptic fibrations on the rational elliptic surface},
\newblock Periodica Mathematica Hungarica Vol. 54 (2), 2007, pp. 137-162


\bibitem[Sy]{Sy}
\newblock M. Symington, 
\newblock \emph{Symplectic rational blowdowns},
\newblock Jour. Diff. Geom. $\mathbf{50}$ (1998), 505-518

\bibitem[We]{We}
\newblock Weisstein, Eric W.,
\newblock \emph{Cubic Formula} from Mathworld - A Wolfram Web Resource,
\newblock http://mathworld.wolfram.com/CubicFormula.html


\bibitem[Wu]{Wu}
\newblock W.-T. Wu,
\newblock \emph{Sur le classes caracteristique des structures fibr\'ees sph\'eriques},
\newblock Actualit\'es Sci. Industr. $\mathbf{1183}$ (1952)

\bibitem[Ya]{Ya}
\newblock K. Yasui,
\newblock \emph{Small exotic rational surfaces without 1- or 3-handles}, preprint,
\newblock http://front.math.ucdavis.edu/0807.0373



\end{thebibliography}
\end{document}